\documentclass[12pt]{article}

\usepackage{epsfig,color,amsmath,amssymb,euscript}
 \usepackage{mathabx}
\usepackage{subfigure}
\usepackage{fancyhdr}
\usepackage{natbib}
\usepackage{rotating}
\usepackage{amsmath}
\usepackage{amsfonts}
\usepackage{setspace}
\usepackage{fancyhdr}
\usepackage{lscape}
\usepackage{amsthm}
\usepackage{enumerate} 
\usepackage{threeparttable}
\usepackage{diagbox}
\usepackage{multirow}
\usepackage{makecell}
\usepackage{xcolor}
\usepackage{colortbl}
\usepackage[colorlinks,linkcolor = blue,anchorcolor=blue, citecolor= blue]{hyperref}
\usepackage{url} 

\usepackage{graphicx}
\usepackage{enumerate}
\usepackage{mathrsfs}
\usepackage{verbatim}
\usepackage{blkarray}

 \usepackage{graphicx}

\usepackage{authblk}
\usepackage{setspace}
\usepackage{threeparttable}
\usepackage{caption}

\usepackage[resetlabels]{multibib}

\allowdisplaybreaks

\newcites{S}{References} 

\usepackage{algorithm,algcompatible}
\algnewcommand\INPUT{\item[\textbf{Input:}]}%
\algnewcommand\OUTPUT{\item[\textbf{Output:}]}%

\makeatletter
\newcommand{\fdsy@scale}{1.0}
\newcommand\fdsy@mweight@normal{Book}
\newcommand\fdsy@mweight@small{Book}
\newcommand\fdsy@bweight@normal{Medium}
\newcommand\fdsy@bweight@small{Medium}

\DeclareFontFamily{U}{FdSymbolC}{}

\DeclareFontShape{U}{FdSymbolC}{m}{n}{
	<-7.1> s * [\fdsy@scale] FdSymbolC-\fdsy@mweight@small
	<7.1-> s * [\fdsy@scale] FdSymbolC-\fdsy@mweight@normal
}{}
\DeclareFontShape{U}{FdSymbolC}{b}{n}{
	<-7.1> s * [\fdsy@scale] FdSymbolC-\fdsy@bweight@small
	<7.1-> s * [\fdsy@scale] FdSymbolC-\fdsy@bweight@normal
}{}

\DeclareSymbolFont{arrows}{U}{FdSymbolC}{m}{n}
\SetSymbolFont{arrows}{bold}{U}{FdSymbolC}{b}{n}

\DeclareMathSymbol{\upvDash}{\mathrel}{arrows}{233}

\DeclareMathSymbol{\upmodels}{\mathrel}{arrows}{237}
\makeatother

\usepackage{enumerate}
\usepackage{natbib}
\usepackage{url} 
\usepackage{epsfig,color,amsmath,amssymb,euscript}
\usepackage{fancyhdr}
\usepackage{rotating}
\usepackage{amsthm, mathrsfs}
\usepackage{amsfonts}
\usepackage{setspace}
\usepackage{lscape}
\usepackage{threeparttable}
\usepackage{diagbox}
\usepackage{multirow}
\usepackage{makecell}
\usepackage{xcolor}
\usepackage{enumitem}
\usepackage{graphicx}
\usepackage{chngcntr}

\newtheorem{theorem}{Theorem}

\newtheorem{lemma}{Lemma}
\newtheorem{proposition}{Proposition}

\theoremstyle{definition}
\newtheorem{as}{Condition}

\newtheorem{remark}{Remark}

\allowdisplaybreaks[3]

\newcommand{\cov}{{\rm Cov}}
\newcommand{\diag}{{\rm diag}}

\newcommand{\var}{\mbox{Var}}

\newcommand{\bE}{\mathbb {E}}
\newcommand{\bP}{{\mathbb{P}}}
\newcommand{\bR}{{\mathbb R}}

\def\T{{\mathrm{\scriptscriptstyle \top} }}

\newcommand{\bA}{{\mathbf A}}
\newcommand{\bB}{{\mathbf B}}
\newcommand{\bF}{{\mathbf F}}

\newcommand{\bH}{{\mathbf H}}
\newcommand{\bI}{{\mathbf I}}

\newcommand{\bM}{{\mathbf M}}

\newcommand{\bW}{{\mathbf W}}

\newcommand{\ba}{{\mathbf a}}
\newcommand{\bb}{{\mathbf b}}
\newcommand{\bfd}{{\mathbf d}}
\newcommand{\bc}{{\mathbf c}}

\newcommand{\bff}{{\mathbf f}}
\newcommand{\bg}{{\mathbf g}}
\newcommand{\bh}{\mathbf h}

\newcommand{\br}{{\mathbf r}}

\newcommand{\bs}{{\mathbf s}}
\newcommand{\bu}{{\mathbf u}}
\newcommand{\bv}{{\mathbf v}}
\newcommand{\bw}{{\mathbf w}}
\newcommand{\bx}{{\mathbf x}}
\newcommand{\by}{{\mathbf y}}

\newcommand{\bfeta}  {\boldsymbol{\eta}}

\newcommand{\bdelta} {\boldsymbol{\delta}}

\newcommand{\bSigma}{\boldsymbol{\Sigma}}

\newcommand{\bgamma}{\boldsymbol{\gamma}}

\newcommand{\bTheta} {\boldsymbol{\Theta}}

\newcommand{\bPhi} {\boldsymbol{\Phi}}
\newcommand{\bPsi} {\boldsymbol{\Psi}}

\newcommand{\bmu} {\boldsymbol{\mu}}

\newcommand{\bGamma} {\boldsymbol{\Gamma}}

\newcommand{\bD}{{\mathbf D}}
\newcommand{\bzero}{{\mathbf 0}}

\newcommand{\bPi}{\boldsymbol{\Pi}}
\newcommand{\bXi}{\boldsymbol{\Xi}}
\newcommand{\bchi}{\boldsymbol{\chi}}

\newcommand{\cA}{{\mathcal A}}

\newcommand{\cF}{{\mathcal F}}

\newcommand{\cH}{{\mathcal H}}
\newcommand{\cI}{{\mathcal I}}
\newcommand{\cJ}{{\mathcal J}}
\newcommand{\cK}{{\mathcal K}}

\newcommand{\cS}{{\mathcal S}}

\newcommand{\cV}{{\mathcal V}}
\newcommand{\cW}{{\mathcal W}}
\newcommand{\cX}{{\mathcal X}}

\newcommand{\ext}{{\rm ext}}

\def\spacingset#1{\renewcommand{\baselinestretch}%
	{#1}\small\normalsize} \spacingset{1}

\newcommand{\blind}{1}
\makeatletter

\def\singlespace{\def\baselinestretch{1}\@normalsize}

\begin{document}

\if1\blind
{
\spacingset{1.25}
  \title{\bf Statistical Inference for  High-Dimensional Spectral Density Matrix}
  \author[1,2]{Jinyuan Chang}
  \author[3]{Qing Jiang}
  \author[4]{Tucker McElroy}
  \author[5]{Xiaofeng Shao}

\affil[1]{\it \small Joint Laboratory of Data Science and Business Intelligence, Southwestern University of Finance and Economics, Chengdu, China}
\affil[2]{\it \small Academy of Mathematics and Systems Science, Chinese Academy of Sciences, Beijing, China}
\affil[3]{\it \small Faculty of Arts and Sciences, Beijing Normal University, Zhuhai, China}
\affil[4]{\it \small Research and Methodology Directorate, U.S. Census Bureau} 
\affil[5]{\it \small Department of Statistics and Data Science, and Department of Economics, Washington University in St Louis}

\setcounter{Maxaffil}{0}
		
		\renewcommand\Affilfont{\itshape\small}
		\date{\vspace{-5ex}}
		\maketitle
} \fi

\if0\blind
{
  \bigskip
  \bigskip
  \bigskip
  \begin{center}
    {\LARGE\bf Statistical Inference for  High-Dimensional Spectral Density Matrix}
\end{center}
  \medskip
} \fi

\spacingset{1.4}
\begin{abstract}
The spectral density matrix is a fundamental object of interest in time series analysis, and it encodes both contemporary and dynamic linear relationships between   component processes of the multivariate system. In this paper we develop novel inference procedures for the spectral density matrix in the high-dimensional setting. Specifically, we introduce a new global testing procedure to test the nullity of the cross-spectral density for a given set of frequencies and across pairs of component indices. For the first time,  both Gaussian approximation and  parametric bootstrap methodologies are employed to conduct inference for a high-dimensional parameter formulated in the frequency domain, and  new technical tools are developed to  provide asymptotic guarantees of the size accuracy and power for global testing.
We further propose a multiple testing procedure for simultaneously testing the nullity of the cross-spectral density at a given set of frequencies. The method is shown to control the false discovery rate. Both numerical simulations and a real data illustration demonstrate the usefulness of the proposed testing methods.
\end{abstract}

\noindent%

{
\spacingset{1.69}
\setlength{\abovedisplayskip}{0.2\baselineskip}
\setlength{\belowdisplayskip}{0.2\baselineskip}
\setlength{\abovedisplayshortskip}{0.2\baselineskip}
\setlength{\belowdisplayshortskip}{0.2\baselineskip}


\section{Introduction}\label{sec:introduction}

The spectral density matrix plays an important role in time series analysis, as it completely characterizes the second order properties of a multivariate stationary time series; moreover, it is a crucial quantity in the  modeling,  inference, and prediction of time series. It has been used in discriminant analysis for time series \citep{TK_2000}, identification and estimation of generalized dynamic factor models \citep{FHLR_2000}, and non-correlation/independence testing for two time series \citep{E_2007,S_2009}, to name a few examples. Estimation and inference for the univariate spectral density and multivariate low-dimensional spectral density matrix have a long history, and classical methods based on the smoothed periodogram or lag-window estimates have been well documented in classical textbooks such as  \cite{H_1970}, \cite{P_1981}, \cite{B_2001} and \cite{BD_1991}. More recent contributions can be found in \cite{SW_2007}, \cite{LW_2009}, and \cite{WZ_2018}, among others.

With the advancement of science and technology, high-dimensional time series have been increasingly collected in many areas, such as neuroimaging studies, finance, and climate science, as well as official statistics.
This situation motivates the development of new statistical methodology and theory to   accommodate the data's high dimensionality. For example,	the classical smoothed periodogram estimate of the multivariate spectral density matrix can be improved by a  shrinkage approach; see  \cite{BV_2008,BV_2009}, \cite{FO_2011}, and \cite{FV_2014}, among others. The popular regularization approaches used in covariance and precision matrix estimation, such as  graphical LASSO \citep{YL_2007}, thresholding \citep{BL_2008},  and constrained $L_1$ minimization \citep{CLL_2011}, have been extended to estimate either the high-dimensional spectral density matrix or its inverse by  \cite{SLKB_2018}, \cite{FLLY_2019} and \cite{T_2022}. Related work along this line also includes
\cite{ZW_2021}, which established convergence rates of regularized estimates of the spectral density matrix and its inverse under a high-dimensional locally stationary framework. Recently, \cite{MM_2023}
introduced a new estimator of the high-dimensional spectral density matrix --
under the assumptions of low rank and sparse structure -- via minimizing a quadratic loss under a nuclear norm, together with a  $L_1$ norm constraint to control the latent rank and residual sparsity pattern.


While the literature on estimating the high-dimensional spectral density matrix or its inverse has been growing rapidly, there seems to be  relatively less work devoted to inference for the spectral density matrix, which is our focus here.	In this paper, we aim to develop  new theory and methods for the inference of the spectral density matrix of a high-dimensional weakly stationary time series. In particular, we propose a new maximum-type test statistic to test for the joint hypothesis that the cross-spectral density is zero for a given set of frequencies and pairs of indices. We establish a Gaussian approximation result for our maximum-type test statistic, and provide a computationally feasible parametric bootstrap-based approach to approximate its finite sample distribution. Building on  the newly established theory,  we further develop a multiple testing procedure to recover the support of the spectral density matrix for a given set of frequencies. 

Since the seminal work of \cite{CCK_2013}, the technique of Gaussian approximation has undergone rapid developments; see \cite{CCK_2017, CCK_2019},  \cite{Chen2018}, \cite{CK2019}, \cite{FK2020}, and \cite{CCKK_2022a}, among others. A recent review of high-dimensional data bootstrap and Gaussian approximation is provided by \cite{CCKK_2022b}.
The extension of Gaussian approximation to high-dimensional time series was first developed by \cite{ZW2017} and \cite{ZC2018}, but their focus was mainly on  inference for quantities such as means and  autocovariance matrices, which are formulated in the time domain. By contrast, our parameter of interest is the cross-spectral density measured over a set of frequencies, and our technical treatment is rather different from that of \cite{ZW2017} and \cite{ZC2018}.  The closest works to our paper, from a technical perspective, are \cite{CJS2021} and \cite{CCW2021}.

Let $\{\bx_t\}$ be a $p$-dimensional weakly stationary time series. Denote by $\bF(\omega)=\{f_{i,j}(\omega)\}_{p\times p}\in\mathbb{C}^{p\times p}$ the spectral density matrix of $\{\bx_t\}$ at frequency $\omega\in[-\pi,\pi)$. 
Our main goal in this paper is to establish the Gaussian approximation to the distribution of
\begin{align}\label{eq:mathcalT}
	\mathcal{T}=	\sup_{\omega\in\cJ}\max_{(i,j)\in\mathcal{I}}\bigg|\sqrt{\frac{n}{l_n}}\{\hat{f}_{i,j}(\omega)-f_{i,j}(\omega)\}\bigg|^2 \,,
\end{align}
where $\cI\subset\{1,\ldots,p\}^2$, $\cJ\subset[-\pi,\pi)$, and $\hat{f}_{i,j}(\omega)$ is some kernel-based estimate of $f_{i,j}(\omega)$ with bandwidth $l_n$. 
Such a Gaussian approximation result provides a technical tool for the inference of the high-dimensional spectral density matrix. For each given $(i,j)\in\mathcal{I}$ and $\omega\in\mathcal{J}$,  $|\sqrt{n/l_n}\{\hat{f}_{i,j}(\omega)-f_{i,j}(\omega)\}|^2$ usually converges in distribution to the  sum of squares of two correlated normal random variables. Our setting is quite different from those considered in   existing works, and this generates several technical challenges for establishing the Gaussian approximation to the distribution of $\mathcal{T}$. 
See our detailed discussion in Remark \ref{rk:1}(d) in Section \ref{sec:ga}. 
In addition,  we  apply our Gaussian approximation results to obtain   false discovery rate (FDR) control in the multiple testing procedure, which not only expands  the application of our Gaussian approximation results, but also extends the validity of FDR control to the high-dimensional time series setting.

As we mentioned earlier, the literature on the inference for the high-dimensional spectral density matrix is scarce, and we are only aware of two recent papers. Motivated by testing the mutual independence of the component series in a $p$-dimensional complex-valued Gaussian time series, \cite{RLV_2021} investigated the asymptotic distribution for the maximum of smoothing-based estimators of the coherence (the standarized cross-spectral density), and showed that its null-distribution converges to a Gumbel limiting distribution when $p/n=o(1)$, as well as some other conditions on the smoothing span.
\cite{KP_2022} developed new statistical inference procedures for coherences and partial coherences of a $p$-dimensional real-valued time series. They addressed the estimation of partial coherence using a debiased approach, and developed a testing procedure for the null hypothesis that the partial coherences do not exceed some user-specified threshold value within a frequency band of interest.
When $p=o(n^{\tau})$ for some constant $\tau>0$, they showed that the limiting distribution for the maximum of sample partial coherence over frequencies
is the classical Gumbel distribution (or its variant).

In contrast to these two works, we focus on the Gaussian approximation to the distribution of $\mathcal{T}$ defined as \eqref{eq:mathcalT}, which can also be used to conduct inference for coherences (see Section \ref{sec:preliminary} for details).
It is worth noting that the maximum-type test statistics in both \cite{RLV_2021} and \cite{KP_2022} are taken over a set of frequencies that are equally spaced, with the spacing having larger order of magnitude than  $2\pi/n$ (the spacing for consecutive Fourier frequencies). This construction appears to be necessary in order for these authors to obtain a Gumbel limiting distribution. Technically speaking,
the derivation of the  Gumbel limiting distribution relies on the weak dependency among test statistics across frequencies. 
Hence, in \cite{RLV_2021} and \cite{KP_2022}, the spacing between the frequencies cannot be too close. 
In addition, \cite{RLV_2021} requires the time series to be mean-zero stationary and complex-valued Gaussian time series with mutually independent component time series under the null. 
By contrast, we allow some dependence between different components of time series under the null, and we consider   real-valued time series in this paper.  Our proposed method can also be extended to  complex-valued time series but is not pursued here.
 \cite{KP_2022} mainly focuses on the hypothesis testing and support recovery of the partial coherence matrix of real-valued high-dimensional time series.
From the viewpoint of asymptotic approximation, it is a common belief that the  maximum-type statistics usually converge slowly to the Gumbel distributions \citep{Hall_1991}, and the  bootstrap is often employed to provide a better finite sample approximation.
By contrast, our Gaussian approximation theory does not require any dependence structure among 
the quantities $|\sqrt{n/l_n}\{\hat{f}_{i,j}(\omega)-f_{i,j}(\omega)\}|^2$ across frequencies, so
we  have no requirements on the spacing of frequencies $\omega\in\mathcal{J}$ as long as they are distinct and fall into $[-\pi,\pi)$.  Notice that the limiting distribution of $\mathcal{T}$ defined in \eqref{eq:mathcalT} may not admit a closed form (or even does not exist), but its finite sample distribution can nevertheless be well-approximated by its parametric bootstrap counterpart. The technical tools employed in these two papers and ours are very different.


A primary application of our  inference procedure for the high-dimensional spectral density matrix  is the division of large time series databases into batches suitable for joint analysis.  For example,  a database of county-level time series can be organized into batches by state, and it is of interest to know if there is significant content in the cross-spectra of two batches; if so, there may be merit in jointly modeling the batches, but otherwise analysis can proceed upon smaller collections.
A second   motivation for our work comes from the processing and analysis of seasonal time series at statistical agencies such as the U.S. Census Bureau (USCB). At USCB, many weekly, monthly, and quarterly economic time series are published, all of which exhibit seasonality of varying types and degrees. The USCB performs seasonal adjustment of such time series in order to remove  seasonality; both the original and seasonally adjusted data are published for public use.
Seasonal adjustment is a vast world-wide undertaking, with the statistical agencies of all developed countries (as well as many private companies) adjusting thousands or millions of time series every month (or quarter).
One important task in practice is to evaluate the effectiveness of seasonal adjustment, that is,  whether the strong seasonality has been properly removed or whether there is an issue of over-adjustment \citep{mcelroy_2021,mcelroy_2022}.
For example, seasonal adjustment methods that are intended to remove seasonality sometimes --
when applied to a time series with weak seasonality -- result in zero power spectrum at a ``seasonal frequency'', e.g.,  $f_{i,i}(\pi/2)=0$ for a quarterly time series.  
This motivates a joint testing problem over many time series, where we wish to identify which series have zero power at the seasonal frequencies (i.e., have been over-adjusted).

The rest of the paper is organized as follows. Section~\ref{main} presents the spectral density matrix estimate, 
a general Gaussian approximation procedure, and the related theory.  Section~\ref{applications} contains two applications, including global testing and support recovery via a multiple testing procedure. Section~\ref{numerical} investigates the finite sample performance of the proposed testing procedures via numerical simulation, and Section~\ref{data}  provides an illustration based on  county-level  quarterly time series
of new hires. Section~\ref{discussion} concludes the paper. All technical details are relegated to the supplementary material. 
For practical convenience  we have developed -- in the R package {\verb"HDTSA"} \citep{CHLY2024} -- two R-functions {\verb"SpecTest"} and {\verb"SpecMulTest"}  that implement the global testing and the multiple testing procedures in an automatic manner, respectively.

{\it Notation.} Denote by $I(\cdot)$ the indicator function. For any positive integer $q\geq2$, we write $[q]=\{1,\ldots,q\}$, and let
$[q]^2=[q]\times [q]$ denote the Cartesian product
of $[q]$. Let $|\cF|$ be the cardinality of  a countable set  $\cF$.
For two positive real-valued sequences  $\{a_n\}$ and $\{b_n\}$, we write $a_n\lesssim b_n$  if $\limsup_{n\rightarrow\infty}a_n/b_n\leq c_0$ for some positive constant $c_0$, $a_n\asymp b_n$ if $a_n\lesssim b_n$ and $b_n\lesssim a_n$ hold simultaneously, and $a_n\ll b_n$  if $\limsup_{n\rightarrow\infty}a_n/b_n=0$.
For a complex-valued number $x$,  denote by $|x|$  its modulus. The operator $\otimes$ denotes the Kronecker product.
For any real-valued numbers $x$ and $y$, we write $|x|_{+}=\max(0,x)$ and $x\vee y=\max(x,y)$. 
Denote by $\mathbb{S}^{q-1}$ the $q$-dimensional  unit sphere. For a $q$-dimensional vector $\ba$, denote by $\ba_{\mathcal{L}}$ the subvector of $\ba$ consisting of the components indexed by a given index set $\mathcal{L}\subset[q]$. For any $q_1\times q_2$ matrix ${\bf M}=(m_{i,j})_{q_1\times q_2}$, let $|{\bf M}|_\infty=\max_{i\in[q_1],j\in[q_2]}|m_{i,j}|$. Let ${\bf 1}_d$ and ${\bf I}_d$ be, respectively, a $d$-dimensional vector with all components being $1$, and a $d\times d$  identity matrix.

\section{Some technical results} \label{main}  
\subsection{Preliminary}\label{sec:preliminary}
Let $\bx_t=(x_{1,t},\ldots,x_{p,t})^\T$ be a $p$-dimensional weakly stationary time series with mean vector $\bmu=\mathbb{E}(\bx_t)$ and autocovariance matrix $\bGamma(k)\equiv\{\gamma_{i,j}(k)\}_{p\times p}=\cov(\bx_{t+k},\bx_t)$. When $\sum_{k=-\infty}^{\infty}|\gamma_{i,j}(k)|<\infty$ for each $i,j\in[p]$, we can  define the spectral density matrix $\bF(\omega)$ for $\omega\in[-\pi,\pi)$ as
\begin{align*}
	\bF(\omega)\equiv\{f_{i,j}(\omega)\}_{p\times p}=\frac{1}{2\pi}\sum_{k=-\infty}^\infty\bGamma(k)e^{-\iota k\omega}\,,
\end{align*}
where $\iota=\sqrt{-1}$. 
Given the observations $\mathcal{X}_n=\{\bx_1,\ldots,\bx_n\}$, we can estimate $\bF(\omega)$ by 
\begin{equation}\label{eq:spectralest}
	\widehat{\bF}(\omega)\equiv\{\hat{f}_{i,j}(\omega)\}_{p\times p}=\frac{1}{2\pi}\sum_{k=-l_n}^{l_n} \mathcal{W}\bigg(\frac{k}{l_n}\bigg)\widehat{\bGamma}(k)e^{-\iota k\omega}\,,
\end{equation}
where $\mathcal{W}(\cdot)$ is a symmetric kernel function, $l_n=o(n)$ is the bandwidth, and
\begin{align}\label{eq:hatGammak}
	\widehat{\bGamma}(k)=\frac{1}{n}\sum_{t=\max(1,-k+1)}^{\min(n,n-k)} (\bx_{t+k}-\bar{\bx})(\bx_t-\bar{\bx})^\T 
\end{align}
with $\bar{\bx}=n^{-1}\sum_{t=1}^n\bx_t$. To reduce the bias involved in \eqref{eq:spectralest}, for some constant $c\in(0,1]$, we adopt  the flat-top kernel suggested by \cite{Politis_2011}:
\begin{align} \label{eq:Wu}
	\mathcal{W}(u)= I(|u|\leq c) + \frac{|u|-1}{c-1} I(c<|u|\leq 1) \,.
\end{align}

For given $(i,j)\in [p]^2$ and $i<j$, 
define the coherence spectrum at frequency $\omega$ by
\[
{\rm coh}_{i,j}(\omega)=\frac{f_{i,j}(\omega)}{\{f_{i,i}(\omega) f_{j,j}(\omega)\}^{1/2}}\,.
\]
As mentioned in \cite{P_1981}, the coherence may be interpreted as the correlation coefficient between the random coefficients in the spectral representations of the components in $x_{i,t}$ and $x_{j,t}$ at frequency $\omega$.
Thus ${\rm coh}_{i,j}(\omega)=0$ for all $\omega\in [-\pi,\pi)$ is equivalent to $f_{i,j}(\omega)=0$
for all $\omega\in [-\pi,\pi)$, which implies the  two processes
$\{x_{i,t}\}$ and $\{x_{j,t}\}$ are linearly unrelated at all lags. Within the scope of linear time series models,  the joint modeling of $\{ (x_{i,t},x_{j,t}) \}$ can be simplified by modeling the linear serial dependence of $\{ x_{i,t} \}$ and $\{ x_{j,t} \}$ separately, provided that $\max_{\omega\in[-\pi,\pi) }|f_{i,j}(\omega)|=0$. 
For a $p$-dimensional weakly stationary time series $\{ \bx_t \}$, it is thus of great importance to recover the support
of nonzero coherence, i.e.,
\begin{align}\label{eq:sf}
	\mathscr{S}_f:=\bigg\{(i,j) \in [p]^{2}\setminus \{(1,1),\ldots,(p,p)\} \, : \, \max_{  \omega\in [-\pi,\pi)} |f_{i,j}(\omega)| \not= 0\bigg\}\,.
\end{align}
It is worth noting that the coherence also  plays  an important role in characterizing the functional connectivity  between neural regions within the brain based on functional magnetic resonance imaging  data; see \cite{SMD_2004} and \cite{B_2016}.
Additionally, in terms of joint modeling in the frequency domain, it is of interest to understand the behavior of the cross-spectrum (or coherence) for a specific set of frequencies. For example, for a quarterly time series, we are interested in the coherence at the ``seasonal frequency'' $\pi/2$ and the ``trend frequency'' $0$; if joint  modeling of low-frequency fluctuation is of particular interest, we may want to focus on a pre-specified frequency interval such as $[-\omega_0,\omega_0]$, where $\omega_0$ is determined by the user based on the characteristics of time series, such as sample size and expected seasonal behavior (e.g., monthly or quarterly).  For the application of seasonal over-adjustment mentioned in Section \ref{sec:introduction}, we also wish to   consider the joint testing problem  
$H_0: f_{i,i}(\omega_0)=0$
for all $i\in[p]$, and the recovery of the support of indices that correspond to nonzero power auto-spectrum at $\omega_0$, where $\omega_0\in\{-\pi,-\pi/2,0,\pi/2\}$.

To study the above-mentioned testing and support recovery problems,  we define
\begin{align*}
	T_{n}(\omega;\mathcal{I})=\max_{(i,j)\in\mathcal{I}}\bigg|\sqrt{\frac{n}{l_n}}\{\hat{f}_{i,j}(\omega)-f_{i,j}(\omega)\}\bigg|^2\,
\end{align*}
for a given index set $\mathcal{I}\subset[p]^2$ and $\omega\in[-\pi,\pi)$.
Given a subset $\mathcal{J}\subset[-\pi,\pi)$, we will first establish in Section \ref{sec:ga} the Gaussian approximation to the distribution of $\sup_{\omega\in\mathcal{J}}T_n(\omega;\mathcal{I})$. 
In practical problems, we mainly focus on two kinds of configurations for $\mathcal{J}$, viz. (a) $\mathcal{J}=\{\omega_1,\ldots,\omega_K\}$ is a set with $K$ distinct frequencies $-\pi\le \omega_1<\omega_2<\cdots<\omega_K < \pi$, where $K$ may grow with the sample size $n$; (b) $\mathcal{J}=[\omega_{L},\omega_U]$ is an interval with $-\pi\leq \omega_L<\omega_U\leq \pi$, where $\mathcal{J}=[\omega_L,\pi)$ if $\omega_U=\pi$.  Based on the established Gaussian approximation theory, we can address the following inference problems of a high-dimensional time series:
\begin{itemize}
	\item (Global hypothesis testing). For given $\mathcal{I}$ and $\mathcal{J}$, consider the testing problem
	$
	H_{0}:f_{i,j}(\omega)=0$ for any $(i,j)\in\mathcal{I}$ and $\omega\in\mathcal{J}$ versus $H_{1}:H_{0}$ is not true.
	Here $\mathcal{I}$ and $\mathcal{J}$ can be chosen according to the user's interest. If we are interested in simultaneously testing the zero power auto-spectrum at a seasonal frequency $\omega_0\in  \{-\pi, -\pi/2, 0, \pi/2\}$, we can select $\mathcal{I}=\{(1,1),\ldots,(p,p)\}$ and $\mathcal{J}=\{\omega_0\}$. If testing for the zero cross-spectrum at all frequencies is of interest, we can set $\mathcal{I}=[p]^{2} \setminus \{(1,1),\ldots,(p,p)\}$ and $\mathcal{J}=[-\pi,\pi)$. See Section \ref{subsec:global} for details.
	
	\item (Support recovery). In the event that the global null is rejected, we are interested in the support of nonzero elements, that is $\{(i,j)\in \mathcal{I}:\sup_{\omega\in \mathcal{J}} |f_{i,j}(\omega)|\not=0\}$. For example, to recover $\mathscr{S}_f$ defined as \eqref{eq:sf}, we can consider a multiple testing problem with $(p^2-p)/2$ marginal hypotheses
	$
	H_{0,i,j}: \max_{\omega\in [-\pi,\pi) } |f_{i,j}(\omega)|=0$ versus $H_{1,i,j}: \max_{\omega\in [-\pi,\pi) } |f_{i,j}(\omega)|\not=0$,
	with $(i,j)\in [p]^{2}$ and $i<j$. For each given $(i,j)$, we can obtain the p-value for the marginal null hypothesis $H_{0,i,j}$ by our established Gaussian approximation result for the distribution of $\sup_{\omega\in\mathcal{J}}T_n(\omega;\mathcal{I})$ with selecting $\mathcal{I}=\{(i,j)\}$ and $\mathcal{J}=[-\pi,\pi)$. Based on the $(p^2-p)/2$ obtained  p-values, we can propose a FDR control procedure to estimate $\mathscr{S}_f$ as in \eqref{eq:sf}. See Section \ref{subsec:multiple} for details.
\end{itemize}

\subsection{A general Gaussian approximation procedure}\label{sec:ga}
Write $r=|\mathcal{I}|$, $\tilde{n}=n-2l_n$, $\bGamma(k)\equiv\{\gamma_{i,j}(k)\}_{p\times p}$ and $\widehat{\bGamma}(k)\equiv\{\hat{\gamma}_{i,j}(k)\}_{p\times p}$. Define $\mathring{\bx}_t\equiv(\mathring{x}_{1,t},\ldots,\mathring{x}_{p,t})^{\T}=\bx_t-\bmu$ for any $t\in[n]$. Let $\bchi(\cdot)=\{\chi_1(\cdot),\chi_2(\cdot)\}$ be a given bijective mapping from $[r]$ to $\mathcal{I}$ such that for any $(i,j)\in\mathcal{I}$, there exists a unique $\ell\in[r]$ satisfying $(i,j)=\bchi(\ell)$.  
For each $t\in[\tilde{n}]$ and $\ell\in[r]$, we define a $(2l_n+1)$-dimensional vector
\begin{align}\label{eq:cjt}
	\bc_{\ell,t}=\frac{1}{2\pi}\{\mathring{x}_{\chi_1(\ell),t}\mathring{x}_{\chi_2(\ell),t+l_n}-{\gamma}_{\bchi(\ell)}(-l_n),\ldots,\mathring{x}_{\chi_1(\ell),t+2l_n}\mathring{x}_{\chi_2(\ell),t+l_n}-{\gamma}_{\bchi(\ell)}(l_n)\}^\T\,.
\end{align}
Notice that $\bE(\bc_{\ell,t})=\bzero$. 
Let $\bc_t=(\bc_{1,t}^\T,\ldots,\bc_{r,t}^\T)^\T$
with $\bc_{\ell,t}$ defined in \eqref{eq:cjt}. Define
\begin{align}\label{eq:etacheckext}
	\check{\bfeta}^{\ext}(\omega)\equiv\{\check{\eta}^{\ext}_1(\omega),\ldots,\check{\eta}_{2r}^{\ext}(\omega)\}^\T=\{\bI_r\otimes \bA(\omega) \}\frac{1}{\sqrt{n}}\sum_{t=1}^{\tilde{n}}\bc_{t}  \,,
\end{align}
where
\begin{equation}\label{eq:A}
	\bA(\omega)=\frac{1}{\sqrt{l_n}}\left(
	\begin{array}{ccc}
		\cos(-l_n\omega) & \cdots & \cos(l_n\omega) \\
		-\sin(-l_n\omega) & \cdots & -\sin(l_n\omega) \\
	\end{array}
	\right)\textrm{diag}\{
	\mathcal{W}(-l_n/l_n),\ldots,\mathcal{W}(l_n/l_n)\}\,.
\end{equation}
Denote the long-run covariance of the sequence $\{\bc_t\}_{t=1}^{\tilde{n}}$ by \begin{equation}\label{eq:Xi}
	\bXi={\rm Var}\bigg(\frac{1}{\sqrt{\tilde{n}}} \sum_{t=1}^{\tilde{n}} \bc_{t}\bigg)\,.
\end{equation}
For any $\omega_1,\omega_2\in[-\pi,\pi)$, we define 
\begin{align}\label{eq:long_run}
	\bSigma(\omega_1,\omega_2)
	= \{\bI_r\otimes \bA(\omega_1)\} \bXi \{\bI_r\otimes \bA^\T(\omega_2)\} \,.
\end{align}
Then $\cov\{\check{\bfeta}^{\ext}(\omega_1),\check{\bfeta}^{\ext}(\omega_2)\}=(\tilde{n}/n)\bSigma(\omega_1,\omega_2)$ for any $\omega_1,\omega_2\in[-\pi,\pi)$ with $\check{\bfeta}^{\ext}(\omega)$ defined as \eqref{eq:etacheckext}.
To investigate the limiting distribution of $\sup_{\omega\in\mathcal{J}}T_n(\omega;\mathcal{I})$, we need the following regularity conditions. The validity of these conditions are discussed in Section A of the supplementary material. 

\begin{as}\label{as:moment}
	There exist some universal constants $C_1>0$ and $C_2>1$ such that $
	\mathbb{E}\{\exp(C_1|x_{j,t}|^{2})\}\leq C_2
	$ for any $t\in[n]$ and $j\in[p]$.
\end{as}

\begin{as}\label{as:betamixing}
	Let $\mathcal{F}_{-\infty}^u$ and $\mathcal{F}_{u+k}^{+\infty}$ be the $\sigma$-fields generated respectively by $\{\bx_{t}\}_{t\leq u}$ and $\{\bx_t\}_{t\geq u+k}$. Define
	$$
		 \alpha_n(k) :=\sup_t\sup_{(A,B)\in\mathcal{F}_{-\infty}^t\times\mathcal{F}_{t+k}^{+\infty}}|\bP(A\cap B)-\bP(A)\bP(B)|\,.$$
	 There exist some universal constants $C_3>0$ and $C_4>0$ such that  $\alpha_n(k) \leq C_3\exp(-C_4k)$ for any positive $k$.
\end{as}

\begin{as}\label{as:eigen}
	There exists a universal constant $C_5>0$ such that
	$ \inf_{\omega\in\mathcal{J} } \bfd^{\T}\bSigma(\omega,\omega)\bfd\geq C_5 $ for any $\bfd\in\bigcup_{j=1}^{r}\{\bfd\in\mathbb{S}^{2r-1}: \bfd_{S_j}\in\mathbb{S}^{1}\}$ with $S_j=\{2j-1,2j\}$.
\end{as}

Notice that in these conditions there are no explicit requirements on the cross-series dependence, and both weak and strong cross-series dependence are allowed by our theory. In particular, the marginal covariance matrix can be banded or AR(1)-type, representing weak cross-series dependence. Or it can be a compound symmetric matrix, which implies strong cross-series dependence.   As a result, we can establish
the following  Gaussian approximation result.

\begin{proposition}\label{tm:1}
	Assume $r\geq n^\kappa$ for some sufficiently small constant $\kappa>0$, and let Conditions {\rm\ref{as:moment}--\ref{as:eigen}} hold. As $n\rightarrow\infty$, the following two assertions are valid.
	
	{\rm (i)} If $\mathcal{J}=\{\omega_1,\ldots,\omega_K\}$ and $\log(Kr)\ll n^{1/9}l_n^{-1}\log^{-8/3}(l_n)$, with the bandwidth $l_n$ in \eqref{eq:spectralest} satisfying $l_n\log^{8/3}(l_n)\ll n^{1/9}$ and $l_n\geq \max\{2,C\log(Kr)\}$ for some sufficiently large constant $C>0$, then
	$$		\sup_{u\geq 0}\bigg|\mathbb{P}\bigg\{\sup_{\omega\in\mathcal{J}}{T}_n(\omega;\mathcal{I})\leq u\bigg\} -\mathbb{P}\bigg\{\max_{j\in [Kr]} (s_{n,2j-1}^2+s_{n,2j}^2)\leq u\bigg\}\bigg|\lesssim \frac{l_n\log^{2/3}(l_n)\log(Kr)}{n^{1/9}} $$
	for a $(2Kr)$-dimensional normally distributed random vector $\bs_{n,\by}=(s_{n,1},\ldots,s_{n,2Kr})^\T\sim \mathcal{N}(\bzero,\bH\bXi\bH^\T)$,
	where  $\bH=\{\bI_r\otimes \bA^\T(\omega_1),\ldots,\bI_r\otimes \bA^\T(\omega_K)\}^\T$
	with $\bA(\omega)$ defined in {\rm\eqref{eq:A}}, and $\bXi$ is defined in {\rm\eqref{eq:Xi}}.

	{\rm (ii)} If $\mathcal{J}=[\omega_L,\omega_U]$  and $\log r\ll n^{1/9}l_n^{-1}\log^{-8/3}(l_n)$, with the bandwidth $l_n$ in \eqref{eq:spectralest} satisfying $l_n\log^{8/3}(l_n)\ll n^{1/9}$ and $l_n\geq \max( 2,C\log r )$ for some sufficiently large constant $C>0$, then
	$$ 	\sup_{u\geq 0} \bigg|\mathbb{P}\bigg\{\sup_{\omega\in\mathcal{J}}{T}_n(\omega;\mathcal{I})\leq u\bigg\} -\mathbb{P}\bigg[\sup_{\omega\in\mathcal{J}}\max_{j\in[r]} \{g_{n,2j-1}^2(\omega)+g_{n,2j}^2(\omega)\}\leq u\bigg]\bigg| \lesssim \frac{l_n\log^{2/3}(l_n)\log r}{n^{1/9}} $$
	for a $(2r)$-dimensional Gaussian process $\bg_{n}(\omega)=\{g_{n,1}(\omega),\ldots,g_{n,2r}(\omega)\}^{\T}$ with mean zero and covariance function $\bSigma(\omega_1,\omega_2)$ defined as \eqref{eq:long_run}.
\end{proposition}

\begin{remark}\label{rk:1}
	(a) In Proposition \ref{tm:1} and other theoretical results of this paper, we focus on the high-dimensional scenario by assuming $r\geq n^{\kappa}$ for some sufficiently small constant $\kappa>0$. Such an assumption is quite mild in the literature of high-dimensional data analysis and it is not necessary for our theory, which is just used to simplify the presentation. In our theoretical proofs, we need to  compare $\log n$ with  $\log(Kr)$ or $\log r$ in many places. Without such a restriction, the proof of Proposition \ref{tm:1} will become much lengthier and some  $\log(Kr)$ and $\log r$ terms in the theoretical results should be replaced by  $\log(nKr)$ and $\log(nr)$, respectively.  Our proposed Gaussian approximation procedure also works for the scenario with fixed $r$.

 (b) If the bandwidth $l_n\asymp n^{\delta}$ for some constant $0<\delta<1/9$, then Proposition \ref{tm:1}(i) holds provided that $\log(Kr)\ll \min\{n^{\delta}, n^{1/9-\delta}\log^{-8/3}(n)\}$, and Proposition \ref{tm:1}(ii) holds provided that $\log r\ll \min\{n^{\delta}, n^{1/9-\delta}\log^{-8/3}(n)\}$. 
 
(c)  In practice, we can select the bandwidth $l_n$ by adapting the simple rule suggested in Section 2.1 of \cite{Politis_2003}. More specifically, let $l_n=2\hat{m}$, where $\hat{m}$ is the smallest positive integer such that $p^{-2}\sum_{i=1}^p\sum_{j=1}^p|\hat{\rho}_{i,j}(\hat{m}+k)|<2\sqrt{n^{-1}\log n}$ for $k=1,\ldots,5$, and $\hat{\boldsymbol{\rho}}(k)\equiv \{\hat{\rho}_{i,j}(k)\}_{p\times p} = \diag\{\widehat{\bGamma}(0)\}^{-1/2} \widehat{\bGamma}(k) \diag\{\widehat{\bGamma}(0)\}^{-1/2}$ with $\widehat{\bGamma}(k)$ specified in \eqref{eq:hatGammak}. Our numerical results in Section \ref{numerical} verify its good finite-sample performance. 

(d) Since $\{\bx_t\}$ is an $\alpha$-mixing sequence, we know $\{\bc_t\}$ is also an $\alpha$-mixing sequence. Denote by $\tilde{\alpha}_n(k)$ the $\alpha$-mixing coefficient of the sequence  $\{\bc_t\}$. By Condition \ref{as:betamixing},  it holds that 
$
 \tilde{\alpha}_n(k)\leq \alpha_n(|k-2l_n|_+)\leq C_3\exp(-C_4|k-2l_n|_+)$. When $l_n$ diverges with $n$, the $\alpha$-mixing coefficients of the sequence $\{\bc_t\}$  vary with $n$. For $\cJ=\{\omega_1,\ldots,\omega_K\}$, it is essential to establish the Gaussian approximation for $\bP(n^{-1/2}\sum_{t=1}^{n-2l_n}\bH\bc_t\in A)$, where $A$ is a 2-sparsely convex set, and 
 $\bH$ involves the frequency domain content. The most related works are \cite{CJS2021} and \cite{CCW2021}. By comparison, \cite{CJS2021} only considered the Gaussian approximation for   hyperrectangle sets, not for sparsely convex sets, and the theoretical results for sparsely convex sets in \cite{CCW2021} did not allow the $\alpha$-mixing coefficients to vary with $n$. 
	Therefore, the existing theoretical results cannot be applied to establish Proposition \ref{tm:1}(i) when $l_n$ diverges with $n$. For $\cJ=[\omega_L,\omega_U]$, we need to construct the Gaussian approximation for the distribution of the supremum of a stochastic process.
	The most related work is \cite{CCK_2014}. With i.i.d. observations, they developed a new direct approach to approximate the supremum of general empirical processes by a sequence of supremums of Gaussian processes. However, our setting focuses on dependent observations and $\mathcal{T}$ is not the supremum of an empirical process.
	Hence, the results of \cite{CCK_2014} cannot be applied to our setting. 
 \end{remark}

\subsection{Parametric bootstrap procedure}\label{sec:parboot}
To apply Proposition \ref{tm:1} to approximate the distribution of $\sup_{\omega\in\mathcal{J}}T_n(\omega;\mathcal{I})$, we need to propose an estimate of the long-run covariance $\bXi$ given in \eqref{eq:Xi}.
Recall $\widehat{\bGamma}(k)=\{\hat\gamma_{i,j}(k)\}_{p\times p}$ and $\bar\bx=(\bar{x}_{1},\ldots,\bar{x}_p)^{\T}=n^{-1}\sum_{t=1}^{n}\bx_t$.  Let $\hat{\mathring{\bx}}_t=(\hat{\mathring{x}}_{1,t},\ldots,\hat{\mathring{x}}_{p,t})^{\T}= \bx_t-\bar\bx$.
For each $\ell\in[r]$, define a vector
\begin{align}\label{eq:hatcl}
	\hat{\bc}_{\ell,t}=\frac{1}{2\pi}\big\{\hat{\mathring{x}}_{\chi_1(\ell),t}\hat{\mathring{x}}_{\chi_2(\ell),t+l_n} -{\hat\gamma}_{\bchi(\ell)}(-l_n),\ldots, \hat{\mathring{x}}_{\chi_1(\ell),t+2l_n}\hat{\mathring{x}}_{\chi_2(\ell),t+l_n} -{\hat\gamma}_{\bchi(\ell)}(l_n)\big\}^\T  \,,
\end{align}
which provides an approximation to $\bc_{\ell,t}$ defined in \eqref{eq:cjt}.
Write
$
\hat{\bc}_t=(\hat{\bc}_{1,t}^{\T},\ldots,\hat{\bc}_{r,t}^{\T})^{\T}$. 
Based on such defined $\hat{\bc}_t$, we propose a kernel-type estimator suggested by \cite{Andrews_1991} for the long-run covariance matrix $\bXi$ as follows:
\begin{equation}\label{eq:Xihat}
	\widehat{\bXi}=\sum_{q=-\tilde{n}+1}^{\tilde{n}-1}\mathcal{K}\bigg(\frac{q}{b_n}\bigg)\widehat{\bPi}(q)  \,,
\end{equation}
where 
$\widehat{\bPi}(q)=
     \tilde{n}^{-1} \sum_{t=\max(1,-q+1)}^{\min(\tilde{n}, \tilde{n}-q)} \hat{\bc}_{t+q}\hat{\bc}_{t}^\T$, $b_n$ is the bandwidth, and $\mathcal{K}(\cdot)$ is a symmetric kernel function.
When $r$ is fixed, \cite{Andrews_1991} systematically investigated  the  theoretical properties of
such an estimator for the long-run covariance matrix, and shows that the Quadratic Spectral kernel
$	\mathcal{K}_{\rm QS}(u)={25}/(12\pi^2u^2)\{{\sin(6\pi u/5)}/(6\pi u/5)-\cos(6\pi u/5)\} 
$
is the optimal kernel 
in the sense of minimizing the asymptotic truncated mean square error.
In our numerical work, we adopt this Quadratic Spectral kernel with the
data-driven selected bandwidth suggested in Section 6 of \cite{Andrews_1991}, 
i.e. $b_n=1.3221(\hat{a}\tilde{n})^{1/5}$, where $\hat a = \{ \sum_{s=1}^{r(2l_n+1)} 4\hat\rho_s^2 \hat\sigma_s^4 (1-\hat\rho_s)^{-8}\} / \{ \sum_{s=1}^{r(2l_n+1)}\hat\sigma_s^4(1-\hat\rho_s)^{-4}\}$ with $\hat\rho_s$ and $\hat\sigma_s^2$ being, respectively, the estimated autoregressive coefficient and innovation variance from fitting an AR(1) model to time series $\{\hat{c}_{s,t}\}_{t=1}^{\tilde{n}}$, the $s$-th component sequence of $\{\hat\bc_t\}_{t=1}^{\tilde n}$, where $\hat\bc_t=\{\hat{c}_{1,t},\ldots,\hat{c}_{r(2l_n+1),t}\}^{\T}$. 
Although Andrews' method is developed for low-dimensional data, both our theoretical and simulation results show that this estimator  works reasonably well when $r$ is large in relation to $n$. The performance of such kernel-type estimator with different kernels and choices of  bandwidth has been studied in \cite{CJS2021}; their numerical results show that this method is robust for different kernels and bandwidths, and still works even in the high-dimensional case.

\begin{as}\label{as:kernel}
	The symmetric kernel function ${\cal K}(\cdot)$ is continuously differentiable with bounded derivative on $\bR$ satisfying  {\rm (i)} $\mathcal{K}(0)=1$  and {\rm (ii)} $|{\cal K}(x)|\leq C_6|x|^{-\vartheta}$ as $|x|\rightarrow\infty$, for some universal constants $C_6>0$ and $\vartheta>1$.
\end{as}
Condition \ref{as:kernel} is commonly used for   nonparametric estimation of the long-run covariance matrix; see \cite{NW_1987} and \cite{Andrews_1991}. For   kernel functions with bounded support -- such as the Parzen kernel and the Bartlett kernel -- we have $\vartheta=\infty$ in Condition \ref{as:kernel}.
As indicated in \cite{Andrews_1991}, to enforce that $\widehat{\bXi}$ given in \eqref{eq:Xihat} be positive semi-definite we can require the kernel function $\mathcal{K}(\cdot)$ to satisfy $\int_{-\infty}^\infty \mathcal{K}(x)e^{-\iota x\lambda}\,{\rm d}x\geq0$ for any $\lambda\in\mathbb{R}$ with $\iota=\sqrt{-1}$. The Quadratic Spectral, Bartlett, and Parzen kernels   all satisfy this requirement.

To construct the parametric bootstrap procedure, let $(\epsilon_1,\ldots,\epsilon_{\tilde{n}})^{\T} \sim \mathcal{N}(\bzero,\boldsymbol{\Theta})$ be independent of $\mathcal{X}_n=\{\bx_1,\ldots,\bx_n\}$, where $\boldsymbol{\Theta}$ is an  $\tilde{n}\times\tilde{n}$ matrix with $(i,j)$-th element $\mathcal K\{(i-j)/b_n\}$.
Following the same arguments in \cite{CYZ_2017}, conditionally on $\cX_n$, we have
\begin{align*}
	\frac{1}{\sqrt{\tilde{n}}}	\sum_{t=1}^{\tilde{n}}\epsilon_t\hat{\bc}_t \sim \mathcal{N}(\bzero,\widehat{\bXi})
\end{align*}
with $\widehat{\bXi}$ given in \eqref{eq:Xihat}. Hence, conditionally on $\cX_n$,
\begin{align}\label{eq:hatetaext}
	\hat{\bfeta}^{\ext}(\omega):=\{\hat{\eta}^{\ext}_1(\omega),\ldots,\hat{\eta}^{\ext}_{2r}(\omega)\}^\T
	=\{\bI_r\otimes \bA(\omega) \} \bigg(\frac{1}{\sqrt{\tilde{n}}}\sum_{t=1}^{\tilde{n}}\epsilon_t\hat{\bc}_t \bigg)
\end{align}
is a $(2r)$-dimensional Gaussian process with mean zero and covariance function $\{\bI_r\otimes \bA(\omega_1) \} \widehat{\bXi} \{\bI_r\otimes \bA^{\T}(\omega_2) \}$. Letting
\begin{align}
	\xi_{\mathcal{J}} := \sup_{\omega\in\mathcal{J}}\max_{\ell\in [r]} \big\{|\hat{\eta}^{\ext}_{2\ell-1}(\omega)|^2+|\hat{\eta}^{\ext}_{2\ell}(\omega)|^2\big\}\,, \label{eq:xi}
\end{align}
our next result
shows that  the distribution of $\sup_{\omega\in\cJ}T_n(\omega;\cI)$ can be approximated by the  distribution of $\xi_{\cJ}$ conditional on $\cX_n$.

\begin{proposition}\label{tm:2}
	Assume Conditions {\rm \ref{as:moment}}--{\rm \ref{as:kernel}} hold and $r\geq n^\kappa$ for some sufficiently small constant $\kappa>0$. Let the bandwidth $b_n$ in \eqref{eq:Xihat}
	satisfy $b_n\asymp n^{\rho}$    for some   constant $\rho$ satisfying $0<\rho< (\vartheta-1)/(3\vartheta-2)$ with $\vartheta$ specified in Condition {\rm\ref{as:kernel}}.
	For $\xi_{\mathcal{J}}$ defined in \eqref{eq:xi},
	the following two assertions are true as $n\rightarrow\infty$.
	
	{\rm (i)} If $\mathcal{J}=\{\omega_1,\ldots,\omega_K\}$,   then
	$$
		\sup_{u\geq 0}\bigg|\mathbb{P}\bigg\{\sup_{\omega\in\mathcal{J}}{T}_n(\omega;\mathcal{I})\leq u\bigg\} - \mathbb{P}(\xi_{\mathcal{J}}\leq u\,|\,\cX_n)\bigg|=o_{\rm p}(1)$$
	provided that $\log(Kr)\ll f_1(l_n,n;\vartheta,\rho)$,
	with the bandwidth $l_n$ in \eqref{eq:spectralest} satisfying the restriction $\max\{2,C'\log(Kr)\}\leq l_n\ll n^{f_2(\vartheta,\rho)}$ for some sufficiently large constant $C'>0$, where $f_1(l_n,n;\vartheta,\rho)$ and $f_2(\vartheta,\rho)$ are defined by \eqref{eq:f1.f2} and \eqref{eq:f1.f2.2} in the Appendix, respectively.

	{\rm (ii)} If $\mathcal{J}=[\omega_L,\omega_U]$, then
	$$
		\sup_{u\geq 0}\bigg|\mathbb{P}\bigg\{\sup_{\omega\in\mathcal{J}}{T}_n(\omega;\mathcal{I})\leq u\bigg\} - \mathbb{P}(\xi_{\mathcal{J}}\leq u\,|\,\cX_n)\bigg|=o_{\rm p}(1)$$
	provided that $\log r\ll f_1(l_n,n;\vartheta,\rho)$,
	with the bandwidth $l_n$ in \eqref{eq:spectralest} satisfying the restriction $\max ( 2,C'\log r ) \leq l_n\ll n^{f_2(\vartheta,\rho)}$ for some sufficiently large constant $C'>0$.
\end{proposition}

\begin{remark}
	Proposition \ref{tm:2} requires the bandwidth $l_n$ involved in \eqref{eq:spectralest} for the estimation of the high-dimensional spectral density matrix $\bF(\omega)$ to satisfy the restriction $l_n\ll n^{f_2(\vartheta,\rho)}$. Such  a restriction together with $b_n\asymp n^\rho$ is applied to guarantee that the long-run covariance matrix estimate $\widehat{\bXi}$ has a suitable convergence rate to $\bXi$ under the loss $|\cdot|_\infty$. If we select the kernel function $\cK(\cdot)$ involved in \eqref{eq:Xihat} with bounded support such as the Parzen kernel and the Bartlett kernel, then $\vartheta=\infty$ in Condition \ref{as:kernel}, which implies $l_n$ should satisfy the restriction $l_n\ll n^{f_2(\infty,\rho)}$ with $f_2(\infty,\rho)=\min\{\rho/3, (1-3\rho)/2\}$. Furthermore, letting the bandwidth $l_n\asymp n^{\delta}$ for some constant $0<\delta<f_2(\infty,\rho)$ and writing $\tilde\delta=\min\{\delta, \, (\rho-3\delta)/2,\, (1-3\rho-2\delta)/8\}  $, then Proposition \ref{tm:2}(i) holds provided that $\log(Kr)\ll \min\{n^{\tilde{\delta}}, n^{1/9-\delta}\log^{-8/3}(n)\}$, and Proposition \ref{tm:2}(ii) holds provided that $\log r\ll \min\{n^{\tilde{\delta}},n^{1/9-\delta}\log^{-8/3}(n)\}$.
	
	
\end{remark}



	

\vspace{-1em}

\section{Applications} 	\label{applications}

In this section, we present two applications of our established Gaussian approximation theory to   inference for the high-dimensional spectral density matrix, including the global hypothesis testing in Section \ref{subsec:global} and the multiple testing with FDR control in Section \ref{subsec:multiple}.

\subsection{Global hypothesis testing}
\label{subsec:global}

Recall $\bF(\cdot)=\{f_{i,j}(\cdot)\}_{p\times p}$. Given $(\mathcal{I},\mathcal{J})$ such that $\mathcal{I}\subset[p]^2$ with $|\mathcal{I}|=r$ and $\mathcal{J}\subset[-\pi,\pi)$, we consider the following hypothesis testing problem:
\begin{align}\label{eq:H0}
	H_{0}:f_{i,j}(\omega)=0~\textrm{for any}~(i,j)\in\mathcal{I}~\textrm{and}~\omega\in\mathcal{J}~~~~\textrm{versus}~~~~H_{1}:H_{0}~\textrm{is not true}\,.
\end{align}
We propose a test statistic for the hypothesis testing problem \eqref{eq:H0} as follows:
\begin{align}\label{eq:T2}
	T_{n}=\sup_{\omega\in\mathcal{J}}\max_{(i,j)\in\mathcal{I}}\bigg|\sqrt{\frac{n}{l_n}}\hat f_{i,j}(\omega)\bigg|^2 \,.
\end{align}
Notice that $T_n=\sup_{\omega\in\cJ}T_n(\omega;\cI)$ under the null hypothesis $H_0$. For given significance level $\alpha\in(0,1)$, based on the theoretical result of Proposition \ref{tm:2} we define the critical value
\begin{align*}
	\hat{\rm cv}_{\alpha} := \inf\{u>0: \,\bP(\xi_{\cJ}\leq u\,|\, \cX_n)\geq 1-\alpha\} \,.
\end{align*}
Then we reject $H_0$ specified in \eqref{eq:H0} at nominal level $\alpha$ if $T_n>\hat{\rm cv}_{\alpha}$. Practically, we can always draw $\xi_{\cJ,1},\ldots,\xi_{\cJ,B}$ independently by \eqref{eq:xi} for some large integer $B$, and  select the $\lfloor B\alpha\rfloor$-th largest value among them as the critical value $\hat{\rm cv}_{\alpha}$.

Based on the selection of $\cJ$, we define
\begin{equation}\label{eq:MJ}
	M_{\cJ} =\left\{
	\begin{aligned}
		K\,,~~~&\textrm{if}~\cJ=\{\omega_1,\ldots,\omega_K\} \,,\\
		n\,,\,~~~&\textrm{if}~\cJ=[\omega_L,\omega_U]\,.
	\end{aligned}
	\right.
\end{equation}
The next theorem states the theoretical guarantee of our proposed global test.

\begin{theorem}\label{tm:H0}
	Assume Conditions {\rm \ref{as:moment}}--{\rm \ref{as:kernel}} hold and $r\geq n^\kappa$ for some sufficiently small constant $\kappa>0$.
	Let  the bandwidth $b_n$ in \eqref{eq:Xihat} satisfy $b_n\asymp n^{\rho}$    for some   constant $\rho$ satisfying $0<\rho< (\vartheta-1)/(3\vartheta-2)$ with $\vartheta$ specified in Condition {\rm\ref{as:kernel}}, and $\log(M_{\cJ}r)\ll f_1(l_n,n;\vartheta,\rho)$
	with the bandwidth $l_n$ in \eqref{eq:spectralest} satisfying $\max\{2,C'\log(M_{\cJ}r)\}\leq l_n\ll n^{f_2(\vartheta,\rho)}$ for some sufficiently large constant $C'>0$, where $f_1(l_n,n;\vartheta,\rho)$ and $f_2(\vartheta,\rho)$ are defined as \eqref{eq:f1.f2} and \eqref{eq:f1.f2.2} in the Appendix, respectively. As $n\rightarrow\infty$, the following two assertions are true.
	
	{\rm (i)} Under the null hypothesis $H_0$,  then
	$\bP(T_n>\hat{\rm cv}_\alpha)\rightarrow\alpha$.
	
	{\rm (ii)} Write $\lambda(M_{\cJ},r,\alpha)=\{2\log(2M_{\cJ}r)\}^{1/2}+\{2\log(4/\alpha)\}^{1/2}$ for given $(M_{\cJ},r,\alpha)$, and $\varrho=\sup_{\omega\in\mathcal{J}}\max_{\ell\in[2r]}\sigma^2_\ell(\omega)$
	with $\sigma^2_\ell(\omega)$ being the $\ell$-th element in the main diagonal of $\bSigma(\omega,\omega)$ defined in \eqref{eq:long_run}. Under the alternative hypothesis $H_1$, then $\bP(T_n>\hat{\rm cv}_\alpha)\rightarrow 1$ provided that $
	\sup_{\omega\in\mathcal{J}}\max_{(i,j)\in\mathcal{I}} |f_{i,j}(\omega)|
	\geq  2n^{-1/2}l_n^{1/2}\varrho^{1/2}\lambda(M_{\cJ},r,\alpha) (1+\epsilon_n) $
	for some positive $\epsilon_n$ satisfying $\epsilon_n\rightarrow 0$ and $\epsilon_n^2\varrho l_n^{-2}\log^{-2}(l_n)\log^{-1}(n)
	\lambda(M_{\cJ},r,\alpha) \rightarrow\infty$.
\end{theorem}

Theorem \ref{tm:H0}(i) shows that the size of our proposed global test can maintain  the nominal level $\alpha$ asymptotically.
Theorem \ref{tm:H0}(ii) indicates that our proposed global test is consistent under certain local alternatives.

\subsection{Multiple testing with FDR control}
\label{subsec:multiple}
If the global null is rejected, it is important to recover the pairs of indices that correspond to non-zero cross-spectrum (or coherence). That is, we wish to estimate the support of $\mathscr{S}_f$ as in \eqref{eq:sf}. It turns out that  support recovery can be formulated as a simultaneous testing problem of $(p^2-p)/2$
hypotheses $H_{0,i,j}: \max_{\omega\in [-\pi,\pi) } |f_{i,j}(\omega)|=0$ versus $H_{1,i,j}: \max_{\omega\in  [-\pi,\pi) } |f_{i,j}(\omega)|\not=0$,
for $(i,j)\in[p]^{2}$ and $i<j$.  Below we will present a more general version of the above multiple testing problem. This generality is needed in our real data analysis, where the interest is to recover the support set at the state level based on county-level time series. See Section \ref{data} for details.

Given $\{\mathcal{I}^{(q)},\mathcal{J}^{(q)}\}$ with $\mathcal{I}^{(q)}\subset[p]^2$ and $\mathcal{J}^{(q)}\subset[-\pi,\pi)$, we consider $Q$ hypothesis testing problems:
\begin{align*}
	H_{0,q}:f_{i,j}(\omega)=0~\textrm{for any}~ (i,j)\in\mathcal{I}^{(q)}~\textrm{and}~\omega\in\mathcal{J}^{(q)}~~~~\textrm{versus}~~~~H_{1,q}:H_{0,q}~\textrm{is not true}
\end{align*}
for $q\in[Q]$. Similar to \eqref{eq:T2}, we propose the test statistic for $H_{0,q}$ as follows:
\begin{align*}
	T_n^{(q)}=\sup_{\omega\in\mathcal{J}^{(q)}}\max_{(i,j)\in\mathcal{I}^{(q)}}\bigg|\sqrt{\frac{n}{l_n}} \hat{f}_{i,j}(\omega)\bigg|^2\,,
\end{align*}
and	reject $H_{0,q}$ when $T_n^{(q)}$ takes some large values. Let $\mathcal{H}_0=\{q\in [Q]:H_{0,q}~\textrm{is true}\}$ and $\mathcal{H}_1=\mathcal{H}\setminus\mathcal{H}_0$ denote the sets of true nulls and true alternatives, respectively. Write $Q_0=|\mathcal{H}_0|$. For each $q\in[Q]$, let $r_q=|\mathcal{I}^{(q)}|$ and  $\bchi^{(q)}(\cdot)=\{\chi_1^{(q)}(\cdot),\chi_2^{(q)}(\cdot)\}$ be a given bijective mapping from $[r_q]$ to $\mathcal{I}^{(q)}$ such that for any $(i,j)\in\mathcal{I}^{(q)}$, there exists a unique $\ell\in[r_q]$ satisfying $(i,j)=\bchi^{(q)}(\ell)$.
Analogously, we define a $(2l_n+1)$-dimensional vector
\begin{align*} 
	\hat{\bc}_{\ell,t}^{(q)}=\frac{1}{2\pi}\big\{\hat{\mathring{x}}_{\chi_1^{(q)}(\ell),t}\hat{\mathring{x}}_{\chi_2^{(q)}(\ell),t+l_n} -{\hat\gamma}_{\bchi^{(q)}(\ell)}(-l_n),\ldots,\hat{\mathring{x}}_{\chi_1^{(q)}(\ell),t+2l_n}\hat{\mathring{x}}_{\chi_2^{(q)}(\ell),t+l_n} -{\hat\gamma}_{\bchi^{(q)}(\ell)}(l_n)\big\}^\T
\end{align*}
for $\ell\in[r_q]$, where $\hat{\mathring{\bx}}_t=(\hat{\mathring{x}}_{1,t},\ldots,\hat{\mathring{x}}_{p,t})^{\T}= \bx_t-\bar\bx$. 
Let $\hat{\bc}_t^{(q)}=\{\hat{\bc}_{1,t}^{(q),\T},\ldots,\hat{\bc}_{r_q,t}^{(q),\T}\}^\T$ and
\begin{align*} 
	\hat{\bfeta}^{\ext,(q)}(\omega)
	:=\big\{\hat{\eta}^{\ext,(q)}_1(\omega), \ldots,\hat{\eta}^{\ext,(q)}_{2r_q}(\omega)\big\}^\T
	=\big\{\bI_{r_q}\otimes \bA(\omega) \big\} \frac{1}{\sqrt{{\tilde n}}} \sum_{t=1}^{\tilde{n}}\epsilon_t^{(q)}\hat{\bc}_{t}^{(q)} \,,
\end{align*}
where $\{\epsilon_1^{(q)},\ldots,\epsilon_{\tilde{n}}^{(q)}\}^\T\sim \mathcal{N}(\bzero,\bTheta)$ and $\bA(\omega)$ defined as \eqref{eq:A}. Identical to $\hat{\bfeta}^{\ext}(\omega)$ defined in \eqref{eq:hatetaext}, $\hat{\bfeta}^{\ext,(q)}(\omega)$ is a $(2r_q)$-dimensional Gaussian process with mean zero and covariance function $\{\bI_r\otimes \bA(\omega_1)\}\widehat{\bXi}^{(q)}\{\bI_r\otimes \bA^{\T}(\omega_2)\}$, where $\widehat{\bXi}^{(q)}$ is defined in the same manner of \eqref{eq:Xihat} but with replacing $\hat{\bc}_t$ by $\hat{\bc}_t^{(q)}$. Letting
\begin{align*} 
	\xi_{\mathcal{J}^{(q)}}^{(q)} :=  \sup_{\omega\in\mathcal{J}^{(q)}}\max_{\ell\in [r_q]} \big\{|\hat{\eta}^{\ext,(q)}_{2\ell-1}(\omega)|^2+|\hat{\eta}^{\ext,(q)}_{2\ell}(\omega)|^2\big\} \,,
\end{align*}
we can  show -- identical to Proposition \ref{tm:2} -- that
\begin{align}\label{eq:diffTnq}
	\max_{q\in\mathcal{H}_0}\sup_{u\geq 0} \big|\bP\big\{T_n^{(q)}>u\big\} - \bP\big\{\xi_{\mathcal{J}^{(q)}}^{(q)}>u\,|\,\mathcal{X}_n\big\}\big| = o_{\rm p}(1) \,.
\end{align}
Denote by ${\rm pv}^{(q)}=\bP\{\xi_{\mathcal{J}^{(q)}}^{(q)}\geq T_n^{(q)}\,|\,\mathcal{X}_n\}$ and $V_n^{(q)}=\Phi^{-1}\{1-{\rm pv}^{(q)}\}$ the p-value of $H_{0,q}$ and its normal quantile transformation, respectively. For the threshold value $t$ such that $H_{0,q}$ is rejected if $V_n^{(q)}\geq t$, denote the total number of false positives by $R_0(t)=\sum_{q\in\mathcal{H}_0}I\{V_n^{(q)}\geq t\}$, and the total number of rejections by $R(t)=\sum_{q\in\mathcal{H}}I\{V_n^{(q)}\geq t\}$. The false discovery proportion (FDP) and false discovery rate (FDR) are defined, respectively,  as
\begin{align*}
	{\rm FDP}(t)=\frac{R_0(t)}{1\vee R(t)}~~~\textrm{and}~~~{\rm FDR}(t)=\mathbb{E}\{{\rm FDP}(t)\}.
\end{align*}

Given a prescribed level $\alpha\in(0,1)$, the key objective for FDR control is to find the smallest $\hat{t}$ such that ${\rm FDR}(\hat{t})\leq \alpha$. 
To do this, 
we first consider ${\rm FDP}(t)$.
Since the true null hypotheses set $\mathcal{H}_0$ is unknown, we need to estimate $R_0(t)$, i.e., the numerator of ${\rm FDP}(t)$. By \eqref{eq:diffTnq}, it holds that $\bP\{V_n^{(q)}\geq t\}=1-\Phi(t)+o(1)$ for any $q\in\cH_0$. An ideal estimate of ${\rm FDP}(t)$ is
$
\widetilde{\rm FDP}(t) = \{Q_0\{1-\Phi(t)\}/\{1\vee R(t)\}$.
Since $Q_0$ is unknown, $\widetilde{\rm FDP}(t)$ is infeasible in practice and we can only estimate ${\rm FDP}(t)$ via a more conservative way:
\begin{align*}
	\widehat{\rm FDP}(t)=\frac{Q\{1-\Phi(t)\}}{1\vee R(t)} \,.
\end{align*}
For given $\alpha\in(0,1)$, we choose
\begin{align}\label{eq:t.hat}
	\hat{t} = \inf\big\{0< t\leq (2\log  Q-2\log \log Q)^{1/2}: \, \widehat{\rm FDP}(t)\leq \alpha\big\} \,.
\end{align}
If $\hat{t}$ defined in \eqref{eq:t.hat} does not exist, let $\hat{t}=(2\log Q)^{1/2}$. We reject all $H_{0,q}$'s with $V_{n}^{(q)}\geq\hat{t}$.

	To analyze the theoretical properties of our proposed multiple testing procedure, we need to measure the  dependency among the marginal test statistics $\{T_n^{(q)}\}_{q\in[Q]}$. Since the limiting distribution of $T_n^{(q)}$ is not pivotal and does not admit an explicit form (or even does not exist), characterization of  the dependency among $\{T_n^{(q)}\}_{q\in[Q]}$ is nontrivial. To overcome this difficulty, we consider a transformation of the test statistics $\{T_n^{(q)}\}_{q\in[Q]}$, that is,
	$
		\zeta^{(q)}=\Phi^{-1}[F_q\{T_n^{(q)}\}]$,
	where $\Phi(\cdot)$ and $F_q(\cdot)$ are, respectively, the cumulative distribution functions of the standard normal distribution $\mathcal{N}(0,1)$, and $T_n^{(q)}$. Due to $\zeta^{(q)}\sim \mathcal{N}(0,1)$ for each $q\in[Q]$, 
 following \cite{CHKW_2024},
 we can measure the dependency between $T_n^{(q)}$ and $T_n^{(q')}$ by the correlation between $\zeta^{(q)}$ and $\zeta^{(q')}$. It is obvious that the independence between $T_n^{(q)}$ and $T_n^{(q')}$ is equivalent to ${\rm Corr}\{\zeta^{(q)},\zeta^{(q')}\}=0$. For some constant $\gamma>0$ and any $q\in[Q]$, define the set
	\begin{align*}
		\mathcal{S}_q(\gamma) = \big\{q'\in[Q]: q'\neq q, |{\rm Corr}\{\zeta^{(q)},\zeta^{(q')}\}| \geq \log^{-2-\gamma} (Q)\big\}\,.
	\end{align*}
	For given $q\in[Q]$, the other $Q-1$ test statistics $\{T_n^{(q')}\}_{q'\in[Q]\setminus\{q\}}$ can be considered in two scenarios: (i) if $q'\in\cS_q(\gamma)$, the test statistic $T_n^{(q')}$ has relatively strong dependence with $T_n^{(q)}$, and (ii) if $q'\notin\cS_q(\gamma)$, the test statistic $T_n^{(q')}$ has quite weak dependence with $T_n^{(q)}$.
	To construct the theoretical guarantee of the proposed multiple testing procedure, it is common practice to analyze these two scenarios separately with different technical tools. See also \cite{Liu_2013} and \cite{CSZ_2016}. 
	For each $q\in[Q]$, similar to \eqref{eq:MJ}, we define
	\begin{align}\label{eq:MJq}
		M_{\cJ^{(q)}} = \left\{
		\begin{aligned}
			K_q\,,~~~&\textrm{if}~\cJ^{(q)}=\{\omega_1^{(q)},\ldots,\omega_{K_q}^{(q)}\} \,,\\
			n\,,\,~~~&\textrm{if}~\cJ^{(q)}=[\omega_L^{(q)},\omega_U^{(q)}]\,.
		\end{aligned}
		\right.
	\end{align}
	Write $M_{\max}=\max_{q\in\cH_0}M_{\cJ^{(q)}}$, 
	$r_{\max}=\max_{q\in\mathcal{H}_0}r_q$ 
	and $r_{\min}=\min_{q\in\mathcal{H}_0}r_q$.
	Theorem \ref{tm:FDR} provides the theoretical guarantee of our proposed multiple testing procedure.

	\begin{theorem}\label{tm:FDR}
		Assume Conditions {\rm \ref{as:moment}}--{\rm \ref{as:kernel}} hold,  $r_{\min}\geq n^{\kappa}$ for some sufficiently small constant $\kappa>0$, and
		$\max_{1\leq q\neq q'\leq Q}|{\rm Corr}\{\zeta^{(q)}, \zeta^{(q')}\}|\leq r_\zeta$ for some constant $r_\zeta\in(0,1)$; also assume that $\max_{q\in[Q]}|\mathcal{S}_q(\gamma)|=o(Q^{\nu})$ for some constants $\gamma>0$ and $0<\nu<(1-r_\zeta)/(1+r_\zeta)$.
		Let the bandwidth $b_n$ in \eqref{eq:Xihat} satisfy $b_n\asymp n^{\rho}$
		for some  constant $\rho$ satisfying $0<\rho< (\vartheta-1)/(3\vartheta-2)$ with $\vartheta$ specified in Condition {\rm \ref{as:kernel}}.
		If  $Q\ll n^{f_3(\vartheta,\rho)}$ for $f_3(\vartheta,\rho)$ defined as \eqref{eq:f3.f4.f5} in the Appendix, then
		$
			\limsup_{n,Q\rightarrow\infty}{\rm FDR}(\hat t)\leq \alpha Q_0/Q$ and
			$\lim_{n,Q\rightarrow\infty}\bP\{{\rm FDP}(\hat t)\leq \alpha Q_0/Q+\varepsilon\}=1$
		for any $\varepsilon>0$, provided that
		$\log(M_{\max}r_{\max})\ll f_4(l_n,n,Q;\vartheta,\rho)$ with
		the bandwidth $l_n$ in \eqref{eq:spectralest}
		satisfying $\max\{2,C'\log(M_{\max}r_{\max})\}\leq l_n\ll  f_5(n,Q;\vartheta,\rho)$ for some sufficiently large constant $C'>0$, where $f_4(l_n,n,Q;\vartheta,\rho)$ and $f_5(n,Q;\vartheta,\rho)$ are defined as \eqref{eq:f3.f4.f5.2} and \eqref{eq:f3.f4.f5.3} in the Appendix, respectively.
	\end{theorem}

	\section{Numerical simulations}
	\label{numerical}

	In this section, we  investigate the finite sample performance of our proposed methods. The number of parametric bootstrap replications used to determine the critical value $\hat{\rm cv}_\alpha$ in the global testing problem and the p-values in the multiple testing problem is selected as $B=1000$. We set the sample size to be $n\in\{300, 600\}$,  and the dimension to be $p\in\{50, 100,  200\}$, which covers low-, moderate- and high-dimensional scenarios. All reported simulation results in this section are based on 1000 replications.
	
	
	\subsection{Global hypothesis testing}\label{sec:4.1} Let $\mathcal{I}=\{(i,j)\in[p]^2\,:\,i>j\}$ with $r=|\mathcal{I}|=p(p-1)/2$.
	Three types of $\mathcal{J}$ are considered, viz. (i) quarterly seasonal frequencies ($K=4$) such that
	$\mathcal{J}=\mathcal{J}^{(4)}=\{-\pi, -\pi/2, 0, \pi/2\}$; (ii) monthly seasonal frequencies ($K=12$) such that $\mathcal{J}=\mathcal{J}^{(12)}=\{-\pi, -5\pi/6,\ldots,5\pi/6\}$; (iii) Fourier frequencies ($K=n$) such that $\mathcal{J}=\mathcal{J}^{(n)}=\{-\pi,-(n-2)\pi/n,-(n-4)\pi/n,\ldots,(n-2)\pi/n\}$.
	Notice that the effective dimension of the parameter we are testing is $p(p-1)|\mathcal{J}|/2$, which ranges from 4900 (i.e., when $p=50$ and $\mathcal{J}=\mathcal{J}^{(4)}$) to 
 11940000 (corresponding to the case $p=200$ and $\mathcal{J}=\mathcal{J}^{(n)}$ with $n=600$).
 For the flat-top kernel \eqref{eq:Wu} involved in \eqref{eq:spectralest} for the estimate of the spectral density matrix, we set the constant $c\in\{0.5,0.8\}$. 	
 The associated bandwidths $l_n$ and $b_n$ are determined as stated in Section \ref{main}.
 To examine the empirical size, we consider the following models:
	\begin{itemize}[leftmargin=1.62cm]
		\item[Model 1.] Cross-sectionally uncorrelated but dependent sequence: $x_{j,t}=|y_{j-p/10,t}|$ if $j\in\{p/10+1,\ldots,p/5\}$ and $x_{j,t}=y_{j,t}$ otherwise, where  $\by_t=(y_{1,t},\ldots,y_{p,t})^{\T} \overset{{\rm i.i.d.}}{\sim} \mathcal{N}(\bzero,\, a^2 {\bf I}_p)$ with $a\in\{0.2,0.4,0.6\}$. 
		\item[Model 2.] VAR(1) model: $\bx_t=-a \bx_{t-1}+\boldsymbol{\varepsilon}_t$ with $\boldsymbol\varepsilon_t\overset{{\rm i.i.d.}}{\sim} \mathcal{N}\{\bzero,\, (1-a^2){\bf I}_p\}$ and $a\in\{0.05,0.1,0.2\}$.
		\item[Model 3.] VMA(1) model: $\bx_t=\boldsymbol\varepsilon_t-a\boldsymbol\varepsilon_{t-1}$ with $\boldsymbol\varepsilon_t\overset{{\rm i.i.d.}}{\sim}\mathcal{N}(\bzero,{\bf I}_p)$ and $a\in\{0.15,0.2,0.25\}$.
            \item[Model 4.] VARMA(2,2) model: $\bx_t=\bPhi_1\bx_{t-1} +\bPhi_2\bx_{t-2} +\boldsymbol{\varepsilon}_t +\bTheta_1\boldsymbol{\varepsilon}_{t-1} +\bTheta_2\boldsymbol{\varepsilon}_{t-2}$, where $\boldsymbol{\varepsilon}_t=(\varepsilon_{1,t},\ldots,\varepsilon_{p,t})^{\T}$ with $\varepsilon_{j,t}\overset{{\rm i.i.d.}}{\sim}t_5$,    $\bPhi_1=\diag(0.6{\bf 1}_{p/2}^{\T},0.4{\bf 1}_{p/2}^{\T})$,
    $\bPhi_2=0.15 \bI_p$,  
    $\bTheta_1= -a\cdot \diag(0.5{\bf 1}_{p/2}^{\T},0.25{\bf 1}_{p/2}^{\T})$ with $a\in\{0.2,0.25,0.3\}$, and
    $\bTheta_2 = -0.05\bI_p$. 
	\end{itemize}	
 Due to the lack of competing methods, we only focus on the examination of the performance of our test. As seen from Table \ref{tab:H0}, our proposed test has relatively accurate sizes when the dimension $p$ is low for Models 1--4. 
 When the sample size $n$ is fixed, the empirical sizes  tend to decrease as the dimension $p$ increases, which shows the impact on the parametric bootstrap-based approximation from the dimension $p$. When the  dimension $p$ is fixed,  the empirical sizes are closer to the nominal level as the sample size  increases from $n=300$ to $n=600$.  For most settings, the size is below the nominal level and our test is conservative. In settings where our test is over-sized, the amount of over-rejection appears quite mild.  So overall the Type-I error is well controlled. 
 
 Based on the above size results, we can see that the two choices for $c$ (i.e., $c=0.5, 0.8$) deliver very similar results, and also our test seems insensitive to  $K$, since setting $K=4, 12, n$ does not have much impact on the rejection rates.

	\begin{table}[htbp]
		\scriptsize
		\centering
		\caption{Empirical sizes of the proposed global tests for Models 1--4 at the $5\%$ nominal level based on 1000 repetitions. All numbers reported below are multiplied by 100.}
		\vspace{0.7em}
		\resizebox{!}{5cm}{
    \begin{tabular}{ccc|cccc|cccc|cccc|cccc}
          &       &       & \multicolumn{4}{c|}{Model 1}  & \multicolumn{4}{c|}{Model 2}  & \multicolumn{4}{c|}{Model 3}  & \multicolumn{4}{c}{Model 4} \\[0.4em]
    $n$     & $p$     & $c$     & $a$     & $K=4$   & $K=12$  & $K=n$   & $a$     & $K=4$   & $K=12$  & $K=n$   & $a$     & $K=4$   & $K=12$  & $K=n$   & $a$     & $K=4$   & $K=12$  & $K=n$ \\[0.4em]
    300   & 50    & 0.5   & 0.20   & 2.8   & 2.8   & 2.8   & 0.05 & 3.5   & 3.6   & 3.6   & 0.15  & 3.9   & 3.9   & 4.0   & 0.20   & 4.3   & 4.8   & 4.2  \\
          &       &       & 0.40   & 2.9   & 2.9   & 2.7   & 0.10 & 4.1   & 4.1   & 3.9   & 0.20  & 3.8   & 3.8   & 3.7   & 0.25  & 4.6   & 5.1   & 5.0  \\
          &       &       & 0.60   & 3.0   & 3.0   & 2.9   & 0.20 & 4.7   & 4.5   & 4.2   & 0.25  & 3.5  & 3.8   & 4.0   & 0.30   & 5.5   & 5.9   & 5.6  \\[0.2em]
          &       & 0.8   & 0.20  & 2.9   & 2.9   & 2.9   & 0.05 & 3.3   & 3.8   & 3.3   & 0.15  & 4.0   & 3.9   & 3.8   & 0.20   & 4.9   & 4.7   & 4.3  \\
          &       &       & 0.40  & 2.7   & 2.8   & 3.0   & 0.10 & 3.8   & 3.9   & 4.0   & 0.20  & 3.9   & 3.8   & 3.8   & 0.25  & 4.8   & 4.8   & 5.0  \\
          &       &       & 0.60  & 2.7   & 3.0   & 2.9   & 0.20 & 4.5   & 4.7   & 4.5   & 0.25  & 3.6   & 3.4   & 3.7   & 0.30   & 5.7   & 5.5   & 6.0  \\[0.2em]
          & 100   & 0.5   & 0.20  & 1.7   & 1.6   & 1.5   & 0.05 & 1.8   & 2.1   & 1.9   & 0.15  & 2.8   & 2.8   & 2.9   & 0.20   & 1.9   & 1.9   & 1.7  \\
          &       &       & 0.40   & 1.9   & 1.8   & 1.9   & 0.10 & 2.7   & 3.0  & 3.2   & 0.20  & 2.9   & 2.9   & 3.1   & 0.25  & 2.3   & 2.4   & 1.9  \\
          &       &       & 0.60   & 1.7   & 1.9   & 1.9   & 0.20 & 3.9   & 3.2   & 3.8   & 0.25  & 2.6   & 2.7   & 2.8   & 0.30   & 2.3   & 2.7   & 2.8  \\[0.2em]
          &       & 0.8   & 0.20   & 1.8   & 1.7   & 1.6   & 0.05 & 1.9   & 1.9   & 1.9   & 0.15  & 2.7   & 2.8   &  2.8  & 0.20   & 1.9   & 2.2   & 1.9  \\
          &       &       & 0.40   & 1.7   & 1.9   & 2.1   & 0.10 & 3.4   & 3.2   & 3.4   & 0.20  & 2.7   & 2.7   &  3.1  & 0.25  & 2.0   & 2.4   & 2.3  \\
          &       &       & 0.60   & 1.6   & 2.1   & 2.2   & 0.20 & 3.5   & 3.7   & 3.6   & 0.25  & 3.0   & 3.0   &  2.5  & 0.30   & 2.9   & 2.7   & 2.8  \\[0.2em]
          & 200   & 0.5   & 0.20   & 1.3   & 1.4   & 1.4   & 0.05 & 1.5   & 1.9   & 1.5   & 0.15  & 2.0   & 1.9   &  1.8    & 0.20   & 1.3   & 1.3   & 1.2  \\
          &       &       & 0.40   & 1.1   & 1.2   & 1.4   & 0.10 & 2.1   & 1.9   & 2.0   & 0.20  & 1.9   & 1.8   &  2.1    & 0.25  & 1.5   & 1.2   & 1.5  \\
          &       &       & 0.60   & 1.4   & 1.1   & 1.2   & 0.20 & 2.3   & 2.4   & 2.3   & 0.25  & 1.8   & 1.6   &  1.8    & 0.30   & 1.6   & 1.6   & 1.6  \\[0.2em]
          &       & 0.8   & 0.20   & 1.3   & 1.5   & 1.5   & 0.05 & 1.5   & 1.7   & 1.6   & 0.15  & 1.8   & 1.7   & 1.8   & 0.20   & 1.1   & 1.2   & 1.5  \\
          &       &       & 0.40   & 1.4   & 1.4   & 1.5   & 0.10 & 2.1   & 2.0   & 1.9   & 0.20  & 2.2   & 1.8   & 1.9   & 0.25  & 1.4   & 1.4   & 1.5  \\
          &       &       & 0.60   & 1.3   & 1.1   & 1.2   & 0.20 & 2.7   & 2.3   &  2.5  & 0.25  & 1.8   & 2.0   & 1.9   & 0.30   & 1.5   & 1.6   & 1.6  \\[0.2em]
    600   & 50    & 0.5   & 0.20   & 2.7   & 3.0   & 3.0   & 0.05 & 3.8   & 3.7   & 4.1   & 0.15  & 4.3   & 4.3   & 4.4   & 0.20   & 5.9   & 5.5   & 5.5  \\
          &       &       & 0.40   & 2.9   & 3.1   & 2.9   & 0.10 & 4.4   & 4.6   &  4.2  & 0.20  & 4.4   & 4.6   & 4.6   & 0.25  & 6.6   & 6.2   & 6.2  \\
          &       &       & 0.60   & 2.9   & 3.1   & 3.0   & 0.20 & 4.6   & 5.3   & 5.2   & 0.25  & 3.8   & 4.2   & 4.5   & 0.30   & 7.1   & 7.2   & 7.2  \\[0.2em]
          &       & 0.8   & 0.20   & 3.0   & 3.1   & 2.9   & 0.05 & 3.5   & 3.8   & 3.6   & 0.15  & 4.4   & 4.5   & 4.4   & 0.20   & 5.6   & 5.8   & 5.6  \\
          &       &       & 0.40   & 3.1   & 2.9   & 3.1   & 0.10 & 4.3   & 4.4   & 4.3   & 0.20  & 4.3   & 4.7   & 4.5   & 0.25  & 6.2   & 6.3   & 6.6  \\
          &       &       & 0.60   & 2.7   & 3.2   & 2.9   & 0.20 & 5.1   & 4.9   & 5.2   & 0.25  & 4.3   & 3.9   & 4.4   & 0.30   & 7.2   & 7.0   & 6.8  \\[0.2em]
          & 100   & 0.5   & 0.20   & 2.8   & 2.6   & 2.6   & 0.05 & 3.0   & 2.9   & 3.0   & 0.15  & 4.7   & 4.9   & 4.3   & 0.20   & 3.1   & 3.1   & 3.3  \\
          &       &       & 0.40   & 2.9   & 2.8   & 2.6   & 0.10 & 4.8   & 4.3   & 4.4   & 0.20  & 4.7   & 4.6   & 4.5   & 0.25  & 3.9   & 4.2   & 3.8  \\
          &       &       & 0.60   & 2.6   & 2.5   & 2.7   & 0.20 & 5.5   & 5.7   & 5.5   & 0.25  & 4.5   & 4.4   & 4.2   & 0.30   & 4.6   & 4.7   & 5.1  \\[0.2em]
          &       & 0.8   & 0.20   & 2.9   & 2.7   & 2.8   & 0.05 & 3.2   & 2.9   & 3.3   & 0.15  & 4.7  & 4.5   & 4.5   & 0.20   & 3.4   & 3.5   & 3.2  \\
          &       &       & 0.40   & 2.8   & 2.7   & 3.0   & 0.10 & 4.4   & 4.3   & 4.5   & 0.20  & 4.5   & 4.6   & 4.6   & 0.25  & 3.8   & 3.9   & 3.9  \\
          &       &       & 0.60   & 2.6   & 2.6   & 2.5   & 0.20 & 5.5   & 5.7   & 5.6   & 0.25  & 4.8   & 4.0   & 4.2   & 0.30   & 5.0   & 4.8   & 4.8  \\[0.2em]
          & 200   & 0.5   & 0.20   & 2.8   & 2.7   & 2.5   & 0.05 & 2.9   & 2.9   & 3.1   & 0.15  & 3.7   & 4.2   &  3.6  & 0.20   & 2.2   & 1.9   & 2.1 \\
          &       &       & 0.40   & 2.4   & 2.8   & 2.5   & 0.10 & 3.8   & 3.7   & 3.6   & 0.20  & 3.3   & 4.0   &  3.7  & 0.25  & 2.5   & 2.0   &  2.2 \\
          &       &       & 0.60   & 2.4   & 2.6   & 2.7   & 0.20 & 4.9   & 4.6   & 4.5   & 0.25  & 3.9   & 3.6   & 3.0   & 0.30   & 2.9   & 2.8   &  2.7 \\[0.2em]
          &       & 0.8   & 0.20   & 2.7   & 2.4   & 2.8   & 0.05 & 3.0   & 2.6   & 3.0   & 0.15  & 3.8   & 3.5   & 4.2   & 0.20   & 1.9   & 1.8   & 1.9 \\
          &       &       & 0.40   & 2.5   & 2.8   &   3.1 & 0.10 & 3.7   & 3.6   & 3.9   & 0.20  & 3.6   & 3.5   & 4.1   & 0.25  & 2.4   & 2.3   & 2.7 \\
          &       &       & 0.60   & 2.7   & 2.9   &  3.2  & 0.20 & 4.6   & 4.6   & 4.2   & 0.25  & 3.4   & 3.5   & 3.8   & 0.30   & 2.9   & 2.7   & 2.8 \\
    \end{tabular}%
			\label{tab:H0}%
		}
	\end{table}%

	To study the empirical power of the proposed method, we consider the following models:
 \vspace{-2em}
	\begin{itemize}[leftmargin=1.62cm]
		\item[Model 5.] 
 $\bx_t=\boldsymbol\Psi\boldsymbol{\varepsilon}_t$, where $\boldsymbol{\varepsilon}_t\overset{{\rm i.i.d.}}{\sim} \mathcal{N}(\bzero,\,  {\bf I}_p)$ and $\boldsymbol\Psi=(\psi_{k,l})_{p\times p}$ for $\psi_{k,l}=0.4 I(k=l)+aI(|k-l|=1)$ with $a\in\{0.05,0.1,0.15\}$.
		\item[Model 6.]  $\bx_t=\boldsymbol\Psi \bx_{t-1}+\boldsymbol{\varepsilon}_t$, where  $\boldsymbol\varepsilon_t\overset{{\rm i.i.d.}}{\sim} \mathcal{N}\{\bzero,\, (1-0.1^2){\bf I}_p\}$ and $\boldsymbol\Psi=(\psi_{k,l})_{p\times p}$ for $\psi_{k,l}=-0.1 I(k=l)+aI(|k-l|=1)$ with  $a\in\{0.15,0.2,0.25\}$.
		\item[Model 7.]  $\bx_t=\boldsymbol\varepsilon_t-\boldsymbol\Psi\boldsymbol\varepsilon_{t-1}$,  where $\boldsymbol\varepsilon_t\overset{{\rm i.i.d.}}{\sim}\mathcal{N}(\bzero,{\bf I}_p)$  and $\boldsymbol\Psi=(\psi_{k,l})_{p\times p}$ for $\psi_{k,l}=0.2 I(k=l)+aI(|k-l|=1)$ with $a\in\{0.15,0.2,0.25\}$.
  		\item[Model 8.]  $\bx_t=\bPhi_1\bx_{t-1} +\bPhi_2\bx_{t-2} +\boldsymbol{\varepsilon}_t +\bTheta_1\boldsymbol{\varepsilon}_{t-1} +\bTheta_2\boldsymbol{\varepsilon}_{t-2}$, where $\boldsymbol{\varepsilon}_t=(\varepsilon_{1,t},\ldots,\varepsilon_{p,t})^{\T}$ with $\varepsilon_{j,t}\overset{{\rm i.i.d.}}{\sim}t_5$, $\bPhi_1=(\phi_{1,k,l})_{p\times p}$ with $\diag(\bPhi_1)=(0.6{\bf 1}_{p/2}^{\T},0.4{\bf 1}_{p/2}^{\T})^{\T}$ and $\phi_{1,k,l}=a I(k-l=1)$ for $k\neq l$ and $a\in\{0.15,0.2,0.25\}$,
    $\bPhi_2=0.15 \bI_p$, $\bTheta_1= -0.25\cdot  \diag(0.5{\bf 1}_{p/2}^{\T},0.25{\bf 1}_{p/2}^{\T})$, and $\bTheta_2 = -0.05\bI_p$.
	\end{itemize}
	
	Table \ref{tab:H1} shows that the power also appears insensitive to the choices of $c$ and $K$.  As the distance from the null hypothesis increases (e.g., $a$ increases), the empirical power of our proposed test grows rapidly to 1. Moreover, enlarging the sample size $n$ helps to increase the empirical power.
	Overall, the power performance of our proposed test is consistent with our theory under the alternative.

 \begin{table}[htbp]
		\scriptsize
		\centering
		\caption{Empirical powers of the proposed global tests for Models 5--8 at the $5\%$ nominal level based on 1000 repetitions. All numbers reported below are multiplied by 100.}
		\resizebox{!}{5cm}{
    \begin{tabular}{ccc|cccc|cccc|cccc|cccc}
          &       &       & \multicolumn{4}{c|}{Model 5}  & \multicolumn{4}{c|}{Model 6}  & \multicolumn{4}{c|}{Model 7}  & \multicolumn{4}{c}{Model 8} \\[0.4em]
    $n$     & $p$     & $c$     & \multicolumn{1}{c}{$a$} & $K=4$ & $K=12$ & $K=n$ & $a$    & $K=4$   & $K=12$  & $K=n$   & $a$    & $K=4$   & $K=12$  & $K=n$   & $a$    & $K=4$   & $K=12$  & $K=n$ \\[0.4em]
    300   & 50    & 0.5   & 0.05  &  88.3 &  87.3 &  87.9 & 0.15  & 99.6  & 99.5  & 99.4  & 0.15  & 94.4  & 94.3  & 99.4  & 0.15  & 74.9  & 73.2  & 73.6 \\
          &       &       & 0.10  & 100   &  100  &  100  & 0.20  & 100   & 100   & 100   & 0.20  & 100   & 100   & 100   & 0.20  & 96.6  & 95.7  & 96.1 \\
          &       &       & 0.15  & 100   &  100  &  100  & 0.25  & 100   & 100   & 100   & 0.25  & 100   & 100   & 100   & 0.25  & 98.1  & 97.9  & 98.3 \\[0.2em]
          &       & 0.8   & 0.05  &  87.7 & 88.1  & 80.6  & 0.15  & 99.4  & 99.3  & 99.5  & 0.15  & 93.9  & 94.1  &  94.1 & 0.15  & 74.4  & 72.9  & 74.0 \\
          &       &       & 0.10  & 100   &  100  & 100   & 0.20  & 100   & 100   & 100   & 0.20  & 100   & 100   & 100   & 0.20  & 96.0    & 95.9  & 95.9 \\
          &       &       & 0.15  & 100   & 100   & 100   & 0.25  & 100   & 100   & 100   & 0.25  & 100   & 100   & 100   & 0.25  & 97.8  & 98.2  & 98.2 \\[0.2em]
          & 100   & 0.5   & 0.05  & 80.9  & 81.0  & 81.1  & 0.15  & 99.4  & 99.3  & 99.6  & 0.15  & 91.9  & 91.9  & 91.9  & 0.15  & 64.8  & 63.4  & 63.7 \\
          &       &       & 0.10  & 100   &  100  & 100   & 0.20  & 100   & 100   &  100  & 0.20  & 100   & 100   & 100   & 0.20  & 96.0    & 95.5  & 95.6  \\
          &       &       & 0.15  & 100   &  100  & 100   & 0.25  & 100   & 100   & 100   & 0.25  & 100   & 100   & 100   & 0.25  & 97.0    & 97.2  & 97.4  \\[0.2em]
          &       & 0.8   & 0.05  & 80.2  & 80.1  & 80.6  & 0.15  & 99.2  & 99.3  & 99.4  & 0.15  & 92.1  & 91.9  &  91.9 & 0.15  & 64.2  & 63.7  & 65.0 \\
          &       &       & 0.10  & 100   &  100  & 100   & 0.20  & 100   & 100   & 100   & 0.20  & 100   & 100   & 100   & 0.20  & 96.0  & 95.8  & 96.1 \\
          &       &       & 0.15  & 100   &  100  & 100   & 0.25  & 100   & 100   & 100   & 0.25  & 100   & 100   & 100   & 0.25  & 96.9  & 97.0    & 97.0 \\[0.2em]
          & 200   & 0.5   & 0.05  & 74.1  & 74.4  & 75.1  & 0.15  & 99.1  & 98.9  & 99.0  & 0.15  & 86.3  & 86.5  & 85.9  & 0.15  & 46.2  & 45.3  &  46.6 \\
          &       &       & 0.10  & 100   & 100   & 100   & 0.20  & 100   & 100   & 100   & 0.20  & 100   & 100   &   100 & 0.20  & 91.3  & 91.6  &  91.5 \\
          &       &       & 0.15  & 100   & 100   & 100   & 0.25  & 100   & 100   & 100   & 0.25  & 100   & 100   &  100   & 0.25  & 96.2  & 95.8  & 95.3 \\[0.2em]
          &       & 0.8   & 0.05  & 73.8  & 73.8  & 74.6  & 0.15  & 99.1  & 98.8  & 99.0  & 0.15  & 86.6  & 86.1  &  86.3   & 0.15  & 45.5  & 46.1  & 46.4 \\
          &       &       & 0.10  & 100   & 100   & 100   & 0.20  & 100   & 100   & 100   & 0.20  & 100   & 100   & 100  & 0.20  & 91.7  & 91.7  & 92.2 \\
          &       &       & 0.15  & 100   &  100  & 100   & 0.25  & 100   & 100   & 100   & 0.25  & 100   & 100   &  100  & 0.25  & 95.8  & 95.9  & 96.4 \\
    600   & 50    & 0.5   & 0.05  &  100  & 99.9  & 100   & 0.15  & 100   & 100   & 100   & 0.15  & 100   & 100   & 100   & 0.15  & 99.3  & 99.3  & 99.5 \\
          &       &       & 0.10  & 100   & 100   &   100 & 0.20  & 100   & 100   & 100   & 0.20  & 100   & 100   & 100   & 0.20  & 99.8  & 99.8  & 99.8 \\
          &       &       & 0.15  & 100   & 100   &   100 & 0.25  & 100   & 100   & 100   & 0.25  & 100   & 100   & 100   & 0.25  & 100   & 100   & 100 \\[0.2em]
          &       & 0.8   & 0.05  & 100   & 100   & 100   & 0.15  & 100   & 100   & 100   & 0.15  & 100   & 100   &  100  & 0.15  & 99.3  & 99.6  & 99.4 \\
          &       &       & 0.10  & 100   & 100   & 100   & 0.20  & 100   & 100   & 100   & 0.20  & 100   & 100   & 100   & 0.20  & 99.7  & 99.8  & 99.8 \\
          &       &       & 0.15  & 100   & 100   & 100   & 0.25  & 100   & 100   & 100   & 0.25  & 100   & 100   & 100   & 0.25  & 100   & 100   &  100 \\[0.2em]
          & 100   & 0.5   & 0.05  & 100   & 100   & 100   & 0.15  & 100   & 100   & 100   & 0.15  & 100   & 100   & 100   & 0.15  & 98.5  & 98.4  & 98.1 \\
          &       &       & 0.10  & 100   &  100  & 100      & 0.20  & 100   & 100   & 100   & 0.20  & 100   & 100   & 100   & 0.20  & 99.4  & 99.5  & 99.5 \\
          &       &       & 0.15  & 100   &  100  & 100   & 0.25  & 100   & 100   & 100   & 0.25  & 100   & 100   & 100   & 0.25  & 99.8  & 99.7  & 99.9 \\[0.2em]
          &       & 0.8   & 0.05  & 100   &  100  & 100   & 0.15  & 100   & 100   & 100   & 0.15  & 100   & 100   & 100   & 0.15  & 98.4  & 98.6  & 98.3 \\
          &       &       & 0.10  & 100   &  100  & 100   & 0.20  & 100   & 100   & 100   & 0.20  & 100   & 100   & 100   & 0.20  & 99.6  & 99.4  & 99.5 \\
          &       &       & 0.15  & 100   &  100  & 100   & 0.25  & 100   & 100   & 100   & 0.25  & 100   & 100   & 100   & 0.25  & 99.7  & 99.9  & 99.9 \\[0.2em]
          & 200   & 0.5   & 0.05  & 100   &  100  & 100   & 0.15  & 100   & 100   & 100   & 0.15  & 100   & 100   & 100   & 0.15  & 97.3  & 97.1  & 97.2 \\
          &       &       & 0.10  & 100   &  100  & 100   & 0.20  & 100   & 100   & 100   & 0.20  & 100   & 100   & 100   & 0.20  & 99.6  & 99.5  & 99.5 \\
          &       &       & 0.15  & 100   &  100  & 100   & 0.25  & 100   & 100   & 100   & 0.25  & 100   & 100   & 100   & 0.25  & 98.6  & 98.5  & 98.7  \\[0.2em]
          &       & 0.8   & 0.05  & 100   & 100   & 100   & 0.15  & 100   & 100   & 100   & 0.15  & 100   & 100   & 100   & 0.15  & 97.2  & 97.3  & 97.2 \\
          &       &       & 0.10  & 100   & 100   & 100   & 0.20  & 100   & 100   & 100   & 0.20  & 100   & 100   & 100   & 0.20  & 99.4  & 99.5  & 99.4 \\
          &       &       & 0.15  & 100   & 100   & 100   & 0.25  & 100   & 100   & 100   & 0.25  & 100   & 100   & 100   & 0.25  & 98.7  & 98.7  & 98.8 \\
    \end{tabular}%
			\label{tab:H1}%
		}
	\end{table}%

	\subsection{Multiple testing}
	\label{sec:mult-test}
	
	We divide the $p\times p$ matrix ${\bf F}(\omega)$ into a $10\times 10$ block submatrix, i.e., $\bF(\omega)=\{\bF_{i,j}(\omega)\}_{i,j\in[10]}$ where each $\bF_{i,j}(\omega)$ is a submatrix of size $(0.1p)\times (0.1p)$.  The groups of $\{\cJ_k\}_{k=1}^{K}$ are set as follows: (i) when $K=4$, let $\cJ_k=\{-\pi+(k-1)\pi/2\}$ be a set with single element for each $k\in[4]$; (ii) when $K=12$, let $\cJ_k=\{-\pi+(k-1)\pi/6\}$ be a set with single element for each $k\in[12]$; (iii) when $K=n$, let $\cJ_k=\{ -\pi+2(k-1)\pi/n\}$ be a set with single element for each $k\in[n]$.
 For each given $(i,j)\in[10]^2$ with $i\leq j$,  we consider $K$ marginal hypothesis testing problems
	$
	H_{0,i,j,k}: \bF_{i,j}(\omega)={\bf 0}$ for any $\omega\in\cJ_k$ versus $H_{1,i,j,k}:H_{0,i,j,k}$ is not true,
	for $k\in[K]$. Hence, the total number of the marginal hypothesis testing is $Q=55K$.
	Due to the fact that ${\bf F}(\omega)=(2\pi)^{-1}\boldsymbol\Psi^2$ and $(2\pi)^{-1}({\bf I}_p+\boldsymbol\Psi^2-2\boldsymbol\Psi\cos\omega)$, respectively, in Models 5 and 7, then $Q_0=|\cH_0|=36K$ in these two model settings.
	For Model 6,  since  $\bx_t=\sum_{j=0}^{\infty}\boldsymbol\Psi^j\boldsymbol{\varepsilon}_{t-j}$ with $\boldsymbol\varepsilon_t\overset{{\rm i.i.d.}}{\sim} \mathcal{N}\{\bzero,\, (1-0.1^2){\bf I}_p\}$, by Example 11.8.1 in \cite{BD_1991}, we have $\bF(\omega)
=0.495\pi^{-1}({\bf I}_p - \bPsi e^{-\iota\omega})^{-1}\{({\bf I}_p - \bPsi e^{\iota\omega})^{-1}\}^{\T}=0.495\pi^{-1}(\sum_{j=0}^\infty\bPsi^{j}e^{-\iota j\omega})(\sum_{k=0}^{\infty}\bPsi^k e^{\iota k\omega})^{\T}$. For Model 8, letting $\bPhi(x)={\bf I}_p - \bPhi_1 x-\bPhi_2 x^2$ and $\bTheta(x)={\bf I}_p+\bTheta_1 x+\bTheta_2 x^2$, again by Example 11.8.1 in \cite{BD_1991} we have $\bF(\omega) = (2\pi)^{-1}\sigma_\varepsilon^2\bPhi^{-1}(e^{-\iota\omega})\bTheta(e^{-\iota\omega}) \bTheta^{\T}(e^{\iota\omega})\{\bPhi^{-1}(e^{\iota\omega})\}^{\T} $, where $\sigma_\varepsilon^2$ denotes the variance of $\varepsilon_{j,t}$. Note that $\bPsi$ and $\bPhi_1$ are banded matrices. We know all the sub-nulls in Models 6 and 8 are false and thus  $Q_0=0$.} Hence, we only consider the performance of our multiple testing procedure in Model 5 with $a\in\{0.05,0.1,0.15\}$ and Model 7 with $a\in\{0.2,0.4,0.6\}$. As with the global testing procedure, we also use the flat-top kernel and the Quadratic Spectral kernel, respectively, in \eqref{eq:spectralest} and \eqref{eq:Xihat} with the associated bandwidths $l_n$ and $b_n$ determined in the same manner as those in Section \ref{sec:4.1}.

	\begin{table}[t!]
 \vspace{-0.5em}
		\scriptsize
		\centering
		\caption{Empirical FDRs and powers of the proposed multiple testing procedure based on 1000 repetitions. All numbers reported below are multiplied by 100.}
		\resizebox{!}{5cm}{
    \begin{tabular}{rcc|ccccccc|ccccccc}
          &       &       &       & \multicolumn{6}{c|}{Model 5}                  &       & \multicolumn{6}{c}{Model 7} \\[0.2em]
          &       &       &       & \multicolumn{2}{c}{$K=4$} & \multicolumn{2}{c}{$K=12$} & \multicolumn{2}{c|}{$K=n$} &       & \multicolumn{2}{c}{$K=4$} & \multicolumn{2}{c}{$K=12$} & \multicolumn{2}{c}{$K=n$} \\[0.2em]
          & $n$     & $p$     & $a$    & $c=0.5$ & $c=0.8$ & $c=0.5$ & $c=0.8$ & $c=0.5$ & $c=0.8$ & $a$    & $c=0.5$ & $c=0.8$ & $c=0.5$ & $c=0.8$ & $c=0.5$ & $c=0.8$ \\[0.4em]
    \multicolumn{1}{l}{Empirical FDRs} & 300   & 50    & 0.05  & 1.8   & 1.8   & 2.0   & 1.8   & 1.9   & 1.9   & 0.20   & 1.7   & 1.8   & 1.9   & 2.0   & 1.9   & 2.1  \\
          &       &       & 0.10   & 2.1   & 2.1   & 2.3   & 2.2   & 2.2   & 2.2   & 0.40   & 1.8   & 1.8   & 2.0   & 2.1   & 2.2   & 2.1  \\
          &       &       & 0.15  & 2.1   & 2.1   & 2.2   & 2.2   & 2.2   & 2.2   & 0.60   & 1.9   & 2.0   & 2.0   & 2.2   & 2.2   & 2.1  \\[0.2em]
          &       & 100   & 0.05  & 1.5   & 1.4   & 1.3   & 1.4   & 1.3   & 1.3   & 0.20   & 1.4   & 1.4   & 1.5   & 1.5   & 1.5   & 1.6  \\
          &       &       & 0.10   & 1.7   & 1.7   & 1.6   & 1.7   & 1.7   & 1.6   & 0.40   & 1.4   & 1.6   & 1.6   & 1.6   & 1.6   & 1.6  \\
          &       &       & 0.15  & 1.8   & 1.8   & 1.7   & 1.8   & 1.7   & 1.7   & 0.60   & 1.7   & 1.6   & 1.7   & 1.8   & 1.7   & 1.8  \\[0.2em]
          &       & 200   & 0.05  & 1.2   & 1.0   & 1.1   & 1.0   & 1.1   & 1.2   & 0.20   & 1.0   & 0.9   & 1.0   & 1.0   & 1.2   & 0.9  \\
          &       &       & 0.10   & 1.6   & 1.4   & 1.4   & 1.3   & 1.5   & 1.4   & 0.40   & 1.3   & 1.2   & 1.3   & 1.3   & 1.3   & 1.2  \\
          &       &       & 0.15  & 1.4   & 1.4   & 1.4   & 1.3   & 1.3   & 1.4   & 0.60   & 1.2   & 1.3   & 1.4   & 1.4   & 1.3   & 1.4  \\[0.2em]
          & 600   & 50    & 0.05  & 2.3   & 2.6   & 2.4   & 2.4   & 2.5   & 2.4   & 0.20   & 2.7   & 2.7   & 2.8   & 2.7   & 2.7   & 2.7  \\
          &       &       & 0.10   & 2.8   & 2.6   & 2.6   & 2.6   & 2.6   & 2.6   & 0.40   & 2.6   & 2.7   & 2.8   & 2.9   & 2.9   & 2.8  \\
          &       &       & 0.15  & 2.9   & 2.7   & 2.8   & 2.8   & 2.8   & 2.7   & 0.60   & 2.7   & 2.7   & 2.8   & 2.8   & 2.7   & 2.9  \\[0.2em]
          &       & 100   & 0.05  & 2.5   & 2.4   & 2.3   & 2.3   & 2.5   & 2.4   & 0.20   & 2.2   & 2.4   & 2.4   & 2.4   & 2.5   & 2.3  \\
          &       &       & 0.10   & 2.4   & 2.4   & 2.4   & 2.5   & 2.6   & 2.3   & 0.40   & 2.4   & 2.3   & 2.4   & 2.4   & 2.4   & 2.3  \\
          &       &       & 0.15  & 2.5   & 2.4   & 2.3   & 2.3   & 2.4   & 2.5   & 0.60   & 2.4   & 2.5   & 2.6   & 2.6   & 2.6   & 2.7  \\[0.2em]
          &       & 200   & 0.05  & 2.0   & 2.1   & 2.1   & 2.2   & 2.1   & 2.0   & 0.20   & 1.8   & 2.0   & 2.0   & 1.9   & 2.0   & 2.1  \\
          &       &       & 0.10   & 2.1   & 2.1   & 2.2   & 2.1   & 2.0   & 2.1   & 0.40   & 1.9   & 1.9   & 1.9   & 1.9   & 1.8   & 1.9  \\
          &       &       & 0.15  & 2.0   & 2.0   & 2.1   & 2.2   & 2.2   & 2.0   & 0.60   & 2.2   & 2.2   & 2.1   & 2.1   & 2.2   & 2.1  \\[0.4em]
    \multicolumn{1}{l}{Empirical powers} & 300   & 50    & 0.05  & 32.4  & 32.3  & 42.8  & 43.1  & 42.7  & 42.8  & 0.20   & 32.9  & 32.9  & 38.7  & 38.5  & 38.4  & 38.4  \\
          &       &       & 0.10   & 95.5  & 95.6  & 96.7  & 96.8  & 96.8  & 96.8  & 0.40   & 50.9  & 51.0  & 69.8  & 69.7  & 72.2  & 72.0  \\
          &       &       & 0.15  & 99.6  & 99.7  & 99.8  & 99.9  & 99.8  & 99.8  & 0.60   & 62.6  & 62.7  & 81.2  & 81.2  & 83.0  & 83.1  \\[0.2em]
          &       & 100   & 0.05  & 21.8  & 21.8  & 38.2  & 38.2  & 38.1  & 38.0  & 0.20   & 30.5  & 30.4  & 36.2  & 36.3  & 35.7  & 36.2  \\
          &       &       & 0.10   & 91.5  & 91.5  & 93.3  & 93.5  & 93.4  & 93.3  & 0.40   & 49.2  & 49.1  & 66.2  & 66.3  & 69.0  & 68.9  \\
          &       &       & 0.15  & 99.1  & 99.1  & 99.5  & 99.5  & 99.6  & 99.5  & 0.60   & 58.3  & 58.1  & 77.9  & 78.0  & 79.4  & 79.4  \\[0.2em]
          &       & 200   & 0.05  & 11.0  & 10.8  & 33.1  & 33.0  & 33.5  & 33.2  & 0.20   & 28.2  & 28.2  & 33.3  & 33.4  & 34.1  & 33.5  \\
          &       &       & 0.10   & 86.0  & 86.2  & 88.5  & 88.7  & 88.6  & 88.4  & 0.40   & 46.8  & 46.8  & 61.7  & 61.7  & 62.0  & 62.2  \\
          &       &       & 0.15  & 97.8  & 97.9  & 98.8  & 98.7  & 98.6  & 98.7  & 0.60   & 53.5  & 53.5  & 73.9  & 73.8  & 74.2  & 73.9  \\[0.2em]
          & 600   & 50    & 0.05  & 83.0  & 82.8  & 85.4  & 85.4  & 85.1  & 85.3  & 0.20   & 48.6  & 48.6  & 64.0  & 64.1  & 66.0  & 66.4  \\
          &       &       & 0.10   & 100  & 100  & 100  & 100  & 100  & 100  & 0.4   & 60.3  & 60.4  & 83.2  & 83.2  & 84.1  & 84.0  \\
          &       &       & 0.15  & 100  & 100  & 100  & 100  & 100  & 100  & 0.60   & 94.0  & 94.0  & 97.5  & 97.5  & 98.0  & 98.1  \\[0.2em]
          &       & 100   & 0.05  & 75.5  & 75.3  & 78.8  & 78.8  & 79.8  & 79.3  & 0.20   & 46.2  & 46.1  & 59.4  & 59.5  & 60.4  & 60.2  \\
          &       &       & 0.10   & 100  & 100  & 100  & 100  & 100  & 100  & 0.40   & 56.7  & 56.7  & 80.6  & 80.5  & 81.2  & 81.4  \\
          &       &       & 0.15  & 100  & 100  & 100  & 100  & 100  & 100  & 0.60   & 88.5  & 88.5  & 94.9  & 94.8  & 95.4  & 95.1  \\[0.2em]
          &       & 200   & 0.05  & 67.6  & 67.7  & 71.6  & 71.7  & 71.8  & 71.8  & 0.20   & 43.4  & 43.5  & 54.9  & 55.0  & 56.1  & 56.3  \\
          &       &       & 0.10   & 99.3  & 99.5  & 99.8  & 99.8  & 99.8  & 99.9  & 0.40   & 53.8  & 53.8  & 78.3  & 78.3  & 79.6  & 79.5  \\
          &       &       & 0.15  & 100  & 100  & 100  & 100  & 100  & 100  & 0.60   & 83.2  & 83.0  & 91.9  & 91.8  & 93.0  & 92.5  \\
    \end{tabular}%
		}
		\label{tab:fdr}%
  \vspace{-0.5em}
	\end{table}%
	
		Theorem \ref{tm:FDR} implies that when $n$ and $Q$ grow to infinity, the FDR
	should be controlled at the level of  $\alpha Q_0/Q$ with high probability, which equals  3.27\% in our settings with $\alpha=5\%$.
	In the $v$-th simulation replication, we can obtain $\hat{t}_v$ defined as \eqref{eq:t.hat}. For each $q\in[Q]$, denote by $V_{n,v}^{(q)}$ the normal quantile transformation of the p-value for $H_{0,q}$ in the $v$-th simulation replication. See its definition below \eqref{eq:diffTnq}. Besides the empirical FDR, we also consider the empirical power of the proposed multiple testing procedure defined as
	\begin{align*}
		\frac{1}{1000}\sum_{v=1}^{1000}\frac{1}{Q-Q_0}\sum_{q\in\mathcal{H}_1} I\{V_{n,v}^{(q)}>\hat{t}_v\} \,.
	\end{align*}
	As shown in Table \ref{tab:fdr}, the empirical FDR
	becomes closer to the limit rate 3.27\% as $n$ increases, and the proposed multiple testing procedure tends to be more conservative when $p$ becomes larger. On the other hand, the corresponding  empirical powers grow quickly  as $n$ increases, regardless of the smaller empirical FDR.
	The results also show that different  choices of   $c$ have little influence on  the empirical FDR and power in the models being examined. The choice of $K$ has little impact on the empirical FDR but appears to have some impact on the power. 
 In particular, when $K$ changes from 4 to 12, there is a notable increase in the empirical powers, while  the increase becomes insignificant when comparing $K=12$ with $K=n$.

	\section{Real Data Analysis}
	\label{data}
	
	\subsection{Batching county-level hires data}
	
	As an illustration of the techniques of this paper, we study data on new hires
	at a national level, obtained from  the Quarterly Workforce Indicators (QWI)
	of the Longitudinal Employer-Household Dynamics 
program at the U.S. Census Bureau \citep{abowd2009lehd}.\footnote{Data (Hires All:Counts) was extracted from {\it https://ledextract.ces.census.gov/qwi/all} on  October 5, 2022; all counties
		in each state were selected, with All NAICS
		and All Ownership (Firm Characteristics),
		No Worker Characteristics, and all available
		quarters.}
	The quarterly data is available for all counties, and we wish to
	obtain a classification of the database whereby we associate clusters of time series
	pertaining to various states, such that they are suitable for joint analysis.
	The national QWI hires data covers a variable number of years, with some states
	providing time series going back to 1990 (e.g., Washington), and others (e.g., Massachusetts) only commencing
	at 2010.   For each of 51 states (excluding D.C. but including Puerto Rico) there is a  new hires time series for each county.
	Additional description of the data, along with its relevancy to labor economics, can be found in
	\cite{hyatt2019labor}.
	
	Given QWI county-level data on new hires, we want to know whether we may analyze the
	data  state-by-state, or whether there is additional time series information to be gleaned by examining relationships across states.      For any two states $i$ and $j$, with $1 \leq i\neq j \leq 51$,
	let $B_i$ and $B_j$ denote batches of time series indices corresponding to the counties within a state.  Between each pair $(i,j)$, we test whether   the cross-spectrum between $B_i$ and $B_j$  is not identically zero,
	with rejection of the null hypothesis indicating there may be merit in considering both batches in a joint time series model.
	We restrict ourselves to examining the cross-spectrum between distinct batches of series $B_i$ and $B_j$,
	which is in contrast to the simulation of Section \ref{sec:mult-test}, where one can also test
	$B_i$ with itself.  This latter procedure would investigate each state, inquiring whether the county-level
	time series of that state should be jointly modeled; instead we focus on whether the county-level
	time series of two distinct states should be jointly modeled.
	Since the time series are quarterly, we assess the cross-spectrum at the seasonal frequencies, viz.
	$\mathcal{J}= \{-\pi, -\pi/2, 0, \pi/2\}$.
	We   apply  the global test by  taking a supremum over the four seasonal frequencies;
	we also consider multiple testing with FDR control.
	
	Many of the series exhibit strong trend and seasonal effects, and it is important to
	ensure that the data is stationary;  therefore, we apply either regular differencing or seasonal differencing.
	Although these operations may over-difference certain series,  our theory has no
	requirement that the spectral density be non-zero, so there is no impediment to analysis with this approach.
	In addition to this data pre-processing, it is   necessary to find common sample sizes
	for each pair of batches $B_i$ and $B_j$ since the start dates differ greatly.  In fact,
	if we were to  consider all $p = 3218$ county-level time series in one huge batch,
	the maximal common sample size is $n = 25$, which is clearly too tiny for such a huge $p$.
	Instead,  for each of the  $\binom{51}{2} =1275$ possible pairings $(B_i, B_j)$, we determine the common sample available,
	and the dimension of the resulting paired data set is $p = | B_i \cup B_j |$.
	The sample size $n$ in each case   is defined to be the most recent contiguous block of times where both batches are fully observed (no missing values);
 by excluding the sparsely measured
	counties of Kalawao, HI and McPherson, NE,
	we ensure there is a common sample for
	every pairing.
Our methodology requires  that   $\tilde{n} = n - 2  l_n  $  (with $l_n$ selected by the data-driven method discussed in Remark \ref{rk:1}(c))
 is positive; we impose $n > 2 l_n + 1$ and $n-d > 7$
 (where $d=1$ for regular differencing and $d=4$
 for seasonal differencing), the latter condition
 ensuring that the data-driven $l_n$ can be 
 calculated, which ensures $\tilde{n} \geq 2$.
 For $284$ state pairs (or $22.3 \%$)
 where regular differencing is used
 the common sample does  not meet these requirements,
 and these cases are skipped over;
 for the case of seasonal differencing, only
 $1$ state pair violates the requirements.

	\begin{figure}[h!]
		\centering
		\includegraphics[width = 0.6 \textwidth]{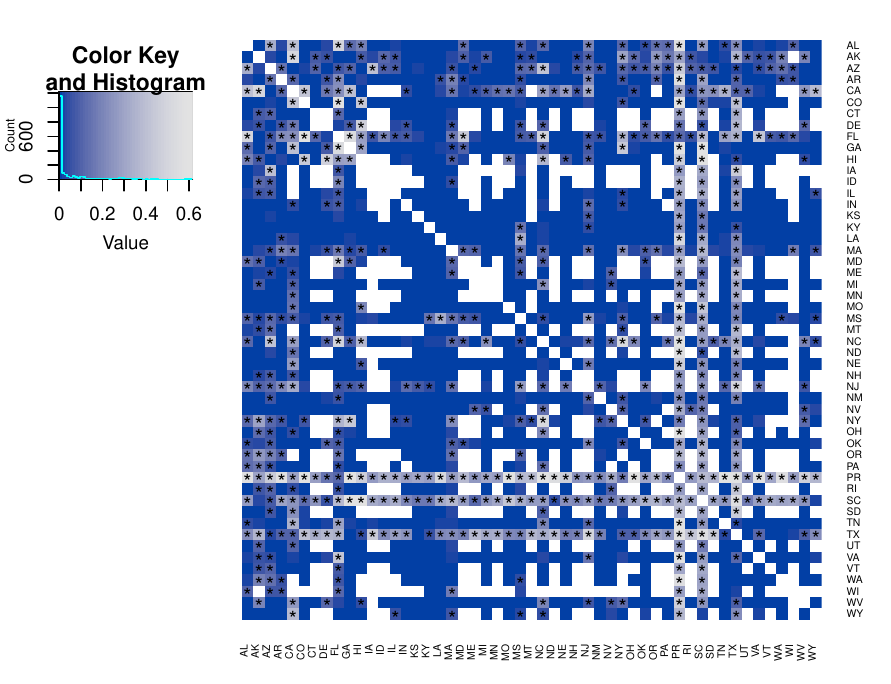}
        \vspace{-0.8em}
		\caption{Heatmap of p-values for $51$ state pairs, testing whether the cross-spectrum of each pair	is zero at seasonal frequencies   $\mathcal{J}= \{-\pi, -\pi/2, 0, \pi/2\}$.  A white box indicates
     a pair for which no test is computed; a star 
     marks pairs that are not significant using FDR control, where $\alpha = 5\%$.   Each series has been differenced. }
		\label{fig:heatmapGR}
	\end{figure}

	\begin{figure}[h!]
		\centering
		\includegraphics[width = 0.6 \textwidth]{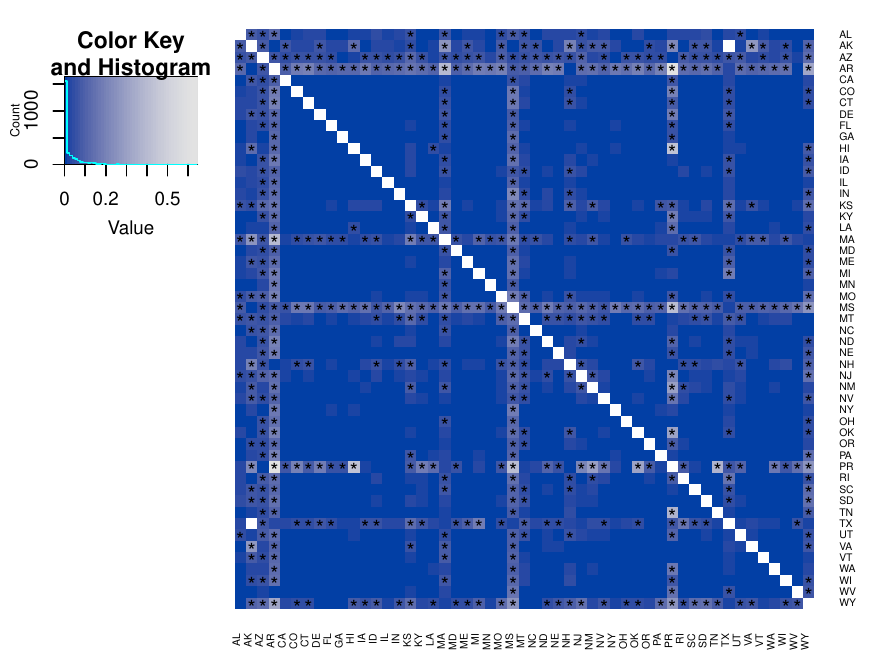}
  \vspace{-0.8em}
		\caption{Heatmap of p-values for $51$ state pairs, testing whether the cross-spectrum of each pair		is zero at seasonal frequencies   $\mathcal{J}= \{-\pi, -\pi/2, 0, \pi/2\}$.  A white box indicates
     a pair for which no test is computed; a star 
     marks pairs that are not significant using FDR control, where $\alpha = 5\%$.  Each series has been seasonally differenced. }
		\label{fig:heatmapAR}
	\end{figure}

	Our final output is summarized with two heat maps, displayed in Figures \ref{fig:heatmapGR} and
	\ref{fig:heatmapAR}, for the cases of regular differencing and seasonal differencing respectively.
	Low p-values have a darker color, and most of the pairs end up rejecting the null hypothesis;
	we have marked with a star those pairs that do not reject the null hypothesis when using the FDR control	method with $\alpha =5\%$.   The diagonal,
 as well as any state pairs with 
 insufficient common sample,  is marked white, indicating that no test is conducted in such cases.
	For the regular differencing there are  $338$  pairs with no association, and  $295$  pairs of no association for the seasonal differencing;
 hence there are $653$ pairs
 under regular differencing for which a joint
 analysis may be useful, and $979$ such pairs
 under seasonal differencing.

\subsection{Detecting seasonal over-adjustment in Texas hires}
	
As a second application of our methods,
we focus on the county hires data for Texas,
and wish to detect over-adjustment in any
of the seasonal adjustments.
First we seasonally adjust the data using
the automatic methods of the X-13ARIMA-SEATS
software \citep{us2020x}, excluding Loving
county (due to many missing values); note
that there is no official seasonal adjustment
for these series, and we use default settings
in the software so as to eliminate the impact
of human intervention in the analysis.
We then analyze the $p=253$ county-level
seasonally adjusted time series, testing
$H_0: f_{i,i}(\omega)=0$ for all $i\in[p]$
with $\mathcal{J}= \{-\pi, -\pi/2,  \pi/2\}$,
i.e., the seasonal frequencies (excluding
the trend frequency).

We consider two approaches to the testing:
first, we can examine the entire batch of
Texan county data by setting
$\mathcal{I}=[p]$, thereby obtaining a single
test statistic that measures the over-adjustment
problem's prevalence for all series.
This yields a p-value of $0.203$, indicating	that at a $5 \%$ level there is a failure
to reject the null hypothesis of over-adjustment,
i.e., there is over-adjustment.
 Secondly, we can conduct a univariate test for each of the $p$ series, and use the FDR control to manage the
multiple testing results.  In this case
(setting the nominal level at $5 \%$)
the support recovery yields $49$ counties
where the null hypothesis is not rejected;
these can then be scrutinized by a human analyst.
	
	This example showcases how the practical
	problem of seasonally adjusting thousands
	or millions of time series can be managed
	with limited computational and human resources:
	the automatic methodology of the X-13ARIMA-SEATS
	software can be applied with default settings,
	and the output can be quickly assessed
	for defective adjustments -- using FDR control and support recovery, the (hopefully small)
	subset of problematic seasonal adjustments
	can then be examined by a seasonal adjustment
	expert.  However, we have not here addressed
	the more subtle problem of seasonal
	under-adjustment (where the spectral density
	of a seasonally adjusted time series still
	has a local peak at some of the seasonal frequencies), which we leave for future research.

	\section{Discussion}
	\label{discussion}
	
	Motivated by the increasing availability of high-dimensional time series, we develop new  inference methodology and theory for the spectral density matrix in the high-dimensional setting, which has not yet
	been fully explored in the literature. We overcome  both methodological and theoretical challenges that high dimensionality induces by extending the celebrated Gaussian approximation
	and multiplier bootstrap to the testing of a high-dimensional parameter formulated in the frequency domain, which seems to be the first such effort in the literature. In particular, we develop a  maximum-type test statistic and  computationally feasible parametric bootstrap approximation to test the nullity of coherence at a pre-specified set of frequencies for given  component pairs.
	The theoretical justification is established under a setting that allows weak temporal dependence,  flexible contemporary dependence across $p$ components, and exponential rate growth for the dimension. In addition,
	we  develop
	a multiple testing procedure to  recover the support set for the nullity of coherence at a given set of frequencies. A rigorous theory for the FDR
	control is also provided.  Finally, we illustrate the size and power of the proposed tests through simulations and real data analysis.

Given $\mathcal{I}\subset[p]^2$ and $\mathcal{J}\subset[-\pi,\pi)$, we can also establish the Gaussian approximation theory for
 \begin{align*}
\mathcal{T}^{{\rm s}}=\sup_{\omega\in\mathcal{J}}\max_{(i,j)\in\cI}
 \bigg|\sqrt{\frac{n}{l_n}}\frac{\hat{f}_{i,j}(\omega)-f_{i,j}(\omega)}{\{\hat{f}_{i,i}(\omega)\hat{f}_{j,j}(\omega)\}^{1/2}}\bigg|^2\,.
 \end{align*}
 Recall $r=|\mathcal{I}|$ and $\tilde{n}=n-2l_n$. For the bijective mapping $\bchi(\cdot)=\{ \chi_1(\cdot),\chi_2(\cdot)\}$ specified in Section \ref{sec:ga}, write ${\bf W}(\omega)={\rm diag}\{ \hat{f}_{\chi_1(1),\chi_1(1)}^{-1/2}(\omega) \hat{f}_{\chi_2(1),\chi_2(1)}^{-1/2}(\omega), \ldots, \hat{f}_{\chi_1(r),\chi_1(r)}^{-1/2}(\omega) \hat{f}_{\chi_2(r),\chi_2(r)}^{-1/2}(\omega)\}$. Letting 
$$
 \xi_{\cJ}^{\rm s}=\sup_{\omega\in\cJ}\max_{\ell\in[r]} \{|\hat{\eta}_{2\ell-1}^{\rm ext, {\rm s}}(\omega)|^2 +|\hat{\eta}_{2\ell}^{\rm ext, {\rm s}}(\omega)|^2\}$$
with $\{\hat{\eta}_1^{\rm ext, s}(\omega),\ldots,\hat{\eta}_{2r}^{\rm ext, s}(\omega) \}^{\T}=\{{\bf W}(\omega)\otimes \bA(\omega) \} (\tilde{n}^{-1/2}\sum_{t=1}^{\tilde{n}}\epsilon_t \hat{\bc}_t)$ for $\bA(\omega)$ defined in \eqref{eq:A}, $\hat{\bc}_t=(\hat{\bc}_{1,t}^\T,\ldots,\hat{\bc}_{r,t}^\T)^\T$
with $\hat{\bc}_{\ell,t}$ defined in \eqref{eq:hatcl}, and
$(\epsilon_1,\ldots,\epsilon_{\tilde{n}})^{\T}\sim \mathcal{N}(\bzero, \bTheta)$ with $\bTheta$ defined in Section \ref{sec:ga},
 we can show $$\sup_{u\geq0}|\bP(\mathcal{T}^{\rm s}\leq u) - \bP(\xi_{\cJ}^{\rm s}\leq u\,|\,\cX_n)|=o_{\rm p}(1)\,.$$ Then we can also use the studentized test statistics  
 \begin{equation*}
 \begin{split}
 T_n^{\rm s} =&~ \sup_{\omega\in\cJ}\max_{(i,j)\in\cI}
 \bigg|\sqrt{\frac{n}{l_n}}\frac{\hat{f}_{i,j}(\omega)}{\{ \hat{f}_{i,i}(\omega)\hat{f}_{j,j}(\omega)\}^{1/2}}\bigg|^2\,,\\
 T_n^{(q),{\rm s}}=&~\sup_{\omega\in\mathcal{J}^{(q)}}\max_{(i,j)\in\mathcal{I}^{(q)}}\bigg|\sqrt{\frac{n}{l_n}}\frac{\hat{f}_{i,j}(\omega)}{\{ \hat{f}_{i,i}(\omega)\hat{f}_{j,j}(\omega)\}^{1/2}}\bigg|^2
 \end{split}
 \end{equation*}
 in Sections \ref{subsec:global} and \ref{subsec:multiple}, respectively, for the associated inference problems. 
Although $T_n^{\rm s}$ and $T_n^{(q),{\rm s}}$ can correct for heterogeneity, their numerical performance is not robust when the sample size $n$ is small due to the fact that the estimation of the denominator in $T_n^{\rm s}$ and $T_n^{(q),{\rm s}}$ may
	not be accurate enough. Such
	phenomenon has been empirically observed in \cite{CQYZ_2018}. When the sample size $n$ is small, we suggest to use nonstudentized statistics $T_n$ and $T_n^{(q)}$, respectively, for the inference problems considered in Sections \ref{subsec:global} and \ref{subsec:multiple}. 
 We can also use the Gaussian approximation technique to construct the simultaneous inference for the coherence matrix. See Section G in the supplementary material for details.

	To conclude, we mention several topics that are worth future investigation. Firstly, there are two tuning parameters involved in our procedure, i.e., $l_n$ and $b_n$, the choices of which can have an impact on the finite sample performance. Some theoretical investigation that can lead to a data-driven formula in the high-dimensional setting would be desirable.  Andrews' rule for the choice of $b_n$ is used here, but there is no good theoretical justification for it at this moment, and there might be better formulas in practice. Secondly, as the precision matrix plays an important role for high-dimensional independent and identically distributed data, the partial coherence matrix is the analogue in the spectral domain and its inference in the high-dimensional setting would be of great importance; see \cite{KP_2022} for a recent effort.
	Thirdly, one of the limitations of the spectral density matrix is that it can only characterize second-order properties. Recently, \cite{BK_2019} proposed quantile coherency to characterize the cross-series nonlinear dependence in the frequency domain, as an extension of quantile spectrum for univariate time series developed in \cite{H_2013} and \cite{KVDH_2016}. It would be interesting to extend our results to detect nonlinear dependence in the frequency domain for high-dimensional time series.
	We leave these topics for future research.

	\section*{Appendix}
	For given $(\vartheta,\rho)$ such that $\vartheta>1$ and $0<\rho<(\vartheta-1)/(3\vartheta-2)$, let
	\begin{align}
		f_1(l_n,n;\vartheta,\rho)  :=&\,  \min\bigg\{ \frac{n^{\rho/2}}{l_n^{3/2}}, \,
		\frac{n^{(1-2\rho)/5}}{l_n^{2/5}}, \, \frac{n^{(\vartheta+2\rho-3\rho\vartheta-1)/(8\vartheta-4)}}{ l_n^{1/4}}, \label{eq:f1.f2}  \\
		&~~~~~~~~~~~~~~~~~~ \frac{n^{(2\vartheta+3\rho-4\rho\vartheta-2)/(12\vartheta-6)}}{l_n^{1/6}}, \,
		\frac{n^{1/9}}{l_n\log^{8/3}(l_n)} \bigg\}  \,,\notag \\
		f_2(\vartheta,\rho):=&\, \min\bigg(\frac{\rho}{3}, \, \frac{\vartheta+2\rho-3\rho\vartheta-1}{2\vartheta-1}\bigg) \,, \label{eq:f1.f2.2} \\
		f_3(\vartheta,\rho)  :=&\,  \min\bigg(\frac{\rho}{6},\,
		\frac{1-2\rho}{12},\,
		\frac{\vartheta+2\rho-3\rho\vartheta-1}{12\vartheta-6}\bigg)  \,, \label{eq:f3.f4.f5}\\
		f_4(l_n,n,Q;\vartheta,\rho):=&\, \min\bigg\{
		\frac{n^{1/9}}{l_n\log^{2/3}(l_n)Q^2},\,
		\frac{n^{(\vartheta+2\rho-3\rho\vartheta-1)/(8\vartheta-4)}}{l_n^{1/4}Q^{3/2}}, \frac{n^{(1-2\rho)/5}}{l_n^{2/5}Q^{12/5}},\, \label{eq:f3.f4.f5.2}\\
		&~~~~~~~~~~~~~~~\frac{n^{(2\vartheta+3\rho-4\rho\vartheta-2)/(12\vartheta-6)}}{l_n^{1/6}Q}, \,
		\frac{n^{\rho/2}}{l_n^{3/2}Q^3},\,
		\frac{n^{1/9}}{l_n\log^{8/3}(l_n)},\,
		\frac{n^{\rho/3}}{Q^2} \bigg\}  \,, \notag\\
		f_5(n,Q;\vartheta,\rho):=&\, \min\bigg\{\frac{n^{(\vartheta+2\rho-3\rho\vartheta-1)/(2\vartheta-1)}}{Q^6},\,
		\frac{n^{\rho/3}}{Q^2}\bigg\} \,. \label{eq:f3.f4.f5.3}
	\end{align}

\section*{Supplementary Materials}

The supplementary material contains detailed discussion for Conditions \ref{as:moment}--\ref{as:eigen}, all technical proofs of the main results, and the procedure for the statistical inference of high-dimensional coherence matrix.





 







\setlength{\bibsep}{0.2pt plus 1ex}
\begin{spacing}{0.95}
\bibliographystyle{jasa}

\bibliography{HDspec_bib_s}
\end{spacing}

\newpage
\spacingset{1.7} 

\newpage
\appendix

\setcounter{equation}{0}
\setcounter{table}{0}
\setcounter{lemma}{0}
\newtheorem{props}{Proposition}
\setcounter{props}{0}
\setcounter{page}{1}

\renewcommand{\thelemma}{L\arabic{lemma}}
\numberwithin{equation}{section}
\renewcommand{\theprops}{P\arabic{props}}
\renewcommand{\thepage}{S\arabic{page}}
\renewcommand{\thetable}{S\arabic{table}}

\newpage
\setcounter{page}{1}

\begin{center}
{\large\bf SUPPLEMENTARY MATERIAL} \medskip\\
{\bf Jinyuan Chang, ~\, Qing Jiang, ~\,Tucker McElroy, ~\,Xiaofeng Shao$^*$}
\end{center}

Throughout the supplementary material, we use
$C$ to denote a generic universal positive finite constant that may be different in different uses.   For any positive integer $q\geq2$, we write $[q]=\{1,\ldots,q\}$. For any real-valued numbers $x$ and $y$, we write $|x|_{+}=\max(0,x)$ and $x\vee y=\max(x,y)$. 
For two positive real-valued sequences $\{a_n\}$ and $\{b_n\}$, we write $a_n\lesssim b_n$ or $b_n\gtrsim a_n$ if there is a universal constant $C>0$ such that $\limsup_{n\rightarrow\infty}a_n/b_n\leq C$, write $a_n\asymp b_n$ if $a_n\lesssim b_n$ and $b_n\lesssim a_n$ hold simultaneously, and write $a_n\ll b_n$ or $b_n\gg a_n $ if $\limsup_{n\rightarrow\infty}a_n/b_n=0$.
Denote by $\mathbb{S}^{q-1}$ the $q$-dimensional unit sphere. For a complex-valued number $x$,  denote by $|x|$, $\Re(x)$ and $\Im(x)$ the modulus, the real part and the imaginary part of $x$, respectively. For any $q_1\times q_2$ complex-valued matrix $\bM=(m_{i,j})_{q_1\times q_2}$, let $|\bM|_{\infty}=\max_{i\in[q_1],j\in[q_2]}|m_{i,j}|$, $|\bM|_1=\sum_{i=1}^{q_1}\sum_{j=1}^{q_2}|m_{i,j}|$, and $|\bM|_0=\sum_{i=1}^{q_1}\sum_{j=1}^{q_2}I(m_{i,j}\neq0)$, where $I(\cdot)$ denotes the indicator function. Specifically, if $q_2=1$, we use $|\bM|_\infty=\max_{i\in[q_1]}|m_{i,1}|$, $|\bM|_1=\sum_{i=1}^{q_1}|m_{i,1}|$, and $|\bM|_0=\sum_{i=1}^{q_1}I(m_{i,1}\neq 0)$ to denote the $L_\infty$-norm, $L_1$-norm and $L_0$-norm of the  $q_1$-dimensional complex-valued vector $\bM$, respectively.
For two $p$-dimensional real-valued vectors $\ba=(a_1,\dots,a_p)^{\T}$ and $\bw=(w_1,\ldots,w_p)^{\T}$, we say $\ba\leq \bw$ if $a_j\leq w_j$ for any $j\in[p]$.
Let $D$ be a subset of a Euclidean space. For given $\epsilon>0$ and some metric $d(\cdot,\cdot)$, a subset of $D$ denoted by $D_{\epsilon}$ is called an $\epsilon$-net of $D$ if for any $u\in D$ there exists some $\tilde{u}\in D_{\epsilon}$ such that $d(u,\tilde{u})\leq \epsilon$. 
For a countable set $\mathcal{F}$, we use $|\mathcal{F}|$ to denote the cardinality of $\mathcal{F}$.
For a $q$-dimensional vector $\ba$, denote by $\ba_{\mathcal{L}}$ the subvector of $\ba$ collecting the components indexed by a given index set $\mathcal{L}\subset[q]$.

\section{Discussion for Conditions \ref{as:moment}--\ref{as:eigen}}

Condition \ref{as:moment} is equivalent to $\max_{j\in[p]}\max_{t\in[n]}\bP(|x_{j,t}|>u)\leq C_*\exp(-C_{**} u^2)$ for any $u>0$ with some positive constants $C_*$ and $C_{**}$, which is a common assumption in the literature of ultra high-dimensional data analysis for deriving exponential-type upper bounds for the tail probabilities of certain statistics.  See \cite{CHYY_2023_supp}  and reference therein. Our technical proofs indeed allow $\max_{j\in[p]}\max_{t\in[n]}\bP(|x_{j,t}|>u)\leq C_2\exp(-C_1u^{\tau_1})$ for any $u>0$ with some constant $\tau_1\in(0,2]$,
but we set $\tau_1=2$ in Condition \ref{as:moment} to simplify the presentation.

Condition \ref{as:betamixing} provides a simple upper bound for the $\alpha$-mixing coefficient $\alpha_n(k)$, which is mild in the literature. Under certain conditions, VAR processes, multivariate ARCH processes, and multivariate GARCH processes all satisfy Condition \ref{as:betamixing}; see \cite{Boussama_etal:2011_supp}, 
\cite{HP_2009_supp}, 
and \cite{Wong_etal_2020_supp}. 
\cite{CCW_2021_supp} also gives two high-dimensional time series models satisfying Condition \ref{as:betamixing}. More generally, repeating our current proof of Proposition \ref{tm:1} with more lengthy arguments, we can show our procedure still works when Condition \ref{as:betamixing} is replaced by a more general upper bound
$
\alpha_n(k)\leq C_3\exp\{-C_4(L_n^{-1}k)^{\tau_2}\}$, where $L_n$ is allowed to diverge with $n$, and $\tau_2\in(0,1]$ is a constant. With this new upper bound, the presentation of Proposition 1 will be more complicated, which still allows $K$ and $r$ to diverge exponentially fast if $L_n$ diverges slower than $n^\upsilon$ for some constant $\upsilon\in(0,1)$.  Such general upper bound for high-dimensional dependent sequence was considered in \cite{CHLT_2022_supp}. In our paper, we just use Condition \ref{as:betamixing} to simplify our presentation. Recall $\bx_t=(x_{1,t},\ldots,x_{p,t})^\T$. When $\{x_{1,t}\}_{t\geq1},\ldots,\{x_{p,t}\}_{t\geq1}$ are $p$ independent ARMA processes with continuous innovation distributions, then each $\{x_{j,t}\}_{t\geq1}$ is $\alpha$-mixing with exponential decay rates. See Section 2.6.1 of 
\cite{FY_2003_supp}. 
Applying Theorem 5.1 of \cite{Bradley_2005_supp}, 
we know the $\alpha$-mixing coefficient $\alpha_n(k)$ of $\{\bx_t\}_{t\geq1}$ satisfies $\alpha_n(k)\leq p\exp(-ck)$ for some universal constant $c>0$, which implies the general upper bound holds for $\tau_2=1$ and $L_n\asymp \log p$. If we just require $\max_{j\in[p]}\max_{t\in[n]} \bP(|x_{j,t}|>u)=O\{u^{-(\nu+\epsilon)}\}$ for any $u>0$ in Condition \ref{as:moment} and $ \alpha_n(k)=O\{k^{-\nu(\nu+\epsilon)/(2\epsilon)}\}$ for all $k\geq 1$ in Condition \ref{as:betamixing} with some constants $\nu>2$ and $\epsilon>0$, 
we can apply the Fuk-Nagaev-type inequalities to construct the upper bounds for the tail probabilities of certain statistics for which our procedure still works for $K$ and $r$ diverging at some polynomial rates of $n$.

Condition \ref{as:eigen} is a mild technical assumption for the validity of the Gaussian approximation, 
where the constant $C_5$ can be replaced by some $\varphi_n=o(1)$ at the expenses of lengthier proofs. Recall 
$T_n(\omega;\cI)=\max_{(i,j)\in\cI}|\sqrt{nl_n^{-1}}\{\hat{f}_{i,j}(\omega)-f_{i,j}(\omega)\}|^2$. Define $\check{T}_n(\omega;\cI)=\max_{\ell\in[r]}\{|\check{\eta}^{\rm ext}_{2\ell-1}(\omega)|^2 + |\check{\eta}^{\rm ext}_{2\ell}(\omega)|^2\}$ with
    \begin{align*}
	\check{\bfeta}^{\ext}(\omega)\equiv\{\check{\eta}^{\ext}_1(\omega),\ldots,\check{\eta}_{2r}^{\ext}(\omega)\}^\T=\{\bI_r\otimes \bA(\omega) \}\frac{1}{\sqrt{n}}\sum_{t=1}^{\tilde{n}}\bc_{t}  \,.
\end{align*}    
    Under some regularity conditions, Lemma L5 in Section \ref{sec:auxlemma} indicates that $|\sup_{\omega\in\cJ}T_n(\omega;\cI) - \sup_{\omega\in\cJ}\check{T}_n(\omega;\cI)|$ is asymptotically negligible.     
    Hence, the key step in our theoretical analysis is to characterize the distribution of $\sup_{\omega\in\cJ}\check{T}_n(\omega;\cI)$ by Gaussian approximation.
   Notice that
    \begin{align*}
        \bigg\{\sup_{\omega\in\cJ}\check{T}_n(\omega;\cI)\leq u\bigg\} 
        =\bigcap_{\omega\in\cJ}\{\check{T}_n(\omega;\cI)\leq u\} =\bigcap_{\omega\in\cJ}\{\check{\bfeta}^{\rm ext}(\omega)\in B(u)\} \,,
    \end{align*}
    where $B(u)=\bigcap_{j=1}^{r}B_j(u)$ for $B_{j}(u)=\{\bb\in\bR^{2r}:\bb_{S_j}^{\T}\bb_{S_j}\leq u\}$ with $S_j=\{2j-1,2j\}$. Here the set $B_j(u)$ is convex in $\mathbb{R}^{2r}$ that only depends on the components of the $(2r)$-dimensional vectors indexed by $S_j$. We can reformulate $B_j(u)$ as  $$B_j(u)=\bigcap_{\bfd\in\mathbb{S}^{2r-1}:\,\bfd_{S_j}\in\mathbb{S}^{1}} \{\bb\in \mathbb{R}^{2r}:\bfd^{\T}\bb\leq u^{1/2}\}\,.$$ 
    Write $\mathcal{F} = \bigcup_{j=1}^{r}\{\bfd\in\mathbb{S}^{2r-1}: \bfd_{S_j}\in\mathbb{S}^{1}\}$ with $S_j=\{2j-1,2j\}$. 
    Then $B(u)=\bigcap_{\bfd\in\mathcal{F}}\{\bb\in \mathbb{R}^{2r}:\bfd^{\T}\bb\leq u^{1/2}\}$, which implies 
     \begin{align*}
        \mathbb{P}\bigg\{\sup_{\omega\in\cJ}\check{T}_n(\omega;\cI)\leq u\bigg\} 
         =\mathbb{P}\bigg\{ \frac{1}{\sqrt{n}}\sum_{t=1}^{\tilde{n}}  {\bf d}^{\T}\{{\bf I}_r\otimes \bA(\omega)\}\bc_t \leq u^{1/2} \mbox{ for any ${\bf d}\in \mathcal{F}$ and $\omega\in\cJ$}\bigg\}  
    \end{align*}
for any $u>0$. Applying the Gaussian approximation technique to approximate  $\mathbb{P}\{\sup_{\omega\in\cJ}\check{T}_n(\omega;\cI)\leq u\}$, we need to use the anti-concentration inequality (Lemma A.1, \citealp{CCK_2017_supp}) for Gaussian random vector which requires
 \begin{align}\label{eq:covbdd1}
{\rm Var}\bigg[\frac{1}{\sqrt{n}}\sum_{t=1}^{\tilde{n}}  {\bf d}^{\T}\{{\bf I}_r\otimes \bA(\omega)\}\bc_t\bigg]\geq C_*
 \end{align}
 for any $\omega\in\mathcal{J}$ and ${\bf d}\in\mathcal{F}$, where $C_*>0$ is a universal constant. See also \cite{CCK_2013_supp} and \cite{CCK_2022_supp}.  Due to $\tilde{n}=n-2l_n$ with $l_n=o(n)$, we have  
 \[
 {\rm Var}\bigg[\frac{1}{\sqrt{n}}\sum_{t=1}^{\tilde{n}}  {\bf d}^{\T}\{{\bf I}_r\otimes \bA(\omega)\}\bc_t\bigg]={\bf d}^{\T}\bSigma(\omega,\omega){\bf d}\cdot\{1+o(1)\}
\]
 for any $\omega\in\mathcal{J}$ and ${\bf d}\in\mathcal{F}$. Hence, Condition \ref{as:eigen} is equivalent to \eqref{eq:covbdd1}.

\section{Proof of Proposition \ref{tm:1}}
Recall $\bmu=\bE(\bx_t)$ and $\mathring{\bx}_t=(\mathring{x}_{1,t},\ldots,\mathring{x}_{p,t})^{\T}=\bx_t-\bmu$. Write $\bar{\mathring{\bx}}=n^{-1}\sum_{t=1}^n\mathring{\bx}_t$. Then
\begin{align*}
	\widehat{\bGamma}(k)=\left\{
	\begin{aligned}
		\frac{1}{n}\sum_{t=1}^{n-k}(\mathring{\bx}_{t+k}-\bar{\mathring{\bx}})(\mathring{\bx}_t-\bar{\mathring{\bx}})^{\T} \,,~~~~~~&\textrm{if}~k\geq0\,,\\
		\frac{1}{n}\sum_{t=-k+1}^{n}(\mathring{\bx}_{t+k}-\bar{\mathring{\bx}})(\mathring{\bx}_t-\bar{\mathring{\bx}})^{\T} \,,~~~~&\textrm{if}~k<0\,.
	\end{aligned}
	\right.
\end{align*}
For any $\omega\in[-\pi,\pi]$, it holds that
\begin{align}\label{eq:expansion}
	\widehat{\bF}(\omega)-\bF(\omega)=&~\underbrace{\frac{1}{2\pi}\sum_{k=-l_n}^{l_n}\mathcal{W}\bigg(\frac{k}{l_n}\bigg)\bigg\{\widetilde{\bGamma}(k)-\frac{n-|k|}{n}\bGamma(k)\bigg\}e^{-\iota k\omega}}_{\textrm{I}(\omega)}\notag\\
	&~ +\underbrace{\frac{1}{2\pi}\sum_{k=-l_n}^{l_n}\mathcal{W}\bigg(\frac{k}{l_n}\bigg)\frac{n-|k|}{n}\bGamma(k)e^{-\iota k\omega}-\bF(\omega)}_{\textrm{II}(\omega)}\notag\\
	&~ -\underbrace{\frac{1}{2\pi}\sum_{k=1}^{l_n} \mathcal{W}\bigg(\frac{k}{l_n}\bigg)\bar{\mathring{\bx}}\bigg\{\bigg(\frac{1}{n}\sum_{t=1}^{n-k}\mathring{\bx}_t\bigg)e^{-\iota k\omega}+\bigg(\frac{1}{n}\sum_{t=k+1}^n\mathring{\bx}_t\bigg)e^{\iota k\omega}\bigg\}^\T}_{\textrm{III}(\omega)}\\
	&~ -\underbrace{\frac{1}{2\pi}\sum_{k=1}^{l_n} \mathcal{W}\bigg(\frac{k}{l_n}\bigg)\bigg\{\bigg(\frac{1}{n}\sum_{t=1}^{n-k}\mathring{\bx}_{t+k}\bigg)e^{-\iota k\omega}+\bigg(\frac{1}{n}\sum_{t=k+1}^n\mathring{\bx}_{t-k}\bigg)e^{\iota k\omega}\bigg\}\bar{\mathring{\bx}}^\T}_{\textrm{IV}(\omega)}\notag\\
	&~ +\underbrace{\frac{1}{2\pi}\sum_{k=-l_n}^{l_n}\mathcal{W}\bigg(\frac{k}{l_n}\bigg)\bigg\{\frac{n-|k|}{n}-2I(k=0)\bigg\}\bar{\mathring{\bx}}\bar{\mathring{\bx}}^\T e^{-\iota k\omega}}_{\textrm{V}(\omega)}\notag\,,
\end{align}
where
\begin{align*}
	\widetilde{\bGamma}(k)=\left\{
	\begin{aligned}
		\frac{1}{n}\sum_{t=1}^{n-k}\mathring{\bx}_{t+k}\mathring{\bx}_t^\T \,,~~~~~~&\textrm{if}~k\geq0\,,\\
		\frac{1}{n}\sum_{t=-k+1}^{n}\mathring{\bx}_{t+k}\mathring{\bx}_t^\T\,,~~~~&\textrm{if}~k<0\,.
	\end{aligned}
	\right.
\end{align*}
Recall $r=|\mathcal{I}|$ and $\bchi(\cdot)=\{\chi_1(\cdot),\chi_2(\cdot)\}$ is a given bijective mapping from $[r]$ to $\mathcal{I}$ such that for any $(i,j)\in\cI$ there exists a unique $\ell\in[r]$ satisfying $(i,j)=\bchi(\ell)$. We define
\begin{align*}
	\bfeta(\omega)\equiv\{\eta_1(\omega),\ldots,\eta_r(\omega)\}^\T=n^{1/2}l_n^{-1/2}\{\hat{f}_{\bchi(1)}(\omega)-f_{\bchi(1)}(\omega), \ldots,\hat{f}_{\bchi(r)}(\omega)-f_{\bchi(r)}(\omega)\}^\T \,.
\end{align*}
Based on such defined $\bfeta(\omega)$, we consider a $(2r)$-dimensional real-valued vector
\[
\bfeta^{\ext}(\omega)=\{\eta^{\ext}_1(\omega),\ldots,\eta^{\ext}_{2r}(\omega)\}^{\T}=[\Re\{\eta_1(\omega)\}, \Im\{\eta_1(\omega)\},\ldots, \Re\{\eta_r(\omega)\},\Im\{\eta_r(\omega)\}]^\T\,.
\]
Thus, we have
\begin{equation}\label{eq:Tn}
	T_n(\omega;\cI)=\max_{\ell\in[r]}\{|\eta^{\rm ext}_{2\ell-1}(\omega)|^2+|\eta^{\rm ext}_{2\ell}(\omega)|^2\}\,.
\end{equation}
Let
\[
\check{\bfeta}(\omega)=\{\check{\eta}_1(\omega),\ldots,\check{\eta}_r(\omega)\}^\T=n^{1/2}l_n^{-1/2}\{\zeta_{\bchi(1)}(\omega),\ldots,\zeta_{\bchi(r)}(\omega)\}^\T
\] with
\begin{align}\label{eq:zetaij}
	\zeta_{i,j}(\omega)=\frac{1}{2\pi n}\sum_{t=l_n+1}^{n-l_n}\sum_{k=-l_n}^{l_n}\mathcal{W}\bigg(\frac{k}{l_n}\bigg)\{\mathring{x}_{i,t+k}\mathring{x}_{j,t}-\gamma_{i,j}(k)\}e^{-\iota k\omega}\,.
\end{align}
Then $\check{\bfeta}^{\rm ext}(\omega)$ defined in \eqref{eq:etacheckext} can be also reformulated as
\[
\check{\bfeta}^{\ext}(\omega):=\{\check{\eta}_1^{\ext}(\omega),\ldots,\check{\eta}_{2r}^{\ext}(\omega)\}^{\T}=[\Re\{\check{\eta}_1(\omega)\}, \Im\{\check\eta_1(\omega)\}, \ldots,\Re\{\check\eta_r(\omega)\},\Im\{\check{\eta}_r(\omega)\}]^{\T}\,.
\]
Define
\begin{equation}\label{eq:checkt}
	\check{T}_n(\omega;\mathcal{I})=\max_{\ell\in[r]}\{|\check{\eta}_{2\ell-1}^{\ext}(\omega)|^2+|\check{\eta}_{2\ell}^{\ext}(\omega)|^2\}\,.
\end{equation}
Write ${\rm I}(\omega)=\{{\rm I}_{i,j}(\omega)\}_{p\times p}$, ${\rm II}(\omega)=\{{\rm II}_{i,j}(\omega)\}_{p\times p}$, ${\rm III}(\omega)=\{{\rm III}_{i,j}(\omega)\}_{p\times p}$, ${\rm IV}(\omega)=\{{\rm IV}_{i,j}(\omega)\}_{p\times p}$ and ${\rm V}(\omega)=\{{\rm V}_{i,j}(\omega)\}_{p\times p}$.

\subsection{Auxillary lemmas}\label{sec:auxlemma}
To construct Proposition  \ref{tm:1}, we  present the following useful lemmas whose proofs are given in Section \ref{sec:pflems}.
\begin{lemma}\label{la:remainder}
	Assume $r\geq n^\kappa$ for some sufficiently small constant $\kappa>0$. Under Conditions {\rm\ref{as:moment}} and {\rm\ref{as:betamixing}}, if $\log r=O(n^{1/2})$ and $l_n=o(n)$, it holds that
	\begin{align*}
		&\sup_{ \omega\in[-\pi,\pi] }\max_{(i,j)\in\cI}\{|{\rm II}_{i,j}(\omega)|+|{\rm III}_{i,j}(\omega)|+|{\rm IV}_{i,j}(\omega)|+|{\rm V}_{i,j}(\omega)|\}\\
		&~~~~~~~~~~~~~~~~~\lesssim \exp(-Cl_n)+O_{\rm  p}(l_nn^{-1}\log r)\,.
	\end{align*}
\end{lemma}

\begin{lemma}\label{la:leading}
	Assume $r\geq n^\kappa$ for some sufficiently small constant $\kappa>0$.
	Under Condition {\rm\ref{as:moment}}, if $l_n=o(n)$, it  holds that
	\begin{align*}
		\sup_{ \omega\in[-\pi,\pi] }\max_{(i,j)\in\mathcal{I}}|{\rm I}_{i,j}(\omega)-\zeta_{i,j}(\omega)|=O_{\rm p}(l_n^2n^{-1}\log r)
	\end{align*}
	for $\zeta_{i,j}(\omega)$  defined as {\rm\eqref{eq:zetaij}}.
\end{lemma}

Let
\begin{equation}\label{eq:zijt}
	\begin{split}
		z_{i,j,t}^{(1)}(\omega)=&~ \frac{1}{l_n}\sum_{k=-l_n}^{l_n}\mathcal{W}\bigg(\frac{k}{l_n}\bigg)\{\mathring{x}_{i,t+k}\mathring{x}_{j,t}-\gamma_{i,j}(k)\}\cos(k\omega)\,, \\
		z_{i,j,t}^{(2)}(\omega)=&~ \frac{1}{l_n}\sum_{k=-l_n}^{l_n}\mathcal{W}\bigg(\frac{k}{l_n}\bigg)\{\mathring{x}_{i,t+k}\mathring{x}_{j,t}-\gamma_{i,j}(k)\}\sin(k\omega)\,.
	\end{split}
\end{equation}
\begin{lemma}\label{la:tailprob}
	Under Conditions {\rm \ref{as:moment}} and {\rm \ref{as:betamixing}}, if $l_n\geq 2$, it holds that
	\begin{align*}
		&\sup_{ \omega\in[-\pi,\pi] }\max_{i,j\in[p]}\mathbb{P}\bigg\{\bigg|\frac{l_n}{s} \sum_{t=l_n+s_1}^{l_n+s_2}z_{i,j,t}^{(1)}(\omega)\bigg|>u\bigg\} +\sup_{ \omega\in[-\pi,\pi]} \max_{i,j\in[p]}\mathbb{P}\bigg\{\bigg|\frac{l_n}{s} \sum_{t=l_n+s_1}^{l_n+s_2}z_{i,j,t}^{(2)}(\omega)\bigg|>u\bigg\}\\
		&~~~~~~~~~~~~
		\lesssim \exp\{-Csl_n^{-3}\log^{-2} (l_n)u^2\} +\exp\{-Cs^{1/3}l_n^{-2/3}\log^{-1/3} (l_n)u^{1/3}\}
	\end{align*}
	for any $u>0$ and $1\leq s_1<s_2\leq n-2l_n$, where $s=s_2-s_1$, $z_{i,j,t}^{(1)}(\omega)$ and $z_{i,j,t}^{(2)}(\omega)$ are specified in \eqref{eq:zijt}.
\end{lemma}

\begin{lemma}\label{la:uniform}
	Assume $r\geq n^\kappa$ for some sufficiently small constant $\kappa>0$.
	Under Conditions {\rm\ref{as:moment}} and {\rm\ref{as:betamixing}}, if $\log r=O(n^{1/5}l_n^{-1/5})$, $l_n\log l_n=o(n) $
	and  $l_n\geq \max(2,C\log n)$ for some sufficiently large constant $C>0$,
	it holds that
	$$\sup_{ \omega\in[-\pi,\pi] }\max_{(i,j)\in\mathcal{I}}|\hat{f}_{i,j}(\omega) - f_{i,j}(\omega)| =O_{\rm p}\{n^{-1/2}l_n^{3/2}(\log l_n)\log^{1/2}(r)\}\,. $$
\end{lemma}

\begin{lemma}\label{la:app1}
	Assume $r\geq n^\kappa$ for some sufficiently small constant $\kappa>0$.
	Under Conditions {\rm\ref{as:moment}} and {\rm\ref{as:betamixing}},
	if
	$\log r=O(n^{1/5}l_n^{-1/5})$, $l_n\log l_n=o(n) $
	and
	$ l_n\geq \max(2,C'\log n)$ for some sufficiently large constant $C'>0$, it holds that
	\begin{align*}
		&	\bP\bigg\{\bigg|\sup_{\omega\in\mathcal{J}}T_n(\omega;\mathcal{I}) -\sup_{\omega\in\mathcal{J}}\check{T}_n(\omega;\mathcal{I})\bigg|\geq u\bigg\} 		 \\
		&~~~~~~~~~~	 \lesssim rnl_n\exp\bigg(-\frac{Cn^{1/2}u^{1/6}}{l_n^{1/6}}\bigg)+rnl_n\exp\bigg\{-\frac{Cn^{1/2}u^{1/3}}{l_n^{1/2}\log^{1/3}(l_n)\log^{1/6}(r)}\bigg\} \\
		&~~~~~~~~~~~~~	+rl_n^2\exp\bigg(-\frac{Cn^{1/2}u^{1/2}}{l_n^{3/2}}\bigg) + rl_n^2\exp\bigg\{-\frac{Cn^{1/2}u}{l_n^{5/2}(\log l_n)\log^{1/2}(r)}\bigg\} + n^{-1}
	\end{align*}
	for $u\gg n^{-1/2}l_n^{1/2}(\log l_n)\log^{1/2}(r)$ and
	\begin{align*}
		\bigg|\sup_{\omega\in\mathcal{J}}T_n(\omega;\mathcal{I}) -\sup_{\omega\in\mathcal{J}}\check{T}_n(\omega;\mathcal{I})\bigg|=O_{\rm p} \{n^{-1/2}l_n^{5/2}(\log l_n) \log^{3/2}(r)\}
	\end{align*}
	for any $\cJ\subset[-\pi,\pi]$,
	where $T_n(\omega;\cI)$ and $\check{T}_n(\omega;\mathcal{I})$ are defined as {\rm\eqref{eq:Tn}} and {\rm\eqref{eq:checkt}}, respectively.
\end{lemma}

\begin{lemma}\label{la.cont}
	Let $\{X(t):t\in\bR\}$ be a real-valued centered Gaussian random process and define $d^2(s,t)=\bE\{|X(s)-X(t)|^2\}$.
	Given a fixed compact interval $I\subset\bR$,
	if $d(s,t)\leq c_1|s-t|^{\lambda}$ for all $s,t\in I$ with some universal constants $c_1>0$ and $\lambda\in(0,1]$,  then there exist positive and finite universal constants $u_0$ and $c_2$ such that for $u\geq u_0$,
	\begin{align*}
		\sup_{t_0\in I} \bP\bigg\{\sup_{|s|\leq a}|X(t_0+s)-X(t_0)|\geq ua^\lambda\bigg\} \leq  \exp(-c_2 u^2)   
	\end{align*}
	for all $a\in(0,1]$ satisfying $t_0-a\in I$ and $t_0+a\in I$.
\end{lemma}

\subsection{Proof of Proposition \ref{tm:1}(i)}\label{sec:pf:tm1:1}
Recall  $r=|\cI|$ and $\bchi(\cdot)=\{\chi_1(\cdot),\chi_2(\cdot)\}$ is a given bijective mapping from $[r]$ to $\cI$ such that for any $(i,j)\in\cI$ there exists a unique $\ell\in[r]$ satisfying $(i,j)=\bchi(\ell)$. Notice that
$\check{\bfeta}^{\ext}(\omega)=\big\{\check{\eta}^{\ext}_1(\omega),\ldots,  \check{\eta}^{\ext}_{2r}(\omega)\big\}^\T=\{\bI_r\otimes \bA(\omega) \} n^{-1/2}\sum_{t=1}^{\tilde{n}}\bc_{t}$,
where $\bc_t=(\bc_{1,t}^\T,\ldots,\bc_{r,t}^\T)^\T$ with $\bc_{\ell,t}=(2\pi)^{-1}\{\mathring{x}_{\chi_1(\ell),t}\mathring{x}_{\chi_2(\ell),t+l_n}-{\gamma}_{\bchi(\ell)}(-l_n),\ldots,
\mathring{x}_{\chi_1(\ell),t+2l_n}\mathring{x}_{\chi_2(\ell),t+l_n}-{\gamma}_{\bchi(\ell)}(l_n)\}^\T$, and
\begin{align*}
	\bA(\omega)=&~\frac{1}{\sqrt{l_n}}\left(
	\begin{array}{ccc}
		\cos(-l_n\omega) & \cdots & \cos(l_n\omega) \\
		-\sin(-l_n\omega) & \cdots & -\sin(l_n\omega) \\
	\end{array}
	\right)
	\diag\{\mathcal{W}(-l_n/l_n),\ldots,\mathcal{W}(l_n/l_n)\}\,.
\end{align*}
Since $\cJ=\{\omega_1,\ldots,\omega_K\}$, then for given $\omega_1,\ldots,\omega_K\in[-\pi,\pi)$, we define a $(2Kr)$-dimensional vector
$
\mathring{\bfeta}^{\ext}=(\mathring{\eta}^{\ext}_1,\ldots,\mathring{\eta}^{\ext}_{2Kr})^\T=[\{\check{\bfeta}^{\ext}(\omega_1)\}^\T,\ldots,\{\check{\bfeta}^{\ext}(\omega_K)\}^\T]^\T$.
Furthermore, we have
\begin{align*}
	\mathring{\bfeta}^{\ext}=\left(
	\begin{array}{c}
		\bI_r\otimes \bA(\omega_1)  \\
		\vdots \\
		\bI_r\otimes \bA(\omega_K)  \\
	\end{array}
	\right)
	\frac{1}{\sqrt{{n}}}\sum_{t=1}^{\tilde{n}}\bc_{t}=:\bH\bigg(\frac{1}{\sqrt{{n}}}\sum_{t=1}^{\tilde{n}}\bc_{t}\bigg)\,.
\end{align*}
Then $\max_{k\in[K]}\check{T}_n(\omega_k;\mathcal{I})=\max_{j\in[Kr]}(|\mathring{\eta}_{2j-1}^{{\rm ext}}|^2+|\mathring{\eta}_{2j}^{{\rm ext}}|^2)$. Let $\ba_t=\bH\bc_t=:(a_{1,t},\ldots,a_{2Kr,t})^{\T}$ and $\bs_{n,\ba}=\tilde{n}^{-1/2}\sum_{t=1}^{\tilde{n}}\ba_t$. Recall $\bXi={\rm Var}(\tilde{n}^{-1/2}\sum_{t=1}^{\tilde{n}}\bc_t)$. Let $\bs_{n,\by}=(s_{n,1},\ldots,s_{n,2Kr})^{\T}\sim \mathcal{N}(\bzero,\bSigma)$ with $\bSigma=\bH\bXi\bH^\T$.
Define
\begin{align}
	\tilde{\varrho}_n:=&~ \sup_{\bu\in\mathbb{R}^{2Kr},\nu\in[0,1]} \big|\mathbb{P}(\sqrt{\nu}\mathring{\bfeta}^{\ext}+\sqrt{1-\nu}\bs_{n,\by}\leq \bu)-\mathbb{P}(\bs_{n,\by}\leq \bu)\big| \,, \label{eq:tildevarrho}\\
	\check{\varrho}_n:=&~\sup_{u\geq 0}\bigg|\mathbb{P}\bigg\{\max_{k\in[K]} \check{T}_n(\omega_k;\mathcal{I})\leq u\bigg\}-\mathbb{P}\bigg\{\max_{j\in[Kr]}(|s_{n,2j-1}|^2+|s_{n,2j}|^2)\leq u\bigg\}\bigg| \,, \label{eq:checkvarrho}\\
	\varrho_n^*:=&~\sup_{u\geq 0}\bigg|\mathbb{P}\bigg\{\max_{k\in[K]}{T}_n(\omega_k;\mathcal{I})\leq u\bigg\} -\mathbb{P}\bigg\{\max_{j\in[Kr]}(|s_{n,2j-1}|^2+|s_{n,2j}|^2)\leq u\bigg\}\bigg|\,. \label{eq:starvarrho}
\end{align}
The proof of Proposition \ref{tm:1}(i) includes three steps:

Step 1. To show
$
\tilde{\varrho}_n\lesssim n^{-1/9}l_n\log^{2/3}(l_n)\log(Kr)    =o(1) $.

Step 2. To show
$
\check{\varrho}_n\lesssim n^{-1/9}l_n\log^{2/3}(l_n)\log(Kr)    =o(1)$ based on the  result of Step 1.

Step 3. To show
$
\varrho_n^* 	\lesssim n^{-1/9}l_n\log^{2/3}(l_n)\log(Kr)
=o(1)$ based on  Step 2 and Lemma \ref{la:app1}.

The proofs of these three steps are given in Sections \ref{subsec:s1}--\ref{subsec:s3}, respectively, which all require the restrictions $\log(Kr)\ll n^{1/9}l_n^{-1}\log^{-8/3}(l_n)$, $l_n\log^{8/3}(l_n)\ll n^{1/9}$ and $l_n\geq \max\{2,C'\log(Kr)\}$ for some sufficiently large constant $C'>0$. $\hfill\Box$

\subsubsection{Step 1: convergence rate of $\tilde{\varrho}_n$}\label{subsec:s1}
Define
\begin{align*}
	\varrho_n=\sup_{\bu\in\mathbb{R}^{2Kr},\nu\in[0,1]}\big|\mathbb{P}(\sqrt{\nu}\bs_{n,\ba}+\sqrt{1-\nu}\bs_{n,\by}\leq \bu)-\mathbb{P}(\bs_{n,\by}\leq \bu)\big|\,.
\end{align*}
As we will show later,
$\varrho_n \lesssim n^{-1/9}l_n\log^{2/3}(l_n)\log(Kr)$
provided that  $\log(Kr)\ll n^{1/9}l_n^{-1}\log^{-8/3}(l_n)$, $l_n\log^{8/3}(l_n)\ll n^{1/9}$ and $l_n\geq \max\{2,C'\log(Kr)\}$ for some sufficiently large constant $C'>0$.
Recall  $\mathring{\bfeta}^{\ext}=\bH(n^{-1/2}\sum_{t=1}^{\tilde n}\bc_t)=n^{-1/2}\tilde{n}^{1/2}\bs_{n,\ba}$.
For some   $D>0$, define an event $\mathcal{E}_D=\{|\mathring{\bfeta}^{\ext}-\bs_{n,\ba}|_\infty\leq D\}$.
Parallel to Equation \eqref{eq:bb2} in Section \ref{sec:pflem.pn2} for the proof of Lemma \ref{pn:2}, 
then $\tilde{\varrho}_n$ defined as \eqref{eq:tildevarrho} satisfies
\begin{align}\label{eq:tilderhon}
	\tilde{\varrho}_n \lesssim \varrho_n+D\log^{1/2}(Kr)+\mathbb{P}(\mathcal{E}_D^{\rm c}) \
	\lesssim  \frac{l_n\log^{2/3}(l_n)\log(Kr)}{n^{1/9}} +D\log^{1/2}(Kr)+ \mathbb{P}(\mathcal{E}_D^{\rm c}) \,.
\end{align}
Recall $\ba_t=(a_{1,t},\ldots,a_{2Kr,t})^{\T}$. With selecting $D=C''n^{-1/9}l_n\log^{2/3}(l_n)\log^{1/2}(Kr)$,
if $l_n\geq \max\{2,C'\log(Kr)\}$ for some sufficiently large constant $C', C''>0$, by Lemma  \ref{la:tailprob} with $s=\tilde{n}$,
\begin{align*}
	\mathbb{P}(\mathcal{E}_D^{\rm c})
	\leq&~ \sum_{j=1}^{2Kr}\mathbb{P}\bigg(\frac{\sqrt{n}-\sqrt{\tilde{n}}}{\sqrt{n\tilde{n}}} \bigg|\sum_{t=1}^{\tilde{n}}a_{j,t}\bigg|>D\bigg)   \\
	\lesssim&~ Kr\exp\bigg\{-\frac{Cn^{16/9}\log(Kr)}{l_n^{2}\log^{2/3}(l_n)}\bigg\} + Kr\exp\bigg\{-\frac{Cn^{25/54}\log^{1/6}(Kr)}{l_n^{1/2}\log^{1/9} (l_n)}\bigg\}  \\
	\lesssim &~ \frac{l_n\log^{2/3}(l_n)\log(Kr)}{n^{1/9}}
\end{align*}
provided that $\log(Kr)\ll n^{1/9}l_n^{-1}\log^{-8/3}(l_n)$  and $l_n\log^{8/3}(l_n)\ll n^{1/9}$.
Together with  \eqref{eq:tilderhon},
\begin{align}\label{eq:tildevarrho.tail}
	\tilde{\varrho}_n\lesssim \frac{l_n\log^{2/3}(l_n)\log(Kr)}{n^{1/9}}    =o(1) \,.
\end{align}

Now, we show
$\varrho_n \lesssim n^{-1/9}l_n\log^{2/3}(l_n)\log(Kr)$.
Let $B=o(\tilde{n}^{1/2})$ be a positive integer that will diverge with $\tilde{n}$. We first decompose the sequence $[\tilde{n}]$ to the following $L+1$ blocks with $L=\lfloor \tilde{n}/B\rfloor$: $\mathcal{G}_\ell=\{(\ell-1)B+1,\ldots,\ell B\}$ for $\ell\in[L]$ and $\mathcal{G}_{L+1}=\{LB+1,\ldots,\tilde{n}\}$, where $\lfloor \cdot \rfloor$ is the integer truncation operator.
Let $b> h$ be two nonnegative integers such that $B=b+h$, $h=o(b)$ and $h>2l_n$. We then decompose each $\mathcal{G}_{\ell}$ $(\ell\in[L])$ to a `large' block $\cI_\ell$ with length $b$ and a `small' block $\cJ_\ell$ with length $h$. More specifically, $\mathcal{I}_{\ell}=\{(\ell-1)B+1,\ldots,(\ell-1)B+b\}$ and $\mathcal{J}_{\ell}=\{(\ell-1)B+b+1,\ldots,\ell B\}$ for any $\ell\in[L]$, and
$\mathcal{J}_{L+1}=\mathcal{G}_{L+1}$. Define $\tilde{\ba}_\ell=b^{-1/2}\sum_{t\in\mathcal{I}_\ell}\ba_t$ and $\check{\ba}_\ell=h^{-1/2}\sum_{t\in\mathcal{J}_\ell}\ba_t$ for each $\ell\in[L]$, and $\check{\ba}_{L+1}=(\tilde{n}-LB)^{-1/2}\sum_{t\in\mathcal{J}_{L+1}}\ba_t$.
Let $\{ \by_t\}_{t=1}^{\tilde{n}}$ be a sequence of independent normal random vectors with mean zero, where the covariance of $\by_t$ $(t\in\cI_\ell)$ is $\bE(\tilde \ba_\ell \tilde \ba_\ell^\T)$. For each $\ell\in[L]$, define $\tilde \by_\ell=b^{-1/2}\sum_{t\in\cI_\ell} \by_t$.
Write $\bs_{n,\ba}^{(1)}=L^{-1/2}\sum_{\ell=1}^L\tilde{\ba}_\ell$ and $\bs_{n,\by}^{(1)}=L^{-1/2}\sum_{\ell=1}^L\tilde{\by}_\ell$. Define
\begin{align*} 
	\varrho_n^{(1)}:=&~\sup_{\bu\in\mathbb{R}^{2Kr},\nu\in[0,1]}|\mathbb{P}\{\sqrt{\nu}\bs_{n,\ba}^{(1)}+\sqrt{1-\nu}\bs_{n,\by}^{(1)}\leq \bu\}-\mathbb{P}\{\bs_{n,\by}^{(1)}\leq \bu\}|\,,  \\
	\varrho_n^{(2)}:=&~\sup_{\bu\in\mathbb{R}^{2Kr},\nu\in[0,1]}|\mathbb{P}\{\sqrt{\nu}\bs_{n,\ba}+\sqrt{1-\nu}\bs_{n,\by}^{(1)}\leq \bu\}-\mathbb{P}\{\bs_{n,\by}^{(1)}\leq \bu\}|\,. 
\end{align*}
Similarly we  also define $\tilde{\bc}_\ell=b^{-1/2}\sum_{t\in\cI_{\ell}}\bc_t$ and $\check{\bc}_\ell=h^{-1/2}\sum_{t\in\cJ_t}\bc_t$ for each $\ell\in[L]$, and $\check{\bc}_{L+1}=(\tilde{n}-LB)^{-1/2}\sum_{t\in\cJ_{L+1}}\bc_t$.
Write $\tilde{\ba}_\ell=(\tilde{a}_{1,\ell},\ldots,\tilde{a}_{2Kr,\ell})^\T$ and $\tilde{\by}_\ell=(\tilde{y}_{1,\ell},\ldots,\tilde{y}_{2Kr,\ell})^\T$. Let $\widetilde{\bSigma}=L^{-1}\sum_{\ell=1}^L\mathbb{E}(\tilde{\ba}_\ell\tilde{\ba}_\ell^\T)$. Since $\ba_t=\bH\bc_t$, then $\widetilde\bSigma=\bH\widetilde{\bXi}\bH^{\T}$ with $\widetilde{\bXi}=L^{-1}\sum_{\ell=1}^{L}\bE(\tilde\bc_\ell\tilde{\bc}_\ell^{\T})$.
The proof for  Proposition \ref{tm:1}(i) also needs the following lemmas  whose proofs are given in Sections \ref{sec:pflem.xtildetail}--\ref{sec:pflem.pn2}, respectively.

\begin{lemma}\label{la:xtildetail}
	Under Conditions {\rm\ref{as:moment}} and {\rm\ref{as:betamixing}}, if $l_n\geq 2$, it holds that
	\begin{align*}
		\max_{\ell\in [L]}\max_{j\in[2Kr]}\mathbb{P}(|\tilde{a}_{j,\ell}|>\lambda)
		\lesssim \exp\{-Cl_n^{-2}\log^{-2}(l_n)\lambda^2\} +\exp\{-Cl_n^{-1/2}\log^{-1/3}(l_n)b^{1/6}\lambda^{1/3}\}
	\end{align*}
	for any $\lambda>0$.
\end{lemma}

\begin{lemma}\label{la:sigmatilde}
	Under Conditions {\rm\ref{as:moment}} and {\rm\ref{as:betamixing}}, it holds that $|\widetilde{\bSigma}-\bSigma|_\infty\lesssim l_n^2(hb^{-1}+bn^{-1})$.
\end{lemma}

\begin{lemma}\label{pn:1}
	Assume $r\geq n^\kappa$ for some sufficiently small constant $\kappa>0$. Let $h=2l_n+C\log(Kr)$ for some sufficiently large constant $C>0$. Under Conditions {\rm\ref{as:moment}}--{\rm\ref{as:eigen}}, if $l_n\geq 2$, it holds that
	\begin{align*}
		\varrho_n^{(1)}\lesssim L^{-1/6}l_n(\log l_n) \log^{7/6}(Kr)
	\end{align*}
	provided that $\log(Kr)\ll \min(b^{3/20}L^{1/10}l_n^{-3/20}, L^{2/5})$ and $l_n^2 (hb^{-1}+bn^{-1})=o(1)$.
\end{lemma}

\begin{lemma}\label{pn:2}
	Assume $r\geq n^\kappa$ for some sufficiently small constant $\kappa>0$. Let $\min(n^{1/2}, \, nl_n^{-2})\gg b\gg
	\max\{hl_n^2, \, h^{3/4}n^{1/4}\log^{-1/4}(Kr)\}$ and $h=2l_n+C\log(Kr)$ for some sufficiently large constant $C>0$. Under Conditions {\rm\ref{as:moment}}--{\rm\ref{as:eigen}}, if $l_n\geq 2$, it holds that
	\begin{align*}
		\varrho_n^{(2)}\lesssim L^{-1/6}l_n(\log l_n)\log^{7/6}(Kr)
	\end{align*}
	provided that $\log(Kr)\ll\min(b^{3/20}L^{1/10}l_n^{-3/20}, L^{2/5})$.
\end{lemma}

By the  triangle inequality, 
\begin{align*}
	&\big|\mathbb{P}(\sqrt{\nu}\bs_{n,\ba}+\sqrt{1-\nu}\bs_{n,\by}\leq \bu)-\mathbb{P}(\bs_{n,\by}\leq \bu)\big|\\
	&~~~~~~~\leq\big|\mathbb{P}(\sqrt{\nu}\bs_{n,\ba}+\sqrt{1-\nu}\bs_{n,\by}\leq \bu)-\mathbb{P}\{\sqrt{\nu}\bs_{n,\ba}+\sqrt{1-\nu}\bs_{n,\by}^{(1)}\leq \bu\}\big|\\
	&~~~~~~~~~~~+\big|\mathbb{P}\{\bs_{n,\by}^{(1)}\leq \bu\}-\mathbb{P}(\bs_{n,\by}\leq \bu)\big|+\big|\mathbb{P}\{\sqrt{\nu}\bs_{n,\ba}+\sqrt{1-\nu}\bs_{n,\by}^{(1)}\leq \bu\}-\mathbb{P}\{\bs_{n,\by}^{(1)}\leq \bu\}\big|\,,
\end{align*}
which implies
$
\varrho_n\leq \varrho_n^{(2)}+2\sup_{\bu\in\mathbb{R}^{2Kr}} |\mathbb{P}\{\bs_{n,\by}^{(1)}\leq \bu\}-\mathbb{P}(\bs_{n,\by}\leq \bu)| $.
Recall $\bs_{n,\by}^{(1)}\sim \mathcal{N}(\bzero,\widetilde{\bSigma})$ and $\bs_{n,\by}\sim \mathcal{N}(\bzero,\bSigma)$.
By Lemma 1 of Chang et al. (2023a),
\begin{align*}
	\sup_{\bu\in\mathbb{R}^{2Kr}}\big|\mathbb{P}\{\bs_{n,\by}^{(1)}\leq \bu\}-\mathbb{P}(\bs_{n,\by}\leq \bu)\big|\lesssim |\widetilde{\bSigma}-\bSigma|_\infty^{1/3}\log^{2/3}(Kr)\,.
\end{align*}
If  $\min(n^{1/2},  nl_n^{-2})\gg b\gg
\max\{h^{3/4}n^{1/4}\log^{-1/4}(Kr), hl_n^2\}$ with $l_n\geq 2$, by Lemmas \ref{la:sigmatilde} and  \ref{pn:2},
\begin{align*}
	\varrho_n
	\lesssim \frac{b^{1/6}l_n(\log l_n)\log^{7/6}(Kr)}{n^{1/6}} +\frac{l_n^{2/3}h^{1/3}\log^{2/3}(Kr)}{b^{1/3}}
\end{align*}
provided that $\log(Kr)\ll\min(n^{1/10}b^{1/20}l_n^{-3/20}, L^{2/5})$ and  $h=2l_n+C'''\log(Kr)$ for some sufficiently large constant $C'''>0$.
Letting $ l_n \gtrsim \log(Kr)$, then $h\asymp l_n$, which implies
\begin{align*}
	\varrho_n
	\lesssim \frac{b^{1/6}l_n(\log l_n)\log^{7/6}(Kr)}{n^{1/6}} + \frac{l_n\log^{2/3}(Kr)}{b^{1/3}}  \,.
\end{align*}
Selecting $b\asymp n^{1/3}\log^{-2}(l_n)\log^{-1}(Kr)$, if  $\log(Kr)\ll n^{1/9}l_n^{-1}\log^{-8/3}(l_n)$, $\log(Kr)\lesssim l_n$ and $l_n\log^{8/3}(l_n)\ll n^{1/9}$, then it holds that
\begin{align*} 
	\varrho_n \lesssim \frac{l_n\log^{2/3}(l_n)\log(Kr)}{n^{1/9}}=o(1)\,.
\end{align*}

\subsubsection{Step 2: convergence rate of $\check{\varrho}_n$}\label{subsec:s2}
Recall $\mathring{\bfeta}^{\ext}=(\mathring{\eta}_{1}^{\ext},\ldots, \mathring{\eta}_{2Kr}^{\ext})^{\T}$ 
and $\{\max_{k\in[K]}\check{T}_n(\omega_k;\cI)\leq u\} =\{\max_{j\in[Kr]}(|\mathring{\eta}_{2j-1}^{\ext}|^2+|\mathring{\eta}_{2j}^{\ext}|^2)\leq u\}$. For any $j\in[Kr]$ and $u>0$, define
\begin{align}\label{eq:Aju}
	A_j(u)=\{\bb=(b_1,\ldots,b_{2Kr})^{\T}\in\bR^{2Kr}:b_{2j-1}^2+b_{2j}^2\leq u\} \,.
\end{align}
Let $A(u)=\bigcap_{j=1}^{Kr}A_j(u)$, which is a 2-sparsely convex set. See Definition 3.1 of \cite{CCK_2017_supp} 
and Definition 3 of 
\cite{CCW_2021_supp} 
for the definition of $s$-sparsely convex set.
Recall $\bs_{n,\by}=(s_{n,1},\ldots,s_{n,2Kr})^{\T}$.
Then $\{\max_{k\in[K]}\check{T}_n(\omega_k;\cI)\leq u\}=\{\mathring{\bfeta}^{\ext}\in A(u)\}$ and $\{\max_{j\in[Kr]}(|s_{n,2j-1}|^2+|s_{n,2j}|^2)\leq u\}= \{\bs_{n,\by}\in A(u)\}$.
We can reformulate $\check{\varrho}_n$ defined in \eqref{eq:checkvarrho} as follows:
\begin{align*}
	\check{\varrho}_n
	= \sup_{u\geq 0}\big|\bP\{\mathring{\bfeta}^{\ext}\in A(u)\}-\bP\{\bs_{n,\by}\in A(u)\}\big| \,.
\end{align*}
Assume $l_n\log^{8/3}(l_n)\log(Kr)= o(n^{1/9})$. In the sequel, we will consider the convergence rate of $|\bP\{\mathring{\bfeta}^{\ext}\in A(u)\}-\bP\{\bs_{n,\by}\in A(u)\}|$ in two scenarios: (i) $0\leq u\leq n^{-1}$ and (ii) $u>n^{-1}$.

\underline{{\it Scenario 1: } $0\leq u\leq n^{-1}$.} 
By the triangle inequality, it holds that
\begin{align*}
	&\sup_{0\leq u\leq n^{-1}}\big|\bP\{\mathring{\bfeta}^{\ext}\in A(u)\}-\bP\{\bs_{n,\by}\in A(u)\}\big| \\
	&~~~~~~~~~~ \leq \sup_{0\leq u\leq n^{-1}}\bP\bigg\{\max_{k\in[K]} \check{T}_n(\omega_k;\mathcal{I})\leq u\bigg\}
	+ \sup_{0\leq u\leq n^{-1}}\bP\bigg\{\max_{j\in[Kr]}(|s_{n,2j-1}|^2+|s_{n,2j}|^2)\leq u\bigg\} \,.
\end{align*}
For any $0\leq u\leq n^{-1}$, Nazarov's inequality yields that
\begin{align*}
	\bP\bigg\{\max_{j\in[Kr]}(|s_{n,2j-1}|^2+|s_{n,2j}|^2)\leq u\bigg\}
	\leq&~ \bP\bigg\{\max_{j\in[Kr]}(|s_{n,2j-1}|^2+|s_{n,2j}|^2)\leq n^{-1}\bigg\}  \\
	\leq&~ \bP\bigg(\max_{j\in[2Kr]}|s_{n,j}|\leq n^{-1/2}\bigg)\lesssim n^{-1/2}\log^{1/2}(Kr) \,.
\end{align*}
Write  ${\bf 1}_{2Kr}$ as a  $(2Kr)$-dimensional vector whose elements are all equal to 1.
By the triangle inequality, Nazarov's inequality and \eqref{eq:tildevarrho.tail}, for any $0\leq u\leq n^{-1}$,
\begin{align*}
	\bP\bigg\{\max_{k\in[K]} \check{T}_n(\omega_k;\mathcal{I})\leq u\bigg\}
	\leq&~ \bP\bigg(\max_{j\in[2Kr]}|\mathring{\eta}_j^{\ext}|\leq n^{-1/2}\bigg) = \bP\big(-n^{-1/2}{\bf 1}_{2Kr}\leq \mathring{\bfeta}^{\ext}\leq n^{-1/2}{\bf 1}_{2Kr}\big) \\
	\leq&~ \bP\big(\mathring{\bfeta}^{\ext}\leq n^{-1/2}{\bf 1}_{2Kr}\big) - \bP\big(\mathring{\bfeta}^{\ext}\leq -n^{-1/2}{\bf 1}_{2Kr}\big) \\
	\leq&~ 2\tilde{\varrho}_n + \bP\big(\bs_{n,\by}\leq n^{-1/2}{\bf 1}_{2Kr}\big) - \bP\big(\bs_{n,\by}\leq -n^{-1/2}{\bf 1}_{2Kr}\big)\\
	\lesssim&~ n^{-1/9}l_n\log^{2/3}(l_n)\log(Kr)+n^{-1/2}\log^{1/2}(Kr) \\
	\lesssim&~ n^{-1/9}l_n\log^{2/3}(l_n)\log(Kr)
\end{align*}
provided that  $\log(Kr)\ll n^{1/9}l_n^{-1}\log^{-8/3}(l_n)$, $l_n\log^{8/3}(l_n)\ll n^{1/9}$ and $l_n\geq \max\{2,C'\log(Kr)\}$ for some sufficiently large constant $C'>0$, where $\tilde{\varrho}_n$ is defined in \eqref{eq:tildevarrho}.
Hence,
\begin{align}\label{eq:diff.smallu}
	\sup_{0\leq u\leq n^{-1}}\big|\bP\{\mathring{\bfeta}^{\ext}\in A(u)\}-\bP\{\bs_{n,\by}\in A(u)\}\big|  
	\lesssim n^{-1/9}l_n\log^{2/3}(l_n)\log(Kr) \,.
\end{align}

\underline{{\it Scenario 2: }  $u>n^{-1}$.}
Let $B=\{\bb=(b_1,\ldots,b_{2Kr})^{\T}\in\bR^{2Kr}:\max_{j\in[2Kr]}|b_j|\leq 2Krn^{5/2}\}$.  Recall  $\mathring{\bfeta}^{\ext}=n^{-1/2}\sum_{t=1}^{\tilde n}\ba_t=n^{-1/2}\sum_{t=1}^{\tilde n}\bH\bc_t$ and $\ba_t=(a_{1,t},\ldots,a_{2Kr,t})^{\T}$. By Markov inequality,
\begin{align*}
	\bP\{\mathring{\bfeta}^{\ext}\in A(u)\cap B^{\rm c}\}
	\leq&~ \bP(\mathring{\bfeta}^{\ext}\in B^{\rm c})=\bP\bigg(\max_{j\in[2Kr]}\bigg|\frac{1}{\sqrt{n}}\sum_{t=1}^{\tilde n} a_{j,t} \bigg|> 2Krn^{5/2}\bigg)  \\
	\lesssim&~ (Kr)^{-1}n^{-2}\bE\bigg(\max_{t\in[\tilde n]}\max_{j\in[2Kr]}|a_{j,t}|\bigg)
	\lesssim n^{-1}\max_{t\in[\tilde n]}\max_{j\in[2Kr]}\bE(|a_{j,t}|)
\end{align*}
for any $u>0$. Similar to \eqref{eq:zijtail} in Section \ref{subsec:pfLem3} for the proof of Lemma \ref{la:tailprob}, we have
\begin{align*} 
	\max_{t\in[\tilde n]}\max_{j\in[2Kr]}\bP(|a_{j,t}|>u)\lesssim l_n\exp(-Cl_n^{-1/2}u)
\end{align*}
for any $u>0$, which implies  $\max_{t\in[\tilde n]}\max_{j\in[2Kr]}\bE(|a_{j,t}|)\lesssim l_n^{3/2}$ and $\max_{t\in[\tilde n]}\max_{j\in[2Kr]}\bE(|a_{j,t}|^2)\lesssim l_n^{2}$. Hence,
$\bP\{\mathring{\bfeta}^{\ext}\in A(u)\cap B^{\rm c}\}\lesssim n^{-1}l_n^{3/2}\lesssim n^{-1/2}$ for any $u>0$ provided that $l_n\ll n^{1/3}$.
Let $\{{\bff}_t\}_{t=1}^{\tilde{n}}$, independent of $\{\ba_t\}_{t=1}^{\tilde{n}}$, be a centered Gaussian sequence such that $\cov(\bff_t,\bff_s)=\cov(\ba_t,\ba_s)$ for all $t,s\in[\tilde{n}]$. Then $\tilde{n}^{-1/2}\sum_{t=1}^{\tilde{n}}\bff_t\overset{d}{=}\bs_{n,\by}$. Write $\bff_t=(f_{1,t},\ldots,f_{2Kr,t})^{\T}$. By Markov inequality again, if $l_n\ll n^{1/2}$, 
\begin{align*}
	\bP\{\bs_{n,\by}\in A(u)\cap B^{\rm c}\}
	\leq&~ \bP(\bs_{n,\by}\in B^{\rm c})
	=\bP\bigg(\max_{j\in[2Kr]}\bigg|\frac{1}{\sqrt{\tilde{n}}} \sum_{t=1}^{\tilde{n}}f_{j,t}\bigg|> 2Krn^{5/2}\bigg) \\
	\lesssim&~ n^{-1}\max_{t\in[\tilde{n}]}\max_{j\in[2Kr]}\bE(|f_{j,t}|) \leq n^{-1}\max_{t\in[\tilde n]}\max_{j\in[2Kr]} \big\{\bE(|f_{j,t}|^2)\big\}^{1/2}\\
	=&~ n^{-1}\max_{t\in[\tilde n]}\max_{j\in[2Kr]} \big\{\bE(|a_{j,t}|^2)\big\}^{1/2}\lesssim n^{-1}l_n\lesssim n^{-1/2}
\end{align*}
for any $u>0$.
Therefore, if $l_n\ll n^{1/3}$,
\begin{align}\label{eq:diff.largeu}
	&  \sup_{u>n^{-1}} \big|\bP\{\mathring{\bfeta}^{\ext}\in A(u)\}-\bP\{\bs_{n,\by}\in A(u)\}\big|  \notag\\
	&~~~~~~~~~~~ \leq\sup_{u> n^{-1}} \big|\bP\{\mathring{\bfeta}^{\ext}\in A(u)\cap B\}-\bP\{\bs_{n,\by}\in A(u)\cap B\}\big|+Cn^{-1/2}  \,.
\end{align}

In the sequel, we will bound $\sup_{u>n^{-1}}|\bP\{\mathring{\bfeta}^{\ext}\in A(u)\cap B\}-\bP\{\bs_{n,\by}\in A(u)\cap B\}|$.
Let $A^*(u)=A(u)\cap B$ and  $A_j^*(u)=\{\bb=(b_1,\ldots,b_{2Kr})^{\T}: |b_{2j-1}|\leq 2Krn^{5/2}, |b_{2j}|\leq 2Krn^{5/2}, b_{2j-1}^2+b_{2j}^{2}\leq u\}$. Then  $A^*(u)=\bigcap_{j=1}^{Kr}A_j^*(u)$.
For any given $u>n^{-1}$, to simplify the notation, we write $A^{*}(u)$ and $A^{*}_j(u)$ as $A^{*}$ and $A^{*}_j$, respectively. Notice that
$A^*$ contains a ball with radius $\epsilon=n^{-1}$ and center at $\bzero\in A^*$.
Following the identical arguments for Case 1 in the proof of Theorem 7 of 
\cite{CCW_2021_supp} 
with $w^*=\bzero$,  $s_q=2=s$ and $p=2Kr$, then $A^*\in \cA^{\rm si}(1,\tilde{C})$  is a simple convex set for some positive constant $\tilde{C}$,
and the $K_*$-generated set
(generated by the intersection of $K_*$ half-spaces)
for $A^*$ denoted by $A^{K_*}$ satisfies $K_*\leq (2Krn)^C$ and $\max_{\bv\in\cV(A^{K_*})}|\bv|_0\leq 2$,
where $\cV(A^{K_*})$ denotes the set that consists of $K_*$ unit normal vectors outward to the facets of $A^{K_*}$.
See  Definition 2 of 
\cite{CCW_2021_supp} 
for the definition of simple convex set. Define
\begin{align*}
	\rho^{(K_*)} =&~ \big|\bP(\mathring{\bfeta}^{\ext}\in A^{K_*} ) - \bP(\bs_{n,\by}\in A^{K_*})\big| \,,  \\
	\rho^{(K_*),\epsilon} =&~ \big|\bP(\mathring{\bfeta}^{\ext}\in A^{K_*,\epsilon})  - \bP(\bs_{n,\by}\in A^{K_*,\epsilon})\big| \,,
\end{align*}
where
\begin{align*}
	A^{K_*} =&~ \bigcap_{\bv\in\cV(A^{K_*})}\{\bw\in\bR^{2Kr}:\bw^{\T}\bv\leq\cS_{A^{K_*}}(\bv)\} \,,  \\
	A^{K_*,\epsilon}=&~
	\bigcap_{\bv\in\cV(A^{K_*})}\{\bw\in\bR^{2Kr}:\bw^{\T}\bv\leq \cS_{A^{K_*}}(\bv)+\epsilon\} \,,
\end{align*}
with $\cS_{A^{K_*}}(\bv)=\sup\{\bw^{\T}\bv:\bw\in A^{K_*}\}$.
Write $\cV(A^{K_*})=\{\bv_1,\ldots,\bv_{K_*}\}$
and define $\bs_{n,\by}^{\rm si}:=(s_{n,1}^{\rm si}, \ldots, s_{n,K_*}^{\rm si})^{\T} =(\bv_1^{\T}\bs_{n,\by},\ldots,\bv_{K_*}^{\T}\bs_{n,\by})^{\T}$.
By Condition \ref{as:eigen},  we have
\begin{align*}
	\min_{i\in[K_*]}\var(s_{n,i}^{\rm si}) = \min_{\bv\in\cV(A^{K_*})}\var\bigg(\frac{1}{\sqrt{\tilde{n}}}\sum_{t=1}^{\tilde{n}}\bv^{\T}\ba_t\bigg) \geq C\,.
\end{align*}
Parallel to (S.30) in the supplementary material of \cite{CCW_2021_supp} 
with $\epsilon=n^{-1}$, $K=K_*$ and $A=A^*$, we have
\begin{align} \label{eq:gaussAppSi}
	\big|\bP(\mathring{\bfeta}^{\ext}\in A^{*} ) - \bP(\bs_{n,\by}\in A^{*})\big|
	\lesssim n^{-1}\log^{1/2}(Kr)+\rho^{(K_*)} + \rho^{(K_*),\epsilon} \,.
\end{align}	
Let $\ba^{\rm si}_t :=(a_{1,t}^{\rm si}, \ldots, a_{K_*,t}^{\rm si})^{\T}=(\bv_1^{\T}\ba_t, \ldots, \bv_{K_*}^{\T}\ba_t)^{\T}$. 
Notice that the sequence $\{\ba_t^{\rm si}\}_{t=1}^{\tilde{n}}$ is an $\alpha$-mixing sequence with $\alpha$-mixing coefficients $\alpha_{\ba^{\rm si}}(k)\lesssim \exp(-C|k-2l_n|_+)$ for any integer $k$,
$\max_{t\in[\tilde{n}]}\max_{j\in[K_*]}\bP(|a_{j,t}^{\rm si}|>x)\lesssim l_n\exp(-Cl_n^{-1/2}x)$ for any $x>0$,
and $\cov(\bs_{n,\by}^{\rm si})=\cov(\tilde{n}^{-1/2}\sum_{t=1}^{\tilde{n}} \ba_t^{\rm si})$.
Write $\cS_{A^{K_*}}=\{\cS_{A^{K_*}}(\bv_1),\ldots, \cS_{A^{K_*}}(\bv_{K_*})\}^{\T}\in\bR^{K_*}$.
Recall we have shown in Section  \ref{subsec:s1} that
\begin{align*}
	\sup_{\bu\in\bR^{2Kr}}\bigg| \bP\bigg(
	\frac{1}{\sqrt{n}}\sum_{t=1}^{\tilde{n}}\ba_t
	\leq\bu\bigg) - \bP(\bs_{n,\by}\leq \bu)\bigg|
	\lesssim n^{-1/9}l_n\log^{2/3}(l_n)\log(Kr) \,.
\end{align*}
Due to $r\geq n^{\kappa}$ for some sufficiently small constant $\kappa>0$, applying the same arguments  to $\{\ba_t^{\rm si}\}_{t=1}^{\tilde n}$, we can also show
\begin{align*}
	\rho^{(K_*)}  = \bigg| \bP\bigg(
	\frac{1}{\sqrt{n}}\sum_{t=1}^{\tilde{n}}\ba_t^{\rm si}
	\leq\cS_{A^{K_*}}\bigg) - \bP\big(\bs_{n,\by}^{\rm si}\leq \cS_{A^{K_*}}\big)\bigg|
	\lesssim n^{-1/9}l_n\log^{2/3}(l_n)\log(Kr)
\end{align*}
provided that  $\log(Kr)\ll n^{1/9}l_n^{-1}\log^{-8/3}(l_n)$, $l_n\log^{8/3}(l_n)\ll n^{1/9}$ and $l_n\geq \max\{2,C'\log(Kr)\}$ for some sufficiently large constant $C'>0$.
Analogously, we also have
$
\rho^{(K_*),\epsilon} \lesssim n^{-1/9}l_n\log^{2/3}(l_n)\\ \log(Kr) $.
Together with \eqref{eq:gaussAppSi},
$|\bP(\mathring{\bfeta}^{\ext}\in A^{*} ) - \bP(\bs_{n,\by}\in A^{*})|
\lesssim n^{-1/9}l_n\log^{2/3}(l_n)\log(Kr) $. Since the upper bound $n^{-1/9}l_n\log^{2/3}(l_n)\log(Kr) $ holds uniformly over $u>n^{-1}$, then
\begin{align*}
	\sup_{u> n^{-1}} \big|\bP\{\mathring{\bfeta}^{\ext}\in A(u)\cap B\}-\bP\{\bs_{n,\by}\in A(u)\cap B\}\big|
	\lesssim&~ n^{-1/9}l_n\log^{2/3}(l_n)\log(Kr)
\end{align*}
provided that  $\log(Kr)\ll n^{1/9}l_n^{-1}\log^{-8/3}(l_n)$, $l_n\log^{8/3}(l_n)\ll n^{1/9}$ and $l_n\geq \max\{2,C'\log(Kr)\}$ for some sufficiently large constant $C'>0$. Together with \eqref{eq:diff.smallu} and \eqref{eq:diff.largeu}, it holds that
\begin{align}\label{eq:checkrhon}
	\check{\varrho}_n	= \sup_{u\geq 0}\big|\bP\{\mathring{\bfeta}^{\ext}\in A(u)\}-\bP\{\bs_{n,\by}\in A(u)\}\big| \lesssim n^{-1/9}l_n\log^{2/3}(l_n)\log(Kr) \,.
\end{align}	

\subsubsection{Step 3: convergence rate of $\varrho_n^*$}\label{subsec:s3}
For any $\varepsilon>0$, it holds that
\begin{align}
	\varrho_n^* \leq&~\check{\varrho}_n +\mathbb{P}\bigg\{\bigg|\max_{k\in[K]}T_n(\omega_k;\mathcal{I})- \max_{k\in[K]}\check{T}_n(\omega_k;\mathcal{I})\bigg|>\varepsilon\bigg\}  \notag\\ &+\sup_{u\geq0}\mathbb{P}\bigg\{u-\varepsilon<\max_{j\in[Kr]}(|s_{n,2j-1}|^2+|s_{n,2j}|^2)\leq u+\varepsilon\bigg\} \,,  \label{eq:rhonstar}
\end{align}
where $\check{\varrho}_n$ and $\varrho_n^{*}$ are defined in \eqref{eq:checkvarrho} and \eqref{eq:starvarrho}, respectively.
Notice that $r\geq n^{\kappa}$ for some sufficiently small constant $\kappa>0$.
Selecting $\varepsilon = C''''n^{-1/2}l_n^{3}\log^{2}(Kr)$ for some sufficiently large constant $C''''>0$,
Lemma \ref{la:app1} implies
\begin{align}\label{eq:diff.Tn.checkTn}
	\mathbb{P}\bigg\{\bigg|\max_{k\in[K]}T_n(\omega_k;\mathcal{I})-\max_{k\in[K]}\check{T}_n(\omega_k;\mathcal{I})\bigg|>\varepsilon\bigg\}
	\lesssim n^{-1}
	=o(1)
\end{align}
provided that $\log r=O(n^{1/5}l_n^{-1/5})$, $l_n\log l_n=o(n) $
and
$ l_n\geq \max(2,\tilde{C}\log n)$ for some sufficiently large constant $\tilde{C}>0$.
For any $u>0$, we  reformulate $A_j(u)$ defined in \eqref{eq:Aju} as follows:
\begin{align*}
	{A}_j(u) = \bigcap_{\bfd\in\{\bfd\in\mathbb{S}^{2Kr-1}:\,\bfd_{S_j}\in\mathbb{S}^1\}} \{\bb\in\mathbb{R}^{2Kr}:\bfd^{\T}\bb\leq\sqrt{u}\}
\end{align*}
with $S_j=\{2j-1,2j\}$. Let $\mathcal{F}=\bigcup_{j=1}^{Kr}\{\bfd\in\mathbb{S}^{2Kr-1}:\bfd_{S_j}\in\mathbb{S}^1\}$. Recall $A(u)=\bigcap_{j=1}^{Kr}A_j(u)$. Then
\begin{align}\label{eq:Adef}
	{A}(u)=\bigcap_{\bfd\in\mathcal{F}}\{\bb\in\mathbb{R}^{2Kr}:\bfd^{\T}\bb\leq\sqrt{u}\}\,.
\end{align}
Following the arguments of Step 1 in the proof of Proposition 1 in 
\cite{CJS_2021_supp},
there exists an $\epsilon$-net of $\cF$, denoted by $\cF_\epsilon$, satisfying (i) $\cF_\epsilon\subset\cF$,
(ii) $\epsilon^{-2} \leqslant |\cF_\epsilon| \leqslant Kr\{(2+\epsilon)\epsilon^{-1}\}^2$, and (iii) $A_1(u;\mathcal{F}_\epsilon)\subset {A}(u)\subset A_2(u;\mathcal{F}_\epsilon)$
with $A_1(u;\mathcal{F}_\epsilon)= \bigcap_{\bfd\in\mathcal{F}_\epsilon} \{\bb\in\mathbb{R}^{2Kr}:\bfd^{\T}\bb\leq(1-\epsilon)\sqrt{u}\}$ and $A_2(u;\mathcal{F}_\epsilon)=\bigcap_{\bfd\in\mathcal{F}_\epsilon} \{\bb\in\mathbb{R}^{2Kr}:\bfd^{\T}\bb\leq\sqrt{u}\}$.
Recall $\varepsilon = C''''n^{-1/2}l_n^{3}\log^{2}(Kr)$ for some sufficiently large constant $C''''>0$.  On one hand,		
\begin{align}\label{eq:antis1}
	&\sup_{0\leq u\leq\varepsilon}\mathbb{P}\bigg\{u-\varepsilon<\max_{j\in[Kr]}(|s_{n,2j-1}|^2+|s_{n,2j}|^2)\leq u+\varepsilon\bigg\} \notag\\
	&~~~~~~~~~~~~~~~~ \leq \sup_{0\leq u\leq\varepsilon}\mathbb{P}\big(|s_{n,1}|^2 \leq u+\varepsilon\big)
	\lesssim \sqrt{\varepsilon}
	\lesssim \frac{l_n\log^{2/3}(l_n)\log(Kr)}{n^{1/9}}
\end{align}
provided that $l_n\ll n^{5/18}$, where the second step is based on the anti-concentration inequality of normal random variable. On the other hand, for any $u>\varepsilon$,
\begin{align*}
	&\mathbb{P}\bigg\{u-\varepsilon<\max_{j\in[Kr]}(|s_{n,2j-1}|^2+|s_{n,2j}|^2)\leq u+\varepsilon\bigg\}
	= \bP\{\bs_{n,\by}\in A(u+\varepsilon)\} - \bP\{\bs_{n,\by}\in A(u-\varepsilon)\}  \notag\\
	&~~~~~~~~~~~\leq \bP\{\bs_{n,\by}\in A_2(u+\varepsilon;\cF_\epsilon)\}  -  \bP\{\bs_{n,\by}\in A_1(u-\varepsilon;\cF_\epsilon)\} \notag\\
	&~~~~~~~~~~~ = \bP\bigg(\max_{\bfd\in\cF_\epsilon}\bfd^{\T}\bs_{n,\by}\leq \sqrt{u+\varepsilon}\bigg)-\bP\bigg\{\max_{\bfd\in\cF_\epsilon}\bfd^{\T}\bs_{n,\by}\leq(1-\epsilon)\sqrt{u-\varepsilon}\bigg\}  \,.
	\notag\\
	&~~~~~~~~~~~\leq  \mathbb{P}\bigg\{(1-\epsilon)\sqrt{u} < \max_{\bfd\in\mathcal{F}_{\epsilon}} \bfd^{\T}\bs_{n,\by}\leq \sqrt{u} + \sqrt{\varepsilon}\bigg\}  \notag\\
	&~~~~~~~~~~~~~~ + \mathbb{P}\bigg\{(1-\epsilon)(\sqrt{u}-\sqrt{\varepsilon})< \max_{\bfd\in\mathcal{F}_{\epsilon}}\bfd^{\T}\bs_{n,\by} \leq (1-\epsilon)\sqrt{u}\bigg\} \,.
\end{align*}
Write $\rho_{\bs}(u) = \bP\{(1-\epsilon)\sqrt{u}<\max_{\bfd\in\cF_\epsilon}\bfd^{\T}\bs_{n,\by} \leq \sqrt{u}\}$. 
Following  the same arguments  below Equation (28) of \cite{CJS_2021_supp}, 
we have
\begin{align} \label{eq:antis2}
	\sup_{u>\varepsilon}\mathbb{P}\bigg\{u-\varepsilon<\max_{j\in[Kr]}(|s_{n,2j-1}|^2+|s_{n,2j}|^2)\leq u+\varepsilon\bigg\}
	\lesssim \sqrt{\varepsilon\log(|\cF_\epsilon|)}
	+\sup_{u>\varepsilon}\rho_{\bs}(u) \,.
\end{align}
Selecting $\epsilon=n^{-1}$, then $|\mathcal{F}_\epsilon|\leq Kr(2n+1)^2$. Recall  $\varepsilon = C''''n^{-1/2}l_n^{3}\log^{2}(Kr)$ for some sufficiently large constant $C''''>0$ and $r\geq n^{\kappa}$ for some sufficiently small constant $\kappa>0$.
Then $ \sqrt{\varepsilon\log(|\cF_\epsilon|)}\lesssim n^{-1/4}l_n^{3/2}\log^{3/2}(Kr)\lesssim  n^{-1/9}l_n\log^{2/3}(l_n)\log(Kr)$ provided that $\log(Kr)\ll n^{5/18}l_n^{-1}\log^{4/3}(l_n)$  and $l_n\log^{-4/3}(l_n)\ll n^{5/18}$.
In addition,
if $u\leq nl_n^{3/2}$, by Nazarov's inequality,  $\rho_{\bs}(u)\lesssim \epsilon\sqrt{u\log(|\mathcal{F}_\epsilon|)}\lesssim n^{-1/2}l_n^{3/4}\log^{1/2}(Kr)\lesssim  n^{-1/9}l_n\log^{2/3}(l_n)\log(Kr)$.
If $u>nl_n^{3/2}$, by Markov inequality and Lemma 7.4 of \cite{FSZ_2018_supp}, 
\begin{align*}
	\rho_{\bs}(u)
	\leq&~  \mathbb{P}\bigg\{(1-\epsilon)\sqrt{u} < \max_{\bfd\in\mathcal{F}_{\epsilon}} \bfd^{\T}\bs_{n,\by}\bigg\}
	\leq\frac{\mathbb{E}(\max_{\bfd\in\mathcal{F}_\epsilon}|\bfd^{\T}\bs_{n,\by}|)}{(1-\epsilon)\sqrt{u}} \\
	\lesssim&~ \frac{l_n^{1/4}\log^{1/2}(Kr)}{n^{1/2}}
	\lesssim  \frac{l_n\log^{2/3}(l_n)\log(Kr)}{n^{1/9}} \,.
\end{align*}
Hence,
\begin{align*}
	\sup_{u>\varepsilon}\rho_{\bs}(u)\lesssim n^{-1/9}l_n\log^{2/3}(l_n)\log(Kr)  \,.
\end{align*}
Together with \eqref{eq:antis2},
\begin{align*}
	\sup_{u>\varepsilon}\mathbb{P}\bigg\{u-\varepsilon<\max_{j\in[Kr]}(|s_{n,2j-1}|^2+|s_{n,2j}|^2)\leq u+\varepsilon\bigg\} \lesssim  \frac{l_n\log^{2/3}(l_n)\log(Kr)}{n^{1/9}}
\end{align*}
provided that $\log(Kr)\ll n^{5/18}l_n^{-1}\log^{4/3}(l_n)$  and $l_n\log^{-4/3}(l_n)\ll n^{5/18}$. Combining with \eqref{eq:antis1}, it holds that
\begin{align}\label{eq:diff.epsilon}
	\sup_{u\geq0}\mathbb{P}\bigg\{u-\varepsilon<\max_{j\in[Kr]}(|s_{n,2j-1}|^2+|s_{n,2j}|^2)\leq u+\varepsilon\bigg\}
	\lesssim  \frac{l_n\log^{2/3}(l_n)\log(Kr)}{n^{1/9}}  \,.
\end{align}
If $l_n\geq \max\{2,C'\log(Kr)\}$ for some sufficiently large constant $C'>0$, by \eqref{eq:checkrhon}, \eqref{eq:rhonstar} and \eqref{eq:diff.Tn.checkTn},
\begin{align}\label{eq:gass1}
	\varrho_n^*
	\lesssim  n^{-1/9}l_n\log^{2/3}(l_n)\log(Kr)
	= o(1)
\end{align}
provided that $\log(Kr)\ll n^{1/9}l_n^{-1}\log^{-8/3}(l_n)$ and $l_n\log^{8/3}(l_n)\ll n^{1/9}$.
Therefore the proof of Proposition \ref{tm:1}(i) is completed. $\hfill\Box$

\subsection{Proof of Proposition \ref{tm:1}(ii)}
In this part, we always assume $\log r\ll n^{1/9}l_n^{-1}\log^{-8/3}(l_n)$, $l_n\log^{8/3}(l_n)\ll n^{1/9}$ and $l_n\geq \max(2,C'\log r)$ for some sufficiently large constant $C'>0$.
Recall $\mathcal{J}=[\omega_L,\omega_U]$. Given $\tilde{K}\asymp n$, let $\omega_L=\omega_1^*<\cdots<\omega_{\tilde{K}}^*=\omega_U$ be the isometric partition of $[\omega_L,\omega_U]$ satisfying $\max_{k\in [\tilde{K}-1]}(\omega_{k+1}^*-\omega_k^*) =O(\tilde{K}^{-1}) \rightarrow0$. Let $B_k'=[\omega_{k}^*,\omega_{k+1}^*)$ for each $k\in[\tilde{K}-2]$ and $B_{\tilde{K}-1}'=[\omega_{\tilde{K}-1}^*,\omega_{\tilde{K}}^*]$. Write $\mathcal{J}_{\tilde{K}}=\{\omega_1^*,\ldots,\omega_{\tilde{K}}^*\}$.
For any $\omega\in[\omega_L,\omega_U]$, there exists $k_\omega\in[\tilde{K}-1]$ such that $\omega\in B_{k_\omega}'$ and $|\omega-\omega_{k_\omega}^*|\leq (\omega_U-\omega_L)/(\tilde{K}-1)$. 	
As given in \eqref{eq:checkt}, $\check{T}_n(\omega;\mathcal{I})=\max_{s\in[r]} \{|\check{\eta}_{2s-1}^{\rm ext}(\omega)|^2+|\check{\eta}_{2s}^{\rm ext}(\omega)|^2\}$.
Recall $\bg_{n}(\omega)=\{g_{n,1}(\omega),\ldots,g_{n,2r}(\omega)\}^{\T}$
is a $(2r)$-dimensional  Gaussian process with mean zero and covariance function $\bSigma(\omega_1,\omega_2)=\{{\bf I}_r\otimes \bA(\omega_1)\}\bXi\{{\bf I}_r\otimes \bA^{\T}(\omega_2)\}$.
Our proof of Proposition \ref{tm:1}(ii) requires the following two inequalities:

\underline{\emph{Inequality 1.}} 	For any $u\gg n^{-1/2}l_n^{1/2}(\log l_n)\log^{1/2}(r)$, it holds that
\begin{align}\label{eq:diff.Tn.cTnK}
	&  \mathbb{P}\bigg\{\bigg|\sup_{\omega\in\mathcal{J}}T_n(\omega;\mathcal{I}) - \max_{\omega\in\mathcal{J}_{\tilde{K}}}\check{T}_n(\omega;\mathcal{I})\bigg|>u\bigg\}   \notag\\
	&~~~~~~~~~~ \lesssim rnl_n\exp\bigg(-\frac{Cn^{1/2}u^{1/6}}{l_n^{1/6}}\bigg)+rnl_n\exp\bigg\{-\frac{Cn^{1/2}u^{1/3}}{l_n^{1/2}\log^{1/3}(l_n)\log^{1/6}(r)}\bigg\}+\frac{1}{n}	 \notag\\
	&~~~~~~~~~~~~~   +rnl_n\exp\bigg(-\frac{Cn^{1/2}u^{1/2}}{l_n^{3/2}}\bigg) + rnl_n\exp\bigg\{-\frac{Cn^{1/2}u}{l_n^{5/2}(\log l_n) \log^{1/2}(r)}\bigg\}   \,.
\end{align}

\underline{\emph{Inequality 2.}} 		For any $u\geq C''l_nn^{-1/2}$ for some constant $C''>0$, it holds that
\begin{align} \label{eq:gcont2}
	&\bP\bigg[\bigg|\sup_{\omega\in\mathcal{J}}\max_{j\in[r]}\big\{|g_{n,2j-1}(\omega)|^2+|g_{n,2j}(\omega)|^2\big\}
	-  \sup_{\omega\in\mathcal{J}_{\tilde{K}}}\max_{j\in[r]}\big\{|g_{n,2j-1}(\omega)|^2+|g_{n,2j}(\omega)|^2\big\}\bigg| >u\bigg]  \notag \\
	&~~~~~~~~~~~~~~~\lesssim nr\exp\bigg(-\frac{Cn^2u}{l_n^{4}}\bigg) + nr \exp\bigg(-\frac{Cnu^2}{l_n^2}\bigg)
	+ nr\exp\bigg(-\frac{Cn}{l_n^4}\bigg) \,.
\end{align}
The proofs for \eqref{eq:diff.Tn.cTnK} and \eqref{eq:gcont2} are given in Sections \ref{subsec:p12.inq1} and \ref{subsec:p12.inq2}, respectively.

Write $\tilde\bs_{n,\tilde{K}}:=(\tilde{s}_{n,1},\ldots,\tilde{s}_{n,2\tilde{K}r})^{\T} =\{\bg_{n}^{\T}(\omega_1^*), \ldots, \bg_{n}^{\T}(\omega_{\tilde{K}}^*)\}^{\T}$.
We know $\tilde\bs_{n,\tilde{K}}$ is a $(2\tilde{K}r)$-dimensional Gaussian random vector with mean zero and covariance matrix  $\widetilde{\bH}\bXi\widetilde{\bH}^{\T}$, where 	  $\widetilde{\bH}=\{\bI_r\otimes \bA^\T(\omega_1^*),\ldots,\bI_r\otimes \bA^\T(\omega_{\tilde{K}}^*)\}^\T$.
Note that $\tilde{K}\asymp n$ and $r\geq n^{\kappa}$ for some sufficiently small constant $\kappa>0$.
Since $\sup_{\omega\in\mathcal{J}_{\tilde{K}}}\max_{j\in[r]}\{|g_{n,2j-1}(\omega)|^2+|g_{n,2j}(\omega)|^2\} = \max_{j\in[\tilde{K}r]}(|\tilde{s}_{n,2j-1}|^2+|\tilde{s}_{n,2j}|^2)$,
following the same arguments in Section \ref{subsec:s2} for deriving the convergence rate of $\check{\varrho}_n$ given in \eqref{eq:checkvarrho}, we have
\begin{align}\label{eq:p12.s3}
	\hbar_n:=&~	\sup_{u\geq 0}\bigg|\bP\bigg\{\sup_{\omega\in\mathcal{J}_{\tilde{K}}}{\check{T}}_n(\omega;\mathcal{I})\leq u\bigg\}	-\bP\bigg[\sup_{\omega\in\mathcal{J}_{\tilde{K}}}\max_{j\in[r]}\big\{|g_{n,2j-1}(\omega)|^2+|g_{n,2j}(\omega)|^2\big\}\leq u\bigg]\bigg|  \notag\\
	\lesssim&~ \frac{l_n\log^{2/3}(l_n)\log r}{n^{1/9}}=o(1)
	\,.
\end{align}
Define two events
\begin{align*}
	&~~~~~~~~~~~~~~~ \mathcal{D}_1= \bigg\{\bigg|\sup_{\omega\in\mathcal{J}}{T}_n(\omega;\mathcal{I})- \sup_{\omega\in\mathcal{J}_{\tilde{K}}}{\check{T}}_n(\omega;\mathcal{I})\bigg|\leq D\bigg\}  \,,\\
	&\mathcal{D}_2= \bigg\{\bigg| \sup_{\omega\in\mathcal{J}}\max_{j\in[r]}  \big\{|g_{n,2j-1}(\omega)|^2+|g_{n,2j}(\omega)|^2\big\} \\
	&~~~~~~~~~~~~~~~~~~~~~~~~~~
	-\sup_{\omega\in\mathcal{J}_{\tilde{K}}} \max_{j\in[r]}  \big\{|g_{n,2j-1}(\omega)|^2+|g_{n,2j}(\omega)|^2\big\} \bigg|\leq D\bigg\} \,,
\end{align*}
where $D=Cn^{-1/2}l_n^3\log^{2}(r)$ for some sufficiently large constant $C>0$.
Notice that $r\geq n^{\kappa}$ for some sufficiently small constant $\kappa>0$.
It follows from \eqref{eq:diff.Tn.cTnK} and \eqref{eq:gcont2}  that  $\bP(\mathcal{D}_1^{\rm c})
\lesssim n^{-1}
=o(1)$ and $\bP(\mathcal{D}_2^{\rm c})
\lesssim n^{-1}=o(1)$.
By \eqref{eq:p12.s3}, for any $u\geq 0$,
\begin{align*}
	&\bP\bigg\{\sup_{\omega\in\mathcal{J}}{T}_n(\omega;\mathcal{I})\leq u\bigg\}
	\leq \bP\bigg\{\sup_{\omega\in\mathcal{J}}{T}_n(\omega;\mathcal{I})\leq u,\,\mathcal{D}_1\bigg\} + \bP(\mathcal{D}_1^{\rm c}) \\
	&~~~~~~~~~~ \leq \bP\bigg[\sup_{\omega\in\mathcal{J}_{\tilde{K}}}\max_{j\in[r]} \big\{|g_{n,2j-1}(\omega)|^2+|g_{n,2j}(\omega)|^2\big\}\leq u+D\bigg] +\bP(\mathcal{D}_1^{\rm c}) + \hbar_n\\
	&~~~~~~~~~~ \leq \bP\bigg[\sup_{\omega\in\mathcal{J}}\max_{j\in[r]} \big\{|g_{n,2j-1}(\omega)|^2+|g_{n,2j}(\omega)|^2\big\}\leq u\bigg] + \bP(\mathcal{D}_2^{\rm c}) +\bP(\mathcal{D}_1^{\rm c}) + \hbar_n  \\
	&~~~~~~~~~~~~~ + \bP\bigg[u-D<\sup_{\omega\in\mathcal{J}_{\tilde{K}}}\max_{j\in[r]}  \big\{|g_{n,2j-1}(\omega)|^2+|g_{n,2j}(\omega)|^2\big\} \leq u+D\bigg]  \,.
\end{align*}
Likewise, we also have
\begin{align*}
	\bP\bigg\{\sup_{\omega\in\mathcal{J}}{T}_n(\omega;\mathcal{I})\leq u\bigg\}  \geq&~ \bP\bigg[\sup_{\omega\in\mathcal{J}}\max_{j\in[r]} \big\{|g_{n,2j-1}(\omega)|^2+|g_{n,2j}(\omega)|^2\big\}\leq u\bigg] - \bP(\mathcal{D}_2^{\rm c}) -\bP(\mathcal{D}_1^{\rm c}) - \hbar_n\\
	& - \bP\bigg[u-D<\sup_{\omega\in\mathcal{J}_{\tilde{K}}}\max_{j\in[r]}  \big\{|g_{n,2j-1}(\omega)|^2+|g_{n,2j}(\omega)|^2\big\} \leq u+D\bigg]
\end{align*}
for any $u\geq 0$. Therefore,
\begin{align}\label{eq:diff.Tn.g}
	& \sup_{u\geq 0}\bigg|\bP\bigg\{\sup_{\omega\in\mathcal{J}}{T}_n(\omega;\mathcal{I})\leq u\bigg\} - \bP\bigg[\sup_{\omega\in\mathcal{J}} \max_{j\in[r]}\big\{|g_{n,2j-1}(\omega)|^2+|g_{n,2j}(\omega)|^2\big\}\leq u\bigg]\bigg| \\
	&~~~~~~~
	\lesssim \sup_{u\geq 0}\bP\bigg[u-D<\sup_{\omega\in\mathcal{J}_{\tilde{K}}}\max_{j\in[r]}  \big\{|g_{n,2j-1}(\omega)|^2+|g_{n,2j}(\omega)|^2\big\} \leq u+D\bigg]   \notag\\
	&~~~~~~~~~~~ + \frac{l_n\log^{2/3}(l_n)\log r}{n^{1/9}}  + \frac{1}{n}   \notag\\
	&~~~~~~~ \lesssim  \frac{l_n\log^{2/3}(l_n)\log r}{n^{1/9}} + \sup_{u\geq 0}\bP\bigg[u-D<\sup_{\omega\in\mathcal{J}_{\tilde{K}}}\max_{j\in[r]}  \big\{|g_{n,2j-1}(\omega)|^2+|g_{n,2j}(\omega)|^2\big\} \leq u+D\bigg]   \,.  \notag
\end{align}
Recall $\sup_{\omega\in\mathcal{J}_{\tilde{K}}}\max_{j\in[r]}\{|g_{n,2j-1}(\omega)|^2+|g_{n,2j}(\omega)|^2\} = \max_{j\in[\tilde{K}r]}(|\tilde{s}_{n,2j-1}|^2+|\tilde{s}_{n,2j}|^2)$.	
Parallel to \eqref{eq:diff.epsilon} in Section \ref{subsec:s3},
we also have
\begin{align}\label{eq:gInv}
	&\sup_{u\geq 0}\bP\bigg[u-D<\sup_{\omega\in\mathcal{J}_{\tilde{K}}}\max_{j\in[r]}  \big\{|g_{n,2j-1}(\omega)|^2+|g_{n,2j}(\omega)|^2\big\} \leq u+D\bigg]  \notag\\
 &~~~~~~~~~~~~~~	 \lesssim \frac{l_n\log^{2/3}(l_n)\log r}{n^{1/9}}  \,,
\end{align}
which implies
\begin{align*} 
	&	\sup_{u\geq 0}\bigg|\bP\bigg\{\sup_{\omega\in\mathcal{J}}{T}_n(\omega;\mathcal{I})\leq u\bigg\} - \bP\bigg[\sup_{\omega\in\mathcal{J}}\max_{j\in[r]} \big\{|g_{n,2j-1}(\omega)|^2+|g_{n,2j}(\omega)|^2\big\}\leq u\bigg]\bigg|  \notag\\
	&~~~~~~~~~~~~~~	
	\lesssim  \frac{l_n\log^{2/3}(l_n)\log r}{n^{1/9}}
	=o(1) \,.
\end{align*}
Therefore  the proof of Proposition \ref{tm:1}(ii) is completed. $\hfill\Box$

\subsubsection{Proof of \eqref{eq:diff.Tn.cTnK}} \label{subsec:p12.inq1}
By the triangle inequality,
\begin{align}
	& \bigg|\sup_{\omega\in\mathcal{J}}\check{T}_n(\omega;\mathcal{I}) -\sup_{\omega\in\mathcal{J}_{\tilde{K}}}\check{T}_n(\omega;\mathcal{I})\bigg|
	\leq \sup_{\omega\in\cJ}|\check{T}_n(\omega;\mathcal{I}) - \check{T}_n(\omega_{k_\omega}^*;\mathcal{I})| \notag\\
	&~~~~~~~~~~ \lesssim \sup_{\omega\in\cJ}\max_{s\in[2r]}|\check{\eta}_s^{\ext}(\omega)-\check{\eta}_{s}^{\ext}(\omega_{k_\omega}^*)|^2
	+ \sup_{\omega\in\cJ}\max_{s\in[2r]}|\check{\eta}_{s}^{\ext}(\omega_{k_\omega}^*)||\check{\eta}_s^{\ext}(\omega)-\check{\eta}_{s}^{\ext}(\omega_{k_\omega}^*)|  \,.  \label{eq:diff.checkwk}
\end{align}
Since $\tilde{K}\asymp n$, $\sum_{k=-l_n}^{l_n}|k|\mathcal{W}(k/l_n)\asymp l_n^2$ and $|\omega-\omega_{k_\omega}^*|\leq 2\pi/(\tilde{K}-1)$, by mean value theorem and the triangle inequality, it holds that
\begin{align*}
	& \sup_{\omega\in\cJ}\max_{s\in[2r]}|\check{\eta}_{s}^{\rm ext}(\omega)-\check{\eta}_{s}^{\rm ext}(\omega_{k_\omega}^*)| \\
	&~~~~~~~~~~ \leq  \frac{1}{2\pi\sqrt{nl_n}}\sup_{\omega\in\cJ}\max_{s\in[r]}\sum_{t=l_n+1}^{n-l_n}\sum_{k=-l_n}^{l_n} |k|\mathcal{W}\bigg(\frac{k}{l_n}\bigg)|\mathring{x}_{\chi_1(s),t+k}\mathring{x}_{\chi_2(s),t}-\gamma_{\bchi(s)}(k)| |\omega-\omega_{k_\omega}^*|  \\
	&~~~~~~~~~~\lesssim \frac{l_n^{3/2}}{n^{3/2}}\max_{s\in[r]}\max_{-l_n\leq k\leq l_n}\sum_{t=l_n+1}^{n-l_n} |\mathring{x}_{\chi_1(s),t+k}\mathring{x}_{\chi_2(s),t}-\gamma_{\bchi(s)}(k)| \,.
\end{align*}
By the Bonferroni inequality, for any $u>0$,
\begin{align*}
	&\bP\bigg\{\sup_{\omega\in\cJ}\max_{s\in[2r]}|\check{\eta}_s^{\ext}(\omega)-\check{\eta}_s^{\ext}(\omega_{k_\omega}^*)| > u\bigg\} \\
	&~~~~~~~~~~~~ \leq \sum_{s=1}^{r}\sum_{k=-l_n}^{l_n}\sum_{t=l_n+1}^{n-l_n} \bP\bigg\{|\mathring{x}_{\chi_1(s),t+k}\mathring{x}_{\chi_2(s),t}-\gamma_{\bchi(s)}(k)| >\frac{C\sqrt{n}u}{l_n^{3/2}}\bigg\}  \\
	&~~~~~~~~~~~~ \lesssim rnl_n\exp(-Cn^{1/2}l_n^{-3/2}u)  \,,
\end{align*}
where the last step is based on \eqref{eq:t1} in Section \ref{subsec:pfLem2} for the proof of Lemma \ref{la:leading}.
By \eqref{eq:zetaij},  $\sup_{\omega\in\cJ}\max_{s\in[2r]}|\check{\eta}_s^{\ext}(\omega)|\leq n^{1/2}l_n^{-1/2}\sup_{\omega\in\cJ}\max_{(i,j)\in\cI}|\zeta_{i,j}(\omega)|$. By   \eqref{eq:zetatail} in Section \ref{subsec:pfLem4} for the proof of Lemma \ref{la:uniform}, if $l_n\geq 2$,
\begin{align*}
	&\bP\bigg\{\sup_{\omega\in\cJ}\max_{s\in[2r]}|\check{\eta}^{\ext}_s(\omega)|>u\bigg\}
	\leq\bP\bigg\{\sup_{\omega\in\cJ}\max_{(i,j)\in\cI}|\zeta_{i,j}(\omega)|>\frac{l_n^{1/2}u}{n^{1/2}}\bigg\} \\
	&~~~~~~~~~~ \lesssim rn\exp\bigg\{-\frac{Cu^2}{l_n^2\log^2 (l_n)}\bigg\} +rn\exp\bigg\{-\frac{Cn^{1/6}u^{1/3}}{l_n^{1/2}\log^{1/3} (l_n)}\bigg\}
	+ rn^2l_n\exp\bigg(-\frac{Cn^{1/2}u}{l_n^{3/2}}\bigg)
\end{align*}
for any $u>0$.
If  $\log r\ll n^{1/5}l_n^{-1/5}$ and $l_n=o(n)$, by the Bonferroni inequality and   \eqref{eq:diff.checkwk},
for some sufficiently large constant $C_*>0$, it holds that
\begin{align}\label{eq:diff.checkTndis}
	& \mathbb{P}\bigg\{\bigg|\sup_{\omega\in\mathcal{J}}\check{T}_n(\omega;\mathcal{I}) -\sup_{\omega\in\mathcal{J}_{\tilde{K}}}\check{T}_n(\omega;\mathcal{I})\bigg|>u\bigg\}   \notag\\
	&~~~~~~~~~\leq \mathbb{P}\bigg\{\sup_{\omega\in\cJ }\max_{s\in[2r]} |\check{\eta}_s^{\rm ext}(\omega) -\check{\eta}_s^{\rm ext} (\omega_{k_\omega}^*)|>\sqrt{\frac{Cu}{2}}\bigg\} \notag\\
	&~~~~~~~~~~~~ + \mathbb{P}\bigg\{\sup_{\omega\in \cJ} \max_{s\in[2r]} |\check{\eta}_s^{\rm ext}(\omega) -\check{\eta}_s^{\rm ext} (\omega_{k_\omega}^*)|>\frac{Cu}{C_{*}l_n(\log l_n)\log^{1/2}(r)}\bigg\}   \notag\\
	&~~~~~~~~~~~~   + \mathbb{P}\bigg\{\sup_{\omega\in\cJ}\max_{s\in[2r]} |\check{\eta}_s^{\rm ext} (\omega)|>\frac{C_{*}l_n(\log l_n) \log^{1/2}(r)}{2}\bigg\}  \notag\\
	&~~~~~~~~~ \lesssim rnl_n\exp\bigg(-\frac{Cn^{1/2}u^{1/2}}{l_n^{3/2}}\bigg) + rnl_n\exp\bigg\{-\frac{Cn^{1/2}u}{l_n^{5/2}(\log l_n) \log^{1/2}(r)}\bigg\} + n^{-C_{**} }
\end{align}
for any $u>0$, where  $C_{**}>0$ is a  constant only depending on $C_*$  such that  $C_{**}\rightarrow\infty$ as $C_*\rightarrow\infty$. We can select a specified $C_*>0$ such that $C_{**}=1$.
If  $\log r\ll n^{1/5}l_n^{-1/5}$, $l_n\log l_n=o(n) $
and $l_n\geq \max(2,\tilde{C}\log n)$ for some sufficiently large constant $\tilde{C}>0$, by the Bonferroni inequality,  \eqref{eq:Tn.Tnc.tail} in Section \ref{subsec:pfLem5} for the proof of Lemma \ref{la:app1}, and \eqref{eq:diff.checkTndis},
\begin{align*}
	&  \mathbb{P}\bigg\{\bigg|\sup_{\omega\in\mathcal{J}}T_n(\omega;\mathcal{I}) - \max_{\omega\in\mathcal{J}_{\tilde{K}}}\check{T}_n(\omega;\mathcal{I})\bigg|>u\bigg\}   \notag\\
	&~~~~~~~~
	\leq \mathbb{P}\bigg\{\bigg|\sup_{\omega\in\mathcal{J}}T_n(\omega;\mathcal{I}) -\sup_{\omega\in\mathcal{J}}\check{T}_n(\omega;\mathcal{I})\bigg| >\frac{u}{2}\bigg\}
	+ \mathbb{P}\bigg\{\bigg|\sup_{\omega\in\mathcal{J}}\check{T}_n(\omega;\mathcal{I}) - \max_{\omega\in\mathcal{J}_{\tilde{K}}}\check{T}_n(\omega;\mathcal{I})\bigg|>\frac{u}{2}\bigg\}  \notag\\
	&~~~~~~~~ \lesssim rnl_n\exp\bigg(-\frac{Cn^{1/2}u^{1/6}}{l_n^{1/6}}\bigg)+rnl_n\exp\bigg\{-\frac{Cn^{1/2}u^{1/3}}{l_n^{1/2}\log^{1/3}(l_n)\log^{1/6}(r)}\bigg\}	 \notag\\
	&~~~~~~~~~~~   +rnl_n\exp\bigg(-\frac{Cn^{1/2}u^{1/2}}{l_n^{3/2}}\bigg) + rnl_n\exp\bigg\{-\frac{Cn^{1/2}u}{l_n^{5/2}(\log l_n) \log^{1/2}(r)}\bigg\} + n^{-1}
\end{align*}
for $u\gg n^{-1/2}l_n^{1/2}(\log l_n)\log^{1/2}(r)$. $\hfill\Box$


\subsubsection{Proof of \eqref{eq:gcont2}}\label{subsec:p12.inq2}
Write $\Delta_{\tilde{K}}=(\omega_U-\omega_L)/(\tilde{K}-1)$.
By the triangle inequality and the Bonferroni inequality, 	for any $u>0$,
\begin{align}\label{eq:gcont}
	& \bP\bigg[\bigg|\sup_{\omega\in\mathcal{J}}\max_{j\in[r]}\{|g_{n,2j-1}(\omega)|^2+|g_{n,2j}(\omega)|^2\}  -  \sup_{\omega\in\mathcal{J}_{\tilde{K}}}\max_{j\in[r]}\{|g_{n,2j-1}(\omega)|^2+|g_{n,2j}(\omega)|^2\}\bigg| >u\bigg]  \notag\\
	&~~~~~~~~~~ \leq  \bP\bigg\{\max_{j\in[r]}\sup_{\omega\in\mathcal{J}}\big||g_{n,2j-1}(\omega)|^2+|g_{n,2j}(\omega)|^2 - |g_{n,2j-1}(\omega_{k_\omega}^*)|^2-|g_{n,2j}(\omega_{k_\omega}^*)|^2\big|>u\bigg\} \notag\\
	&~~~~~~~~~~\leq \bP\bigg\{\max_{j\in[2r]}\sup_{\omega\in\mathcal{J}}\big||g_{n,j}(\omega)|^2 - |g_{n,j}(\omega_{k_\omega}^*)|^2\big|>\frac{u}{2}\bigg\} \notag\\
	&~~~~~~~~~~\lesssim \tilde{K}r\max_{j\in[2r]}\max_{\omega\in\mathcal{J}_{\tilde{K}}} \bP\bigg\{\sup_{|s|\leq \Delta_{\tilde{K}}}\big||g_{n,j}(\omega+s)|^2 - |g_{n,j}(\omega)|^2\big|>\frac{u}{2}\bigg\}  \notag\\
	&~~~~~~~~~~ \leq \tilde{K}r\max_{j\in[2r]}\max_{\omega\in\mathcal{J}_{\tilde{K}}}\bP\bigg\{\sup_{|s|\leq\Delta_{\tilde{K}}}|g_{n,j}(\omega+s)-g_{n,j}(\omega)| >\frac{\sqrt{u}}{2}\bigg\} \notag\\
	&~~~~~~~~~~~~~ +\tilde{K}r\max_{j\in[2r]}\max_{\omega\in\mathcal{J}_{\tilde{K}}}\bP\bigg\{|g_{n,j}(\omega)| \sup_{|s|\leq\Delta_{\tilde{K}}}|g_{n,j}(\omega+s)-g_{n,j}(\omega)|>\frac{u}{8}\bigg\}
\end{align}
Without loss of generality, we focus on the scenario where $j=2j'-1$ with $j'\in[r]$. The arguments for $j=2j'$ with $j'\in[r]$	are identical.
Write $u_{t,j',k}=\mathring{x}_{\chi_1(j'),t+k}\mathring{x}_{\chi_2(j'),t}-\gamma_{\bchi(j')}(k)$.
Since   $\bg_{n}(\omega)=\{g_{n,1}(\omega),\ldots,g_{n,2r}(\omega)\}^{\T}$ is  a $(2r)$-dimensional  Gaussian process with mean zero and covariance function $\bSigma(\omega_1,\omega_2)$ defined as \eqref{eq:long_run}, we know $g_{n,2j'-1}(\omega)$ is a univariate Gaussian process with mean zero and covariance function
\begin{align*} 
	\cov\{g_{n,2j'-1}(\omega_1), g_{n,2j'-1}(\omega_2)\} = \bE\bigg[&\bigg\{\frac{1}{2\pi\sqrt{l_n\tilde{n}}}\sum_{t=l_n+1}^{n-l_n}\sum_{k=-l_n}^{l_n} \mathcal{W}\bigg(\frac{k}{l_n}\bigg)u_{t,j',k} \cos(k\omega_1)\bigg\}  \\
	&~~~~\times \bigg\{\frac{1}{2\pi\sqrt{l_n\tilde{n}}}\sum_{t=l_n+1}^{n-l_n}\sum_{k=-l_n}^{l_n} \mathcal{W}\bigg(\frac{k}{l_n}\bigg)u_{t,j',k} \cos(k\omega_2)\bigg\} \bigg]  \,, \notag
\end{align*}
which implies, for any $\omega_1,\omega_2\in\cJ$,
\begin{align}\label{eq:gdist}
	&\bE\big\{|g_{n,2j'-1}(\omega_1)-g_{n,2j'-1}(\omega_2)|^2\big\}  \notag \\
	&~~~~~~ = \bE\bigg[\bigg|\frac{1}{2\pi\sqrt{l_n\tilde{n}}}\sum_{t=l_n+1}^{n-l_n}\sum_{k=-l_n}^{l_n} \mathcal{W}\bigg(\frac{k}{l_n}\bigg)u_{t,j',k} \{\cos(k\omega_1)-\cos(k\omega_2)\} \bigg|^2\bigg] \notag\\
	&~~~~~~  \lesssim \frac{|\omega_1-\omega_2|^2}{\tilde{n}l_n} \sum_{k_1,k_2=-l_n}^{l_n} |k_1k_2|\mathcal{W}\bigg(\frac{k_1}{l_n}\bigg)\mathcal{W}\bigg(\frac{k_2}{l_n}\bigg)
	\bigg\{\sum_{t_1,t_2=l_n+1}^{n-l_n}\big|\bE(u_{t_1,j',k_1}u_{t_2,j',k_2})\big|\bigg\} \notag\\
	&~~~~~~ \lesssim \frac{|\omega_1-\omega_2|^2l_n^3}{\tilde{n}} \max_{-l_n\leq k_1,k_2\leq l_n}\sum_{t_1,t_2=l_n+1}^{n-l_n}\big|\bE(u_{t_1,j',k_1}u_{t_2,j',k_2})\big| \,.
\end{align}
The last step of \eqref{eq:gdist} is based on the fact $\sum_{k=-l_n}^{l_n}|k|\mathcal{W}(k/l_n)\asymp l_n^2$.
As  shown in  \eqref{eq:t1} in Section \ref{subsec:pfLem2} for the proof of Lemma \ref{la:leading} that	$\max_{t,j',k}\bP(|u_{t,j',k}|>u) \lesssim \exp(-Cu)$ for any $u>0$, it holds that $\bE(u_{t,j',k}^4)\leq C$. 
For any $-l_n\leq k_1,k_2\leq l_n$, by Cauchy-Schwarz inequality and Davydov's inequality,
\begin{align*} 
	& \max_{-l_n\leq k_1, k_2 \leq l_n} \sum_{t_1,t_2=l_n+1}^{n-l_n}\big|\bE(u_{t_1,j',k_1}u_{t_2,j',k_2})\big|  \notag\\
	&~~~~~~~~~~~~ \leq\sum_{t_1=l_n+1}^{n-l_n} \max_{-l_n\leq k \leq l_n}\bE(u_{t_1,j',k}^2)
	+\max_{-l_n\leq k_1, k_2 \leq l_n} \sum_{t_1\neq t_2}\big|\bE(u_{t_1,j',k_1}u_{t_2,j',k_2})\big|  \notag\\
	&~~~~~~~~~~~~\lesssim \tilde{n} + \sum_{t_1> t_2}\exp(-C|t_1-t_2-2l_n|_+)
	\lesssim \tilde{n}l_n \,.
\end{align*}
Together with \eqref{eq:gdist}, we have $\bE\{|g_{n,2j'-1}(\omega_1)-g_{n,2j'-1}(\omega_2)|^2\}\leq C l_n^4|\omega_1-\omega_2|^2 $ for any $\omega_1,\omega_2\in\cJ$.
Recall $\Delta_{\tilde{K}}=(\omega_U-\omega_L)/(\tilde{K}-1)$.
By Lemma \ref{la.cont}, there exist universal constants $\tilde{C}_1,\tilde{C}_2>0$ such that
\begin{align}\label{eq:I1}
	\sup_{\omega\in\cJ_{\tilde{K}}}\bP\bigg\{\sup_{|s|\leq\Delta_{\tilde{K}}}|g_{n,2j'-1}(\omega+s)-g_{n,2j'-1}(\omega)| >\frac{\sqrt{u}}{2}\bigg\} \leq \exp\bigg(-\frac{\tilde{C}_2{\tilde{K}}^2u}{l_n^{4}}\bigg)
\end{align}
for any $u\geq \tilde{C}_1l_n^4\Delta_{\tilde{K}}^{2}$.
Define the event $\mathcal{D}=\{|g_{n,2j'-1}(\omega)|\leq \tilde{K}^{1/2}l_n^{-1}\}$
and write $\sigma_{2j'-1}^2(\omega):=\bE\{g_{n,2j'-1}^2(\omega)\}$.
It holds that $\sup_{\omega\in\cJ}\max_{j'\in[r]}\sigma_{2j'-1}(\omega)\leq \check{C}l_n$ for some universal constant $\check{C}>0$.
Notice that $\tilde{K}\asymp n$ and $g_{n,2j'-1}(\omega)\sim\mathcal{N}\{0,\sigma^2_{2j'-1}(\omega)\}$ for any given $\omega\in\cJ$.
By \eqref{eq:I1},  if $l_n\ll n^{1/4}$,
\begin{align}\label{eq:I2}
	&\sup_{\omega\in\cJ_{\tilde{K}}}\bP\bigg\{|g_{n,{2j'-1}}(\omega)| \sup_{|s|\leq\Delta_{\tilde{K}}}|g_{n,2j'-1}(\omega+s)-g_{n,2j'-1}(\omega)|>\frac{u}{8}\bigg\} \notag\\
	&~~~~~~~~~~=
	\sup_{\omega\in\cJ_{\tilde{K}}}\bP\bigg\{|g_{n,,2j'-1}(\omega)| \sup_{|s|\leq\Delta_{\tilde{K}}}|g_{n,2j'-1}(\omega+s)-g_{n,2j'-1}(\omega)|>\frac{u}{8},\, \mathcal{D}\bigg\} \notag\\
	&~~~~~~~~~~~~~ + \sup_{\omega\in\cJ_{\tilde{K}}}\bP\bigg\{|g_{n,2j'-1}(\omega)| \sup_{|s|\leq\Delta_{\tilde{K}}}|g_{n,2j'-1}(\omega+s)-g_{n,2j'-1}(\omega)|>\frac{u}{8},\, \mathcal{D}^{\rm c}\bigg\} \notag\\
	&~~~~~~~~~~\leq\sup_{\omega\in\cJ_{\tilde{K}}}\bP\bigg\{ \sup_{|s|\leq\Delta_{\tilde{K}}}|g_{n,2j'-1}(\omega+s)-g_{n,2j'-1}(\omega)|>\frac{l_nu}{8\tilde{K}^{1/2}}\bigg\}  \notag\\
	&~~~~~~~~~~~~~ + \sup_{\omega\in\cJ_{\tilde{K}}}\bP\big\{|g_{n,2j'-1}(\omega)|> \tilde{K}^{1/2}l_n^{-1}\big\} \notag\\
	&~~~~~~~~~~\lesssim \exp\bigg(-\frac{Cnu^2}{l_n^2}\bigg) + \exp\bigg(-\frac{Cn}{l_n^4}\bigg)
\end{align}
for any $u\geq 4\tilde{C}_1^{1/2}\tilde{K}^{1/2}l_n\Delta_{\tilde{K}}$.
Notice that $\tilde{K}\asymp n$.
By \eqref{eq:gcont}, \eqref{eq:I1} and \eqref{eq:I2}, if $l_n\ll n^{1/4}$, we have
\begin{align*}
	&\bP\bigg[\bigg|\sup_{\omega\in\mathcal{J}}\max_{j\in[r]}\big\{|g_{n,2j-1}(\omega)|^2+|g_{n,2j}(\omega)|^2\big\}
	-  \sup_{\omega\in\mathcal{J}_{\tilde{K}}}\max_{j\in[r]}\big\{|g_{n,2j-1}(\omega)|^2+|g_{n,2j}(\omega)|^2\big\}\bigg| >u\bigg]  \notag \\
	&~~~~~~~~~~~~\lesssim nr\exp\bigg(-\frac{Cn^2u}{l_n^{4}}\bigg) + nr \exp\bigg(-\frac{Cnu^2}{l_n^2}\bigg)
	+ nr\exp\bigg(-\frac{Cn}{l_n^4}\bigg) 
\end{align*}
for any $u\geq \max(\tilde{C}_1l_n^4\Delta_{\tilde{K}}^{2}, 4\tilde{C}_1^{1/2}\tilde{K}^{1/2}l_n\Delta_{\tilde{K}})=Cl_nn^{-1/2}$.  $\hfill\Box$

\section{Proof of Proposition \ref{tm:2}}
To construct the proof of Proposition \ref{tm:2}, we need the following lemma whose proof is given in Section \ref{sec:pfLem11}.
\begin{lemma}\label{lem.Xi}
	Assume $r\geq n^{\kappa}$ for some sufficiently small constant $\kappa>0$ and that  Conditions {\rm \ref{as:moment}}, {\rm \ref{as:betamixing}} and {\rm \ref{as:kernel}} hold. Let $b_n\asymp n^\rho$ for some  constant $\rho\in(0,1)$. Then it holds that
	\begin{align*}
		|\widehat{\bXi}-\bXi|_\infty =&~ O_{\rm p}\bigg\{\frac{\log^{2} (r)} {n^{(\vartheta+2\rho-3\rho\vartheta-1)/(2\vartheta-1)}}\bigg\}
		+O_{\rm p}\bigg\{\frac{\log^{4}(r)} {n^{(2\vartheta+3\rho-4\rho\vartheta-2)/(2\vartheta-1)}}\bigg\} \\
		&   + O_{\rm p}\bigg\{\frac{\log^{1/2}(r)}{n^{(1-2\rho)/2}}\bigg\}
		+O_{\rm p}\bigg(\frac{l_n\log r}{n^{1-\rho}}\bigg)
		+ O\bigg(\frac{l_n^2}{n^{\rho}}\bigg)
	\end{align*}
	provided that $\log r \ll \min\{n^{(2\vartheta+2\rho-2\rho\vartheta-2)/(7\vartheta-4)}, n^{1/5}l_n^{-1/5}\}$, where $\bXi$ and $\widehat{\bXi}$ are defined as \eqref{eq:Xi} and \eqref{eq:Xihat}, respectively.	
\end{lemma}

\subsection{Proof of  Proposition \ref{tm:2}(i)}
Recall $\mathcal{J}=\{\omega_1,\ldots,\omega_K\}$, $\mathcal{X}_n=\{\bx_1,\ldots,\bx_n\}$
and	$\hat\bfeta^{\rm ext}(\omega)= \{\hat\eta_1^{\rm ext}(\omega),\ldots,\hat\eta_{2r}^{\rm ext}(\omega)\}^\T = \{\bI_r\otimes\bA(\omega)\}{\tilde n}^{-1/2} \sum_{t=1}^{\tilde n}\epsilon_t\hat{\bc}_t $, where $(\epsilon_1,\ldots,\epsilon_{\tilde{n}})^{\T}\sim\mathcal{N}(\bzero,\bTheta)$ independent of $\cX_n$ and the $(i,j)$-th component of $\bTheta$ is $\mathcal{K}\{(i-j)/b_n\}$.
Let
\begin{align*}
	\hat{\mathring\bfeta}^{\rm ext}:=\big(\hat{\mathring{\eta}}^{\rm ext}_1, \ldots, \hat{\mathring{\eta}}^{\rm ext}_{2Kr}\big)^\T =\big[\{\hat{\bfeta}^{\ext}(\omega_1)\}^{\T}, \ldots, \{\hat{\bfeta}^{\ext}(\omega_K)\}^{\T}\big]^{\T} \,.
\end{align*}
Then $\hat{\mathring\bfeta}^{\rm ext}\,|\,\cX_n \sim \mathcal{N}(\bzero,\bH\widehat{\bXi}\bH^{\T})$  and $\xi_{\cJ}$ given in \eqref{eq:xi} can be reformulated as
\begin{align*}
	\xi_{\mathcal{J}}
	=\max_{\ell\in[Kr]}\big(|\hat{\mathring\eta}^{\rm ext}_{2\ell-1}|^2+|\hat{\mathring\eta}^{\rm ext}_{2\ell}|^2\big) \,,
\end{align*}
where $\widehat{\bXi}$ is defined as \eqref{eq:Xihat} and $\bH=\{\bI_r\otimes \bA^{\T}(\omega_1),\ldots, \bI_r\otimes \bA^{\T}(\omega_K)\}^{\T}$. Let $\bs_{n,\by}=(s_{n,1},\ldots,s_{n,2Kr})^{\T} \sim \mathcal{N}(\bzero,\bH\bXi\bH^{\T})$ with $\bXi$ defined as \eqref{eq:Xi}.
For $\varrho_n^*$ given in \eqref{eq:starvarrho}, \eqref{eq:gass1} shows that
\begin{align}\label{eq:diff.Tn.s}
	&~\sup_{u\geq 0}\bigg|\bP\bigg\{\sup_{\omega\in\mathcal{J}}{T}_n(\omega;\mathcal{I})\leq u\bigg\}-\bP\bigg\{\max_{\ell\in[Kr]}(|s_{n,2\ell-1}|^2+|s_{n,2\ell}|^2)\leq u\bigg\}\bigg| \notag\\
	&~~~~~~~~~~~~~~ 	\lesssim
	n^{-1/9}l_n\log^{2/3}(l_n)\log (Kr) =o(1)
\end{align}
provided that $\log(Kr)\ll n^{1/9}l_n^{-1}\log^{-8/3}(l_n)$,  $l_n\log^{8/3}(l_n)\ll n^{1/9}$ and $l_n\geq \max\{2,C'\log(Kr)\}$ for some sufficiently large constant $C'>0$.
In the sequel, we always assume $\log(Kr)\ll n^{1/9}l_n^{-1}\log^{-8/3}(l_n)$,  $l_n\log^{8/3}(l_n)\ll n^{1/9}$ and $l_n\geq \max\{2,C'\log(Kr)\}$ for some sufficiently large constant $C'>0$. Notice that $\bH\bXi\bH^{\T}=\var(\tilde{n}^{-1/2}\sum_{t=1}^{\tilde{n}}\bH\bc_t)$.
Write $\bH=(\bh_1,\ldots,\bh_{2Kr})^{\T}$, where each $\bh_j$ is a $r(2l_n+1)$-dimensional vector. For any $j\in[Kr]$, there exists unique pair $(\ell,s)\in[r]\times[K]$ such that
$\bh_{2j-1}^{\T}\bc_t = (2\pi)^{-1}l_n^{1/2}z_{\bchi(\ell),t}^{(1)}(\omega_s)$ and $\bh_{2j}^{\T}\bc_t  = -(2\pi)^{-1}l_n^{1/2}z_{\bchi(\ell),t}^{(2)}(\omega_s)$
with $z_{i,j,t}^{(1)}(\omega)$ and $z_{i,j,t}^{(2)}(\omega)$  defined in \eqref{eq:zijt}. By  \eqref{eq:zij.tail} in Section \ref{subsec:pfLem3} for the proof of Lemma \ref{la:tailprob}, if $l_n\geq 2$, it holds that for any $u>0$
\begin{align*}
	\max_{j\in[2Kr]}\bP(|\bh_j^{\T}\bc_t|>u) \leq C\exp\bigg( - \frac{Cu}{l_n^{1/2}\log l_n}\bigg) \,.
\end{align*}
By the definition of Orlicz norm, we have $\max_{j\in[2Kr]}|\bh_j^{\T}\bc_t|_{\psi_1} \leq Cl_n^{1/2}\log l_n$. For any $u\geq 0$, it holds that $\{\max_{\ell\in[Kr]}(|s_{n,2\ell-1}|^2+|s_{n,2\ell}|^2)\leq u\}=\{\bs_{n,\by}\in A(u)\}$ and $\{\max_{\ell\in[Kr]}(|\hat{\mathring\eta}^{\rm ext}_{2\ell-1}|^2+|\hat{\mathring\eta}^{\rm ext}_{2\ell}|^2)\leq u\}=\{\hat{\mathring \bfeta}^{\rm ext} \in{A}(u)\}$, where ${A}(u)$  given in \eqref{eq:Adef} is a 2-sparsely convex set.
Recall $\bs_{n,\by}\sim \mathcal{N}(\bzero,\bH\bXi\bH^{\T})$ and  $\hat{\mathring{\bfeta}}^{\ext}\,|\, \cX_n\sim\mathcal{N}(\bzero,\bH\widehat{\bXi}\bH^{\T})$.
As shown in the proof of  Theorem 10 of \cite{CCW_2021_supp} 
with $(s,p,B_n,n)=(2,2Kr,Cl_n^{1/2}\log l_n,\tilde{n})$,
\begin{align} \label{eq:diff.ghat}
	&\sup_{u\geq 0}\bigg|\bP\bigg\{\max_{\ell\in[Kr]}\big(|s_{n,2\ell-1}|^2+|s_{n,2\ell}|^2\big)\leq u\bigg\}
	-\bP\bigg\{\max_{\ell\in[Kr]}\big(|\hat{\mathring\eta}^{\rm ext}_{2\ell-1}|^2+|\hat{\mathring\eta}^{\rm ext}_{2\ell}|^2\big)\leq u\,\bigg|\,\mathcal{X}_n\bigg\}\bigg| \notag\\
	&~~~~~~~~~~~~ \lesssim \big|\bH(\widehat{\bXi}-\bXi)\bH^{\T}\big|_\infty^{1/3} \log^{2/3}(Kr) +n^{-1}l_n^{1/2}\log l_n+ n^{-1}\log^{1/2}(Kr)   \notag\\
	&~~~~~~~~~~~~ \lesssim
	|\widehat{\bXi}-\bXi|_\infty^{1/3}l_n^{1/3} \log^{2/3}(Kr) +n^{-1}l_n^{1/2}\log l_n+ n^{-1}\log^{1/2}(Kr)  \,,
\end{align}
where the last step is	due to $|\bH(\widehat\bXi-\bXi)\bH^{\T}|_\infty\leq\max_{i,j\in [2Kr]}|\bh_i|_1|\widehat{\bXi}-\bXi|_\infty|\bh_j|_1 \lesssim l_n|\widehat{\bXi}-\bXi|_\infty$.
Together with \eqref{eq:diff.Tn.s}, it holds that
\begin{align}  \label{eq:gass3}
	&	\sup_{u\geq 0}\bigg|\mathbb{P}\bigg\{\sup_{\omega\in\mathcal{J}}{T}_n(\omega;\mathcal{I})\leq u\bigg\} - \mathbb{P}(\xi_{\mathcal{J}}\leq u\,|\,\mathcal{X}_n)\bigg| \notag\\
	&~~~~~~~~~~~~~	\lesssim
	n^{-1/9}l_n\log^{2/3}(l_n)\log (Kr)
	+  |\widehat\bXi-{\bXi}|_\infty^{1/3}l_n^{1/3}\log^{2/3}(Kr) \,.
\end{align}
By Lemma \ref{lem.Xi},
\begin{align*}
	&|\widehat\bXi-{\bXi}|_\infty^{1/3}l_n^{1/3}\log^{2/3}(Kr)  \\
	&~~~~~~~~ \lesssim \frac{l_n\log^{2/3}(Kr)}{n^{\rho/3}}
	+O_{\rm p}\bigg\{\frac{l_n^{1/3}\log^{4/3} (Kr)} {n^{(\vartheta+2\rho-3\rho\vartheta-1)/(6\vartheta-3)}}
	\bigg\}
	+O_{\rm p}\bigg\{\frac{l_n^{1/3}\log^{2}(Kr)} {n^{(2\vartheta+3\rho-4\rho\vartheta-2)/(6\vartheta-3)} } \bigg\}  \notag\\
	&~~~~~~~~~~~
	+O_{\rm p}\bigg\{\frac{l_n^{1/3}\log^{5/6}(Kr)}{n^{(1-2\rho)/6}}\bigg\}
	+O_{\rm p}\bigg\{\frac{l_n^{2/3}\log (Kr)}{n^{(1-\rho)/3}}\bigg\}
\end{align*}
provided that $\log r \ll \min\{n^{(2\vartheta+2\rho-2\rho\vartheta-2)/(7\vartheta-4)}, n^{1/5}l_n^{-1/5}\}$.
Together with \eqref{eq:gass3},
\begin{align*} 
	&\sup_{u\geq 0}\bigg|\mathbb{P}\bigg\{\sup_{\omega\in\mathcal{J}}{T}_n(\omega;\mathcal{I})\leq u\bigg\} - \mathbb{P}(\xi_{\mathcal{J}}\leq u\,|\,\mathcal{X}_n)\bigg|  \notag\\
	&~~~~~~~~~~~ \lesssim  
	o(1)
	+ \frac{l_n\log^{2/3}(l_n)\log(Kr)}{n^{1/9}}
	+ \frac{l_n\log^{2/3}(Kr)}{n^{\rho/3}}
	+O_{\rm p}\bigg\{\frac{l_n^{1/3}\log^{4/3} (Kr)} {n^{(\vartheta+2\rho-3\rho\vartheta-1)/(6\vartheta-3)}}
	\bigg\}  \\
	&~~~~~~~~~~~~~~ 		+O_{\rm p}\bigg\{\frac{l_n^{1/3}\log^{2}(Kr)} {n^{(2\vartheta+3\rho-4\rho\vartheta-2)/(6\vartheta-3)} } \bigg\}
	+O_{\rm p}\bigg\{\frac{l_n^{1/3}\log^{5/6}(Kr)}{n^{(1-2\rho)/6}}\bigg\}
	+O_{\rm p}\bigg\{\frac{l_n^{2/3}\log (Kr)}{n^{(1-\rho)/3}}\bigg\}  \notag\\
	&~~~~~~~~~~~ = o_{\rm p}(1)
\end{align*}
provided that
$\log(Kr)\ll f_1(l_n,n;\vartheta,\rho)$
and $\max\{2,C'\log(Kr)\}\leq l_n\ll  n^{f_2(\vartheta,\rho)}$
for some sufficiently large constant $C'>0$, where $f_1(l_n,n;\vartheta,\rho)$ and $f_2(\vartheta,\rho)$ are specified in \eqref{eq:f1.f2} and \eqref{eq:f1.f2.2}, respectively.
We have completed the proof of Proposition \ref{tm:2}(i). $\hfill\Box$


\subsection{Proof of  Propostion \ref{tm:2}(ii)}
In this part, we always assume 	$\log r\ll f_1(l_n,n;\vartheta,\rho)$
and $\max(2,C'\log r)\leq l_n\ll n^{f_2(\vartheta,\rho)}$
for some sufficiently large constant $C'>0$, where $f_1(l_n,n;\vartheta,\rho)$ and $f_2(\vartheta,\rho)$ are specified in \eqref{eq:f1.f2} and \eqref{eq:f1.f2.2}, respectively.
Recall $\mathcal{J}=[\omega_L,\omega_U]$ and $\mathcal{X}_n=\{\bx_1,\ldots,\bx_n\}$. Given $\tilde{K}\asymp n$, write $\omega_k^*= \omega_L+ (k-1)(\omega_U-\omega_L)/(\tilde{K}-1)$ for each $k\in[\tilde{K}]$. Let $B_k'=[\omega_{k}^*,\omega_{k+1}^*)$ for each $k\in[\tilde{K}-2]$ and $B_{\tilde{K}-1}'=[\omega_{\tilde{K}-1}^*,\omega_{\tilde{K}}^*]$.
Recall	$\hat\bfeta^{\rm ext}(\omega)= \{\hat\eta_1^{\rm ext}(\omega),\ldots,\hat\eta_{2r}^{\rm ext}(\omega)\}^\T = \{\bI_r\otimes\bA(\omega)\}{\tilde n}^{-1/2} \sum_{t=1}^{\tilde n}\epsilon_t\hat{\bc}_t $, where $(\epsilon_1,\ldots,\epsilon_{\tilde{n}})^{\T}\sim\mathcal{N}(\bzero,\bTheta)$ independent of $\cX_n$ and the $(i,j)$-th component of $\bTheta$ is $\mathcal{K}\{(i-j)/b_n\}$.
Write $\mathcal{J}_{\tilde{K}}=\{\omega_1^*,\ldots,\omega_{\tilde{K}}^*\}$
and
\begin{align*}
	\xi_{\mathcal{J}_{\tilde{K}}} := \sup_{\omega\in\mathcal{J}_{\tilde{K}}}\max_{j\in[r]}\{|\hat{\eta}^{\ext}_{2j-1}(\omega)|^2+|\hat{\eta}^{\ext}_{2j}(\omega)|^2\} \,.
\end{align*}
Conditionally on $\mathcal{X}_n$,  $\hat{\bfeta}^{\ext}(\omega)$ is a $(2r)$-dimensional   Gaussian process with mean zero and   covariance function  $\{\bI_r\otimes \bA(\omega_1)\} \widehat\bXi \{\bI_r\otimes \bA^{\T}(\omega_2)\}$, where $\widehat\bXi$ is given in {\rm \eqref{eq:Xihat}}. 	Recall that $\xi_{\mathcal{J}} = \sup_{\omega\in\mathcal{J}}\max_{j\in[r]}\{|\hat{\eta}^{\ext}_{2j-1}(\omega)|^2+|\hat{\eta}^{\ext}_{2j}(\omega)|^2\}$.
By the triangle inequality,
\begin{align}\label{eq:diff.Tn.XiJ}
	& \sup_{u\geq 0} \bigg|\bP\bigg\{\sup_{\omega\in\mathcal{J}}{T}_n(\omega;\mathcal{I})\leq u\bigg\} - \bP(\xi_{\mathcal{J}}\leq u\,|\,\mathcal{X}_n)\bigg|  \notag\\
	&~~~~~~~~~~\leq \sup_{u\geq 0}\bigg|\bP\bigg\{\sup_{\omega\in\mathcal{J}}{T}_n(\omega;\mathcal{I}) \leq u\bigg\} - \bP\bigg[\sup_{\omega\in\mathcal{J}_{\tilde{K}}}\max_{j\in[r]} \big\{|g_{n,2j-1}(\omega)|^2+|g_{n,2j}(\omega)|^2\big\} \leq u\bigg]\bigg|  \notag\\
	&~~~~~~~~~~~~~ +\sup_{u\geq 0}\bigg|\bP\bigg[\sup_{\omega\in\mathcal{J}_{\tilde{K}}}\max_{j\in[r]} \big\{|g_{n,2j-1}(\omega)|^2+|g_{n,2j}(\omega)|^2\big\} \leq u\bigg] - \bP(\xi_{\mathcal{J}_{\tilde{K}}}\leq u\,|\,\mathcal{X}_n)\bigg|  \notag\\
	&~~~~~~~~~~~~~ + \sup_{u\geq 0} \big|\bP(\xi_{\mathcal{J}_{\tilde{K}}}\leq u\,|\,\mathcal{X}_n) -\bP(\xi_{\mathcal{J}}\leq u\,|\,\mathcal{X}_n) \big| \,,
\end{align}
where $\bg_{n}(\omega)=\{g_{n,1}(\omega),\ldots,g_{n,2r}(\omega)\}^{\T}$ is a $(2r)$-dimensional Gaussian process with mean zero and covariance function $\bSigma(\omega_1,\omega_2)=\{\bI_r\otimes \bA(\omega_1)\} \bXi \{\bI_r\otimes \bA^{\T}(\omega_2)\}$ with $\bXi$   defined as   {\rm\eqref{eq:Xi}}.

Following the arguments to derive \eqref{eq:diff.Tn.g} with  $D=Cn^{-1/2}l_n^3\log^2(r)$ for some sufficiently large constant $C>0$, we have
\begin{align} \label{eq:diff.Tn.gnK}
	& \sup_{u\geq 0}\bigg|\bP\bigg\{\sup_{\omega\in\mathcal{J}}{T}_n(\omega;\mathcal{I}) \leq u\bigg\} - \bP\bigg[\sup_{\omega\in\mathcal{J}_{\tilde{K}}}\max_{j\in[r]} \big\{|g_{n,2j-1}(\omega)|^2+|g_{n,2j}(\omega)|^2\big\} \leq u\bigg]\bigg|  \\
	&~~~~~~\lesssim
	\frac{l_n\log^{2/3}(l_n)\log r}{n^{1/9}}
	+ \sup_{u\geq 0}\bP\bigg[u-D<\sup_{\omega\in\mathcal{J}_{\tilde{K}}}\max_{j\in[r]}  \big\{|g_{n,2j-1}(\omega)|^2+|g_{n,2j}(\omega)|^2\big\} \leq u+D\bigg]  \notag\\
	&~~~~~~\lesssim  \frac{l_n\log^{2/3}(l_n)\log r}{n^{1/9}} = o(1) \,, \notag
\end{align}
where the last step is based on \eqref{eq:gInv}.
Write $\widetilde{\bH}=\{\bI_r\otimes \bA^\T(\omega_1^*),\ldots,\bI_r\otimes \bA^\T(\omega_{\tilde{K}}^*)\}^\T$.
Recall $	\xi_{\mathcal{J}_{\tilde{K}}} = \sup_{\omega\in\mathcal{J}_{\tilde{K}}}\max_{j\in[r]}\{|\hat{\eta}^{\ext}_{2j-1}(\omega)|^2+|\hat{\eta}^{\ext}_{2j}(\omega)|^2\} $.
Notice that $\tilde{K}\asymp n$ and $r\geq n^{\kappa}$ for some sufficiently small constant $\kappa>0$.
Parallel to \eqref{eq:diff.ghat}, by Lemma \ref{lem.Xi}, it holds that
\begin{align}\label{eq:diff.g}
	&\sup_{u\geq 0}\bigg|\bP\bigg[\sup_{\omega\in\mathcal{J}_{\tilde{K}}}\max_{j\in[r]} \big\{|g_{n,2j-1}(\omega)|^2+|g_{n,2j}(\omega)|^2\big\} \leq u\bigg] - \bP(\xi_{\mathcal{J}_{\tilde{K}}}\leq u\,|\,\mathcal{X}_n)\bigg| \notag\\
	&~~~~~~~~~~~~ \lesssim |\widehat{\bXi}-\bXi|_\infty^{1/3}l_n^{1/3} \log^{2/3}(\tilde{K}r) +n^{-1}l_n^{1/2}\log l_n+ n^{-1}\log^{1/2}(\tilde{K}r)
	= o_{\rm p}(1) \,.
\end{align}
Define the event $\mathcal{D}=\{|\widehat{\bXi}-\bXi|_\infty\leq l_n\}$. Similarly as \eqref{eq:gcont2}, 	for any $u\geq C''l_nn^{-1/2}$ with some constant $C''>0$, we have
\begin{align*}	\bP\big(|\xi_{\mathcal{J}}-\xi_{\mathcal{J}_{\tilde{K}}}|>u, \, \mathcal{D}\,|\,\mathcal{X}_n\big)
	\lesssim
	nr\exp(-Cnl_n^{-4}u)+nr\exp(-Cnl_n^{-2}u^2)
	+nr\exp(-Cnl_n^{-4})  \,.
\end{align*}
The proof  follows the arguments in Section \ref{subsec:p12.inq2} with replacing $\bP(\cdot)$, $\bE(\cdot)$ and  $\bXi$ by  $\bP(\cdot\,|\,\cX_n)$, $\bE(\cdot\,|\,\cX_n)$ and $\widehat\bXi$, respectively. Hence, 	for any $u\geq C''l_nn^{-1/2}$ with some constant $C''>0$,
\begin{align}\label{eq:XiK_Xn}	\bP\big(|\xi_{\mathcal{J}}-\xi_{\mathcal{J}_{\tilde{K}}}|>u \,|\,\mathcal{X}_n\big)
	\lesssim&~
	nr\exp(-Cnl_n^{-4}u)+nr\exp(-Cnl_n^{-2}u^2) \notag\\
	&+nr\exp(-Cnl_n^{-4})
	+ I(|\widehat{\bXi}-\bXi|_\infty>l_n) \,.
\end{align}	
By Lemma \ref{lem.Xi}, we have $\bE\{ I(|\widehat{\bXi}-\bXi|_\infty>l_n)\}=\bP(|\widehat{\bXi}-\bXi|_\infty>l_n)=o(1)$, which implies $I(|\widehat{\bXi}-\bXi|_\infty>l_n)=o_{\rm p}(1)$.
Define the event $\tilde{\mathcal{D}}=\{|\xi_{\mathcal{J}}-\xi_{\mathcal{J}_{\tilde{K}}}|\leq {D}\}$, where 	 ${D}=Cn^{-1/2}l_n^{3}\log^{2}(r)$ for some sufficiently large constant $C>0$.
Recall $r\geq n^{\kappa}$ for some sufficiently small constant $\kappa>0$.
It follows from \eqref{eq:XiK_Xn} that $\bP(\tilde{\mathcal{D}}^{\rm c}\,|\,\cX_n)
\lesssim n^{-1} + I(|\widehat{\bXi}-\bXi|_\infty>l_n)
=o_{\rm p}(1)$. By the triangle inequality,
\begin{align*}
	\bP\big(\xi_{\mathcal{J}}\leq u\,|\,\mathcal{X}_n\big)
	\leq&~ \bP\big(\xi_{\mathcal{J}_{\tilde{K}}}\leq u+{D}\,|\,\mathcal{X}_n\big) + \bP(\tilde{\mathcal{D}}^{\rm c}\,|\,\mathcal{X}_n) \\
	=&~ \bP\big(\xi_{\mathcal{J}_{\tilde{K}}}\leq u\,|\,\mathcal{X}_n\big) +  \bP\big(u<\xi_{\mathcal{J}_{\tilde{K}}}\leq u+{D}\,|\,\mathcal{X}_n\big)
	+ \bP(\tilde{\mathcal{D}}^{\rm c}\,|\,\mathcal{X}_n) \\
	\leq&~    2 \sup_{u\geq 0}\bigg|\bP\bigg[\sup_{\omega\in\mathcal{J}_{\tilde{K}}}\max_{j\in[r]} \big\{|g_{n,2j-1}(\omega)|^2+|g_{n,2j}(\omega)|^2\big\} \leq u\bigg] - \bP\big(\xi_{\mathcal{J}_{\tilde{K}}}\leq u\,|\,\mathcal{X}_n\big)\bigg|  \\
	& +  \bP\bigg[u<\sup_{\omega\in\mathcal{J}_{\tilde{K}}}\max_{j\in[r]} \big\{|g_{n,2j-1}(\omega)|^2+|g_{n,2j}(\omega)|^2\big\} \leq u+{D}\bigg] \\
	& + \bP\big(\xi_{\mathcal{J}_{\tilde{K}}}\leq u\,|\,\mathcal{X}_n\big)
	+ \bP\big(\tilde{\mathcal{D}}^{\rm c}\,|\,\mathcal{X}_n\big)
\end{align*}
for any $u\geq 0$. Likewise, we can obtain the reverse inequality. Since $\tilde{K}\asymp n$, by   \eqref{eq:gInv} and \eqref{eq:diff.g}, 
\begin{align*}
	&\sup_{u\geq 0}\big|\bP\big(\xi_{\mathcal{J}_{\tilde{K}}}\leq u\,|\,\mathcal{X}_n\big) -\bP\big(\xi_{\mathcal{J}}\leq u\,|\,\mathcal{X}_n\big) \big| \\
	&~~~~~~ \lesssim  \sup_{u\geq 0}\bigg|\bP\bigg[\sup_{\omega\in\mathcal{J}_{\tilde{K}}}\max_{j\in[r]} \big\{|g_{n,2j-1}(\omega)|^2+|g_{n,2j}(\omega)|^2\big\} \leq u\bigg] - \bP\big(\xi_{\mathcal{J}_{\tilde{K}}}\leq u\,|\,\mathcal{X}_n\big)\bigg|
	+ \bP(\tilde{\mathcal{D}}^{\rm c}\,|\,\mathcal{X}_n)  \\
	&~~~~~~~~~   +  \sup_{u\geq 0}\bP\bigg[u-{D}<\sup_{\omega\in\mathcal{J}_{\tilde{K}}}\max_{j\in[r]} \big\{|g_{n,2j-1}(\omega)|^2+|g_{n,2j}(\omega)|^2\big\} \leq u+{D}\bigg]  \\
	&~~~~~~ \lesssim   |\widehat{\bXi}-\bXi|_\infty^{1/3}l_n^{1/3} \log^{2/3}(r) +  n^{-1/9}l_n\log^{2/3}(l_n)\log r  + I(|\widehat{\bXi}-\bXi|_\infty>l_n)
	= o_{\rm p}(1)  \,.
\end{align*}
Together with \eqref{eq:diff.Tn.gnK} and \eqref{eq:diff.g},   \eqref{eq:diff.Tn.XiJ} implies
\begin{align}\label{eq:gass4}
	&	\sup_{u\geq 0} \bigg|\bP\bigg\{\sup_{\omega\in\mathcal{J}}{T}_n(\omega;\mathcal{I})\leq u\bigg\} - \bP\big(\xi_{\mathcal{J}}\leq u\,|\,\mathcal{X}_n\big)\bigg|  \\
	&~~~~~~\lesssim   |\widehat{\bXi}-\bXi|_\infty^{1/3} l_n^{1/3}\log^{2/3}(r) +  n^{-1/9}l_n\log^{2/3}(l_n)\log r 	+I(|\widehat{\bXi}-\bXi|_\infty>l_n)
	= o_{\rm p}(1)  \,.  \notag
\end{align}
We have completed the proof of Proposition \ref{tm:2}(ii).
$\hfill\Box$

\section{Proof of Theorem \ref{tm:H0}}
In this part, we always assume 	$\log(M_{\cJ}r)\ll
f_{1}(l_n,n;\vartheta,\rho)
$
and $\max\{2,C'\log(M_{\cJ}r)\}\leq l_n\ll n^{f_2(\vartheta,\rho)}$ 
for some sufficiently large constant $C'>0$, where $f_1(l_n,n;\vartheta,\rho)$ and $f_2(\vartheta,\rho)$ are specified in \eqref{eq:f1.f2} and \eqref{eq:f1.f2.2}, respectively.

\subsection{Proof of Theorem \ref{tm:H0}(i)}
Define 
\begin{align*}
	T_n^{\rm G}=
 \left\{
	\begin{aligned}
	\max_{j\in[Kr]} (|s_{n,2j-1}|^2+|s_{n,2j}|^2)\,, ~~~~~~~~~~~ & \mbox{if~} \mathcal{J}=\{\omega_1,\ldots,\omega_K\}\,, \\
		\sup_{\omega\in\mathcal{J}}\max_{j\in[r]}\{|g_{n,2j-1}(\omega)|^2+|g_{n,2j}(\omega)|^2\}\,, ~~~& \mbox{if~}\mathcal{J}= [\omega_{L},\omega_U] \,,
	\end{aligned}
 \right.
\end{align*}	
where $\bs_{n,\by}=(s_{n,1},\ldots,s_{n,2Kr})^{\T}$ and $\bg_{n}(\omega)=\{g_{n,1}(\omega),\ldots,g_{n,2r}(\omega)\}^{\T}$ are defined in Proposition \ref{tm:1}. Since $T_n=\sup_{\omega\in\mathcal{J}}T_n(\omega;\mathcal{I})$ under $H_0$,   Proposition \ref{tm:1} shows that 
\begin{align}\label{eq:diff.Tn.TnG}
    \sup_{u>0}|\bP(T_n>u)-\mathbb{P}(T_n^{\rm G}>u)|=o(1) \,.
\end{align}
For any  $\epsilon>0$, let ${\rm cv}_\alpha^{(\epsilon)}$ and ${\rm cv}_\alpha^{(-\epsilon)}$ be two constants which satisfy $\bP\{T_n^{\rm G}>{\rm cv}_\alpha^{(\epsilon)}\}=\alpha+\epsilon$ and $\bP\{T_n^{\rm G}>{\rm cv}_\alpha^{(-\epsilon)}\}=\alpha-\epsilon$, respectively. By triangle inequality and Proposition \ref{tm:2}, 
\begin{align}\label{eq:diff.TnG.xi}
	\sup_{u>0}|\bP(T_n^{\rm G}>u) - \bP(\xi_{\cJ}>u\,|\,\mathcal{X}_n)|=o_{\rm p}(1)\,.
\end{align}
Notice that 
$\mathbb{P}(\xi_{\cJ}>\hat{\rm cv}_\alpha\,|\,\mathcal{X}_n)=\alpha$. We  claim that for any $\epsilon>0$, it holds that $\mathbb{P}\{{\rm cv}_\alpha^{(\epsilon)}<\hat{\rm cv}_\alpha<{\rm cv}_\alpha^{(-\epsilon)}\}\rightarrow 1$  as $n\rightarrow\infty$. Otherwise, if $\hat{\rm cv}_\alpha\leq {\rm cv}_\alpha^{(\epsilon)}$, by \eqref{eq:diff.TnG.xi}, we have
\begin{align*}
	\alpha=\bP(\xi_{\cJ}>\hat{\rm cv}_\alpha\,|\,\mathcal{X}_n) 
 \geq&\, \bP\{\xi_{\cJ} > {\rm cv}_\alpha^{(\epsilon)} \,|\,\mathcal{X}_n \}  \\
	=&\, \bP\{T_n^{\rm G}>{\rm cv}_\alpha^{(\epsilon)}\} + o_{\rm p}(1) = \alpha+\epsilon+o_{\rm p}(1) \,,
\end{align*}
which is a contradiction with probability approaching one as $n\rightarrow\infty$. Analogously, if $\hat{\rm cv}_\alpha \geq {\rm cv}_\alpha^{(-\epsilon)}$, by \eqref{eq:diff.TnG.xi}, we have
\begin{align*}
	\alpha=\bP(\xi_{\cJ}>\hat{\rm cv}_\alpha \,|\, \mathcal{X}_n) 
 \leq&\, \bP\{\xi_{\cJ}>{\rm cv}_\alpha^{(-\epsilon)}\,|\,\mathcal{X}_n\} \\
	=&\, \bP\{T_n^{\rm G}>{\rm cv}_\alpha^{(-\epsilon)}\} +o_{\rm p}(1)
	=\alpha-\epsilon+o_{\rm p}(1)\,,
\end{align*}
which is also a contradiction  with probability approaching one as $n\rightarrow\infty$. 

For any $\epsilon>0$, define the event $\mathcal{E}_{\epsilon}=\{{\rm cv}_\alpha^{(\epsilon)} < \hat{\rm cv}_\alpha < {\rm cv}_\alpha^{(-\epsilon)}\}$. Then $\bP(\mathcal{E}_{\epsilon})\rightarrow 1$ as $n\rightarrow\infty$. On the one hand, by \eqref{eq:diff.Tn.TnG},
\begin{align*}
	\bP(T_n>\hat{\rm cv}_\alpha) 
	\leq&~ \bP(T_n>\hat{\rm cv}_\alpha,\mathcal{E}_{\epsilon}) + \bP(\mathcal{E}_{\epsilon}^{\rm c})
	\leq \bP\{T_n>{\rm cv}_\alpha^{(\epsilon)}\} +o(1) \\
	=&~ \bP\{T_n^{\rm G}>{\rm cv}_\alpha^{(\epsilon)}\} +o(1) =\alpha+\epsilon+o(1) \,,
\end{align*}
which implies that $\varlimsup_{n\rightarrow\infty}\bP(T_n>\hat{\rm cv}_\alpha)\leq \alpha+\epsilon$. On the other hand, by \eqref{eq:diff.Tn.TnG},
\begin{align*}
	\bP(T_n>\hat{\rm cv}_\alpha) 
	\geq&~ \bP(T_n>\hat{\rm cv}_\alpha,\mathcal{E}_{\epsilon})
	\geq \bP\{T_n>{\rm cv}_\alpha^{(-\epsilon)}\} - \bP(\mathcal{E}_{\epsilon}^{\rm c})  \\
	\geq&~ \bP\{T_n^{\rm G}>{\rm cv}_\alpha^{(-\epsilon)}\} -o(1) =\alpha-\epsilon-o(1) \,,
\end{align*}
which implies that $\varliminf_{n\rightarrow\infty}\bP(T_n>\hat{\rm cv}_\alpha)\geq \alpha-\epsilon$. Since  $\bP(T_n>\hat{\rm cv}_\alpha)$ does not depend on $\epsilon$, by letting $\epsilon\rightarrow0^+$, we have $\lim_{n\rightarrow\infty}\bP(T_n>\hat{\rm cv}_\alpha)=\alpha$. Hence, we  complete the proof of Theorem \ref{tm:H0}(i). $\hfill\Box$

\subsection{Proof of Theorem \ref{tm:H0}(ii)}
For $M_{\cJ}$ defined in \eqref{eq:MJ},
we define $\mathcal{J}_M=\{\omega_1^*,\ldots,\omega_{M_{\cJ}}^*\}$ with
\begin{equation*}
\omega_m^*=\left\{
\begin{aligned}
	\omega_m\,,~~~~~~~~~~~~~~~~~&\textrm{if}~\cJ=\{\omega_1,\ldots,\omega_K\} \,,\\
	\omega_L+\frac{(m-1)(\omega_U-\omega_L)}{M_{\cJ}-1} \,,~~~&\textrm{if}~\cJ=[\omega_L,\omega_U]\,.
\end{aligned}
\right.
\end{equation*}
For $\hat{\bfeta}^{\ext}(\omega)=\{\hat{\eta}_1^{\ext}(\omega),\ldots,\hat{\eta}_{2r}^{\ext}(\omega)\}^{\T}$ defined in \eqref{eq:hatetaext},  let
\begin{align*}
\xi_{\mathcal{J}_M}:=\sup_{\omega\in\mathcal{J}_M}\max_{ j\in [r]}\big\{|\hat\eta_{2j-1}^{\rm ext}(\omega)|^2+|\hat\eta_{2j}^{\rm ext}(\omega)|^2\big\}\,.
\end{align*}
As mentioned below \eqref{eq:hatetaext}, conditionally on $\mathcal{X}_n$,  $\hat{\bfeta}^{\ext}(\omega)$ is a $(2r)$-dimensional   Gaussian process with mean zero and   covariance function  $\{\bI_r\otimes \bA(\omega_1)\} \widehat\bXi \{\bI_r\otimes \bA^{\T}(\omega_2)\}$ with $\widehat\bXi$ given in {\rm \eqref{eq:Xihat}}.
Write
\begin{align*}
\hat{\mathring\bfeta}^{\rm ext, \natural} :=\big(\hat{\mathring\eta}^{\rm ext, \natural}_1,\ldots, \hat{\mathring\eta}^{\rm ext,  \natural}_{2M_{\cJ}r}\big)^{\T}
=\big[\{\hat\bfeta^{{\rm ext}}(\omega_1^*)\}^{\T}, \ldots,\{\hat\bfeta^{{\rm ext}}(\omega_{M_{\cJ}}^*)\}^{\T}\big]^{\T} \,.
\end{align*}
Then  $\hat{\mathring\bfeta}^{\rm ext, \natural}\,|\,\cX_n\sim \mathcal{N}(\bzero,\widehat\bSigma^{\natural})$ with $\widehat\bSigma^{\natural}:=\{\hat\sigma^{\natural}(\ell_1,\ell_2)\}_{2M_{\cJ}r\times2M_{\cJ}r} =\bH^{\natural}\widehat\bXi\bH^{^{\natural},\T}$, where $\bH^{\natural}=\{\bI_r\otimes\bA^{\T}(\omega_1^*),\ldots,\bI_r\otimes\bA^{\T}(\omega_{M_{\cJ}}^*)\}^{\T}$.
Write $\bSigma^{\natural}:=\{\sigma^{\natural}(\ell_1,\ell_2)\}_{2M_{\cJ}r\times2M_{\cJ}r}=\bH^{\natural}\bXi\bH^{^{\natural},\T}$ and $\tilde{\varrho}=\max_{\ell\in[2M_{\cJ}r]}\sigma^{\natural}(\ell,\ell)$.
Recall  $\xi_{\mathcal{J}}=\sup_{\omega\in\mathcal{J}}\max_{j\in [r]}\{|\hat\eta_{2j-1}^{\rm ext}(\omega)|^2+|\hat\eta_{2j}^{\rm ext}(\omega)|^2\}$. 	
Define the events
\begin{align*}
\mathcal{E}_1(\nu_1)= \big\{|\xi_{\mathcal{J}}-\xi_{\mathcal{J}_M}|\leq \nu_1  \big\} ~~\mbox{ and }~~
\mathcal{E}_2(\nu_2)= \bigg\{\max_{\ell\in[2M_{\cJ}r]}\bigg|\frac{\hat\sigma^{\natural}(\ell,\ell)}{\sigma^{\natural}(\ell,\ell)}-1\bigg|\leq \nu_2 \bigg\} \,,
\end{align*}
where  $\nu_1=\tilde{C}' n^{-1/2}l_n^{3}\log^{2}(M_{\cJ}r)$ and $\nu_2=\tilde{C}''\log^{-2}(M_{\cJ}r)$ for some sufficiently large constants $\tilde{C}',\tilde{C}''>0$.  Following the same arguments as \eqref{eq:XiK_Xn} and Lemma \ref{lem.Xi}, we  have $\bP\{\mathcal{E}_1^{\rm c}(\nu_1)\,|\,\cX_n\}=o_{\rm p}(1)$ and $\max_{\ell\in[2M_{\cJ}r]}|\hat\sigma^{\natural}(\ell,\ell)-\sigma^{\natural}(\ell,\ell)|=o_{\rm p}\{\log^{-2}(M_{\cJ}r)\}$.
By Condition \ref{as:eigen}, we know $\min_{\ell\in[2M_{\cJ}r]}\sigma^{\natural}(\ell,\ell)$ is uniformly bounded away from zero. Hence,
\begin{align*}
\max_{\ell\in[2M_{\cJ}r]}\bigg|\frac{\hat\sigma^{\natural}(\ell,\ell)}{\sigma^{\natural}(\ell,\ell)} -1\bigg|
\leq \frac{\max_{\ell\in[2M_{\cJ}r]}|\hat\sigma^{\natural}(\ell,\ell)-\sigma^{\natural}(\ell,\ell)|} {\min_{\ell\in[2M_{\cJ}r]}\sigma^{\natural}(\ell,\ell)} = o_{\rm p}(\nu_2)\,,
\end{align*}
which implies $\bP\{\mathcal{E}_2^{\rm c}(\nu_2)\}=o(1)$.
Notice that
$\sup_{\omega\in\mathcal{J}_M}\max_{j\in[2r]} |\hat\eta_{j}^{\rm ext}(\omega)|
= \max_{\ell\in[2M_{\cJ}r]} |\hat{\mathring\eta}^{\rm ext, \natural}_{\ell}|  $.
By the triangle inequality and the Bonferroni inequality, for any $u\geq\nu_1$,
\begin{align*}
\bP\big\{
\xi_{\mathcal{J}}>u, \mathcal{E}_1(\nu_1), \mathcal{E}_2(\nu_2) \,|\,\mathcal{X}_n
\big\}
\leq &~ 
\bP\bigg[\sup_{\omega\in\mathcal{J}_M}\max_{ j\in[r]}\big\{|\hat\eta_{2j-1}^{\rm ext}(\omega)|^2+|\hat\eta_{2j}^{\rm ext}(\omega)|^2\big\}>u-\nu_1, \mathcal{E}_2(\nu_2) \,\bigg|\,\mathcal{X}_n\bigg]\\ 
\leq&~ \bP\bigg\{\sup_{\omega\in\mathcal{J}_M}\max_{j\in[2r]}|\hat\eta_{j}^{\rm ext}(\omega)|>\sqrt{\frac{u-\nu_1}{2}}, \mathcal{E}_2(\nu_2) \,\bigg|\,\mathcal{X}_n\bigg\} \\
=&~ \bP\bigg\{\max_{\ell\in[2M_{\cJ}r]} \big|\hat{\mathring\eta}_{\ell}^{\rm ext, \natural}\big| >\sqrt{\frac{u-\nu_1}{2}}, \mathcal{E}_2(\nu_2) \,\bigg|\,\mathcal{X}_n\bigg\}  \,.
\end{align*}
Restricted on $\mathcal{E}_2(\nu_2)$, it holds that
\begin{align*}
\bE\bigg(\max_{\ell\in[2M_{\cJ}r]} \big|\hat{\mathring\eta}^{\rm ext,  \natural}_{\ell}\big| \,\bigg|\, \mathcal{X}_n \bigg)
\leq
&~	[1+\{2\log(2M_{\cJ}r)\}^{-1}] \{2\log(2M_{\cJ}r)\}^{1/2}\max_{\ell\in[2M_{\cJ}r]}\{\hat\sigma^{\natural}(\ell,\ell)\}^{1/2} \\
\leq &~ (1+\nu_2)^{1/2}\tilde{\varrho}^{1/2}[1+\{2\log(2M_{\cJ}r)\}^{-1}]\{2\log(2M_{\cJ}r)\}^{1/2} \,.
\end{align*}
By Borell inequality for Gaussian process,
\begin{align*}
\bP\bigg\{\max_{\ell\in[2M_{\cJ}r]} \big|\hat{\mathring\eta}^{\rm ext, \natural}_{\ell}\big|\geq \bE\bigg(\max_{\ell\in[2M_{\cJ}r]} \big|\hat{\mathring\eta}^{\rm ext,  \natural}_{\ell}\big| \,\bigg|\,\mathcal{X}_n\bigg)+u \,\bigg|\,\mathcal{X}_n\bigg\}
\leq  2\exp\bigg\{-\frac{u^2}{2\max_{\ell\in[2M_{\cJ}r]}\hat\sigma^{\natural}(\ell,\ell)}\bigg\}
\end{align*}
for any $u>0$. Let $u_0=\nu_1+2\tilde{\varrho}(1+\nu_2)([1+\{2\log(2M_{\cJ}r)\}^{-1}]\{2\log(2M_{\cJ}r)\}^{1/2}+\{2\log(4/\alpha)\}^{1/2})^2$.
Restricted on  $\mathcal{E}_2(\nu_2)$, we have
\begin{align*}
\sqrt{\frac{u_0-\nu_1}{2}} \geq \bE\bigg(\max_{\ell\in[2M_{\cJ}r]} \big|\hat{\mathring\eta}^{\rm ext,\natural}_{\ell}\big| \,\bigg|\, \mathcal{X}_n \bigg) + \sqrt{2}\tilde{\varrho}^{1/2}(1+\nu_2)^{1/2}\log^{1/2}\bigg(\frac{4}{\alpha}\bigg) \,,
\end{align*}
which implies
\begin{align*}
\bP\big\{
\xi_{\mathcal{J}}>u_0,  \mathcal{E}_1(\nu_1), \mathcal{E}_2(\nu_2) \,|\,\mathcal{X}_n
\big\}
\leq&~  \bP\bigg\{\max_{\ell\in[2M_{\cJ}r]} \big|\hat{\mathring\eta}_{\ell}^{\rm ext, \natural}\big| >\sqrt{\frac{u_0-\nu_1}{2}}, \mathcal{E}_2(\nu_2) \,\bigg|\,\mathcal{X}_n\bigg\}  \\
\leq&~ 2\exp\bigg\{- \frac{2\tilde{\varrho}(1+\nu_2)\log(4/\alpha)}{2\tilde{\varrho}(1+\nu_2)}\bigg\}
=\frac{\alpha}{2} \,.
\end{align*}
Since $\bP\{\mathcal{E}_1^{\rm c}(\nu_1)\,|\,\cX_n\}+\bP\{\mathcal{E}_2^{\rm c}(\nu_2)\,|\,\cX_n\}=o_{\rm p}(1)$, then
$
\bP\{\mathcal{E}_1^{\rm c}(\nu_1)\,|\,\cX_n\} +\bP\{\mathcal{E}_2^{\rm c}(\nu_2)\,|\,\cX_n\} \leq \alpha/4$ with probability approaching one.
Hence,
$
\bP\big(\xi_{\cJ}>u_0\,|\,\cX_n\big)
\leq 5\alpha/6$
with probability approaching one.
By the definition of $\hat{\rm cv}_\alpha$, it holds with probability approaching one that
\begin{align}\label{eq:cvhat}
\hat{\rm cv}_\alpha \leq &~\nu_1+2\tilde{\varrho}(1+\nu_2)[1+\{2\log(2M_{\cJ}r)\}^{-1}]^2\lambda^2(M_{\cJ},r,\alpha) \notag\\
\leq&~ \nu_1+2{\varrho}(1+\nu_2)[1+\{2\log(2M_{\cJ}r)\}^{-1}]^2\lambda^2(M_{\cJ},r,\alpha) \,,
\end{align}
where $\lambda(M_{\cJ},r,\alpha)=\{2\log(2M_{\cJ}r)\}^{1/2}+\{2\log(4/\alpha)\}^{1/2}$ and $\varrho=\sup_{\omega\in\mathcal{J}}\max_{j\in[2r]}\sigma_j^2(\omega)$ with $\sigma_j^2(\omega)$ being the $j$-th element in the main diagonal of $\bSigma(\omega,\omega)$.

Recall  $  T_{n}=\sup_{\omega\in\mathcal{J}}\max_{(i,j)\in\mathcal{I}}
nl_n^{-1}|\hat f_{i,j}(\omega)|^2$. Let $(\omega_0,i_0,j_0)=\arg\sup_{\omega\in\mathcal{J},(i,j)\in\mathcal{I}}|f_{i,j}(\omega)|^2$.
Then $T_n\geq nl_n^{-1} |\hat{f}_{i_0,j_0}(\omega_{0})|^2$.
Given some $\epsilon_n\rightarrow 0$ satisfying $\epsilon_n^2\varrho l_n^{-2}\log^{-2}(l_n)\log^{-1}(n)\lambda^2(M_{\cJ},r,\alpha) \rightarrow \infty$,
we choose $u_n>0$ such that $(1+\nu_2)[1+\{2\log(2M_{\cJ}r)\}^{-1}+u_n]^2=(1+\epsilon_n)^2$.
Notice that
\begin{align}\label{eq:signal}
nl_n^{-1}|f_{i_0,j_0}(\omega_0)|^2
\geq&~  4\varrho\lambda^2(M_{\cJ},r,\alpha)(1+\epsilon_n)^2 \notag\\
=&~ 4\varrho\lambda^2(M_{\cJ},r,\alpha)(1+\nu_2)[1+\{2\log(2M_{\cJ}r)\}^{-1}+u_n]^2 \,.
\end{align}
Without loss of generality, we  assume $|\Re\{f_{i_0,j_0}(\omega_0)\}| \geq |\Im\{f_{i_0,j_0}(\omega)\}|$. Otherwise, the following arguments can be stated based on $\Im\{f_{i_0,j_0}(\omega)\}$. Also, we assume $\Re\{f_{i_0,j_0}(\omega_0)\}>0$. Otherwise, we can replace $\Re\{f_{i_0,j_0}(\omega_0)\}$ by $-\Re\{f_{i_0,j_0}(\omega_0)\}$ in our proof.
Under these assumptions, by \eqref{eq:signal}, we have $n^{1/2}l_n^{-1/2}\Re\{f_{i_0,j_0}(\omega_0)\}\geq \sqrt{2\varrho(1+\nu_2)}\lambda(M_{\cJ},r,\alpha)[1+\{2\log(2M_{\cJ}r)\}^{-1}+u_n]$.
By \eqref{eq:cvhat},
\begin{align*}
&\bP(T_n>\hat{\rm cv}_\alpha)
\geq  \bP\big\{nl_n^{-1} |\hat{f}_{i_0,j_0}(\omega_{0})|^2>\hat{\rm cv}_\alpha \big\}  \\
&~~	\geq\bP\big\{nl_n^{-1} |\hat{f}_{i_0,j_0}(\omega_{0})|^2 > \nu_1+2\varrho(1+\nu_2)[1+\{2\log(2M_{\cJ}r)\}^{-1}]^2\lambda^2(M_{\cJ},r,\alpha) \big\} -o(1) \\
&~~	\geq \bP\big[n^{1/2}l_n^{-1/2} \Re\{\hat{f}_{i_0,j_0}(\omega_{0}) \}> \sqrt{\nu_1+2\varrho(1+\nu_2)[1+\{2\log(2M_{\cJ}r)\}^{-1}]^2\lambda^2(M_{\cJ},r,\alpha)}\big]-o(1) \\
&~~	\geq 1   -\bP\big(n^{1/2}l_n^{-1/2}[\Re\{\hat{f}_{i_0,j_0}(\omega_0)\}-\Re\{f_{i_0,j_0}(\omega_0)\}] \leq \Delta_n \big)-o(1) \,,
\end{align*}
where
\begin{align*}
\Delta_n
=&~ \sqrt{\nu_1+2\varrho(1+\nu_2)[1+\{2\log(2M_{\cJ}r)\}^{-1}]^2\lambda^2(M_{\cJ},r,\alpha)}   \\
&- \sqrt{2\varrho(1+\nu_2)}\lambda(M_{\cJ},r,\alpha)[1+\{2\log(2M_{\cJ}r)\}^{-1}+u_n] \,.
\end{align*}
As we mentioned below \eqref{eq:I1},  $\varrho\lesssim l_n^2$. Notice that $\nu_2\asymp\log^{-2}(M_{\cJ}r)$, $\lambda(M_{\cJ},r,\alpha)=\{2\log(2M_{\cJ}r)\}^{1/2}+\{2\log(4/\alpha)\}^{1/2}$ and  $\epsilon_n^2\varrho l_n^{-2}\log^{-2}(l_n)\log^{-1}(n)\lambda^2(M_{\cJ},r,\alpha)\rightarrow\infty$, which implies
\begin{align*}
\epsilon_n
\gg \frac{l_n(\log l_n)\log^{1/2}(n)}{\varrho^{1/2}\lambda(M_{\cJ},r,\alpha)}
\gtrsim \frac{(\log l_n)\log^{1/2}(n)}{\lambda(M_{\cJ},r,\alpha)} \gg \frac{1}{\log(2M_{\cJ}r)} \gg \nu_2 \,.
\end{align*}
Since $(1+\epsilon_n)^2 = (1+\nu_2)[1+\{2\log(2M_{\cJ}r)\}^{-1}+u_n]^2$,
we have $u_n\asymp \epsilon_n$,
which implies
$ \varrho\lambda^2(M_{\cJ},r,\alpha)u_n\gg n^{-1/2}l_n^3\log^{2}(M_{\cJ}r)\asymp\nu_1$. Hence,
$
\Delta_n\leq  -C\varrho^{1/2}\lambda(M_{\cJ},r,\alpha)u_n$ for some  constant $C>0$.
Following similar arguments as Lemmas \ref{la:remainder}, \ref{la:leading} and \ref{la:uniform}, we can  also show  $|\Re\{\hat{f}_{i_0,j_0}(\omega_{0})\} - \Re\{f_{i_0,j_0}(\omega_{0})\}|=O_{\rm p} \{ n^{-1/2}l_n^{3/2}(\log l_n)\log^{1/2}(n)\}$.
Since $\varrho^{1/2}\lambda(M_{\cJ},r,\alpha)u_n \gg l_n(\log l_n)\log^{1/2}(n)$, then
\begin{align*}
n^{1/2}l_n^{-1/2}\big|\Re\{\hat{f}_{i_0,j_0}(\omega_{0})\} - \Re\{f_{i_0,j_0}(\omega_{0})\}\big| = O_{\rm p} \{l_n(\log l_n) \log^{1/2}(n)\}=o_{\rm p}\{\varrho^{1/2}\lambda(M_{\cJ},r,\alpha)u_n\} \,,
\end{align*}
which implies
$
\bP\big(n^{1/2}l_n^{-1/2}[\Re\{\hat{f}_{i_0,j_0}(\omega_0)\}-\Re\{f_{i_0,j_0}(\omega_0)\}] \leq \Delta_n 	\big)	
\rightarrow 0$.
Hence,  $\bP(T_n>\hat{\rm cv}_\alpha)\rightarrow 1$. We have completed the proof of Theorem \ref{tm:H0}(ii).
$\hfill\Box$

\section{Proof of Theorem \ref{tm:FDR}}
Notice that $Q\ll n^{f_3(\vartheta,\rho)}$ with $f_3(\vartheta,\rho)$ specified in \eqref{eq:f3.f4.f5}. Without loss of generality, we assume $Q=O(n^{\tilde{\kappa}})$ for $\tilde{\kappa}<f_3(\vartheta,\rho)$.
Recall $V_n^{(q)}=\Phi^{-1}\{1-{\rm pv}^{(q)}\}$ with ${\rm pv}^{(q)}=\bP\{\xi_{\mathcal{J}^{(q)}}^{(q)}\geq T_n^{(q)}\,|\,\mathcal{X}_n\}$
and $F_q(u)=\mathbb{P}\{T_n^{(q)}< u\}$ for any $u\in\mathbb{R}$. Note that $r_{\max}=\max_{q\in\mathcal{H}_0}r_q$, $r_{\min}=\min_{q\in\mathcal{H}_0}r_q$ and $M_{\max}=\max_{q\in\mathcal{H}_0}M_{\cJ^{(q)}}$.

We first prove that $\mathbb{P}\{V_n^{(q)}\geq t\}=1-\Phi(t)+o(1)$ holds uniformly for $q\in\mathcal{H}_0$.
Let $\bXi^{(q)}=\var\{\tilde{n}^{-1/2}\sum_{t=1}^{\tilde n}\bc_t^{(q)}\}$, where $\bc_t^{(q)}=\{\bc_{1,t}^{(q),\T},\ldots, \bc_{r_q,t}^{(q),\T}\}^{\T}$ and $\bc_{\ell,t}^{(q)}$ is defined in the same manner of \eqref{eq:cjt} but with replacing $\bchi(\cdot)$ by $\bchi^{(q)}(\cdot)$. Recall $\widehat{\bXi}^{(q)}$ is defined as \eqref{eq:Xihat} but with replacing $\hat{\bc}_t$ by $\hat{\bc}_t^{(q)}$.
Given $T_n^{(q)}$ with $q\in\mathcal{H}_0$ and some $\epsilon>0$ satisfying $ l_n^2n^{-\rho}\ll\epsilon \ll 1$, define the events
\begin{align*}
\mathcal{E}_{n,1}^{(q)}=&~ \big\{\big|1-F_q\{T_n^{(q)}\}-\bP\{\xi_{\mathcal{J}^{(q)}}^{(q)}\geq T_n^{(q)}\,|\,\mathcal{X}_n\} \big| \leq \tilde{C}_1 A_q\big\}\,, \\
\mathcal{E}_{n,2}= &~ \bigg\{\max_{q\in\mathcal{H}_0}\big|\widehat{\bXi}^{(q)}-\bXi^{(q)}\big|_\infty \leq \epsilon\bigg\}
\end{align*}
for some sufficiently large constant $\tilde{C}_1>0$ independent of $q$, where
\begin{align*}
A_q= \epsilon^{1/3}l_n^{1/3}\log^{2/3}\{M_{\cJ^{(q)}}r_q\}
+ n^{-1/9}l_n\log^{2/3}(l_n)\log\{M_{\cJ^{(q)}}r_q\}  \,.
\end{align*}
When $q\in\cH_0$,
if $\cJ^{(q)}=\{\omega_1^{(q)},\ldots,\omega_{K_q}^{(q)}\}$, parallel to \eqref{eq:gass3},
$
\sup_{u\geq 0 }|\bP\{T_n^{(q)}\leq u\}-\bP\{\xi_{\cJ^{(q)}}^{(q)}\leq u\,|\,\cX_n\}|
\lesssim
n^{-1/9}l_n\log^{2/3}(l_n)\log (K_qr_q)
+  |\widehat\bXi^{(q)}-{\bXi}^{(q)}|_\infty^{1/3}l_n^{1/3}\log^{2/3}(K_qr_q)
$
provided that
$\log(K_qr_q)\ll f_1(l_n,n;\vartheta,\rho)$
and $\max\{2,C'\log(K_qr_q)\}\leq l_n\ll  n^{f_2(\vartheta,\rho)}$
for some sufficiently large constant $C'>0$, where $f_1(l_n,n;\vartheta,\rho)$ and $f_2(\vartheta,\rho)$ are specified in \eqref{eq:f1.f2} and \eqref{eq:f1.f2.2}, respectively.
When $q\in\cH_0$, if $\cJ^{(q)}=[\omega_L^{(q)}, \omega_{U}^{(q)}]$, parallel to \eqref{eq:gass4},
$
\sup_{u\geq 0 }|\bP\{T_n^{(q)}\leq u\}-\bP\{\xi_{\cJ^{(q)}}^{(q)}\leq u\,|\,\cX_n\}|
\lesssim  n^{-1/9}l_n\log^{2/3}(l_n)\log r_q  + |\widehat{\bXi}^{(q)}-\bXi^{(q)}|_\infty^{1/3} l_n^{1/3}\log^{2/3}(r_q) 	+I\{|\widehat{\bXi}^{(q)}-\bXi^{(q)}|_\infty>l_n\}
$
provided that
$\log r_q\ll f_1(l_n,n;\vartheta,\rho)$
and $\max(2,C'\log r_q)\leq l_n\ll  n^{f_2(\vartheta,\rho)}$
for some sufficiently large constant $C'>0$.
Based on the definition of $M_{\cJ^{(q)}}$ given in \eqref{eq:MJq}, regardless of $\cJ^{(q)}=\{\omega_1^{(q)},\ldots,\omega_{K_q}^{(q)}\}$ or $\cJ^{(q)}=[\omega_L^{(q)}, \omega_{U}^{(q)}]$, for any $q\in\cH_0$, we always have
\begin{align*}
\sup_{u\geq 0 }\big|\bP\{T_n^{(q)}\leq u\}-\bP\{\xi_{\cJ^{(q)}}^{(q)}\leq u\,|\,\cX_n\}\big|
\lesssim&~
n^{-1/9}l_n\log^{2/3}(l_n)\log \{M_{\cJ^{(q)}}r_q\} \\
&+I\{|\widehat{\bXi}^{(q)}-\bXi^{(q)}|_\infty>l_n\} \\
&+  |\widehat\bXi^{(q)}-{\bXi}^{(q)}|_\infty^{1/3}l_n^{1/3}\log^{2/3}\{M_{\cJ^{(q)}}r_q\}
\end{align*}
provided that
$\log\{M_{\cJ^{(q)}}r_q\}\ll f_1(l_n,n;\vartheta,\rho)$
and $\max[2,C'\log \{M_{\cJ^{(q)}}r_q\}]\leq l_n\ll  n^{f_2(\vartheta,\rho)}$
for some sufficiently large constant $C'>0$.  Under the event $\mathcal{E}_{n,2}$, $\sup_{u\geq 0 }|\bP\{T_n^{(q)}\leq u\}-\bP\{\xi_{\cJ^{(q)}}^{(q)}\leq u\,|\,\cX_n\}|\lesssim A_q  $ for $q\in\cH_0$. Hence,
$\bP\{\mathcal{E}_{n,1}^{(q),{\rm c}}\cap \mathcal{E}_{n,2}\}=0$ for any $q\in\mathcal{H}_0$.
Define
\begin{align}\label{eq:betau}
\beta(u)=&~ nr_{\max}^2l_n^2 \exp\{-Cn^{(2\vartheta+4\rho-6\rho\vartheta-2)/(2\vartheta-1)}u^2\}
+nr_{\max}^2l_n^2\exp(-Cn^{1-\rho}l_n^{-1}u)  \notag\\
&
+ nr_{\max}^2l_n^2\exp\{-Cn^{(3-2\rho)/3}l_n^{-1}u^{2/3}\}\notag\\
&+ n^2r_{\max}^2l_n^2 \exp\{-Cn^{(\vartheta+2\rho-3\rho\vartheta-1)/(4\vartheta-2)}u^{1/2}\}  \notag\\
&+r_{\max}^2l_n^2\exp\{-Cn^{(1-\rho)/3}u^{1/3}\}
+nr_{\max}^2l_n^2 \exp\{-Cn^{(2\vartheta+3\rho-4\rho\vartheta-2)/(8\vartheta-4)}u^{1/4}\} \notag\\
& +nr_{\max}^2l_n^2\exp\{-Cn^{(2-\rho)/6}l_n^{-1/3}u^{1/6}\}
+ nr_{\max}^2l_n^2\exp\{-Cn^{(3-\rho)/9}l_n^{-1/3}u^{1/9}\}  \notag\\
& + nr_{\max}^2l_n^2\exp\{-Cn^{(4-\rho)/12}u^{1/12}\}
+r_{\max}^2l_n^2\exp\{-Cn^{(6-\rho)/18}u^{1/18}\}  \notag\\
& 	+r_{\max}^2l_n^2\exp\{-Cn^{(2\vartheta+2\rho-2\rho\vartheta-2)/(7\vartheta-4)}\}
\end{align}
for any  $u\gg l_n^2n^{-\rho}$.
Notice that we have shown $|\bXi^*-\bXi|_\infty\lesssim n^{-\rho}l_n^2$ in Section \ref{sec:pfLem11} for the proof of Lemma \ref{lem.Xi}.
Parallel to the upper bound for $\bP(|\widehat{\bXi}-\bXi^*|_\infty>u)$ given in Section \ref{subsec:upperbound},
we  have
\begin{align*}
\max_{q\in\cH_0}\bP\big\{|\widehat{\bXi}^{(q)} -\bXi^{(q)}|_\infty > u\big\}
\lesssim\beta(u)
\end{align*}
for any $u\gg l_n^2n^{-\rho}$.
Recall $V_n^{(q)}=\Phi^{-1}\{1-{\rm pv}^{(q)}\}$ with ${\rm pv}^{(q)}=\bP\{\xi_{\mathcal{J}^{(q)}}^{(q)}\geq T_n^{(q)}\,|\,\mathcal{X}_n\}$.
Notice that $1-F_q\{T_n^{(q)}\}$ follows the uniform distribution on $[0,1]$.
Then
\begin{align} \label{eq:Normal2}
\bP\big\{V_n^{(q)}\geq t\big\}
=&~\bP\big[\bP\{\xi_{\mathcal{J}^{(q)}}^{(q)}\geq T_n^{(q)}\,|\,\mathcal{X}_n\} \leq 1-\Phi(t)\big]   \notag\\
\leq &~ \bP\big[\bP\{\xi_{\mathcal{J}^{(q)}}^{(q)}\geq T_n^{(q)}\,|\,\mathcal{X}_n\}\leq 1-\Phi(t),\,\mathcal{E}_{n,1}^{(q)}\big] + \bP\{\mathcal{E}_{n,1}^{(q),{\rm c}}\cap \mathcal{E}_{n,2}\} + \bP( \mathcal{E}_{n,2}^{\rm c})     \notag\\
\leq&~ \bP\big[1-F_q\{T_n^{(q)}\}\leq 1-\Phi(t)+\tilde{C}_1A_q\big] + \bP( \mathcal{E}_{n,2}^{\rm c})   \notag\\
\leq&~ 1-\Phi(t)+\tilde{C}_1A_q +\tilde{C}_2Q_0\beta(\epsilon) \,,
\end{align}
where $\tilde{C}_2>0$ is a universal  constant independent of $q$ and $t$.
Likewise, we also have
\begin{align*}
\bP\big\{V_n^{(q)}\geq t\big\} =&~\bP\big[\bP\{\xi_{\mathcal{J}^{(q)}}^{(q)}\geq T_n^{(q)}\,|\,\mathcal{X}_n\} \leq 1-\Phi(t)\big]   \notag\\
\geq&~ \bP\big[\bP\{\xi_{\mathcal{J}^{(q)}}^{(q)}\geq T_n^{(q)}\,|\,\mathcal{X}_n\}\leq 1-\Phi(t),\,\mathcal{E}_{n,1}^{(q)}\big]     \notag\\
\geq&~ \bP\big[1-F_q\{T_n^{(q)}\}\leq 1-\Phi(t)-\tilde{C}_1A_q\big]-  \bP( \mathcal{E}_{n,2}^{\rm c})    \notag\\
\geq&~ 1-\Phi(t)-\tilde{C}_1A_q- \tilde{C}_2Q_0\beta(\epsilon)  \,.
\end{align*}
Therefore, for any  $\epsilon>0$  satisfying $ l_n^2n^{-\rho}\ll\epsilon \ll 1$,
\begin{align}\label{eq:tm3delta}
&\max_{q\in\mathcal{H}_0}	\sup_{t\in\mathbb{R}}\big|\bP\{V_n^{(q)}\geq t\} - \{1-\Phi(t)\}\big| \\
&~~~~~~~\lesssim \epsilon^{1/3}l_n^{1/3}\log^{2/3}(M_{\max}r_{\max}) + n^{-1/9}l_n\log^{2/3}(l_n)\log(M_{\max}r_{\max})    +Q_0\beta(\epsilon)
=:  \delta(\epsilon) \,.  \notag
\end{align}

Let $G(t) = 1-\Phi(t)$ and $t_{\max} = (2\log Q - 2\log\log Q)^{1/2}$. We consider two cases: (i) there exists $t\in [0,t_{\max}]$ such that $\widehat{{\rm FDP}}(t)\leq\alpha$, and (ii) $\widehat{{\rm FDP}}(t)>\alpha$ for any $t\in [0,t_{\max}]$.

{\underline{{\emph{Case {\rm (i)}.}}}}  By the definition of $\hat{t}$, it holds that $\widehat{\mathrm{FDP}}(t) > \alpha$ for any $t < \hat{t}$. Notice that $I\{V_{n}^{(q)} \geq t \} \ge I\{V_{n}^{(q)} \geq \hat{t}\}$ for $t < \hat{t}$. Then
\begin{align*}
\frac{Q G(t)}{1 \vee \sum_{q\in [Q]}I\{V_{n}^{(q)} \ge \hat{t}\}} \geq \frac{Q G(t)}{1 \vee \sum_{q\in [Q]}I\{V_{n}^{(q)} \ge t\}}=\widehat{\mathrm{FDP}}(t)> \alpha\,. 
\end{align*}
By letting $t \uparrow \hat{t}$ in the numerator of the first term in above inequality, we have $\widehat{\mathrm{FDP}}(\hat{t}) \ge \alpha$. On the other hand, based on the definition of $\hat{t}$, there exists a sequence $\{ t_i\}$ with $t_i \ge \hat{t}$ and $t_i \downarrow \hat{t}$ such that $\widehat{\mathrm{FDP}}(t_i) \le \alpha$. Thus we have $I\{V_{n}^{(q)} \geq \hat{t}\} \ge I\{V_{n}^{(q)} \geq t_i\}$, which implies that
\begin{align*}
\frac{Q G(t_i)}{1 \vee \sum_{q\in [Q]}I\{V_{n}^{(q)} \geq \hat{t}\}}\leq \frac{Q G(t_i)}{1 \vee \sum_{q\in [Q]}I\{V_{n}^{(q)} \geq t_i\}} \le \alpha\,. 
\end{align*}
Letting $t_i \downarrow \hat{t}$ in the numerator of the first term in above inequality, we have $\widehat{\mathrm{FDP}}(\hat{t}) \le \alpha$. Therefore, we have $\widehat{\mathrm{FDP}}(\hat{t}) = \alpha$ in Case (i).

\underline{\emph{Case {\rm (ii)}.}} We first show that  the threshold of $\hat{t}$ at $(2\log Q)^{1/2}$ leads to no false rejection with probability approaching one.  If $\hat{t} = (2\log Q)^{1/2}$, following from \eqref{eq:tm3delta}, we have
\begin{align*}
\bP\Bigg[ \sum_{q \in \mathcal{H}_0}I\{ V_{n}^{(q)} \ge \hat{t} \}\geq1 \Bigg]
\leq&~ Q_0 \max_{q \in \mathcal{H}_0}\bP\{ V_{n}^{(q)} \ge \hat{t} \} \\
\leq&~ Q_0 G\{(2\log Q)^{1/2} \} +CQ_0\delta(\epsilon)
\lesssim o(1)+Q\delta(\epsilon)
\end{align*}
as $Q\rightarrow\infty$.
Notice that $Q=O(n^{\tilde{\kappa}})$ for some constant $\tilde{\kappa}>0$.
By \eqref{eq:betau}, if
\begin{align*}
\epsilon \gg \max\bigg\{ \frac{\log^{2} (r_{\max})} {n^{(\vartheta+2\rho-3\rho\vartheta-1)/(2\vartheta-1)}}, \,
\frac{\log^{4}(r_{\max})} {n^{(2\vartheta+3\rho-4\rho\vartheta-2)/(2\vartheta-1)}},\,
\frac{\log^{1/2}(r_{\max})}{n^{(1-2\rho)/2}}, \,
\frac{l_n\log r_{\max}}{n^{1-\rho}},\,
\frac{l_n^2}{n^{\rho}}
\bigg\} \,,
\end{align*}
we know $Q^3\beta(\epsilon)=o(1)$. To make $Q^2\delta(\epsilon)=o(1)$, we also need to require
$\epsilon\ll l_n^{-1}\log^{-2}(M_{\max}r_{\max})Q^{-6}$ and $	Q^{18}l_n^{9}\log^6(l_n)\log^9(M_{\max} r_{\max}) \ll n$. If
\begin{align*}
\log(M_{\max}r_{\max}) \ll&~ \min\bigg\{ \frac{n^{(\vartheta+2\rho-3\rho\vartheta-1)/(8\vartheta-4)}}{l_n^{1/4}Q^{3/2}}, \,
\frac{n^{(2\vartheta+3\rho-4\rho\vartheta-2)/(12\vartheta-6)}}{l_n^{1/6}Q},\\
&~~~~~~~~~~~~~~~~~~~~~~~~~~~~~~~~	\frac{n^{(1-2\rho)/5}}{l_n^{2/5}Q^{12/5}},\,
\frac{n^{(1-\rho)/3}}{l_n^{2/3}Q^2},\,
\frac{n^{\rho/2}}{l_n^{3/2}Q^3}\bigg\} \,,
\end{align*}
we can find a suitable $\epsilon=o(1)$ such that $Q^{2}\delta(\epsilon)\rightarrow0$ as $n\rightarrow\infty$,
which implies
$
\mathbb{P}\{{\rm FDP}(\hat{t})=0\}\rightarrow1 $ in Case (ii).

It holds that
\begin{align}\label{eq:tm3p1}
\mathbb{P}\bigg\{{\rm FDP}(\hat{t})\leq\frac{\alpha Q_0}{Q}+\varepsilon\bigg\}=&~\mathbb{P}\bigg\{{\rm FDP}(\hat{t})\leq\frac{\alpha Q_0}{Q}+\varepsilon,~{\rm Case~(i)~holds}\bigg\}\notag\\
&+\mathbb{P}\bigg\{{\rm FDP}(\hat{t})\leq\frac{\alpha Q_0}{Q}+\varepsilon,~{\rm Case~(ii)~holds}\bigg\}\notag\\
=&~\mathbb{P}\big\{{\rm Case~(i)~holds}\big\}+\mathbb{P}\big\{{\rm Case~(ii)~holds}\big\} \notag\\
&-\mathbb{P}\bigg\{{\rm FDP}(\hat{t})>\frac{\alpha Q_0}{Q}+\varepsilon,~{\rm Case~(i)~holds}\bigg\}\notag\\
=&~1-\mathbb{P}\bigg\{{\rm FDP}(\hat{t})>\frac{\alpha Q_0}{Q}+\varepsilon,~{\rm Case~(i)~holds}\bigg\}\,.
\end{align}
Notice that $\widehat{{\rm FDP}}(\hat{t})=\alpha$ and $\hat{t}\in [0,t_{\max}]$ in Case (i).
As we will show in Section \ref{subsec:pftm3p2} that
\begin{align}\label{eq:tm3p2}
\sup_{t \in [0,t_{\max}]}\bigg| \frac{ \mathrm{FDP}(t)}{\widehat{\mathrm{FDP}}(t)} - \frac{Q_0}{Q} \bigg| \to 0~~\mbox{in probability}
\end{align}
as $n,Q\rightarrow\infty$, then we have $\lim_{n,Q\rightarrow\infty}\mathbb{P}\{{\rm FDP}(\hat{t})\leq\alpha Q_0/Q+\varepsilon\}=1$ by \eqref{eq:tm3p1}. Since ${\rm FDR}(t)={\bE}\{{\rm FDP}(t)\}$, then $\limsup_{n,Q\to\infty}{\rm FDR}(\hat{t})\leq\alpha Q_0/Q$. We have completed the proof of Theorem \ref{tm:FDR}. $\hfill\Box$

\subsection{Proof of \eqref{eq:tm3p2}}\label{subsec:pftm3p2}
In the sequel, we always assume $Q^{2}\delta(\epsilon)\rightarrow0$ as $n\rightarrow\infty$.
To prove \eqref{eq:tm3p2}, it suffices to show that
\begin{align*}
\sup_{t \in [0,t_{\max}]}\Bigg| \frac{ \sum_{q\in \mathcal{H}_0} [I\{V_{n}^{(q)} \geq t\} - G(t)] }{Q G(t)} \Bigg|  \to 0~~\mbox{in probability}
\end{align*}
as $n,Q\rightarrow\infty$. Let $0 = t_0 < t_1 < \cdots < t_s = t_{\max}$ such that $t_i - t_{i-1} = \tilde{v}$ for $i \in [s-1]$ and $t_s - t_{s-1} \le \tilde{v}$, where $\tilde{v} = \{(\log Q)(\log\log Q)^{1/2}\}^{-1/2}$. Then   $s \asymp t_{\max}/\tilde{v}$. For any $t\in[t_{i-1},t_i]$,
\begin{align*}
\frac{ \sum_{q\in \mathcal{H}_0} I\{V_{n}^{(q)} \geq t_{i}\} }{Q G(t_{i})}\frac{G(t_{i})}{G(t_{i-1})}  \leq \frac{ \sum_{q\in \mathcal{H}_0} I\{V_{n}^{(q)} \geq t\} }{Q G(t)}
\leq \frac{ \sum_{q\in \mathcal{H}_0} I\{V_{n}^{(q)} \geq t_{i-1}\} }{Q G(t_{i-1})}\frac{G(t_{i-1})}{G(t_i)}\,.
\end{align*}
Notice that there exists a universal constant $\tilde{C}_3>0$ such that $e^{-t^2/2}\leq \max(\tilde{C}_3,2t)\int_t^\infty e^{-x^2/2}\,{\rm d}x$ for any $t>0$. Then
\begin{align*}
0<1-\frac{G(t_i)}{G(t_{i-1})}=\frac{\int_{t_{i-1}}^{t_i}e^{-x^2/2}\,{\rm d}x}{\int_{t_{i-1}}^\infty e^{-x^2/2}\,{\rm d}x}\leq\frac{\tilde{v}e^{-t_{i-1}^2/2}}{\int_{t_{i-1}}^\infty e^{-x^2/2}\,{\rm d}x}\leq \tilde{v}\max(\tilde{C}_3,2t_{i-1})\,,
\end{align*}
which implies that
$\max_{ i\in [s]}|1-G(t_i)/G(t_{i-1})|\leq 2\tilde{v}t_{\max}\rightarrow0$
as $Q\rightarrow\infty$. Thus, to prove \eqref{eq:tm3p2}, it suffices to show that
\begin{align*}
\max_{0\leq i \leq s} \Bigg| \frac{ \sum_{q\in \mathcal{H}_0} [I\{V_{n}^{(q)} \geq t_i\} - G(t_i)] }{Q G(t_i)} \Bigg|  \to 0 ~~\mbox{in probability}
\end{align*}
as $n,Q\rightarrow\infty$. By the Bonferroni inequality and  Markov's inequality, for any $\varepsilon>0$,
\begin{align*} 
& \bP\bigg(\max_{0 \leq i \leq s} \bigg| \frac{ \sum_{q\in \mathcal{H}_0} [I\{V_{n}^{(q)} \geq t_i\} - G(t_i)] }{Q G(t_i)} \bigg| \ge \varepsilon \bigg)\notag \\
&~~~~~~~~~~\le \sum_{i=0}^{s}\bP \bigg( \bigg| \frac{ \sum_{q\in \mathcal{H}_0} [I\{V_{n}^{(q)} \geq t_i\} - G(t_i)] }{Q G(t_i)} \bigg| \ge \varepsilon \bigg) \notag\\
&~~~~~~~~~~\le \frac{\tilde{C}_4}{\tilde{v}}\int_{0}^{t_{\max}} \bP \bigg(\bigg| \frac{ \sum_{q\in \mathcal{H}_0} [I\{V_{n}^{(q)} \geq t\} - G(t)] }{Q G(t)} \bigg| \ge \varepsilon \bigg)\,\mathrm{d}t \notag\\
&~~~~~~~~~~\leq \frac{\tilde{C}_4}{\tilde{v}\varepsilon^2}\int_{0}^{t_{\max}}\mathbb{E}\bigg(\bigg| \frac{ \sum_{q\in \mathcal{H}_0} [I\{V_{n}^{(q)} \geq t\} - G(t)] }{Q G(t)} \bigg|^2\bigg)\,{\rm d}t  \,,
\end{align*}
where $\tilde{C}_4$ is some positive constant, and the second step is based on the relationship between integration and its associated Riemann sum. Thus it suffices to prove
\begin{align*}
\int_{0}^{t_{\max}}\mathbb{E}\bigg(\bigg| \frac{ \sum_{q\in \mathcal{H}_0} [I\{V_{n}^{(q)} \geq t\} - G(t)] }{Q G(t)} \bigg|^2\bigg)\,{\rm d}t = o(\tilde{v}) \,.
\end{align*}
Notice that
\begin{align*} 
&\mathbb{E}\bigg(\bigg| \frac{ \sum_{q\in \mathcal{H}_0} [I\{V_{n}^{(q)} \geq t\} - G(t)] }{Q G(t)} \bigg|^2\bigg)\notag\\
&~~~~~~~~~~\leq2\mathbb{E}\bigg(\bigg| \frac{ \sum_{q\in \mathcal{H}_0} [I\{V_{n}^{(q)} \geq t\} - \mathbb{P}\{V_{n}^{(q)} \geq t\}] }{Q G(t)} \bigg|^2\bigg)\notag\\
&~~~~~~~~~~~~~+2\bigg| \frac{ \sum_{q\in \mathcal{H}_0} [\mathbb{P}\{V_{n}^{(q)}\geq t\}-G(t)] }{Q G(t)} \bigg|^2   \notag\\
&~~~~~~~~~~= 2\sum_{q,q' \in \mathcal{H}_0} \frac{\bP\{ V_{n}^{(q)} \geq t,\, V_{n}^{(q')} \geq t\} - \bP\{ V_{n}^{(q)} \geq t\}\bP\{ V_{n}^{(q')} \geq t\}}{Q^2G^2(t)}\notag\\
&~~~~~~~~~~~~~+2\bigg| \frac{ \sum_{q\in \mathcal{H}_0} [\mathbb{P}\{V_{n}^{(q)} \geq t\}-G(t)] }{Q G(t)} \bigg|^2\,.
\end{align*}
By \eqref{eq:tm3delta}, we have
\begin{align*}
\bigg| \frac{ \sum_{q\in \mathcal{H}_0} [\mathbb{P}\{V_{n}^{(q)} \geq t\}-G(t)] }{Q G(t)} \bigg|^2 \lesssim \bigg\{\frac{Q_0\delta(\epsilon)}{QG(t)}\bigg\}^2\leq\frac{\delta^2(\epsilon)}{G^2(t)}\,.
\end{align*}
Notice that $\int_0^{t_{\max}}\{G(t)\}^{-2}\,{\rm d}t\lesssim t_{\max}\exp(t_{\max}^{2})\lesssim Q^2\log^{-3/2} (Q)$. Recall $\tilde{v} = \{(\log Q)(\log\log Q)^{1/2}\}^{-1/2}$ and $Q^2\delta(\epsilon)\rightarrow 0$. Then
\begin{align*}
\int_0^{t_{\max}}\bigg| \frac{ \sum_{q\in \mathcal{H}_0} [\mathbb{P}\{V_{n}^{(q)} \geq t\}-G(t)] }{Q G(t)} \bigg|^2 \,{\rm d}t \lesssim \int_0^{t_{\max}}\frac{\delta^2(\epsilon)}{G^2(t)}\,{\rm d}t \lesssim \frac{Q^2\delta^2(\epsilon)}{\log^{3/2} (Q)} = o(\tilde{v})\,.
\end{align*}
In the sequel, we focus on proving
\begin{align}\label{eq:toprove5}
\int_0^{t_{\max}}\sum_{q,q' \in \mathcal{H}_0} \frac{\bP\{ V_{n}^{(q)} \geq t, \,V_{n}^{(q')} \geq t\} - \bP\{ V_{n}^{(q)} \geq t\}\bP\{ V_{n}^{(q')} \geq t\}}{Q^2G^2(t)} \,{\rm d}t  = o(\tilde{v})\,.
\end{align}
Define
\begin{align*}
\mathcal{H}_{01} =&~ \{ (q, q'): q,q'\in \mathcal{H}_0,\, q = q' \} \,, \\
\mathcal{H}_{02} =&~ \{ (q,q'):q, q' \in \mathcal{H}_0,\, q \neq q', q \in \mathcal{S}_{q'}(\gamma) \mbox{~or~} q' \in \mathcal{S}_{q}(\gamma) \} \,,\\
\mathcal{H}_{03} =&~ \{(q,q'):q,q'\in\mathcal{H}_{0}\}\backslash (\mathcal{H}_{01} \cup \mathcal{H}_{02}) \,.
\end{align*}
Then $\sum_{q,q'\in\mathcal{H}_0}=\sum_{(q,q')\in\mathcal{H}_{01}}+\sum_{(q,q')\in\mathcal{H}_{02}}+\sum_{(q,q')\in\mathcal{H}_{03}}$.

\underline{\emph{Scenario {\rm 1}: $(q,q')\in\mathcal{H}_{01}$.}}
$\bP\{ V_{n}^{(q)} \geq t,\, V_{n}^{(q')} \geq t\} - \bP\{ V_{n}^{(q)} \geq t\}\bP\{ V_{n}^{(q')} \geq t\}
\leq \bP\{ V_{n}^{(q)} \geq t\}$.  By \eqref{eq:tm3delta}, it holds that
\begin{align*}
\sum_{(q,q')\in \mathcal{H}_{01}}\frac{\bP\{ V_{n}^{(q)} \geq t,\, V_{n}^{(q')} \geq t\} - \bP\{ V_{n}^{(q)} \geq t\}\bP\{ V_{n}^{(q')} \geq t\} }{Q^2G^2(t)} \lesssim \frac{1}{Q G(t)}+\frac{\delta(\epsilon)}{QG^2(t)}\,.
\end{align*}
Since $\int_0^{t_{\max}}\{G(t)\}^{-1}\,{\rm d}t\lesssim \exp(2^{-1}t_{\max}^2)=Q\log^{-1} (Q)$ and $\int_0^{t_{\max}}\{G(t)\}^{-2}\,{\rm d}t\lesssim t_{\max}\exp(t_{\max}^{2})\lesssim Q^2\log^{-3/2} (Q)$, then
\begin{align}
&\int_0^{t_{\max}}\sum_{(q,q')\in \mathcal{H}_{01}}\frac{\bP\{ V_{n}^{(q)} \geq t,\, V_{n}^{(q')} \geq t\} - \bP\{ V_{n}^{(q)} \geq t\}\bP\{ V_{n}^{(q')} \geq t\} }{Q^2G^2(t)}\,{\rm d}t\notag\\
&~~~~~~~~~~\lesssim \int_0^{t_{\max}}\frac{{\rm d}t}{QG(t)}+\delta(\epsilon) \int_0^{t_{\max}}\frac{{\rm d}t}{QG^2(t)}\lesssim \frac{1}{\log Q}+\frac{Q\delta(\epsilon)}{\log^{3/2} (Q)}=o(\tilde{v})\,.\label{eq:toprove51}
\end{align}

\underline{\emph{Scenario {\rm 2}: $(q,q')\in\cH_{02}$.}}
Since $F_q\{T_{n}^{(q)}\}\sim U[0,1]$ and $F_{q'}\{T_{n}^{(q')}\}\sim U[0,1]$,
similar to  \eqref{eq:Normal2}, we  also have
\begin{align}\label{eq:Nquantile}
&	\bP\big\{ V_{n}^{(q)} \geq t, \,  V_{n}^{(q')} \geq t \big\} \notag\\
&~~~~~~~~~	\leq \bP\big[ T_{n}^{(q)} \ge F_q^{-1}\{\Phi(t) - \tilde{C}_1A_q \},\, T_{n}^{(q')} \geq F_{q'}^{-1}\{\Phi(t) - \tilde{C}_1A_{q'} \} \big]
+2\tilde{C}_2Q_0\beta(\epsilon) \notag\\
&~~~~~~~~~	\leq \bP\big[  F_q\{T_{n}^{(q)}\} \geq \Phi(t),\, F_{q'}\{T_{n}^{(q')}\} \geq \Phi(t) \big]
+\tilde{C}_1A_q+\tilde{C}_1A_{q'}
+ 2\tilde{C}_2Q_0\beta(\epsilon)  \notag\\
&~~~~~~~~~	\leq \bP\big[ \zeta^{(q)} \geq t , \,\zeta^{(q')} \geq t \big]+C\delta(\epsilon) \,,
\end{align}
where $\zeta^{(q)} = \Phi^{-1}[F_{q}\{T_{n}^{(q)}\}]\sim \mathcal{N}(0,1)$ and $\zeta^{(q')} = \Phi^{-1}[F_{q'}\{T_{n}^{(q')}\}]\sim \mathcal{N}(0,1)$.
When $(q,q')\in\mathcal{H}_{02}$, since $\max_{q,q' \in [Q]}|\mathrm{Corr}\{\zeta^{(q)}, \zeta^{(q')}\}| \leq r_{\zeta}<1$,   Lemma 2 in \cite{Berman_1962_supp}  
implies that $\bP\{ \zeta^{(q)}\ge t,\, \zeta^{(q')} \ge t \}\lesssim t^{-2}\exp\{-t^2/(1+r_\zeta)\}$ for any $t>\tilde{C}_5$ for some universal constant $\tilde{C}_5>0$. Notice that $e^{-t^2/2}\leq \max(\tilde{C}_3,2t)\int_t^\infty e^{-x^2/2}\,{\rm d}x$ for any $t>0$. By \eqref{eq:Nquantile},
$
\bP\{ V_{n}^{(q)} \geq t,\, V_{n}^{(q')} \geq t\} \lesssim t^{-2} \exp\{-t^2/(1+r_\zeta)\} +\delta(\epsilon)
\lesssim  t^{-2r_{\zeta}/(1+r_\zeta)} \{G(t)\}^{2/(1+r_\zeta)} +\delta(\epsilon)
\lesssim \{G(t)\}^{2/(1+r_\zeta)} +\delta(\epsilon)$ for any $t>\max(1,\tilde{C}_5,\tilde{C}_3/2)$.
Since $\max_{q \in [Q]}|\mathcal{S}_q(\gamma)| = o(Q^{\nu})$, then $|\mathcal{H}_{02}| = O(Q^{1+\nu})$. Due to $\nu<(1-r_\zeta)/(1+r_\zeta)<1$ and $\bP\{ V_{n}^{(q)} \geq t, \, V_{n}^{(q')} \geq t\} \leq1$ for any $0<t<\max(1,\tilde{C}_5, \tilde{C}_3/2)$,
\begin{align}
&\int_0^{t_{\max}}\sum_{(q,q')\in\mathcal{H}_{02}}\frac{\bP\{V_{n}^{(q)}\geq t,V_{n}^{(q')}\geq t\} -\bP\{V_{n}^{(q)}\geq t\}\bP\{V_{n}^{(q')}\geq t\}}{Q^2G^2(t)}\,{\rm d}t\notag \\
&~~~~~~~~ \leq \int_0^{t_{\max}}\sum_{(q,q')\in\mathcal{H}_{02}}\frac{\bP\{V_{n}^{(q)}\geq t,V_{n}^{(q')}\geq t\}} {Q^2G^2(t)}\,{\rm d}t\notag \\
&~~~~~~~~ \lesssim  \int_0^{\max(1,\tilde{C}_5,\tilde{C}_3/2)} \frac{{\rm d}t}{Q^{1-\nu}G^2(t)}  + \int_{\max(1,\tilde{C}_5,\tilde{C}_3/2)}^{t_{\max}} \frac{{\rm d}t}{Q^{1-\nu}\{G(t)\}^{2r_\zeta/(1+r_\zeta)}} \notag\\
&~~~~~~~~~~~+\delta(\epsilon)\int_0^{t_{\max}}\frac{{\rm d}t}{Q^{1-\nu}G^2(t)}\notag\\
&~~~~~~~~ \lesssim \frac{1}{Q^{1-\nu}} + \frac{Q^{(r_\zeta-1)/(1+r_\zeta)+\nu}}{\log^{(1+3r_\zeta)/(2+2r_\zeta)} (Q)}+\frac{Q^{1+\nu}\delta(\epsilon)}{\log^{3/2} (Q)} =o(\tilde{v})\,.  \label{eq:toprove52}
\end{align}

\underline{\emph{Scenario {\rm 3}: $(q,q')\in\mathcal{H}_{03}$.}}
Let $\rho_{q,q'}={\rm Corr}\{\zeta^{(q)},\zeta^{(q')}\}$. By Theorem 2.1.e of 
\cite{LinBai_2010_supp},
\begin{equation}\label{eq:tm3p11}
\bP\big\{ \zeta^{(q)}\geq t,\, \zeta^{(q')}\geq t \big\}
\leq
\left\{\begin{aligned}
	G(t)G\bigg\{\frac{(1-\rho_{q,q'})t}{(1-\rho_{q,q'}^2)^{1/2}}\bigg\}\,, ~~~~~~~~~~&\mbox{if} ~ -1 < \rho_{q,q'} \le 0\,. \\
	(1+\rho_{q,q'})G(t)G\bigg\{\frac{(1-\rho_{q,q'})t}{(1-\rho_{q,q'}^2)^{1/2}}\bigg\}\,, ~~~~&\mbox{if} ~ 0 \le \rho_{q,q'} < 1\,.
\end{aligned}\right.
\end{equation}
Note that $|\rho_{q,q'}|\leq\log^{-2-\gamma} (Q)$ for any $(q,q')\in\mathcal{H}_{03}$. When $-\log^{-2-\gamma} (Q) \le \rho_{q,q'} \le 0$, due to $(1-\rho_{q,q'})/(1-\rho_{q,q'}^2)^{1/2} \ge 1$, we have $G\{(1-\rho_{q,q'})t/(1-\rho_{q,q'}^2)^{1/2}\} \le G(t)$, which implies that $\bP\{ \zeta^{(q)} \geq t,\, \zeta^{(q')} \geq t \} \le G^2(t) \le \{1+\log^{-1-\gamma} (Q)\}G^2(t)$.
When $0 <\rho_{q,q'} \le \log^{-2-\gamma} (Q)$, by the mean-value theorem,
there exists $\tilde{t}$ satisfying $(1-\rho_{q,q'})t/(1-\rho_{q,q'}^2)^{1/2}<\tilde{t}<t$ such that
\begin{align*}
G\bigg\{\frac{(1-\rho_{q,q'})t}{(1-\rho_{q,q'}^2)^{1/2}}\bigg\} = G(t)+\phi(\tilde{t})\bigg\{t-\frac{(1-\rho_{q,q'})t}{(1-\rho_{q,q'}^2)^{1/2}} \bigg\}
\end{align*}
where $\phi(\cdot)$ is the density function of the standard normal distribution $\mathcal{N}(0,1)$. Then if $\rho_{q,q'}>0$,
\begin{align*}
G\bigg\{\frac{(1-\rho_{q,q'})t}{(1-\rho_{q,q'}^2)^{1/2}}\bigg\}\{G(t)\}^{-1} = &~ 1+\frac{\phi(\tilde{t})}{G(t)}\bigg\{t-\frac{(1-\rho_{q,q'})t}{(1-\rho_{q,q'}^2)^{1/2}} \bigg\} \\
\leq &~ 1+\frac{t\phi(\tilde{t})}{t\phi(t)/(1+t^2)}\bigg\{1-\frac{1-\rho_{q,q'}}{(1-\rho_{q,q'}^2)^{1/2}}\bigg\} \\
\leq &~ 1+\frac{\phi(\tilde{t})}{\phi(t)}\cdot2\rho_{q,q'}(1+t^2)
\end{align*}
for any $t>0$. For $0<t<t_{\max}$, it holds that
\begin{align*}
\frac{\phi(\tilde{t})}{\phi(t)}
=&~\exp\bigg\{\frac{(t-\tilde{t})(t+\tilde{t})}{2}\bigg\}
\leq\exp\bigg[t_{\max}^2\bigg\{1-\frac{1-\rho_{q,q'}}{(1-\rho_{q,q'}^2)^{1/2}}\bigg\}\bigg] \\
&~ \leq\exp(2\rho_{q,q'}t_{\max}^2)\leq \exp\{4\log^{-1-\gamma} (Q)\}\lesssim 1 \,,
\end{align*}
which implies that $G\{(1-\rho_{q,q'})t/(1-\rho_{q,q'}^2)^{1/2}\} \le G(t)[1+O\{\log^{-1-\gamma} (Q)\}]$ for any $t\in[0,t_{\max}]$ if $0<\rho_{q,q'}\leq\log^{-2-\gamma} (Q)$, where the term $O\{\log^{-1-\gamma} (Q)\}$ holds uniformly over $t\in[0,t_{\max}]$. Then $(1+\rho_{q,q'})G(t)G\{(1-\rho_{q,q'})t/(1-\rho_{q,q'}^2)^{1/2}\}\leq G^2(t)[1+O\{\log^{-1-\gamma} (Q)\}]$ for any $t\in[0,t_{\max}]$ if $0<\rho_{q,q'}\leq\log^{-2-\gamma} (Q)$. By \eqref{eq:tm3p11},  $\bP\{ \zeta^{(q)} \geq t,\, \zeta^{(q')}\geq t \} \le [1+O\{\log^{-1-\gamma} (Q)\}]G^2(t)$ for $t \in [0,t_{\max}]$. Thus we have
\begin{align*}
\max_{(q,q') \in \mathcal{H}_{03}}\bP\{ V_{n}^{(q)}\geq t,\, V_{n}^{(q')}\geq t\} \le [1+O\{\log^{-1-\gamma} (Q)\}]G^2(t) + C\delta(\epsilon)
\end{align*}
for $t \in [0,t_{\max}]$. Due to $\delta(\epsilon)\int_0^{t_{\max}}\{G(t)\}^{-2}\,{\rm d}t\lesssim Q^2\delta(\epsilon)\log^{-3/2} (Q)=o(\tilde{v})$, it then holds that
\begin{align*}
&\int_0^{t_{\max}}\sum_{(q,q')\in \mathcal{H}_{03}}\frac{\bP\{V_{n}^{(q)}\geq t,\,V_{n}^{(q')}\geq t\} -\bP\{V_{n}^{(q)}\geq t\}\bP\{V_{n}^{(q')}\geq t\} }{Q^2G^2(t)}\,{\rm d}t \\
&~~~~~~~~~~~~~ \leq O\{\log^{-1-\gamma} (Q)\}\cdot\int_0^{t_{\max}}1\,{\rm d}t
+C \int_0^{t_{\max}} \frac{\delta(\epsilon)}{G^2(t)}\,{\rm d}t =o(\tilde{v})\,.
\end{align*}
Together with \eqref{eq:toprove51} and \eqref{eq:toprove52}, we know \eqref{eq:toprove5} holds. Then we can obtain  \eqref{eq:tm3p2} holds.
$\hfill\Box$

\section{Proof of auxillary lemmas}\label{sec:pflems}
\subsection{Proof of Lemma \ref{la:remainder}}
For ${\rm II}(\omega)$ defined in \eqref{eq:expansion}, due to $\bF(\omega)=(2\pi)^{-1}\sum_{k=-\infty}^{\infty}\bGamma(k)e^{-\iota k\omega}$, it holds that
\begin{align*}
{\rm II}(\omega) =\underbrace{\frac{1}{2\pi}\sum_{k=-l_n}^{l_n}\bigg\{\mathcal{W}\bigg(\frac{k}{l_n}\bigg) \frac{n-|k|}{n}-1\bigg\}\bGamma(k)e^{-\iota k\omega}}_{A(\omega)}-\underbrace{\frac{1}{2\pi}\sum_{|k|>l_n}\bGamma(k)e^{-\iota k\omega}}_{B(\omega)}  \,.
\end{align*}
Recall $\gamma_{i,j}(k)=\cov(x_{i,t+k},x_{j,t})$. By Conditions \ref{as:moment} and \ref{as:betamixing}, it follows from Davydov's inequality  that $\max_{i,j\in[p]}|\gamma_{i,j}(k)|\lesssim\exp(-C|k|)$ for any integer $k$. Thus, as $l_n\rightarrow\infty$, we have
\begin{align*}
\sup_{ {\omega\in[-\pi,\pi]} }|B(\omega)|_\infty
=&~\sup_{ {\omega\in[-\pi,\pi]} }\max_{i,j\in[p]}\bigg|\frac{1}{2\pi}\sum_{|k|>l_n}\gamma_{i,j}(k)e^{-\iota k\omega}\bigg|  \\
\lesssim&~\sum_{k>l_n}\exp(-Ck) \lesssim\int_{l_n}^\infty\exp(-Cu)\,{\rm d}u\lesssim \exp(-Cl_n)\,.
\end{align*}
Notice that the flat-top kernel function $\cW(u)$ satisfies that   $\cW(u)=1$ if $|u|\leq c$, $\cW(u)=(|u|-1)/(c-1)$ if $c<|u|\leq 1$ and $\cW(u)=0$ otherwise, where $c\in(0,1]$ is a constant. Then
\begin{align*}
\sup_{ {\omega\in[-\pi,\pi]} }|A(\omega)|_\infty
\leq&~\frac{1}{2\pi}\sum_{k=-\lfloor cl_n\rfloor}^{\lfloor cl_n\rfloor}\bigg\{1-\mathcal{W}\bigg(\frac{k}{l_n}\bigg)\frac{n-|k|}{n}\bigg\} \max_{i,j\in[p]}|\gamma_{i,j}(k)|+\frac{1}{2\pi}\sum_{|k|>\lfloor cl_n\rfloor}\max_{i,j\in[p]}|\gamma_{i,j}(k)|\\
\lesssim&~n^{-1}+\exp(-Cl_n)
\end{align*}
as $l_n\rightarrow\infty$, which implies
$\sup_{ {\omega\in[-\pi,\pi]} }\max_{(i,j)\in\cI}|{\rm II}_{i,j}(\omega)|\lesssim n^{-1}+\exp(-Cl_n)$.

Write  $\bar{\mathring{\bx}}=(\bar{\mathring{x}}_{1}, \ldots, \bar{\mathring{x}}_{p})^{\T}=n^{-1}\sum_{t=1}^n\mathring{\bx}_t$.
Notice that $(2\pi)^{-1}\sum_{k=1}^{l_n}\mathcal{W}({k}/{l_n})\asymp l_n$.
For ${\rm III}(\omega)$ defined in \eqref{eq:expansion},
by the triangle inequality and the Bonferroni inequality, we have
\begin{align*}
&~ \mathbb{P}\bigg[\sup_{ {\omega\in[-\pi,\pi]} }\max_{(i,j)\in\mathcal{I}}|\Re\{{\rm III}_{i,j}(\omega)\}|>\frac{u}{\sqrt{2}} \bigg]
\leq  r\max_{(i,j)\in\mathcal{I}} \mathbb{P}\bigg[\sup_{ {\omega\in[-\pi,\pi]} }|\Re\{{\rm III}_{i,j}(\omega)\}|>\frac{u}{\sqrt{2}} \bigg] \\
\leq &~ r\max_{(i,j)\in\mathcal{I}} \mathbb{P}\bigg\{\sup_{ {\omega\in[-\pi,\pi]} }|\bar{\mathring{x}}_i| \bigg|\frac{1}{2\pi}\sum_{k=1}^{l_n}\mathcal{W}\bigg(\frac{k}{l_n}\bigg) \bigg(\frac{1}{n}\sum_{t=1}^{n-k}\mathring{x}_{j,t}+\frac{1}{n}\sum_{t=k+1}^n\mathring{x}_{j,t}\bigg)\cos(k\omega)\bigg| >\frac{u}{\sqrt{2}}\bigg\} \\
\leq &~ r\max_{(i,j)\in\mathcal{I}} \mathbb{P}\bigg[ |\bar{\mathring{x}}_i| \bigg\{\frac{1}{2\pi}\sum_{k=1}^{l_n}\mathcal{W}\bigg(\frac{k}{l_n}\bigg)\bigg\} \bigg(\max_{k\in[l_n]}\bigg|\frac{1}{n}\sum_{t=1}^{n-k}\mathring{x}_{j,t}\bigg|\bigg)>\frac{u}{2\sqrt{2}}\bigg] \\
&~ + r\max_{(i,j)\in\mathcal{I}} \mathbb{P}\bigg[ |\bar{\mathring{x}}_i| \bigg\{\frac{1}{2\pi}\sum_{k=1}^{l_n}\mathcal{W}\bigg(\frac{k}{l_n}\bigg)\bigg\} \bigg(\max_{k\in[l_n]}\bigg|\frac{1}{n}\sum_{t=k+1}^n\mathring{x}_{j,t}\bigg|\bigg)>\frac{u}{2\sqrt{2}}\bigg] \\
\leq &~ 2r\max_{(i,j)\in\mathcal{I}} \mathbb{P}\bigg(|\bar{\mathring{x}}_i|>\frac{Cu^{1/2}}{l_n^{1/2}}\bigg) +rl_n\max_{(i,j)\in\mathcal{I}}\max_{k\in[l_n]} \mathbb{P}\bigg(\bigg|\frac{1}{n}\sum_{t=1}^{n-k}\mathring{x}_{j,t}\bigg| >\frac{Cu^{1/2}}{l_n^{1/2}}\bigg) \\
&~ +rl_n\max_{(i,j)\in\mathcal{I}}\max_{k\in[l_n]} \mathbb{P}\bigg(\bigg|\frac{1}{n}\sum_{t=k+1}^n\mathring{x}_{j,t}\bigg| >\frac{Cu^{1/2}}{l_n^{1/2}}\bigg)
\end{align*}
for any $u>0$. Since $l_n=o(n)$, Theorem 1 of \cite{MPR_2011_supp} 
yields that
\begin{align*}
& \max_{(i,j)\in\mathcal{I}} \mathbb{P}\bigg(|\bar{\mathring{x}}_i|>\frac{Cu^{1/2}}{l_n^{1/2}}\bigg)
+\max_{(i,j)\in\mathcal{I}}\max_{k\in[l_n]} \mathbb{P}\bigg( \bigg|\frac{1}{n}\sum_{t=1}^{n-k}\mathring{x}_{j,t}\bigg| >\frac{Cu^{1/2}}{l_n^{1/2}}\bigg)  \\
&~~~~~~~~~~~ + \max_{(i,j)\in\mathcal{I}}\max_{k\in[l_n]} \mathbb{P}\bigg( \bigg|\frac{1}{n}\sum_{t=k+1}^n\mathring{x}_{j,t}\bigg| >\frac{Cu^{1/2}}{l_n^{1/2}}\bigg)\\
&~~~~~~\lesssim  n\exp(-Cn^{2/3}l_n^{-1/3}u^{1/3})  + \exp(-Cnl_n^{-1}u)
\end{align*}
for any $u>0$ satisfying $nl_n^{-1/2}u^{1/2}\rightarrow\infty$, which implies
\begin{align*}
\mathbb{P}\bigg[\sup_{ {\omega\in[-\pi,\pi]} }\max_{(i,j)\in\mathcal{I}}|\Re\{{\rm III}_{i,j}(\omega)\}|>\frac{u}{\sqrt{2}} \bigg]
\lesssim rnl_n\exp(-Cn^{2/3}l_n^{-1/3}u^{1/3})  + rl_n\exp(-Cnl_n^{-1}u)
\end{align*}
for any $u>0$ satisfying $nl_n^{-1/2}u^{1/2}\rightarrow\infty$.
Analogously, we also have
\begin{align*}
\mathbb{P}\bigg[\sup_{ {\omega\in[-\pi,\pi]} }\max_{(i,j)\in\mathcal{I}}|\Im\{{\rm III}_{i,j}(\omega)\}|>\frac{u}{\sqrt{2}} \bigg]
\lesssim rnl_n\exp(-Cn^{2/3}l_n^{-1/3}u^{1/3})  + rl_n\exp(-Cnl_n^{-1}u)
\end{align*}
for any $u>0$ satisfying $nl_n^{-1/2}u^{1/2}\rightarrow\infty$.
By the Bonferroni inequality,
\begin{align}\label{eq:IIItail}
\mathbb{P}\bigg\{\sup_{ {\omega\in[-\pi,\pi]} }\max_{(i,j)\in\mathcal{I}}|{\rm III}_{i,j}(\omega)|>u \bigg\} \leq&~   \mathbb{P}\bigg[\sup_{ {\omega\in[-\pi,\pi]} }\max_{(i,j)\in\mathcal{I}}|\Re\{{\rm III}_{i,j}(\omega)\}|>\frac{u}{\sqrt{2}} \bigg]\notag\\
&+ \mathbb{P}\bigg[\sup_{ {\omega\in[-\pi,\pi]} }\max_{(i,j)\in\mathcal{I}}|\Im\{{\rm III}_{i,j}(\omega)\}|>\frac{u}{\sqrt{2}} \bigg] \\
\lesssim&~ rnl_n\exp(-Cn^{2/3}l_n^{-1/3}u^{1/3})  + rl_n\exp(-Cnl_n^{-1}u)\notag
\end{align}
for any $u>0$ satisfying $nl_n^{-1/2}u^{1/2}\rightarrow\infty$, which implies
\begin{align*}
\sup_{ {\omega\in[-\pi,\pi]} }\max_{(i,j)\in\mathcal{I}}|{\rm III}_{i,j}(\omega)| = O_{\rm p}(l_n n^{-1}\log r )
\end{align*}
provided that $\log r=O(n^{1/2})$. Using the same arguments, we also have
\begin{align*}
\sup_{ {\omega\in[-\pi,\pi]} }\max_{(i,j)\in\mathcal{I}}|{\rm IV}_{i,j}(\omega)|
=O_{\rm p} (l_nn^{-1}\log r )
=\sup_{ {\omega\in[-\pi,\pi]} }\max_{(i,j)\in\mathcal{I}}|{\rm V}_{i,j}(\omega)|
\end{align*}
provided that $\log r=O(n^{1/2})$.
We have completed the proof of Lemma \ref{la:remainder}. $\hfill\Box$

\subsection{Proof of Lemma \ref{la:leading}}\label{subsec:pfLem2}
For ${\rm I}(\omega)$ defined in \eqref{eq:expansion},
we can reformulate ${\rm I}_{i,j}(\omega)$ as follows:
\begin{align*} 
{\rm I}_{i,j}(\omega)=\frac{1}{2\pi n} \sum_{t=1}^n\sum_{k=-\min(t-1,l_n)}^{\min(l_n,n-t)} \mathcal{W}\bigg(\frac{k}{l_n}\bigg)\{\mathring{x}_{i,t+k}\mathring{x}_{j,t}-\gamma_{i,j}(k)\}e^{-\iota k\omega} \,.
\end{align*}
For $\zeta_{i,j}(\omega)$ defiined as \eqref{eq:zetaij}, due to $n>2l_n$,
by the triangle inequality,
\begin{align*}
|{\rm I}_{i,j}(\omega)-\zeta_{i,j}(\omega)|
\leq &~ \bigg|\frac{1}{2\pi n} \sum_{t=1}^{l_n}\sum_{k=1-t}^{l_n} \mathcal{W}\bigg(\frac{k}{l_n}\bigg)\{\mathring{x}_{i,t+k}\mathring{x}_{j,t}-\gamma_{i,j}(k)\}e^{-\iota k\omega}\bigg| \\
&~ +\bigg|\frac{1}{2\pi n} \sum_{t=n-l_n+1}^{n}\sum_{k=-l_n}^{n-t} \mathcal{W}\bigg(\frac{k}{l_n}\bigg)\{\mathring{x}_{i,t+k}\mathring{x}_{j,t}-\gamma_{i,j}(k)\}e^{-\iota k\omega}\bigg|  \,.
\end{align*}
Based on  the facts $\sum_{t=1}^{l_n}\sum_{k=1-t}^{l_n}\mathcal{W}(k/l_n)\lesssim l_n^2$ and $\sum_{t=n-l_n+1}^{n}\sum_{k=-l_n}^{n-t}\mathcal{W}(k/l_n)\lesssim l_n^2$, by the Bonferroni inequality,
\begin{align*}
&~ \mathbb{P}\bigg\{\sup_{ {\omega\in[-\pi,\pi]} }\max_{(i,j)\in\mathcal{I}}|{\rm I}_{i,j}(\omega)-\zeta_{i,j}(\omega)|>u\bigg\} \\
\leq &~ \mathbb{P}\bigg\{\sup_{ {\omega\in[-\pi,\pi]} }\max_{(i,j)\in\mathcal{I}}\bigg|\frac{1}{2\pi n} \sum_{t=1}^{l_n}\sum_{k=1-t}^{l_n} \mathcal{W}\bigg(\frac{k}{l_n}\bigg) \{\mathring{x}_{i,t+k}\mathring{x}_{j,t}-\gamma_{i,j}(k)\}e^{-\iota k\omega}\bigg|>\frac{u}{2}\bigg\} \\
&~ + \mathbb{P}\bigg\{\sup_{ {\omega\in[-\pi,\pi]} }\max_{(i,j)\in\mathcal{I}} \bigg|\frac{1}{2\pi n} \sum_{t=n-l_n+1}^{n}\sum_{k=-l_n}^{n-t} \mathcal{W}\bigg(\frac{k}{l_n}\bigg) \{\mathring{x}_{i,t+k}\mathring{x}_{j,t}-\gamma_{i,j}(k)\}e^{-\iota k\omega}\bigg|>\frac{u}{2}\bigg\} \\
\leq &~ \mathbb{P}\bigg\{\max_{(i,j)\in\mathcal{I}}\frac{1}{2\pi n} \sum_{t=1}^{l_n}\sum_{k=1-t}^{l_n} \mathcal{W}\bigg(\frac{k}{l_n}\bigg) |\mathring{x}_{i,t+k}\mathring{x}_{j,t}-\gamma_{i,j}(k)|>\frac{u}{2}\bigg\} \\
&~ + \mathbb{P}\bigg\{\max_{(i,j)\in\mathcal{I}} \frac{1}{2\pi n} \sum_{t=n-l_n+1}^{n}\sum_{k=-l_n}^{n-t} \mathcal{W}\bigg(\frac{k}{l_n}\bigg)|\mathring{x}_{i,t+k}\mathring{x}_{j,t}-\gamma_{i,j}(k)|>\frac{u}{2}\bigg\} \\
\leq &~ r \sum_{t=1}^{l_n}\sum_{k=1-t}^{l_n}\max_{(i,j)\in\mathcal{I}} \mathbb{P}\{|\mathring{x}_{i,t+k}\mathring{x}_{j,t}-\gamma_{i,j}(k)|>Cnl_n^{-2}u\} \\
&~  +r \sum_{t=n-l_n+1}^{n}\sum_{k=-l_n}^{n-t}\max_{(i,j)\in\mathcal{I}} \mathbb{P}\{|\mathring{x}_{i,t+k}\mathring{x}_{j,t}-\gamma_{i,j}(k)|>Cnl_n^{-2}u\}  \\
\lesssim &~ rl_n^2 \max_{t}\max_{-l_n\leq k\leq l_n}\max_{(i,j)\in\mathcal{I}} \mathbb{P}\{|\mathring{x}_{i,t+k}\mathring{x}_{j,t}-\gamma_{i,j}(k)|>Cnl_n^{-2}u\}
\end{align*}
for any $u>0$. By Lemma 2 in the supplementary material of \cite{CTW_2013_supp} 
and Condition \ref{as:moment}, it holds that
\begin{align} \label{eq:t1}
\max_{t}\max_{-l_n\leq k\leq l_n}\max_{(i,j)\in\mathcal{I}} \mathbb{P}\{|\mathring{x}_{i,t+k}\mathring{x}_{j,t}-\gamma_{i,j}(k)|>u\}\leq C\exp(-Cu)
\end{align}
for any $u>0$, which implies
\begin{align}\label{eq:diff.I.zeta}
\mathbb{P}\bigg\{\sup_{ {\omega\in[-\pi,\pi]} }\max_{(i,j)\in\mathcal{I}}|{\rm I}_{i,j}(\omega)-\zeta_{i,j}(\omega)|>u\bigg\}
\lesssim rl_n^2 \exp(-C nl_n^{-2}u)
\end{align}
for any $u>0$.
Recall $r\geq n^{\kappa}$ for some sufficiently small constant $\kappa>0$ and $n>2l_n$.
Then
\[
\sup_{ {\omega\in[-\pi,\pi]} }\max_{(i,j)\in\mathcal{I}}|{\rm I}_{i,j}(\omega)-\zeta_{i,j}(\omega)|=O_{\rm p}(l_n^2n^{-1}\log r) \,.
\]
We have completed the proof of Lemma \ref{la:leading}. $\hfill\Box$

\subsection{Proof of Lemma \ref{la:tailprob}}\label{subsec:pfLem3}
Recall
$z_{i,j,t}^{(1)}(\omega)=l_n^{-1}\sum_{k=-l_n}^{l_n}\mathcal{W}({k}/{l_n})\{\mathring{x}_{i,t+k}\mathring{x}_{j,t}-\gamma_{i,j}(k)\}\cos(k\omega)$ and $\sum_{k=-l_n}^{l_n}\mathcal{W}(k/l_n)\asymp l_n$. By the Bonferroni inequality and \eqref{eq:t1}, we have
\begin{align} \label{eq:zijtail}
\max_{t}\max_{i,j\in[p]}\mathbb{P}\big\{\big|z_{i,j,t}^{(1)}(\omega)\big|>u\big\}
\leq&\, \sum_{k=-l_n}^{l_n}\mathbb{P}\{|\mathring{x}_{i,t+k}\mathring{x}_{j,t}-\gamma_{i,j}(k)|>Cu\}  \notag\\
\leq&\, \tilde{C}_1l_n\exp(-\tilde{C}_2u)
\end{align}
for any $u>0$, where $\tilde{C}_1$ and $\tilde{C}_2$ are two positive constants. Select $\tilde{C}_1^*=\max(2,\tilde{C}_1)$ and $\tilde{C}_2^*=(\tilde{C}_2\log 2)/2$. Notice that $l_n\rightarrow\infty$ as $n\rightarrow\infty$.
When $l_n\geq 2$, we consider two scenarios:

(i) if $u\geq \log^2 (l_n)/(\tilde{C}_2\log l_n - \tilde{C}_2^*)$, then it holds that
\[
\log l_n - \tilde{C}_2u\leq -\frac{\tilde{C}_2^*u}{\log l_n} ~~\Longrightarrow~~ \tilde{C}_1l_n\exp(-\tilde{C}_2u) \leq \tilde{C}_1^* \exp\bigg(-\frac{\tilde{C}_2^*u}{\log l_n }\bigg) \,.
\]

(ii) if $0<u< \log^2 (l_n)/(\tilde{C}_2\log l_n - \tilde{C}_2^*)$, then it holds that
\[
-\frac{\tilde{C}_2^*u}{\log l_n} \geq - \frac{(\log 2)(\log l_n)}{2\log l_n-\log 2} \geq -\log 2
~~\Longrightarrow~~ 1 \leq \tilde{C}_1^* \exp\bigg(-\frac{\tilde{C}_2^*u}{\log l_n }\bigg) \,.
\]
Therefore, we have
\begin{align}\label{eq:zij.tail}
\max_{t}\max_{i,j\in[p]}\mathbb{P}\{|z_{i,j,t}^{(1)}(\omega)|>u\} \leq \tilde{C}_1^* \exp\{-\tilde{C}_2^*(\log l_n )^{-1}u\}
\end{align}
for any $u>0$. Since  $\{z_{i,j,t}^{(1)}(\omega)\}_{t=l_n+s_1}^{l_n+s_2}$ is an $\alpha$-mixing sequence with $\alpha$-mixing coefficients $\{\alpha_{z}(k)\}_{k\geq 1}$ satisfying $\alpha_{z}(k)\leq C\exp(-C|k-2l_n|_+)$ for any integer $k\geq 1$, applying Lemma L1 in the supplementary material of \cite{CCW_2021_supp} 
with $\tilde{B}_{\tilde{n}}=\log l_n$, $\tilde{L}_{\tilde{n}}=1$, $\tilde{j}_{\tilde{n}}=2l_n$, $r_1=1$, $r_2=1$ and $r=1/3$ to obtain that
\begin{align*}
\sup_{ {\omega\in[-\pi,\pi]} }\max_{i,j\in[p]}\bP\bigg(\bigg|\sum_{t=l_n+s_1}^{l_n+s_2}z_{i,j,t}^{(1)}(\omega)\bigg|\geq u\bigg) \lesssim \exp\bigg\{-\frac{Cu^2}{sl_n\log^2 (l_n)}\bigg\} + \exp\bigg\{-\frac{Cu^{1/3}}{l_n^{1/3}\log^{1/3} (l_n)}\bigg\}
\end{align*}
for any $u>0$, which implies
\begin{align*}
\sup_{ {\omega\in[-\pi,\pi]} }\max_{i,j\in[p]}\mathbb{P}\bigg\{\bigg|\frac{l_n}{s} \sum_{t=l_n+s_1}^{l_n+s_2}z_{i,j,t}^{(1)}(\omega)\bigg|>u\bigg\}
\lesssim \exp\bigg\{-\frac{Csu^2}{l_n^3\log^2( l_n)}\bigg\} +\exp\bigg\{-\frac{Cs^{1/3}u^{1/3}}{l_n^{2/3}\log^{1/3} (l_n)}\bigg\}  \,.
\end{align*}
Analogously, we can also show the same tail probability holds for $|l_ns^{-1} \sum_{t=l_n+s_1}^{l_n+s_2}z_{i,j,t}^{(2)}(\omega)|$ with $z_{i,j,t}^{(2)}(\omega)$ specified in \eqref{eq:zijt}. We have completed the proof of Lemma \ref{la:tailprob}. $\hfill\Box$

\subsection{Proof of Lemma \ref{la:uniform}}\label{subsec:pfLem4}
Applying the Bonferroni inequality, by \eqref{eq:IIItail} and \eqref{eq:diff.I.zeta}, it holds that
\begin{align} \label{eq:I-V}
&~ \mathbb{P}\bigg\{\sup_{ {\omega\in[-\pi,\pi]} }\max_{(i,j)\in\mathcal{I}}|{\rm I}_{i,j}(\omega)-\zeta_{i,j}(\omega)+{\rm III}_{i,j}(\omega)+{\rm IV}_{i,j}(\omega)+{\rm V}_{i,j}(\omega)|>\frac{2u}{3}\bigg\}  \notag\\	
&~~~~~~~~~~~~~~~ \lesssim
rnl_n\exp\bigg(-\frac{Cn^{2/3}u^{1/3}}{l_n^{1/3}}\bigg)
+ rl_n^2 \exp\bigg(-\frac{C nu}{l_n^{2}}\bigg)
\end{align}
for any $u\gg n^{-2}l_n$. As shown in the proof of Lemma \ref{la:remainder}, $\sup_{ {\omega\in[-\pi,\pi]} }\max_{(i,j)\in\cI}|{\rm II}_{i,j}(\omega)|\lesssim n^{-1}+\exp(-Cl_n)$. For any $u\gg \max\{n^{-1},\exp(-Cl_n)\}$, we have
\begin{align}\label{eq:II}
\mathbb{P}\bigg\{\sup_{ {\omega\in[-\pi,\pi]} }\max_{(i,j)\in\mathcal{I}} |{\rm II}_{i,j}(\omega)|>\frac{u}{6}\bigg\} = 0\,.
\end{align}
By \eqref{eq:expansion} and the triangle inequality,
\begin{align*}
\sup_{ {\omega\in[-\pi,\pi]} }\max_{i,j\in\cI}|\hat{f}_{i,j}(\omega)-f_{i,j}(\omega)|
\leq&~ \sup_{ {\omega\in[-\pi,\pi]} }\max_{(i,j)\in\mathcal{I}}|{\rm I}_{i,j}(\omega)-\zeta_{i,j}(\omega)+{\rm III}_{i,j}(\omega)+{\rm IV}_{i,j}(\omega)+{\rm V}_{i,j}(\omega)| \\
&~ +\sup_{ {\omega\in[-\pi,\pi]} }\max_{(i,j)\in\mathcal{I}} |{\rm II}_{i,j}(\omega)| +\sup_{ {\omega\in[-\pi,\pi]} }\max_{(i,j)\in\mathcal{I}} |\zeta_{i,j}(\omega)|\,.
\end{align*}
Together with \eqref{eq:I-V} and \eqref{eq:II},
by the Bonferroni inequality,
\begin{align}\label{eq:diff.hatf.f.1}
\mathbb{P}\bigg\{\sup_{ {\omega\in[-\pi,\pi]} }\max_{(i,j)\in\mathcal{I}}|\hat{f}_{i,j}(\omega) - f_{i,j}(\omega)|>u\bigg\}
\lesssim &~ \mathbb{P}\bigg\{\sup_{ {\omega\in[-\pi,\pi]} }\max_{(i,j)\in\mathcal{I}} |\zeta_{i,j}(\omega)|>\frac{u}{6}\bigg\}  \\
& + rnl_n\exp\bigg(-\frac{Cn^{2/3}u^{1/3}}{l_n^{1/3}}\bigg)
+ rl_n^2 \exp\bigg(-\frac{C nu}{l_n^{2}}\bigg)     \notag
\end{align}
for any $u\gg\max\{n^{-1},\exp(-Cl_n)\}$.

Let $-\pi=\omega_1^*<\cdots<\omega_{M+1}^*=\pi$ be the isometric partition of $[-\pi,\pi]$ with $\delta_M=2\pi/M\rightarrow0$. Write $B_1=[\omega_1^*,\omega_2^*]$ and $B_m=(\omega_{m}^*,\omega_{m+1}^*]$ for  $m\in\{2,\ldots,M\}$. For any $\omega\in[-\pi,\pi]$, there exists $m_{\omega}\in[M]$  such that $\omega\in B_{m_\omega}$ and $|\omega-\omega_{m_\omega}^*|\leq2\pi/M$. Recall $r=|\cI|$.

By the triangle inequality and the Bonferroni inequality,
\begin{align} \label{eq:zeta}
\mathbb{P}\bigg\{\sup_{{\omega\in[-\pi,\pi]} }\max_{(i,j)\in\mathcal{I}} |\zeta_{i,j}(\omega)|>\frac{u}{6}\bigg\}
\leq&~  r\max_{(i,j)\in\mathcal{I}}\mathbb{P}\bigg\{\sup_{{\omega\in[-\pi,\pi]} }|\zeta_{i,j}(\omega)| >\frac{u}{6}\bigg\}  \notag\\
\leq &~ r\max_{(i,j)\in\mathcal{I}}\mathbb{P}\bigg\{\sup_{{\omega\in[-\pi,\pi]} }|\zeta_{i,j}(\omega)| -\max_{m\in[M+1]}|\zeta_{i,j}(\omega_m^*)|>\frac{u}{12}\bigg\}  \notag\\
& + r\max_{(i,j)\in\mathcal{I}}\mathbb{P}\bigg\{\max_{m\in[M+1]}|\zeta_{i,j}(\omega_m^*)| >\frac{u}{12}\bigg\}  \notag\\
\lesssim &~ rM\max_{(i,j)\in\mathcal{I}}\max_{m\in[M]}\mathbb{P}\bigg\{\sup_{\omega\in B_m} |\zeta_{i,j}(\omega)-\zeta_{i,j}(\omega_{m_\omega}^*)|>\frac{u}{12}\bigg\}   \notag\\
& +rM\max_{(i,j)\in\mathcal{I}}\max_{m\in[M+1]}\mathbb{P}\bigg\{|\zeta_{i,j}(\omega_m^*)| >\frac{u}{12}\bigg\}
\end{align}
for any $u>0$. Recall $l_n=o(n)$ and
\begin{align*}
\zeta_{i,j}(\omega) = \frac{l_n}{2\pi n}\sum_{t=l_n+1}^{n-l_n}z_{i,j,t}^{(1)}(\omega) - \iota \bigg\{\frac{l_n}{2\pi n}\sum_{t=l_n+1}^{n-l_n}z_{i,j,t}^{(2)}(\omega)\bigg\}
\end{align*}
with $\iota=\sqrt{-1}$, $z_{i,j,t}^{(1)}(\omega)$ and $z_{i,j,t}^{(2)}(\omega)$ specified in \eqref{eq:zijt}.
Applying Lemma \ref{la:tailprob} with $s_1=1$ and $s_2=n-2l_n$, if $l_n\geq 2$, we have
\begin{align}\label{eq:zetam}
&\max_{(i,j)\in\mathcal{I}}\max_{m\in[M+1]}\mathbb{P}\bigg\{|\zeta_{i,j}(\omega_m^*)| >\frac{u}{12}\bigg\} \notag\\
&~~~~~~~~~~~~~~
\lesssim \exp\bigg\{-\frac{Cnu^2}{l_n^3\log^2 (l_n)}\bigg\} +\exp\bigg\{-\frac{Cn^{1/3}u^{1/3}}{l_n^{2/3}
	\log^{1/3} (l_n)}\bigg\}
\end{align}
for any $u>0$. 
For any $\omega\in[-\pi,\pi]$, by the triangle inequality, 
\begin{align*}
|\zeta_{i,j}(\omega)-\zeta_{i,j}(\omega_{m_{\omega}}^*)|	
\leq&~ \frac{1}{2\pi n}\sum_{t=l_n+1}^{n-l_n}\sum_{k=-l_n}^{l_n} \mathcal{W}\bigg(\frac{k}{l_n}\bigg)|k||\omega-\omega_{m_\omega}^*||\mathring{x}_{i,t+k}\mathring{x}_{j,t}-\gamma_{i,j}(k)| \\
\leq&~ \frac{Cl_n^2}{M} \max_{-l_n\leq k\leq l_n}\bigg\{\frac{1}{n}\sum_{t=l_n+1}^{n-l_n} |\mathring{x}_{i,t+k}\mathring{x}_{j,t}-\gamma_{i,j}(k)| \bigg\} \,,
\end{align*}
where the last step is due to $\sum_{k=-l_n}^{l_n}|k|\mathcal{W}(k/l_n)\asymp l_n^2$. By the Bonferroni inequality and \eqref{eq:t1}, 
\begin{align*}
&~ \max_{(i,j)\in\mathcal{I}}\max_{m\in[M]}\mathbb{P}\bigg\{\sup_{\omega\in B_m} |\zeta_{i,j}(\omega)-\zeta_{i,j}(\omega_{m_\omega}^*)|>\frac{u}{12}\bigg\} \\
\leq &~ \max_{(i,j)\in\mathcal{I}}\max_{m\in[M]}\mathbb{P}\bigg[\frac{Cl_n^2}{M} \max_{-l_n\leq k\leq l_n}\bigg\{\frac{1}{n}\sum_{t=l_n+1}^{n-l_n} |\mathring{x}_{i,t+k}\mathring{x}_{j,t}-\gamma_{i,j}(k)| \bigg\} > \frac{u}{12} \bigg] \\
\lesssim &~ l_n\max_{(i,j)\in\mathcal{I}}\max_{-l_n\leq k\leq l_n} \mathbb{P}\bigg\{\frac{1}{n}\sum_{t=l_n+1}^{n-l_n} |\mathring{x}_{i,t+k}\mathring{x}_{j,t}-\gamma_{i,j}(k)| > \frac{CMu}{l_n^2} \bigg\} \\
\lesssim &~ nl_n \max_{t}\max_{(i,j)\in\mathcal{I}}\max_{-l_n\leq k\leq l_n} \mathbb{P}\{|\mathring{x}_{i,t+k}\mathring{x}_{j,t}-\gamma_{i,j}(k)| > CMl_n^{-2}u \}\\
\lesssim &~ nl_n\exp(-CMl_n^{-2}u)
\end{align*}
for any $u>0$. Together with \eqref{eq:zetam}, selecting $M\asymp n$,
\eqref{eq:zeta} implies that  
\begin{align}\label{eq:zetatail}
\mathbb{P}\bigg\{\sup_{{\omega\in[-\pi,\pi]} }\max_{(i,j)\in\mathcal{I}} |\zeta_{i,j}(\omega)|>\frac{u}{6}\bigg\}
\lesssim &~     rn\exp\bigg\{-\frac{Cnu^2}{l_n^3\log^2 (l_n)}\bigg\} +rn\exp\bigg\{-\frac{Cn^{1/3}u^{1/3}}{l_n^{2/3}\log^{1/3} (l_n)}\bigg\}  \notag\\
&  + rn^2l_n\exp\bigg(-\frac{Cnu}{l_n^2}\bigg)
\end{align}
for any $u>0$. Together with \eqref{eq:diff.hatf.f.1},
due to $l_n=o(n)$,
if $l_n\log l_n=o(n) $ and $l_n\geq 2$,
we have
\begin{align*} 
&\mathbb{P}\bigg\{\sup_{ {\omega\in[-\pi,\pi]} }\max_{(i,j)\in\mathcal{I}}|\hat{f}_{i,j}(\omega) - f_{i,j}(\omega)|>u\bigg\}   \notag\\
&~~~~~~~ \lesssim rnl_n\exp\bigg\{-\frac{Cn^{1/3}u^{1/3}}{l_n^{2/3}\log^{1/3} (l_n)}\bigg\}
+  rn\exp\bigg\{-\frac{Cnu^2}{l_n^3\log^2 (l_n)}\bigg\}  + rn^2l_n\exp\bigg(-\frac{Cnu}{l_n^2}\bigg)
\end{align*}
for any $u\gg\max\{ n^{-1},\exp(-Cl_n)\}$.
Notice that $r\geq n^\kappa$ for some sufficiently small constant $\kappa>0$. If $\log r=O(n^{1/5}l_n^{-1/5})$, $l_n\log l_n=o(n) $
and $l_n\geq \tilde{C}\log n$ for some sufficiently large constant $\tilde{C}>0$, we have
\begin{align*}
\sup_{{\omega\in[-\pi,\pi]} }\max_{(i,j)\in\mathcal{I}}|\hat{f}_{i,j}(\omega) - f_{i,j}(\omega)| =O_{\rm p}\bigg\{\frac{l_n^{3/2}(\log l_n) \log^{1/2}(r)}{n^{1/2}}\bigg\} \,.
\end{align*}
We have completed the proof of Lemma \ref{la:uniform}. $\hfill\Box$

\subsection{Proof of Lemma \ref{la:app1}}\label{subsec:pfLem5}
Recall
$
T_n(\omega;\cI)
=nl_n^{-1}\max_{(i,j)\in\cI}|\hat{f}_{i,j}(\omega)-f_{i,j}(\omega)|^2$ and
$	\check{T}_n(\omega;\mathcal{I})
=nl_n^{-1}\max_{(i,j)\in\cI}|\zeta_{i,j}(\omega)|^2 $. By the triangle inequality,
\begin{align}\label{eq:diff.Tn.Tncheck}
\bigg|\sup_{\omega\in\mathcal{J}}T_n(\omega;\mathcal{I}) -\sup_{\omega\in\mathcal{J}}\check{T}_n(\omega;\mathcal{I})\bigg|
\leq&~ \sup_{\omega\in\mathcal{J}}|T_n(\omega;\mathcal{I}) -\check{T}_n(\omega;\mathcal{I})|  \notag\\
\leq&~ nl_n^{-1}\sup_{\omega\in\cJ}\max_{(i,j)\in\cI}|\hat{f}_{i,j}(\omega)-f_{i,j}(\omega)-\zeta_{i,j}(\omega)|^2  \\
& + 2nl_n^{-1}\sup_{\omega\in\cJ}\max_{(i,j)\in\cI}|\zeta_{i,j}(\omega)||\hat{f}_{i,j}(\omega)-f_{i,j}(\omega)-\zeta_{i,j}(\omega)| \,. \notag
\end{align}
By \eqref{eq:expansion} and the triangle inequality,
$|\hat{f}_{i,j}(\omega)-f_{i,j}(\omega)-\zeta_{i,j}(\omega)|
\leq |{\rm I}_{i,j}(\omega)-\zeta_{i,j}(\omega)+{\rm III}_{i,j}(\omega)+{\rm IV}_{i,j}(\omega)+{\rm V}_{i,j}(\omega)|+|{\rm II}_{i,j}(\omega)| $.
By the Bonferroni inequality, it follows from \eqref{eq:I-V} and \eqref{eq:II} that
\begin{align*}
&\bP\bigg\{\sup_{\omega\in\cJ}\max_{(i,j)\in\cI}|\hat{f}_{i,j}(\omega)-f_{i,j}(\omega)-\zeta_{i,j}(\omega)|>u\bigg\}   \\
&~~~~~~~~~~~~~~ \lesssim rnl_n\exp\bigg(-\frac{Cn^{2/3}u^{1/3}}{l_n^{1/3}}\bigg)
+ rl_n^2 \exp\bigg(-\frac{C nu}{l_n^{2}}\bigg)     \notag
\end{align*}
for any $u\gg\max\{n^{-1},\exp(-Cl_n)\}$. Together with \eqref{eq:zetatail}, for any constant $\bar{C}\in(0,1)$, if  $\log r=O(n^{1/5}l_n^{-1/5})$, $l_n\log l_n=o(n) $
and
$l_n\geq \max(2,\tilde{C}\log n)$ for some sufficiently large constant $\tilde{C}>0$, by the Bonferroni inequality,  \eqref{eq:diff.Tn.Tncheck} yields
\begin{align}  \label{eq:Tn.Tnc.tail}
& \mathbb{P}\bigg\{\bigg|\sup_{\omega\in\mathcal{J}}T_n(\omega;\mathcal{I}) -\sup_{\omega\in\mathcal{J}}\check{T}_n(\omega;\mathcal{I})\bigg|>u\bigg\}  \notag\\
&~~~~~~~~~~ \leq
\mathbb{P}\bigg\{\sup_{\omega\in\cJ}\max_{(i,j)\in\cI}|\hat{f}_{i,j}(\omega)-f_{i,j}(\omega)-\zeta_{i,j}(\omega)|>\sqrt{\frac{l_nu}{2n}}\bigg\}   \notag\\
&~~~~~~~~~~~~~ +\mathbb{P}\bigg\{\sup_{\omega\in\mathcal{J}}\sup_{(i,j)\in\mathcal{I}} |\zeta_{i,j}(\omega)|>\frac{l_n^{3/2}(\log l_n)\log^{1/2}(r)}{4\bar{C}n^{1/2}}\bigg\} \notag\\
&~~~~~~~~~~~~~ +\mathbb{P}\bigg\{\sup_{\omega\in\cJ}\max_{(i,j)\in\cI}|\hat{f}_{i,j}(\omega)-f_{i,j}(\omega)-\zeta_{i,j}(\omega)|>\frac{\bar{C}u}{n^{1/2}l_n^{1/2}(\log l_n) \log^{1/2}(r)}\bigg\} \notag\\
&~~~~~~~~~~ \lesssim rnl_n\exp\bigg(-\frac{Cn^{1/2}u^{1/6}}{l_n^{1/6}}\bigg)+rnl_n\exp\bigg\{-\frac{Cn^{1/2}u^{1/3}}{l_n^{1/2}\log^{1/3}(l_n)\log^{1/6}(r)}\bigg\}
\notag\\
&~~~~~~~~~~~~~ +rl_n^2\exp\bigg(-\frac{Cn^{1/2}u^{1/2}}{l_n^{3/2}}\bigg) + rl_n^2\exp\bigg\{-\frac{Cn^{1/2}u}{l_n^{5/2}(\log l_n)\log^{1/2}(r)}\bigg\} + n^{-C_*}
\end{align}
for $u\gg n^{-1/2}l_n^{1/2}(\log l_n)\log^{1/2}(r)$,
where $C_*>0$ is a constant only depending on $\bar{C}$ such that $C_*\rightarrow\infty$ as $\bar{C}\rightarrow0$.  We can select a specified $\bar{C}>0$ such that $C_*>1$. Therefore,
\begin{align*}
\bigg|\sup_{\omega\in\mathcal{J}}T_n(\omega;\mathcal{I}) -\sup_{\omega\in\mathcal{J}}\check{T}_n(\omega;\mathcal{I})\bigg|
=O_{\rm p} \bigg\{\frac{l_n^{5/2}(\log l_n) \log^{3/2}(r)}{n^{1/2}}\bigg\}  \,.
\end{align*}
We have completed the proof of Lemma \ref{la:app1}. $\hfill\Box$

\subsection{Proof of Lemma \ref{la.cont}}
For given $t_0\in I$, we consider  the Gaussian random process $\{Z(t):t\in S\}$ with $Z(t)=X(t_0+t)-X(t_0)$, where $S=\{t\in\bR:|t|\leq a\}$ and $t_0+t\in I$. Since $d(s,t)\leq c_1|s-t|^\lambda$ for all $s,t\in I$, then
\begin{align*}
d_Z^2(s,t)\equiv \bE\{|Z(s)-Z(t)|^2\} =d^2(t_0+s,t_0+t) \leq c_1^2|s-t|^{2\lambda} \,, ~~~\forall s,t\in S \,.
\end{align*}
Let $N_{d_Z}(S,\varepsilon)$ denote the smallest number of (open) $d_Z$-balls of radius $\varepsilon$ needed to cover $S$. Write $D_Z=\sup\{d_Z(s,t):s,t\in S\}$. Notice that $D_Z\leq c_1(2a)^\lambda$ and
$  N_{d_Z}(S,\varepsilon) \leq \max\{1,2ac_1^{1/\lambda} \varepsilon^{-1/\lambda}\}$. Then
\begin{align*}
\int_0^{D_Z}\sqrt{\log N_{d_Z}(S,\varepsilon)}\,{\rm d}\varepsilon
\leq&~ \int_0^{c_1(2 a)^\lambda} \sqrt{\log\bigg\{\max\bigg(1,\, \frac{2ac_1^{1/\lambda}}{\varepsilon^{1/\lambda}}\bigg)\bigg\}}\,{\rm d}\varepsilon \\
=&~ c_12^{\lambda-1}\lambda^{-1/2}\pi^{1/2} a^{\lambda}\,.
\end{align*}
Let $c_*=c_12^{\lambda-1}\lambda^{-1/2}\pi^{1/2}$. By Lemma 5.3 of \cite{MWX_2013_supp}, 
there exists a universal constant $C>0$ such that for any $x>0$,
\begin{align*}
\bP\bigg\{\sup_{s,t\in S}|Z(s)-Z(t)|\geq C(x+c_*a^\lambda)\bigg\} \leq \exp\bigg\{-\frac{x^2}{(c_12^\lambda a^\lambda)^2}\bigg\} \,.
\end{align*}
For any $u>2Cc_*$, we know $a^\lambda u\geq C\{(a^\lambda u)/(2C)+c_*a^\lambda\}$, which implies
\begin{align*}
\bP\bigg\{\sup_{s,t\in S}|Z(s)-Z(t)|\geq a^\lambda u\bigg\}
\leq \exp\bigg(-\frac{u^2}{4^{\lambda+1}C^2c_1^2}\bigg)
\end{align*}
We have completed the proof of Lemma \ref{la.cont}. $\hfill\Box$

\subsection{Proof of Lemma \ref{la:xtildetail}}\label{sec:pflem.xtildetail}
Recall
$\bH=\{\bI_r\otimes \bA^\T(\omega_1),\ldots,\bI_r\otimes \bA^\T(\omega_K)\}^\T\in \bR^{(2Kr)\times (2l_n+1)r}$.
Write $\bH=(\bh_1,\ldots,\bh_{2Kr})^\T$, where each $\bh_j$ is a $r(2l_n+1)$-dimensional vector.
Since $\tilde{\ba}_{\ell}\equiv (\tilde{a}_{1,\ell},\ldots,\tilde{a}_{2Kr,\ell})^{\T} =b^{-1/2}\sum_{t\in\cI_{\ell}}\bH\bc_t$,
then $\tilde{a}_{j,\ell}=b^{-1/2}\sum_{t\in\mathcal{I}_\ell}\bh_j^\T\bc_t$. For $j=1$, we have
$
\tilde{a}_{1,\ell}
= (2\pi)^{-1}b^{-1/2}l_n^{1/2} \sum_{t\in\cI_\ell}z_{\bchi(1),t}^{(1)} (\omega_1)
$
with $z_{i,j,t}^{(1)} (\omega)$ specified in \eqref{eq:zijt}.
By Lemma \ref{la:tailprob} with $s=b$, if  $l_n\geq 2$, then
\begin{align*}
\mathbb{P}(|\tilde{a}_{1,\ell}|>\lambda)
\lesssim
\exp\bigg\{-\frac{C\lambda^2}{l_n^{2}\log^{2}(l_n)}\bigg\}
+\exp\bigg\{-\frac{Cb^{1/6}\lambda^{1/3}}{l_n^{1/2}\log^{1/3}(l_n)}\bigg\}
\end{align*}
for any $\lambda>0$.  Applying the identical arguments, we  know the above inequality also holds for any $j\in[2Kr]$ and $\ell\in[L]$. We have completed the proof of Lemma \ref{la:xtildetail}. $\hfill\Box$

\subsection{Proof of Lemma \ref{la:sigmatilde}}
Write $\bc_t=\{c_{1,t},\ldots,c_{r(2l_n+1),t}\}^{\T}$. By \eqref{eq:t1}, $\max_{t,j}\bP(|c_{j,t}|>u)\leq C\exp(-Cu)$ for any $u>0$.
Recall  $\bXi=\tilde{n}^{-1}\bE\{(\sum_{t=1}^{\tilde n}\bc_t)(\sum_{t=1}^{\tilde n}\bc_t)^{\T}\}$ and $\widetilde{\bXi}=L^{-1}\sum_{\ell=1}^{L}\bE(\tilde\bc_\ell\tilde{\bc}_\ell^{\T})$ with $\tilde{\bc}_\ell=b^{-1/2}\sum_{t\in\cI_\ell}\bc_t$.
By Lemma L3 in  the supplementary material of  \cite{CJS_2021_supp} 
with $r_1=r_2=1$ and $m=2l_n$, we have $|\widetilde{\bXi}-\bXi|_\infty\lesssim l_nhb^{-1}+l_nbn^{-1}$.
Write $\bH=(\bh_1,\ldots,\bh_{2Kr})^\T$, where each $\bh_j$ is a $r(2l_n+1)$-dimensional vector.
Since $\bSigma=\bH\bXi\bH^{\T}$ and $\widetilde\bSigma=\bH\widetilde{\bXi}\bH^{\T}$, due to $	|\bh_j|_1\leq l_n^{-1/2}\sum_{k=-l_n}^{l_n}\mathcal{W}(k/l_n)\asymp l_n^{1/2} $, then
$  |\widetilde{\bSigma}-\bSigma|_\infty\leq |\widetilde{\bXi}-\bXi|_\infty \max_{i\in [2Kr]}|\bh_i|_1^2
\lesssim l_n^2(hb^{-1}+bn^{-1})$.
We have completed the proof of Lemma \ref{la:sigmatilde}. $\hfill\Box$
%

\subsection{Proof of Lemma \ref{pn:1}}
Recall $\{ \by_t\}_{t=1}^{\tilde{n}}$ is a sequence of independent normal random vectors.
Let $\mathcal{W}_n=\{\bw_1,\ldots,\bw_{\tilde{n}}\}$ be a copy of $\mathcal{Y}_n=\{\by_1,\ldots,\by_{\tilde{n}}\}$. Write $\mathcal{A}_n=\{\ba_1,\ldots,\ba_{\tilde{n}}\}$. Assume $\mathcal{A}_n$, $\mathcal{Y}_n$ and $\mathcal{W}_n$ are independent.
Recall $\bs_{n,\by}^{(1)}=L^{-1/2}\sum_{\ell=1}^L \tilde \by_\ell$ with $\tilde \by_\ell=b^{-1/2}\sum_{t\in\cI_\ell} \by_t$, where $\by_t\sim \mathcal{N}\{\bzero,\bE(\tilde \ba_\ell \tilde \ba_\ell^\T)\}$ for any $t\in\cI_\ell$. Let $\bs_{n,\bw}^{(1)}=L^{-1/2}\sum_{\ell=1}^L \tilde \bw_\ell$ with $\tilde \bw_\ell=b^{-1/2}\sum_{t\in\cI_\ell} \bw_t$. Then
\begin{align*}
\varrho_n^{(1)}=\sup_{\bu\in\mathbb{R}^{2Kr},\nu\in[0,1]}\big|\mathbb{P}\{\sqrt{\nu}\bs_{n,\ba}^{(1)}+\sqrt{1-\nu}\bs_{n,\by}^{(1)}\leq \bu\}-\mathbb{P}\{\bs_{n,\bw}^{(1)}\leq \bu\}\big| \,.
\end{align*}
For any $\phi>0$, let $\beta=\phi\log(2Kr)$ and
define
\begin{align*}
M_{\tilde{\ba}}(\phi)=&~\max_{\ell\in[L]}\mathbb{E}(|\tilde{\ba}_\ell |_\infty^3{I}[|\tilde{\ba}_\ell|_\infty>L^{1/2}\{4\phi\log(2Kr)\}^{-1}]) \,,\\
M_{\tilde{\by}}(\phi)=&~\max_{\ell\in[L]}\mathbb{E}(|\tilde{\by}_\ell |_\infty^3{I}[|\tilde{\by}_\ell|_\infty>L^{1/2}\{4\phi\log(2Kr)\}^{-1}]) \,.
\end{align*}
For a given $\bu=(u_1,\ldots,u_{2Kr})^\T\in\mathbb{R}^{2Kr}$, define
$
F_\beta(\bv)=\beta^{-1}\log[\sum_{j=1}^{2Kr}\exp\{\beta(v_j-u_j)\}]$
for any $\bv=(v_1,\ldots,v_{2Kr})^\T\in\mathbb{R}^{2Kr}$. Such defined function $F_{\beta}(\bv)$ satisfies the property $0\leq F_{\beta}(\bv)-\max_{j\in[2Kr]}(v_j-u_j)\leq \beta^{-1}\log(2Kr)=\phi^{-1}$ for any $\bv\in\mathbb{R}^{2Kr}$. Select a thrice continuously differentiable function $g_0:\mathbb{R}\rightarrow[0,1]$ whose derivatives up to the third order are all bounded such that $g_0(t)=1$ for $t\leq 0$ and $g_0(t)=0$ for $t\geq1$. Define $g(t)=g_0(\phi t)$ for any $t\in\mathbb{R}$, and $q(\bv)=g\{F_\beta(\bv)\}$ for any $\bv\in\mathbb{R}^{2Kr}$.
Let
$\mathcal{T}_n=q\{\sqrt{\nu}\bs_{n,\ba}^{(1)}+\sqrt{1-\nu}\bs_{n,\by}^{(1)}\} -q\{\bs_{n,\bw}^{(1)}\}$.
Since $\cov\{\bs_{n,\bw}^{(1)}\}=L^{-1}\sum_{\ell=1}^L\bE(\tilde{\ba}_\ell\tilde{\ba}_\ell^{\T})=\widetilde{\bSigma}$, by Lemma \ref{la:sigmatilde} and Condition \ref{as:eigen}, we know all the elements of the main-diagonal of $\widetilde{\bSigma}$ are uniformly bounded away from 0 provided that $l_n^2(hb^{-1}+bn^{-1})=o(1)$. Hence, following the same arguments in the proof of Lemma 3 in  the supplementary material of \cite{CCW_2021_supp},  
we have
\begin{align}\label{eq:rhon1.exp}
\varrho_n^{(1)}
\lesssim&~\phi^{-1}\log^{1/2}(2Kr)+ \sup_{\bu\in\mathbb{R}^{2Kr},\nu\in[0,1]}|\mathbb{E}(\mathcal{T}_n)|\,.
\end{align}
By Lemma \ref{la:xtildetail}, since $b\gg h>2l_n$,  if $l_n\geq 2$,
$\max_{\ell\in[L]}\max_{j\in[2Kr]}\bE(|\tilde{a}_{j,\ell}|^3)
\lesssim l_n^3\log^3(l_n) $.
Since $\tilde{y}_{j,\ell}$ is normal random variable,
$\mathbb{E}(|\tilde{y}_{j,\ell}|^3)\lesssim \{\mathbb{E}(|\tilde{y}_{j,\ell}|^2)\}^{3/2}=\{\mathbb{E}(|\tilde{a}_{j,\ell}|^2)\}^{3/2}
\leq \mathbb{E}(|\tilde{a}_{j,\ell}|^3)$. Parallel to Equation (S.4) in  the supplementary material of \cite{CCW_2021_supp}, 
\begin{align}\label{eq:Tnbound}
	\sup_{\bu\in\mathbb{R}^{2Kr},\nu\in[0,1]}|\mathbb{E}(\mathcal{T}_n)|
	\lesssim&~\phi L^{1/2}\max_{\ell\in[L]}\mathbb{E}\bigg\{\max_{j\in[2Kr]}| \mathbb{E}(\tilde{a}_{j,\ell}\,|\,\mathcal{F}_{-\ell})|\bigg\} \notag\\
	&+\phi^2\log(2Kr)\max_{\ell\in[L]}\mathbb{E}\bigg[\max_{k,j\in[2Kr]}| \mathbb{E}\{\tilde{a}_{k,\ell}\tilde{a}_{j,\ell}-\mathbb{E}(\tilde{a}_{k,\ell}\tilde{a}_{j,\ell})\,|\,\mathcal{F}_{-\ell}\}|\bigg]  \notag\\
	&+\frac{\phi^3\log^2(2Kr)}{L^{1/2}} \max_{\ell\in[L]}\mathbb{E}\bigg[\max_{j\in[2Kr]} \big|\mathbb{E}\big\{|\tilde{a}_{j,\ell}|^3-\mathbb{E}\big(|\tilde{a}_{j,\ell}|^3\big)\,|\,\mathcal{F}_{-\ell}\big\}\big|\bigg]  \notag\\
	&+\frac{l_n^{3}\log^3(l_n)\phi^3\log^2(2Kr)}{L^{1/2}}\big\{\phi^{-1}\log^{1/2}(2Kr)+\varrho_n^{(1)}\big\} \notag\\
	&+\frac{\phi^3\log^2(2Kr)}{L^{1/2}}\big\{M_{\tilde{\ba}}(\phi)+M_{\tilde{\by}}(\phi)\big\}\,,
\end{align}
where $\mathcal{F}_{-\ell}$ is the $\sigma$-filed generated by $\{\tilde{\ba}_s\}_{s\neq \ell}$.

Selecting $\phi=C'L^{1/6}l_n^{-1}\log^{-1}(l_n)\log^{-2/3}(2Kr)$ for some sufficiently small constant $C'>0$, \eqref{eq:rhon1.exp} and \eqref{eq:Tnbound} imply
\begin{align}
\varrho_n^{(1)}\lesssim&~\frac{l_n(\log l_n)\log^{7/6}(2Kr)}{L^{1/6}}
+\frac{1}{l_n^{3}\log^3(l_n)}M_{\tilde{\ba}}\bigg\{\frac{C'L^{1/6}}{l_n(\log l_n)\log^{2/3}(2Kr)}\bigg\}  \notag\\
&~+\frac{1}{l_n^{3}\log^3(l_n)}M_{\tilde{\by}}\bigg\{\frac{C'L^{1/6}}{l_n(\log l_n)\log^{2/3}(2Kr)}\bigg\} \notag\\
&~+\frac{L^{2/3}}{l_n(\log l_n)\log^{2/3}(2Kr)}\max_{\ell\in[L]}\mathbb{E}\bigg\{\max_{j\in[2Kr]} |\mathbb{E}(\tilde{a}_{j,\ell}\,|\,\mathcal{F}_{-\ell})|\bigg\} \label{eq:rho1}\\
&~+\frac{L^{1/3}}{l_n^2\log^2(l_n)\log^{1/3}(2Kr)}\max_{\ell\in[L]}\mathbb{E}\bigg[\max_{k,j\in[2Kr]}| \mathbb{E}\{\tilde{a}_{k,\ell}\tilde{a}_{j,\ell} -\mathbb{E}(\tilde{a}_{k,\ell}\tilde{a}_{j,\ell})\,|\,\mathcal{F}_{-\ell}\}|\bigg] \notag\\
&~+\frac{1}{l_n^{3}\log^3(l_n)}\max_{\ell\in[L]}\mathbb{E}\bigg[\max_{j\in[2Kr]} \big|\mathbb{E}\big\{|\tilde{a}_{j,\ell}|^3-\mathbb{E}\big(|\tilde{a}_{j,\ell}|^3\big)\,|\,\mathcal{F}_{-\ell}\big\}\big|\bigg]\,. \notag
\end{align}
Denote by $\sigma_{\ell,j,j}^2$ the $(j,j)$-th element of $\mathbb{E}(\tilde{\ba}_{\ell}\tilde{\ba}_\ell^\T)$. By Lemma \ref{la:xtildetail}, we know $\max_{\ell\in[L]}\max_{j\in[2Kr]} \sigma_{\ell,j,j}^2\lesssim l_n^2\log^2 (l_n)$.
It is elementary to verify that
\begin{align*}
\frac{1}{l_n^{3}\log^3(l_n)}M_{\tilde{\by}}\bigg\{\frac{C'L^{1/6}}{l_n(\log l_n)\log^{2/3}(2Kr)}\bigg\}
\lesssim&~ \frac{l_n(\log l_n)\log^{7/6}(Kr)}{L^{1/6}}
\end{align*}
provided that $\log(Kr)=o(L^{2/5})$. By the Bonferroni inequality,  Lemma 2 in the supplementary material of \cite{CTW_2013_supp}  
and Lemma \ref{la:xtildetail}, if $l_n\geq 2$, it holds that
\begin{align*}
&~~~~~~~~~~~~~\max_{\ell\in[L]}\mathbb{P}(|\tilde{\ba}_{\ell}|_\infty>u)
\lesssim Kr\exp\bigg\{-\frac{Cu^2}{l_n^{2}\log^{2}(l_n)}\bigg\} +Kr\exp\bigg\{-\frac{Cb^{1/6}u^{1/3}}{l_n^{1/2}\log^{1/3}(l_n)}\bigg\} \,, \\
&~~\max_{\ell\in[L]}\max_{k,j\in[2Kr]}\mathbb{P}\{|\tilde{a}_{k,\ell}\tilde{a}_{j,\ell}-\mathbb{E}(\tilde{a}_{k,\ell}\tilde{a}_{j,\ell})|>u\}
\lesssim \exp\bigg\{-\frac{Cu}{l_n^{2}\log^{2}(l_n)}\bigg\} +\exp\bigg\{-\frac{Cb^{1/6}u^{1/6}}{l_n^{1/2}\log^{1/3}(l_n)}\bigg\} \,, \\
&~~~~~\max_{\ell\in[L]}\max_{j\in[2Kr]} \mathbb{P}\{|\tilde{a}_{j,\ell}|^3-\mathbb{E}(|\tilde{a}_{j,\ell}|^3)|>u\}
\lesssim \exp\bigg\{-\frac{Cu^{2/3}}{l_n^{2}\log^{2}(l_n)}\bigg\} +\exp\bigg\{-\frac{Cb^{1/6}u^{1/9}}{l_n^{1/2}\log^{1/3}(l_n)}\bigg\}   \,,
\end{align*}
for any $u>0$. Analogously,
if $\log(Kr)\ll \min(L^{2/5}, b^{3/20}L^{1/10}l_n^{-3/20})$ and $h=2l_n+C''\log(Kr)$ for some sufficiently large constant $C''>0$,
using the same arguments in the proof of Lemma 3 in  the supplementary material of \cite{CCW_2021_supp}, 
we have
\begin{align*}
&~~~~~~~~~~~~~~~~~\frac{1}{l_n^{3}\log^3(l_n)}M_{\tilde{\ba}}\bigg\{\frac{C'L^{1/6}}{l_n(\log l_n)\log^{2/3}(2Kr)}\bigg\}
\lesssim \frac{l_n(\log l_n)\log^{7/6}(Kr)}{L^{1/6}}  \,, \\
&~~~~~~~~~~\frac{L^{2/3}}{l_n(\log l_n)\log^{2/3}(Kr)} \max_{\ell\in[L]} \mathbb{E}\bigg\{\max_{j\in[2Kr]}|\mathbb{E}(\tilde{a}_{j,\ell}\,|\,\mathcal{F}_{-\ell})|\bigg\}
\lesssim \frac{l_n(\log l_n)\log^{7/6}(Kr)}{L^{1/6}}\,, \\
&\frac{L^{1/3}}{l_n^2\log^2(l_n)\log^{1/3}(2Kr)}\max_{\ell\in[L]}\mathbb{E}\bigg[\max_{k,j\in[2Kr]}|\mathbb{E}\{\tilde{a}_{k,\ell}\tilde{a}_{j,\ell}-\mathbb{E}(\tilde{a}_{k,\ell}\tilde{a}_{j,\ell})\,|\,\mathcal{F}_{-\ell}\}|\bigg]
\lesssim \frac{l_n(\log l_n)\log^{7/6}(Kr)}{L^{1/6}}  \,, \\
&~~~~~~~~~~~~~\frac{1}{l_n^{3}\log^3(l_n)}\max_{\ell\in[L]}\mathbb{E}\bigg[\max_{j\in[2Kr]}|\mathbb{E}\{|\tilde{a}_{j,\ell}|^3-\mathbb{E}(|\tilde{a}_{j,\ell}|^3)\,|\,\mathcal{F}_{-\ell}\}|\bigg]
\lesssim \frac{l_n(\log l_n)\log^{7/6}(Kr)}{L^{1/6}}  \,.
\end{align*}
Together with \eqref{eq:rho1}, this completes the proof of Lemma \ref{pn:1}. $\hfill\Box$

%

\subsection{Proof of Lemma \ref{pn:2}} \label{sec:pflem.pn2}
The proof of Lemma \ref{pn:2}  follows the idea for   the proof of Lemma 4 in the supplementary material of \cite{CCW_2021_supp}.  
Denote by $\tilde{\alpha}(k)$ the $\alpha$-mixing coefficients of the sequence $\{\ba_t\}_{t=1}^{\tilde{n}}$. Due to $\ba_t=\bH\bc_t$, where $\bc_t=(\bc_{1,t}^{\T},\ldots,\bc_{r,t}^{\T})^{\T}$ with $\bc_{\ell,t}$ defined as \eqref{eq:cjt}, Condition \ref{as:betamixing} implies $\tilde{\alpha}(k)\leq C\exp(-C|k-2l_n|_+)$ for any integer $k\geq 1$.
Different from the setting considered in \cite{CCW_2021_supp}, 
the $\alpha$-mixing coefficients of $\{\ba_t\}_{t=1}^{\tilde{n}}$ vary with the sample size $n$ and thus we
need to refine some upper bounds used there.
In the sequel, we only specify the difference between our proof and that for Lemma 4 in  the supplementary material of \cite{CCW_2021_supp}.  
For some $D_n>0$, consider the event $\mathcal{E}=\{|\bdelta_n|_\infty\leq {D}_{n}\}$ with
\begin{align*}
\bdelta_n=\frac{1}{\sqrt{\tilde{n}}}\sum_{\ell=1}^{L+1}\sum_{t\in\mathcal{J}_\ell}{\ba}_t+\bigg(\frac{1}{\sqrt{\tilde{n}}}-\frac{1}{\sqrt{Lb}}\bigg)\sum_{\ell=1}^L\sum_{t\in\mathcal{I}_\ell}{\ba}_t \,.
\end{align*}
If  $l_n^2 (hb^{-1}+bn^{-1})=o(1)$ with $l_n\geq 2$, parallel to Equation (S.12) in the supplementary material of \cite{CCW_2021_supp},  
by Lemma \ref{pn:1}, we have
\begin{align}
\varrho_n^{(2)}
\lesssim  \frac{l_n(\log l_n)\log^{7/6}(Kr)}{L^{1/6}}+{D}_{n}\log^{1/2}(Kr)+\mathbb{P}(\mathcal{E}^{\rm c}) \label{eq:bb2}   
\end{align}
provided that $\log(Kr)\ll \min(b^{3/20}L^{1/10}l_n^{-3/20}, L^{2/5})$ and $h=2l_n+C'\log(Kr)$ for some sufficiently large constant $C'>0$.
Similar to Lemma \ref{la:xtildetail}, if $l_n\geq 2$, we also have
\begin{align*}
& \max_{j\in[2Kr]} \mathbb{P}\bigg(\bigg|\sum_{\ell=1}^L\sum_{t\in\mathcal{I}_\ell}a_{j,t}\bigg| >\frac{D_n\sqrt{\tilde{n}}b}{2h}\bigg)
\lesssim  \exp\bigg\{-\frac{CD_n^2b^2}{l_n^2\log^2(l_n)h^2}\bigg\} + \exp\bigg\{-\frac{C{n}^{1/6}b^{1/3}D_n^{1/3}}{l_n^{1/2}\log^{1/3}(l_n)h^{1/3}}\bigg\} \,, \\
&~~~ \max_{j\in[2Kr]}\mathbb{P}\bigg(\bigg|\sum_{\ell=1}^{L+1}\sum_{t\in\mathcal{J}_\ell}a_{j,t}\bigg| >\frac{D_n\sqrt{\tilde{n}}}{2}\bigg)
\lesssim \exp\bigg\{-\frac{CD_n^2b}{l_n^2\log^2(l_n)h}\bigg\} +\exp\bigg\{-\frac{Cn^{1/6}D_n^{1/3}}{l_n^{1/2}\log^{1/3}(l_n)}\bigg\}  \,.
\end{align*}
Same as that in the proof of  Lemma 4 in the supplementary material of \cite{CCW_2021_supp},  
due to $h=o(b)$, it holds that
\begin{align*}
\mathbb{P}(\mathcal{E}^{\rm c})
\leq&~ \sum_{j=1}^{2Kr}\mathbb{P}\bigg(\bigg|\sum_{\ell=1}^{L+1}\sum_{t\in\mathcal{J}_\ell}a_{j,t}\bigg|>\frac{D_n\sqrt{\tilde{n}}}{2}\bigg) +\sum_{j=1}^{2Kr}\mathbb{P}\bigg(\bigg|\sum_{\ell=1}^L\sum_{t\in\mathcal{I}_\ell}a_{j,t}\bigg|>\frac{D_n\sqrt{\tilde{n}}b}{2h}\bigg) \\
\lesssim&~  Kr\exp\bigg\{-\frac{CD_n^2b}{l_n^2\log^2(l_n)h}\bigg\} +Kr\exp\bigg\{-\frac{Cn^{1/6}D_n^{1/3}}{l_n^{1/2}\log^{1/3}(l_n)}\bigg\}  \,.
\end{align*}
With selecting $D_n=C''L^{-1/6}l_n(\log l_n)\log^{2/3}(Kr)$ for some sufficiently large constant $C''>0$,
if $\min(n^{1/2}, \, nl_n^{-2})\gg b\gg
\max\{hl_n^2, \, h^{3/4}n^{1/4}\log^{-1/4}(Kr)\}$ and $h=2l_n+C'\log(Kr)$ for some sufficiently large constant $C'>0$, we complete the proof of Lemma \ref{pn:2} by \eqref{eq:bb2}
provided that $\log(Kr)\ll\min(b^{3/20}L^{1/10}l_n^{-3/20}, L^{2/5})$. $\hfill\Box$

\subsection{Proof of Lemma \ref{lem.Xi}}\label{sec:pfLem11}
Define $\bXi^* = \sum_{q=-\tilde n+1}^{\tilde n-1} \mathcal K(q/b_n) \bPi(q) $,
where  $\bPi(q)=\{\Pi_q(\ell_1,\ell_2)\}_{r(2l_n+1)\times r(2l_n+1)}$ with $\bPi(q)=\tilde{n}^{-1} \sum_{t=q+1}^{\tilde n}\bE(\bc_t\bc_{t-q}^{\T})$ if $q\geq0$ and $\bPi(q)=\tilde{n}^{-1} \sum_{t=-q+1}^{\tilde n}\bE(\bc_{t+q}\bc_t^{\T})$ if $q<0$.
By the triangle inequality,
\begin{align*}
|\widehat{\bXi}-\bXi|_\infty \leq
|\widehat{\bXi}-\bXi^*|_\infty + |\bXi^*-\bXi|_\infty \,.
\end{align*}
Recall $\bXi=\var(\tilde{n}^{-1/2}\sum_{t=1}^{\tilde n}\bc_t)$, where $\{\bc_t\}_{t=1}^{\tilde n}$ is an $\alpha$-mixing sequence with $\alpha$-mixing coefficients satisfying $\alpha_{\bc}(k)\lesssim \exp(-C|k-2l_n|_+)$ for any integer $k\geq 1$.
Recall  $\bc_t=\{c_{1,t},\ldots,c_{r(2l_n+1),t}\}^\T=(\bc_{1,t}^{\T},\ldots,\bc_{r,t}^{\T})^{\T}$, where ${\bc}_{j,t}=(2\pi)^{-1}\{\mathring{x}_{\chi_1(j),t}\mathring{x}_{\chi_2(j),t+l_n}-\gamma_{\bchi(j)}(-l_n), \\ \ldots,
\mathring{x}_{\chi_1(j),t+2l_n}\mathring{x}_{\chi_2(j),t+l_n}-\gamma_{\bchi(j)}(l_n)\}^{\T}$ with $\mathring{\bx}_t=(\mathring{x}_{1,t}, \ldots,\mathring{x}_{p,t})^{\T}=\bx_t-\bmu$.
By \eqref{eq:t1}, we know
$\max_{t\in[\tilde n]}\max_{j\in[r(2l_n+1)]}\bP(|c_{j,t}|>u)\lesssim \exp(-Cu)$ for any $u>0$.
Using the same arguments
for deriving  the convergence rate of $|\bSigma_{n,K}^*-\bSigma_{n,K}|_\infty$ in the proof of Proposition 2  in \cite{CJS_2021_supp}, 
we know $|\bXi^*-\bXi|_\infty  \lesssim n^{-\rho}l_n^2$. As we will show in Section \ref{subsec:diff.xi.xistar}, if  $\log r \ll \min\{n^{(2\vartheta+2\rho-2\rho\vartheta-2)/(7\vartheta-4)}, n^{1/5}l_n^{-1/5}\}$,
\begin{align}\label{eq:diff.xi.xistar}
|\widehat{\bXi}-\bXi^*|_\infty	
= &~   O_{\rm p}\bigg\{\frac{\log^{2} (r)} {n^{(\vartheta+2\rho-3\rho\vartheta-1)/(2\vartheta-1)}}\bigg\}
+O_{\rm p}\bigg\{\frac{\log^{4}(r)} {n^{(2\vartheta+3\rho-4\rho\vartheta-2)/(2\vartheta-1)}}\bigg\}  \notag\\
& + O_{\rm p}\bigg\{\frac{\log^{1/2}(r)}{n^{(1-2\rho)/2}}\bigg\} +O_{\rm p}\bigg(\frac{l_n\log r}{n^{1-\rho}}\bigg)
\,,
\end{align}
which implies
\begin{align*}
|\widehat{\bXi}-\bXi|_\infty =&~ O_{\rm p}\bigg\{\frac{\log^{2} (r)} {n^{(\vartheta+2\rho-3\rho\vartheta-1)/(2\vartheta-1)}}\bigg\}
+O_{\rm p}\bigg\{\frac{\log^{4}(r)} {n^{(2\vartheta+3\rho-4\rho\vartheta-2)/(2\vartheta-1)}}\bigg\} \\
&   + O_{\rm p}\bigg\{\frac{\log^{1/2}(r)}{n^{(1-2\rho)/2}}\bigg\}
+O_{\rm p}\bigg(\frac{l_n\log r}{n^{1-\rho}}\bigg)
+ O\bigg(\frac{l_n^2}{n^{\rho}}\bigg)
\,.
\end{align*}
We have completed the proof of Lemma \ref{lem.Xi}.
$\hfill\Box$

\subsubsection{Proof of \eqref{eq:diff.xi.xistar}}	\label{subsec:diff.xi.xistar}
To specify the convergence rate of $|\widehat{\bXi}-\bXi^*|_\infty$, we need the following lemmas whose proofs are given in Sections \ref{subsec:pf:gammatail2}--\ref{subsec:pftildegamma}. 	
\begin{lemma}\label{la:gammatail2}
Under Conditions {\rm \ref{as:moment}} and {\rm \ref{as:betamixing}}, for any $u>0$, it holds that
\begin{align*}
	&\max_{0\leq q\leq M}\max_{ i,j\in[r(2l_n+1)]}\mathbb{P}\bigg[\bigg|\sum_{t=q+1}^{\tilde n}
	\{c_{i,t}c_{j,t-q}-\bE(c_{i,t}c_{j,t-q}) \}\bigg|
	>u\bigg] \\
	&~~~~~~~~~~~~~~~~~~ \lesssim
	\exp\bigg(-\frac{Cu^2}{nM}\bigg)
	+\exp\bigg(-\frac{Cu^{1/4}}{M^{1/4}}\bigg) \,,
\end{align*}
where $2l_n\leq M=o(n)\rightarrow\infty$.
\end{lemma}

\begin{lemma}\label{la:gammatail}
Under Conditions {\rm \ref{as:moment}} and {\rm \ref{as:betamixing}}, for any $1\leq s_1<s_2\leq n-2l_n$, it holds that
\begin{align*}
	\max_{-l_n\leq k\leq l_n}  \max_{i,j\in[p]} \bP\bigg[\bigg|\sum_{t=l_n+s_1}^{l_n+s_2} \{\mathring{x}_{i,t+k}\mathring{x}_{j,t}-\gamma_{i,j}(k)\}\bigg|>u\bigg]  \lesssim
	\exp\bigg(-\frac{Cu^2}{sl_n}\bigg)+\exp\bigg(-\frac{Cu^{1/3}}{l_n^{1/3}}\bigg)
\end{align*}
for any $u>0$, where $s=s_2-s_1$.
\end{lemma}

\begin{lemma} \label{la:tildegammatail}
Under Conditions {\rm\ref{as:moment}} and {\rm\ref{as:betamixing}}, it holds that
\begin{align*}
	&	\max_{(i,j)\in\cI}\max_{-l_n\leq k\leq l_n} \bP\big\{|\hat{\gamma}_{i,j}(k)-\gamma_{i,j}(k)| >u\big\}  \\
	&~~~~~~~~~~~	\lesssim
	\exp(-Cnl_n^{-1}u^2)
	+\exp(-Cn^{1/3}l_n^{-1/3}u^{1/3})
	+ \exp(-Cnu)+ \exp(-Cn^{1/3}u^{1/6})
\end{align*}
for any $u\gg l_nn^{-1}$.
\end{lemma}

Recall $\hat{\bc}_t=\{\hat{c}_{1,t},\ldots,\hat{c}_{r(2l_n+1),t}\}^{\T}=(\hat{\bc}_{1,t}^{\T},\ldots,\hat{\bc}_{r,t}^{\T})^{\T}$, where
$$\hat{\bc}_{\ell,t}=(2\pi)^{-1}\big\{\hat{\mathring{x}}_{\chi_1(\ell),t}\hat{\mathring{x}}_{\chi_2(\ell),t+l_n} -{\hat\gamma}_{\bchi(\ell)}(-l_n),\ldots, \hat{\mathring{x}}_{\chi_1(\ell),t+2l_n}\hat{\mathring{x}}_{\chi_2(\ell),t+l_n} -{\hat\gamma}_{\bchi(\ell)}(l_n)\big\}^\T$$
with $\hat{\mathring{\bx}}_t=(\hat{\mathring{x}}_{1,t},\ldots,\hat{\mathring{x}}_{p,t})^{\T}= \bx_t-\bar\bx$.
Then for any $l\in [r(2l_n+1)]$, there exists a unique triple $(i,j,k)$ such that $\hat{c}_{l,t}=(2\pi)^{-1}\{(x_{i,t+l_n+k}-\bar{x}_i)(x_{j,t+l_n}-\bar{x}_j)-\hat{\gamma}_{i,j}(k)\}$.
Without loss of generality, we  assume $\bmu=\bzero$. Then $\mathring{\bx}_t=\bx_t$, $\gamma_{i,j}(k)=\bE(x_{i,t+k}x_{j,t})$ and
$$\bc_{\ell,t}=(2\pi)^{-1}\big\{{x}_{\chi_1(\ell),t}{x}_{\chi_2(\ell),t+l_n}-\gamma_{\bchi(\ell)}(-l_n), \ldots,{x}_{\chi_1(\ell),t+2l_n}{x}_{\chi_2(\ell),t+l_n}-\gamma_{\bchi(\ell)}(l_n)\big\}^{\T}\,.$$
Define $\tilde{\gamma}_{i,j}(k)=(2\pi)^{-1}\{\hat{\gamma}_{i,j}(k)-\gamma_{i,j}(k)\}$, and write
\begin{align*}
\tilde{\bgamma}:=&~\{\tilde{\gamma}_1,\ldots,\tilde{\gamma}_{r(2l_n+1)}\}^{\T} = \big\{\tilde{\gamma}_{\bchi(1)}(-l_n),\ldots,\tilde{\gamma}_{\bchi(1)}(l_n),\ldots, \tilde{\gamma}_{\bchi(r)}(-l_n),
\ldots,\tilde{\gamma}_{\bchi(r)}(l_n)\big\}^{\T} \,, \\
{\br}_t^{(1)}:=&~\{r_{1,t}^{(1)},\ldots,r_{r(2l_n+1),t}^{(1)}\}^{\T}\\
=&~ (2\pi)^{-1}\big\{x_{\chi_1(1),t}\bar{x}_{\chi_2(1)},\ldots,x_{\chi_1(1),t+2l_n}\bar{x}_{\chi_2(1)},
\ldots,x_{\chi_1(r),t}\bar{x}_{\chi_2(r)},\ldots,
x_{\chi_1(r),t+2l_n}\bar{x}_{\chi_2(r)}\big\}^{\T} \,, \\
{\br}_t^{(2)}:=&~\{r_{1,t}^{(2)},\ldots,r_{r(2l_n+1),t}^{(2)}\}^{\T}
= (2\pi)^{-1}\big\{\bar{x}_{\chi_1(1)}{x}_{\chi_2(1),t+l_n},
\ldots,\bar{x}_{\chi_1(r)}{x}_{\chi_2(r),t+l_n}\big\}^{\T} \otimes{\bf 1}_{2l_n+1} \,,\\
{\br}^{(3)}:=&~\{r_{1}^{(3)},\ldots,r_{r(2l_n+1)}^{(3)}\}^{\T}
=(2\pi)^{-1}\big\{\bar{x}_{\chi_1(1)}\bar{x}_{\chi_2(1)},\ldots,\bar{x}_{\chi_1(r)}\bar{x}_{\chi_2(r)}\big\}^{\T}  \otimes{\bf 1}_{2l_n+1} \,,
\end{align*}
where ${\bf 1}_{2l_n+1}$ is a $(2l_n+1)$-dimensional vector whose elements are all equal to 1.
Thus,
\begin{align*}
\hat{\bc}_t=\bc_t-{\br}_t^{(1)}-{\br}_t^{(2)}+{\br}^{(3)}-\tilde\bgamma\,.
\end{align*}
Then for any $\ell_1,\ell_2\in[r(2l_n+1)]$ and $t_1,t_2\in[\tilde{n}]$,
\begin{align*}
\hat{c}_{\ell_1,t_1}\hat{c}_{\ell_2,t_2}
=&~ c_{\ell_1,t_1}c_{\ell_2,t_2}-c_{\ell_1,t_1}r_{\ell_2,t_2}^{(1)}-c_{\ell_1,t_1}r_{\ell _2,t_2}^{(2)} +c_{\ell_1,t_1}r_{\ell_2}^{(3)} -c_{\ell_1,t_1}\tilde{\gamma}_{\ell_2}  \\
& -r_{\ell_1,t_1}^{(1)}c_{\ell_2,t_2}+r_{\ell_1,t_1}^{(1)}r_{\ell_2,t_2}^{(1)}+r_{\ell_1,t_1}^{(1)}r_{\ell _2,t_2}^{(2)} -r_{\ell_1,t_1}^{(1)}r_{\ell_2}^{(3)} +r_{\ell_1,t_1}^{(1)}\tilde{\gamma}_{\ell_2}  \\
& -r_{\ell_1,t_1}^{(2)}c_{\ell_2,t_2}+r_{\ell_1,t_1}^{(2)}r_{\ell_2,t_2}^{(1)}+r_{\ell_1,t_1}^{(2)}r_{\ell _2,t_2}^{(2)} -r_{\ell_1,t_1}^{(2)}r_{\ell_2}^{(3)} +r_{\ell_1,t_1}^{(2)}\tilde{\gamma}_{\ell_2}  \\
& +r_{\ell_1}^{(3)}c_{\ell_2,t_2}-r_{\ell_1}^{(3)}r_{\ell_2,t_2}^{(1)}-r_{\ell_1}^{(3)}r_{\ell _2,t_2}^{(2)} +r_{\ell_1}^{(3)}r_{\ell_2}^{(3)} -r_{\ell_1}^{(3)}\tilde{\gamma}_{\ell_2}  \\
& -\tilde{\gamma}_{\ell_1}c_{\ell_2,t_2}+\tilde{\gamma}_{\ell_1}r_{\ell_2,t_2}^{(1)}+\tilde{\gamma}_{\ell_1}r_{\ell _2,t_2}^{(2)} -\tilde{\gamma}_{\ell_1}r_{\ell_2}^{(3)} +\tilde{\gamma}_{\ell_1}\tilde{\gamma}_{\ell_2} \\
=&: \sum_{j=1}^{25}C_j(\ell_1,\ell_2,t_1,t_2)\,.
\end{align*}
Recall
$
\widehat\bXi = \sum_{q=-\tilde n+1}^{\tilde n-1} \mathcal K(q/b_n) \widehat\bPi(q)$,
where $\widehat\bPi(q)=\{\widehat{\Pi}_q(\ell_1,\ell_2)\}_{r(2l_1+1)\times r(2l_n+1)}$ with $\widehat{\bPi}(q)=\tilde{n}^{-1} \sum_{t=q+1}^{\tilde n}\hat{\bc}_t\hat{\bc}_{t-q}^{\T}$ if $q\geq0$ and $\widehat\bPi(q)=\tilde{n}^{-1} \sum_{t=-q+1}^{\tilde n}\hat{\bc}_{t+q}\hat{\bc}_t^{\T}$ if $q<0$.
For any $-\tilde{n}+1\leq q\leq \tilde{n}-1$, it holds that
\begin{align*}
\widehat{\bPi}(q)-\bPi(q) = \frac{1}{\tilde{n}}\sum_{t=|q|+1}^{\tilde{n}}\big\{\hat{\bc}_{t-|-q|_+}\hat{\bc}_{t-|q|_+}^{\T} - \bE(\bc_{t-|-q|_+}\bc_{t-|q|_+}^{\T})\big\} \,.
\end{align*}
Write
$$A_1=\max_{\ell_1,\ell_2}\bigg|\sum_{q=-\tilde{n}+1}^{\tilde n-1}\cK\bigg(\frac{q}{b_n}\bigg) \bigg[\frac{1}{\tilde{n}}\sum_{t=|q|+1}^{\tilde{n}} \big\{C_1(\ell_1,\ell_2,t-|-q|_+,t-|q|_+)-\bE(c_{\ell_1,t-|-q|_+}c_{\ell_2,t-|q|_+})\big\}\bigg]\bigg| \,,$$
and
$$A_j=\max_{\ell_1,\ell_2} \bigg|\sum_{q=-\tilde{n}+1}^{\tilde n-1}\cK\bigg(\frac{q}{b_n}\bigg)\bigg\{\frac{1}{\tilde{n}}\sum_{t=|q|+1}^{\tilde{n}} C_j(\ell_1,\ell_2,t-|-q|_+,t-|q|_+)\bigg\}\bigg|$$
for $j\in\{2,\ldots,25\}$.
Recall $\bXi^* = \sum_{q=-\tilde n+1}^{\tilde n-1} \mathcal K(q/b_n) \bPi(q) $.
By the triangle inequality,
\begin{align*}
|\widehat{\bXi}-\bXi^*|_\infty=\Bigg|\sum_{q=-\tilde{n}+1}^{\tilde{n}-1}\mathcal{K}\bigg(\frac{q}{b_n}\bigg) \{\widehat\bPi(q)-\bPi(q)\}\Bigg|_\infty \leq \sum_{j=1}^{25} A_j\,.
\end{align*}
By the symmetrization, we know
$A_2=A_6$, $A_3=A_{11}$, $A_4=A_{16}$, $A_5=A_{21}$, $A_8=A_{12}$, $A_9=A_{17}$, $A_{10}=A_{22}$, $A_{14}=A_{18}$, $A_{15}=A_{23}$ and  $A_{20}=A_{24}$. To specify the convergence rate of $|\widehat{\bXi}-\bXi^*|_\infty$, it suffices to derive the convergence rates of $A_1$, $A_2$, $A_3$, $A_4$, $A_5$, $A_7$, $A_8$, $A_9$, $A_{10}$, $A_{13}$, $A_{14}$, $A_{15}$, $A_{19}$, $A_{20}$ and $A_{25}$, respectively. As we will show in Sections \ref{subsec:A1}--\ref{subsec:A10},
\begin{align*}
&~~~~~~~~~~~~~~~ ~ A_1= O_{\rm p}\bigg\{\frac{\log^{2}(r)}{n^{(\vartheta+2\rho-3\rho\vartheta-1)/(2\vartheta-1)}}\bigg\}  +O_{\rm p}\bigg\{\frac{\log^4(r)}{n^{(2\vartheta+3\rho-4\rho\vartheta-2)/(2\vartheta-1)}}\bigg\}  \,, \\
&~~~~~~~~~~ A_2+A_3=  O_{\rm p}\bigg\{\frac{\log^{1/2}(r)}{n^{(1-2\rho)/2}}\bigg\}   \,, ~~~~
A_4+A_{10}+A_{15}+A_{20} = O_{\rm p}\bigg\{\frac{l_n^{3/2}\log^{3/2}(r)}{n^{(3-2\rho)/2}}\bigg\}  \,,\\
&A_5+A_{25} = O_{\rm p}\bigg(\frac{l_n\log r}{n^{1-\rho}}\bigg) \,,  ~~
A_7+A_8+A_{13}=  O_{\rm p}\bigg(\frac{\log r}{n^{1-\rho}}\bigg)    \,, ~~
A_9+A_{14}+A_{19}=  O_{\rm p}\bigg\{\frac{\log^2(r)}{n^{2-\rho}}\bigg\}
\end{align*}
provided that $\log r \ll \min\{n^{(2\vartheta+2\rho-2\rho\vartheta-2)/(7\vartheta-4)}, n^{1/5}l_n^{-1/5}\}$. Therefore,
\begin{align*}
|\widehat{\bXi}-\bXi^*|_\infty	=
&~ O_{\rm p}\bigg\{\frac{\log^{2}(r)}{n^{(\vartheta+2\rho-3\rho\vartheta-1)/(2\vartheta-1)}}\bigg\}
+O_{\rm p}\bigg\{\frac{\log^4(r)}{n^{(2\vartheta+3\rho-4\rho\vartheta-2)/(2\vartheta-1)}}\bigg\}  \\
&
+ O_{\rm p}\bigg\{\frac{\log^{1/2}(r)}{n^{(1-2\rho)/2}}\bigg\} +O_{\rm p}\bigg(\frac{l_n\log r}{n^{1-\rho}}\bigg)  \,.
\end{align*}
We complete the proof of \eqref{eq:diff.xi.xistar}.
$\hfill\Box$

\subsubsection{Convergence rate of $A_1$}\label{subsec:A1}
Given $\ell_1$ and $\ell_2$,  write $\psi_{t,q}=c_{\ell_1,t+q}c_{\ell_2,t}-\bE(c_{\ell_1,t+q}c_{\ell_2,t})$.  By the triangle inequality and the Bonferrroni inequality,
\begin{align}\label{eq:sigma.trun}
& \bP\bigg(\bigg|\sum_{q=0}^{\tilde{n}-1}\mathcal{K}\bigg(\frac{q}{b_n}\bigg) \bigg[\frac{1}{\tilde{n}}\sum_{t=q+1}^{\tilde{n}} \{c_{\ell_1,t}c_{\ell_2,t-q}-\bE(c_{\ell_1,t}c_{\ell_2,t-q})\}\bigg]\bigg|>u\bigg) \\
&~~~~~~~~~~ \leq \bP\bigg\{ \sum_{q=0}^{M_{1n}} \bigg| \mathcal{K}\bigg(\frac{q}{b_n}\bigg)\bigg| \bigg|\frac{1}{\tilde{n}}\sum_{t=1}^{\tilde{n}-q} \psi_{t,q} \bigg| >\frac{u}{2}\bigg\}
+ \bP\bigg\{ \sum_{q=M_{1n}+1}^{\tilde{n}-1} \bigg| \mathcal{K}\bigg(\frac{q}{b_n}\bigg)\bigg|\bigg|\frac{1}{\tilde{n}}\sum_{t=1}^{\tilde{n}-q} \psi_{t,q} \bigg| >\frac{u}{2}\bigg\}  \notag
\end{align}
for any $u>0$.
Note that $\max_{t\in[\tilde{n}]}\max_{j\in[r(2l_n+1)]}\bP(|c_{j,t}|>u)\leq C\exp(-Cu)$ for any $u>0$.
Lemma 2 in the supplementary material of \cite{CTW_2013_supp} 
yields $\max_{q,t}\bP(|\psi_{t,q}|>u)\leq C\exp(-Cu^{1/2})$ for any $u>0$.
Given $M_{1n}=o(n)\rightarrow\infty$ satisfying $M_{1n}\geq 2l_n$,
by Condition \ref{as:kernel}(ii) and $b_n\asymp n^{\rho}$ for some $\rho\in(0,1)$, we have
$\sum_{q=M_{1n}+1}^{\tilde{n}-1}|\mathcal{K}(q/b_n)| \lesssim\sum_{q=M_{1n}+1}^{\tilde{n}-1}(q/b_n)^{-\vartheta}\lesssim n^{\rho\vartheta}M_{1n}^{1-\vartheta}$. By the triangle inequality and the Bonferroni inequality,
\begin{align}\label{eq:sigma.tail}
&\bP\bigg\{\sum_{q=M_{1n}+1}^{\tilde{n}-1}
\bigg|\mathcal{K}\bigg(\frac{q}{b_n}\bigg)\bigg| \bigg|\frac{1}{\tilde{n}}\sum_{t=1}^{\tilde{n}-q}\psi_{t,q}\bigg|>\frac{u}{2} \bigg\}
\leq  \sum_{q=M_{1n}+1}^{\tilde{n}-1}
\bP\bigg(\bigg|\frac{1}{\tilde{n}}\sum_{t=1}^{\tilde{n}-q}\psi_{t,q}\bigg| >\frac{CM_{1n}^{\vartheta-1}u}{n^{\rho\vartheta}}\bigg) \notag\\
&~~~~~~~~~~~\leq  \sum_{q=M_{1n}+1}^{\tilde{n}-1}\sum_{t=1}^{\tilde{n}-q}
\bP\bigg(|\psi_{t,q}|
>\frac{CM_{1n}^{\vartheta-1}u}{n^{\rho\vartheta}}\bigg)
\lesssim  n^2\exp\bigg\{ -\frac{CM_{1n}^{(\vartheta-1)/2}u^{1/2}}{n^{\rho\vartheta/2}}\bigg\}
\end{align}
for any $u>0$. Write $D_n=\sum_{q=0}^{M_{1n}}|\mathcal{K}(q/b_n)|$.  By Condition \ref{as:kernel}(ii),  $D_n\lesssim b_n\asymp n^\rho$.
By Bonferroni inequality and Lemma \ref{la:gammatail2},
\begin{align*}
& \bP\bigg\{ \sum_{q=0}^{M_{1n}} \bigg| \mathcal{K}\bigg(\frac{q}{b_n}\bigg)\bigg| \bigg|\frac{1}{\tilde{n}}\sum_{t=1}^{\tilde{n}-q} \psi_{t,q} \bigg| >\frac{u}{2}\bigg\}
\leq \sum_{q=0}^{M_{1n}} \bP\bigg(\bigg|\sum_{t=1}^{\tilde{n}-q}\psi_{t,q}\bigg|>\frac{C\tilde{n}u}{D_n}\bigg) \\
&~~~~~~~~~~~~~~~ \lesssim M_{1n}\exp\bigg(-\frac{Cn^{1-2\rho}u^2}{M_{1n}}\bigg)
+ M_{1n} \exp\bigg\{-\frac{Cn^{(1-\rho)/4}u^{1/4}}{M_{1n}^{1/4}}\bigg\}
\end{align*}
for any $u>0$.
Together with \eqref{eq:sigma.trun} and \eqref{eq:sigma.tail}, we have
\begin{align*}
& \max_{\ell_1,\ell_2}\bP\bigg(\bigg|\sum_{q=0}^{\tilde{n}-1}\mathcal{K}\bigg(\frac{q}{b_n}\bigg) \bigg[\frac{1}{\tilde{n}}\sum_{t=q+1}^{\tilde{n}} \{c_{\ell_1,t}c_{\ell_2,t-q}-\bE(c_{\ell_1,t}c_{\ell_2,t-q})\}\bigg]\bigg|>u\bigg) \\
&~~~~~~~~~~~	\lesssim n^2\exp\bigg\{ -\frac{CM_{1n}^{(\vartheta-1)/2}u^{1/2}}{n^{\rho\vartheta/2}}\bigg\}+  M_{1n}\exp\bigg(-\frac{Cn^{1-2\rho}u^2}{M_{1n}}\bigg) \\
&~~~~~~~~~~~~~~~ +M_{1n} \exp\bigg\{-\frac{Cn^{(1-\rho)/4}u^{1/4}}{M_{1n}^{1/4}}\bigg\}
\end{align*}
for any $u>0$, which implies
\begin{align*}
\bP(A_1>u)\lesssim&~ r^2l_n^2\max_{\ell_1,\ell_2} \bP\bigg(\bigg|\sum_{q=0}^{\tilde{n}-1}\mathcal{K}\bigg(\frac{q}{b_n}\bigg) \bigg[{\frac{1}{\tilde{n}}\sum_{t=q+1}^{\tilde{n}} \{c_{\ell_1,t}c_{\ell_2,t-q}}-\bE(c_{\ell_1,t}c_{\ell_2,t-q})\}\bigg]\bigg|>\frac{u}{2}\bigg) \\
\lesssim&~ M_{1n} r^2l_n^2\exp\bigg(-\frac{Cn^{1-2\rho}u^2}{M_{1n}}\bigg) +M_{1n} r^2l_n^2\exp\bigg\{-\frac{Cn^{(1-\rho)/4}u^{1/4}}{M_{1n}^{1/4}}\bigg\}  \\
& + n^2r^2l_n^2\exp\bigg\{ -\frac{CM_{1n}^{(\vartheta-1)/2}u^{1/2}}{n^{\rho\vartheta/2}}\bigg\}
\end{align*}
for any $u>0$. Note that $r\geq n^{\kappa}$ for some sufficiently small constant $\kappa>0$. Hence,
\begin{align*}
A_1=  O_{\rm p}\bigg\{\frac{M_{1n}^{1/2}\log^{1/2}(r)}{n^{(1-2\rho)/2}}\bigg\}
+ O_{\rm p}\bigg\{\frac{M_{1n}\log^4(r)}{n^{1-\rho}}\bigg\}
+ O_{\rm p}\bigg\{\frac{n^{\rho\vartheta}\log^2(r)}{M_{1n}^{\vartheta-1}}\bigg\} \,.
\end{align*}
Selecting $M_{1n}=n^{(1-2\rho+2\rho\vartheta)/(2\vartheta-1)}$, then $n^\rho\ll M_{1n}\ll n$ and
\begin{align*}
A_1 = O_{\rm p}\bigg\{\frac{\log^{2}(r)}{n^{(\vartheta+2\rho-3\rho\vartheta-1)/(2\vartheta-1)}}\bigg\}  +O_{\rm p}\bigg\{\frac{\log^4(r)}{n^{(2\vartheta+3\rho-4\rho\vartheta-2)/(2\vartheta-1)}}\bigg\}  \,.
\end{align*}
Furthermore,
\begin{align}\label{eq:A1tail}
\bP(A_1>u)
\lesssim &~ nr^2l_n^2 \exp\{-Cn^{(2\vartheta+4\rho-6\rho\vartheta-2)/(2\vartheta-1)}u^2\}\notag\\
&+nr^2l_n^2 \exp\{-Cn^{(2\vartheta+3\rho-4\rho\vartheta-2)/(8\vartheta-4)}u^{1/4}\} \notag\\
& + n^2r^2l_n^2 \exp\{-Cn^{(\vartheta+2\rho-3\rho\vartheta-1)/(4\vartheta-2)}u^{1/2}\}
\end{align}
for any $u>0$.
$\hfill\Box$

\subsubsection{Convergence rates of $A_2$ and $A_3$}	
Notice that $A_2\leq \max_{\ell_1,\ell_2}A_2(\ell_1,\ell_2)+\max_{\ell_1,\ell_2}\tilde{A}_2(\ell_1,\ell_2)$, where
\begin{align*} 
A_2(\ell_1,\ell_2)
:=&~ \bigg|\sum_{q=0}^{\tilde{n}-1} \mathcal{K}\bigg(\frac{q}{b_n}\bigg) \bigg\{\frac{1}{\tilde{n}}\sum_{t=q+1}^{\tilde{n}}c_{\ell_1,t}r_{\ell_2,t-q}^{(1)}\bigg\}\bigg|  \,, \\
\tilde{A}_2(\ell_1,\ell_2):=&~ \bigg|\sum_{q=-\tilde{n}+1}^{-1} \mathcal{K}\bigg(\frac{q}{b_n}\bigg) \bigg\{\frac{1}{\tilde{n}}\sum_{t=-q+1}^{\tilde{n}}c_{\ell_1,t+q}r_{\ell_2,t}^{(1)}\bigg\}\bigg|  \,.
\end{align*}
For given $\ell_1$ and $\ell_2$, there exist unique triples $(i_1,j_1,k_1), (i_2,j_2,k_2)\in \cI\times\{-l_n,\ldots,l_n\}$ such that $c_{\ell_1,t}=(2\pi)^{-1}\{x_{i_1,t+l_n+k_1}x_{j_1,t+l_n}-\gamma_{i_1,j_1}(k_1)\}$ and $r_{\ell_2,t}^{(1)}=(2\pi)^{-1}x_{i_2,t+l_n+k_2}\bar{x}_{j_2}$, respectively, which implies
\begin{align}\label{eq:A2.exp}
A_2(\ell_1,\ell_2)
\leq&~ \underbrace{|\bar{x}_{j_2}|\bigg|\sum_{q=0}^{\tilde n-1}\cK\bigg(\frac{q}{b_n}\bigg)\bigg(\frac{1}{4\pi^2\tilde{n}}\sum_{t=q+1}^{\tilde n} x_{i_1,t+l_n+k_1}x_{j_1,t+l_n}x_{i_2,t-q+l_n+k_2}\bigg)\bigg|}_{A_{2,1}(\ell_1,\ell_2)} \notag\\
& +\underbrace{|\gamma_{i_1,j_1}(k_1)||\bar{x}_{j_2}|\bigg|\sum_{q=0}^{\tilde n-1}\cK\bigg(\frac{q}{b_n}\bigg)\bigg(\frac{1}{4\pi^2\tilde{n}}\sum_{t=q+1}^{\tilde n}x_{i_2,t-q+l_n+k_2}\bigg)\bigg|}_{A_{2,2}(\ell_1,\ell_2)} \,.
\end{align}
In the sequel, we will consider the convergence rates of $\max_{\ell_1,\ell_2}A_{2,1}(\ell_1,\ell_2)$ and $\max_{\ell_1,\ell_2}A_{2,2}(\ell_1,\ell_2)$, respectively.

Given $M_{2n}=o(n)\rightarrow\infty$ satisfying $M_{2n}\geq 2l_n$, by the triangle inequality,
\begin{align}\label{eq:A21.exp}
& \bigg|\sum_{q=0}^{\tilde n-1}\cK\bigg(\frac{q}{b_n}\bigg)\bigg(\frac{1}{4\pi^2\tilde{n}}\sum_{t=q+1}^{\tilde n} x_{i_1,t+l_n+k_1}x_{j_1,t+l_n}x_{i_2,t-q+l_n+k_2}\bigg)\bigg|  \notag\\
&~~~~~~~~~~~~ \leq \bigg|\sum_{q=0}^{M_{2n}}\cK\bigg(\frac{q}{b_n}\bigg)\bigg(\frac{1}{4\pi^2\tilde{n}}\sum_{t=q+1}^{\tilde n} x_{i_1,t+l_n+k_1}x_{j_1,t+l_n}x_{i_2,t-q+l_n+k_2}\bigg)\bigg|  \notag\\
&~~~~~~~~~~~~~~~ +\bigg|\sum_{q=M_{2n}+1}^{\tilde n-1}\cK\bigg(\frac{q}{b_n}\bigg)\bigg(\frac{1}{4\pi^2\tilde{n}}\sum_{t=q+1}^{\tilde n} x_{i_1,t+l_n+k_1}x_{j_1,t+l_n}x_{i_2,t-q+l_n+k_2}\bigg)\bigg| \,.
\end{align}
By Lemma 2 in the supplementary material of \cite{CTW_2013_supp} 
and Condition \ref{as:moment},
it holds that $
\max_{t,q,k_1,k_2,i_1,i_2,j_1,j_2}
\bP\{|x_{i_1,t+l_n+k_1}x_{j_1,t+l_n}x_{i_2,t-q+l_n+k_2}-
\bE(x_{i_1,t+l_n+k_1}x_{j_1,t+l_n}x_{i_2,t-q+l_n+k_2})|>u\}\lesssim \exp(-Cu^{2/3})$
for any $u>0$.	
For given $(q,i_1,i_2,j_1,j_2,k_1,k_2)$ such that $0\leq q\leq M_{2n}$ and $-l_n\leq k_1,k_2\leq l_n$, 	
$\{x_{i_1,t+l_n+k_1}x_{j_1,t+l_n}x_{i_2,t-q+l_n+k_2}\}_{t=q+1}^{\tilde{n}}$ is  an $\alpha$-mixing sequence with $\alpha$-mixing coefficients $\tilde{\alpha}_{3}(k)\lesssim \exp(-C|k-q-2l_n|_+)$. Due to
$$\max_{t,q,i_1,i_2,j_1,j_2,k_1,k_2} |\bE(x_{i_1,t+l_n+k_1}x_{j_1,t+l_n}x_{i_2,t-q+l_n+k_2})|\leq C\,,$$
by
Lemma L1 in the supplementary material of \cite{CCW_2021_supp}  
with $(\tilde{B}_{\tilde n},\tilde{L}_{\tilde n},\tilde{j}_{\tilde n}, r_1, r_2,r)=(1,1,q+2l_n,2/3,1,2/7)$, it holds that
\begin{align}\label{eq:x3.tail}
&\max_{0\leq q\leq M_{2n}}\max_{(i_1,j_1),(i_2,j_2)\in\cI}\max_{-l_n\leq k_1,k_2\leq l_n} \bP\bigg(\bigg|\sum_{t=q+1}^{\tilde n}x_{i_1,t+l_n+k_1}x_{j_1,t+l_n}x_{i_2,t-q+l_n+k_2}\bigg|\geq u\bigg)  \notag\\
&~~~~~~~~~~~~~~ \lesssim \exp\bigg(-\frac{Cu^2}{M_{2n}n}\bigg)+\exp\bigg(-\frac{Cu^{2/7}}{M_{2n}^{2/7}}\bigg)
\end{align}
for any $u\geq \tilde{C}_2\tilde{n}$ with some sufficiently large constant $\tilde{C}_2>0$.
Notice that $\sum_{q=0}^{M_{2n}}|\mathcal{K}(q/b_n)|\leq \tilde{C}_2'n^{\rho}$ and $M_{2n}=o(n)$.  By  \eqref{eq:x3.tail}, for some sufficiently large positive constant $\tilde{C}_2''$ satisfying $4\pi^2\tilde{C}_2''\geq \tilde{C}_2\tilde{C}_2'$, it holds that
\begin{align}\label{eq:M2pro}
&\max_{i_1,j_1,i_2,j_2,k_1,k_2}\bP\bigg\{\bigg|\sum_{q=0}^{M_{2n}}\cK\bigg(\frac{q}{b_n}\bigg)\bigg(\frac{1}{4\pi^2\tilde{n}}\sum_{t=q+1}^{\tilde n} x_{i_1,t+l_n+k_1}x_{j_1,t+l_n}x_{i_2,t-q+l_n+k_2}\bigg)\bigg|>\tilde{C}_2''n^{\rho}\bigg\}  \notag\\
&~~~~~~~~~~~~ \leq M_{2n}\max_{0\leq q\leq M_{2n}}\bP\bigg\{\bigg|\frac{1}{\tilde n}\sum_{t=q+1}^{\tilde n}x_{i_1,t+l_n+k_1}x_{j_1,t+l_n}x_{i_2,t-q+l_n+k_2} \bigg| > \tilde{C}_2\bigg\}  \\
&~~~~~~~~~~~~ \lesssim M_{2n}\exp\bigg(-\frac{Cn}{M_{2n}}\bigg) + M_{2n}\exp\bigg(-\frac{Cn^{2/7}}{M_{2n}^{2/7}}\bigg)
\lesssim M_{2n}\exp\bigg(-\frac{Cn^{2/7}}{M_{2n}^{2/7}}\bigg)  \,.  \notag
\end{align}
By Condition \ref{as:kernel}(ii) and $b_n\asymp n^{\rho}$ for some $\rho\in(0,1)$, we have
$\sum_{q=M_{2n}+1}^{\tilde{n}-1}|\mathcal{K}(q/b_n)|
\lesssim n^{\rho\vartheta}M_{2n}^{1-\vartheta}$. By the triangle inequality and the Bonferroni inequality,
\begin{align*} 
& \max_{i_1,j_1,i_2,j_2,k_1,k_2}\bP\bigg\{\bigg|\sum_{q=M_{2n}+1}^{\tilde n-1}\cK\bigg(\frac{q}{b_n}\bigg)\bigg(\frac{1}{4\pi^2\tilde{n}}\sum_{t=q+1}^{\tilde n} x_{i_1,t+l_n+k_1}x_{j_1,t+l_n}x_{i_2,t-q+l_n+k_2}\bigg)\bigg|>\tilde{C}_2''n^{\rho} \bigg\}  \notag\\
&~~~~~~~~~ \leq  \max_{i_1,j_1,i_2,j_2,k_1,k_2} \sum_{q=M_{2n}+1}^{\tilde{n}-1}
\bP\bigg\{\bigg|\frac{1}{\tilde{n}}\sum_{t=q+1}^{\tilde{n}}x_{i_1,t+l_n+k_1}x_{j_1,t+l_n}x_{i_2,t-q+l_n+k_2}\bigg| > \frac{CM_{2n}^{\vartheta-1}}{n^{\rho(\vartheta-1)}} \bigg\}\notag\\
&~~~~~~~~~\lesssim  n^2 \max_{t,q,i_1,j_1,i_2,j_2,k_1,k_2}
\bP\bigg\{ |x_{i_1,t+l_n+k_1}x_{j_1,t+l_n}x_{i_2,t-q+l_n+k_2}|> \frac{CM_{2n}^{\vartheta-1}}{n^{\rho(\vartheta-1)}} \bigg\} \notag\\
&~~~~~~~~~\lesssim  n^2\exp\bigg\{ -\frac{CM_{2n}^{(2\vartheta-2)/3}}{n^{(2\rho\vartheta-2\rho)/3}}  \bigg\} \,.
\end{align*}
Selecting $M_{2n}=n^{(7\rho\vartheta-7\rho+3)/(7\vartheta-4)}$, together with \eqref{eq:M2pro}, by \eqref{eq:A21.exp},
\begin{align}\label{eq:sumx3.tail}
&\max_{i_1,j_1,i_2,j_2,k_1,k_2}\bP\bigg\{\bigg|\sum_{q=0}^{\tilde n-1}\cK\bigg(\frac{q}{b_n}\bigg)\bigg(\frac{1}{4\pi^2\tilde{n}}\sum_{t=q+1}^{\tilde n} x_{i_1,t+l_n+k_1}x_{j_1,t+l_n}x_{i_2,t-q+l_n+k_2}\bigg)\bigg| >2\tilde{C}_2''n^{\rho}\bigg\}  \notag\\
&~~~~~~~~~~ 	\lesssim M_{2n}\exp\bigg(-\frac{Cn^{2/7}}{M_{2n}^{2/7}}\bigg)  + n^2\exp\bigg\{ -\frac{CM_{2n}^{(2\vartheta-2)/3}}{n^{(2\rho\vartheta-2\rho)/3}}  \bigg\}  \\
&~~~~~~~~~~ 	\lesssim	n^2\exp\{-Cn^{(2\vartheta+2\rho-2\rho\vartheta-2)/(7\vartheta-4)}\} \lesssim \exp\{-Cn^{(2\vartheta+2\rho-2\rho\vartheta-2)/(7\vartheta-4)}\}  \,. \notag
\end{align}
Applying Lemma L1 in the supplementary material of \cite{CCW_2021_supp} 
with $(\tilde{B}_{\tilde n}, \tilde{L}_{\tilde n}, \tilde{j}_{\tilde n}, r_1, r_2, r) =(1,1,0,2,1,1/3)$,
\begin{align}\label{eq:sumx.tail}
\max_{j\in[p]}\bP\bigg(\bigg|\sum_{t=s_1}^{s_2} x_{j,t}\bigg|\geq u\bigg)\lesssim \exp\bigg(-\frac{Cu^2}{s_2-s_1}\bigg)+\exp(-Cu^{1/3})
\end{align}
for any  $1\leq s_1\leq s_2\leq n$ and $u>0$.
By the Bonferroni inquality and \eqref{eq:sumx.tail},  we have
\begin{align}\label{eq:A21.tailpro}
&\max_{\ell_1,\ell_2}\bP\{A_{2,1}(\ell_1,\ell_2)>u\}   \notag\\
&~~~~~~	\leq    \max_{i_1,j_1,i_2,k_1,k_2} \bP\bigg\{\bigg|\sum_{q=0}^{\tilde n-1}\cK\bigg(\frac{q}{b_n}\bigg)\bigg(\frac{1}{4\pi^2\tilde{n}}\sum_{t=q+1}^{\tilde n} x_{i_1,t+l_n+k_1}x_{j_1,t+l_n}x_{i_2,t-q+l_n+k_2}\bigg)\bigg|>2\tilde{C}_2''n^{\rho} \bigg\}  \notag\\
&~~~~~~~~~ + \max_{j_2}\bP\bigg(|\bar{x}_{j_2}|\geq \frac{u}{2\tilde{C}_2''n^\rho}\bigg)   \\
&~~~~~~ \lesssim \exp\{-Cn^{(2\vartheta+2\rho-2\rho\vartheta-2)/(7\vartheta-4)}\}   +\exp(-Cn^{1-2\rho}u^2)
+\exp\{-Cn^{(1-\rho)/3}u^{1/3}\}   \notag
\end{align}
for any $u>0$.
Since $\max_{(i,j)\in\cI}\max_{-l_n\leq k\leq l_n}|\gamma_{i,j}(k)|\leq C$,
by the Bonferroni inequality and \eqref{eq:sumx.tail},
\begin{align*}
&\max_{\ell_1,\ell_2}\bP\{A_{2,2}(\ell_1,\ell_2)>u\} \\
&~~~ \leq  \max_{i_2,k_2}\bP\bigg\{\bigg|\sum_{q=0}^{\tilde n-1}\cK\bigg(\frac{q}{b_n}\bigg)\bigg(\frac{1}{4\pi^2\tilde n}\sum_{t=q+1}^{\tilde n}x_{i_2,t-q+l_n+k_2}\bigg)\bigg|>Cb_n^{1/2}u^{1/2}\bigg\} + \max_{j_2}\bP\bigg(|\bar{x}_{j_2}|>\frac{Cu^{1/2}}{b_n^{1/2}}\bigg)  \\
&~~~\leq n\max_{i_2,k_2}\max_{0\leq q\leq \tilde{n}-1}\bP\bigg( \bigg|\frac{1}{\tilde n}\sum_{t=q+1}^{\tilde n}x_{i_2,t-q+l_n+k_2}\bigg|>\frac{Cu^{1/2}}{b_n^{1/2}} \bigg)  + \max_{j_2}\bP\bigg(|\bar{x}_{j_2}|>\frac{Cu^{1/2}}{b_n^{1/2}}\bigg) \\
&~~~\lesssim n\exp(-Cn^{1-\rho}u) + n\exp\{-Cn^{(2-\rho)/6}u^{1/6}\}
\end{align*}
for any $u>0$.
Together with \eqref{eq:A21.tailpro},   \eqref{eq:A2.exp} implies
\begin{align}\label{eq:A2l1l2tail}
&\max_{\ell_1,\ell_2}\bP\{A_2(\ell_1,\ell_2)>u\} \leq \max_{\ell_1,\ell_2}\bP\bigg\{A_{2,1}(\ell_1,\ell_2)>\frac{u}{2}\bigg\} + \max_{\ell_1,\ell_2}\bP\bigg\{A_{2,2}(\ell_1,\ell_2)>\frac{u}{2}\bigg\}   \notag\\
&~~~~~~~~~~~~ \lesssim \exp\{-Cn^{(2\vartheta+2\rho-2\rho\vartheta-2)/(7\vartheta-4)}\} +\exp(-Cn^{1-2\rho}u^2)+\exp\{-Cn^{(1-\rho)/3}u^{1/3}\}    \notag\\
&~~~~~~~~~~~~~~~
+n\exp(-Cn^{1-\rho}u) + n\exp\{-Cn^{(2-\rho)/6}u^{1/6}\}
\end{align}
for any $u>0$. Analogously, we can  show the upper bound given in \eqref{eq:A2l1l2tail}  also holds for $\max_{\ell_1,\ell_2}\bP\{\tilde{A}_2(\ell_1,\ell_2)>u\}$.
Recall $A_2\leq \max_{\ell_1,\ell_2}A_2(\ell_1,\ell_2)+\max_{\ell_1,\ell_2}\tilde{A}_2(\ell_1,\ell_2)$.
By Bonferroni inequality, 	for any $u>0$,
\begin{align}\label{eq:A2tail}
\bP(A_2>u)
\lesssim&~ r^2l_n^2\exp\{-Cn^{(2\vartheta+2\rho-2\rho\vartheta-2)/(7\vartheta-4)}\} +r^2l_n^2\exp(-Cn^{1-2\rho}u^2)  \notag\\
&+r^2l_n^2\exp\{-Cn^{(1-\rho)/3}u^{1/3}\}
+nr^2l_n^2\exp(-Cn^{1-\rho}u) \notag\\
&+ nr^2l_n^2\exp\{-Cn^{(2-\rho)/6}u^{1/6}\} \,.  
\end{align}
Since $l_n=o(n)$ and $r\geq n^{\kappa}$ for some sufficiently small constant $\kappa>0$, then
$
A_2
=  O_{\rm p}\{n^{(2\rho-1)/2}\log^{1/2}(r)\}
$
provided that $\log r\ll  \min\{n^{(2\vartheta+2\rho-2\rho\vartheta-2)/(7\vartheta-4)}, n^{1/5}\}$.
By the definition of $A_3$, using the same arguments, we know the upper bound given in \eqref{eq:A2tail} also holds for $\bP(A_3>u)$ and
$
A_3
=  O_{\rm p}\{n^{(2\rho-1)/2}\log^{1/2}(r)\}
$.   $\hfill\Box$

\subsubsection{Convergence rates of $A_4$ and $A_5$}
Notice that $A_4\leq \max_{\ell_1,\ell_2}A_4(\ell_1,\ell_2)+\max_{\ell_1,\ell_2}\tilde{A}_4(\ell_1,\ell_2)$ and  $A_5\leq \max_{\ell_1,\ell_2}A_5(\ell_1,\ell_2)+\max_{\ell_1,\ell_2}\tilde{A}_5(\ell_1,\ell_2)$, where
\begin{align*} 
A_4(\ell_1,\ell_2)
:=&~ \bigg|\sum_{q=0}^{\tilde{n}-1} \mathcal{K}\bigg(\frac{q}{b_n}\bigg) \bigg\{\frac{1}{\tilde{n}}\sum_{t=q+1}^{\tilde{n}}c_{\ell_1,t}r_{\ell_2}^{(3)}\bigg\}\bigg|  \,, \\
\tilde{A}_4(\ell_1,\ell_2):=&~ \bigg|\sum_{q=-\tilde{n}+1}^{-1} \mathcal{K}\bigg(\frac{q}{b_n}\bigg) \bigg\{\frac{1}{\tilde{n}}\sum_{t=-q+1}^{\tilde{n}}c_{\ell_1,t+q}r_{\ell_2}^{(3)}\bigg\}\bigg|  \,, \\
A_5(\ell_1,\ell_2):=&~\bigg|\sum_{q=0}^{\tilde{n}-1} \mathcal{K}\bigg(\frac{q}{b_n}\bigg) \bigg(\frac{1}{\tilde{n}}\sum_{t=q+1}^{\tilde{n}}c_{\ell_1,t}\tilde{\gamma}_{\ell_2}\bigg)\bigg|  \,, \\
\tilde{A}_5(\ell_1,\ell_2):=&~ \bigg|\sum_{q=-\tilde{n}+1}^{-1} \mathcal{K}\bigg(\frac{q}{b_n}\bigg) \bigg\{\frac{1}{\tilde{n}}\sum_{t=-q+1}^{\tilde{n}}c_{\ell_1,t+q}\tilde{\gamma}_{\ell_2}\bigg\}\bigg|  \,.
\end{align*}
For given $\ell_1$ and $\ell_2$, there exist unique triples $(i_1,j_1,k_1), (i_2,j_2,k_2)\in \cI\times\{-l_n,\ldots,l_n\}$ such that $c_{\ell_1,t}=(2\pi)^{-1}\{x_{i_1,t+l_n+k_1}x_{j_1,t+l_n}-\gamma_{i_1,j_1}(k_1)\}$, $r_{\ell_2}^{(3)}=(2\pi)^{-1}\bar{x}_{i_2}\bar{x}_{j_2}$ and $\tilde{\gamma}_{\ell_2}=(2\pi)^{-1}\{\hat\gamma_{i_2,j_2}(k_2)-\gamma_{i_2,j_2}(k_2)\}$, respectively, which implies
\begin{align*}
A_4(\ell_1,\ell_2)
=&~ |\bar{x}_{i_2}| |\bar{x}_{j_2}| \bigg|\sum_{q=0}^{\tilde n-1}\cK\bigg(\frac{q}{b_n}\bigg)\bigg(\frac{1}{4\pi^2\tilde{n}}\sum_{t=q+1}^{\tilde n} c_{\ell_1,t}\bigg)\bigg| \,, \\
A_5(\ell_1,\ell_2)
=&~ |\hat{\gamma}_{i_2,j_2}(k_2)-\gamma_{i_2,j_2}(k_2)|  \bigg|\sum_{q=0}^{\tilde n-1}\cK\bigg(\frac{q}{b_n}\bigg)\bigg(\frac{1}{4\pi^2\tilde{n}}\sum_{t=q+1}^{\tilde n} c_{\ell_1,t}\bigg)\bigg| \,.
\end{align*}
Due to $b_n\asymp n^\rho$, by \eqref{eq:sumx.tail} and  Lemma \ref{la:gammatail},
\begin{align}\label{eq:A4tail}
\max_{\ell_1,\ell_2}\bP\{A_4(\ell_1,\ell_2)>u\}
\lesssim&~ n \max_{\ell_1}\max_{q}
\bP\bigg(\bigg|\sum_{t=q+1}^{\tilde n} c_{\ell_1,t}\bigg|>\frac{Cnu^{1/3}}{b_n^{1/3}}\bigg)
+ \max_{j}\bP\bigg(|\bar{x}_{j}|>\frac{Cu^{1/3}}{b_n^{1/3}}\bigg) \notag\\
\lesssim&~  n\exp\{-Cn^{(3-2\rho)/3}l_n^{-1}u^{2/3}\} + n\exp\{-Cn^{(3-\rho)/9}l_n^{-1/3}u^{1/9}\}
\end{align}
for any $u>0$. 	By Lemmas \ref{la:gammatail} and \ref{la:tildegammatail},
\begin{align}\label{eq:A51tail}
\max_{\ell_1,\ell_2}\bP\{A_5(\ell_1,\ell_2)>u\}
\leq&~ n \max_{\ell_1}\max_{0\leq q\leq \tilde{n}-1}\bP\bigg(\bigg|\sum_{t=q+1}^{\tilde n} c_{\ell_1,t}\bigg|>\frac{Cnu^{1/2}}{b_n^{1/2}}\bigg) \notag\\
&	+ \max_{i_2,j_2,k_2}\bP\bigg\{|\hat{\gamma}_{i_2,j_2}(k_2)-\gamma_{i_2,j_2}(k_2)|>\frac{Cu^{1/2}}{b_n^{1/2}}\bigg\}  \notag\\
\lesssim&~ n\exp(-Cn^{1-\rho}l_n^{-1}u)
+n\exp\{-Cn^{(2-\rho)/6}l_n^{-1/3}u^{1/6}\} \notag\\
&+ \exp\{-Cn^{(2-\rho)/2}u^{1/2}\}
+   \exp\{-Cn^{(4-\rho)/12}u^{1/12}\}
\end{align}
for any $u\gg l_n^2n^{\rho-2}$.
Analogously, we can  show the upper bounds given in \eqref{eq:A4tail} and \eqref{eq:A51tail}  also hold for $\max_{\ell_1,\ell_2}\bP\{\tilde{A}_4(\ell_1,\ell_2)>u\}$ and $\max_{\ell_1,\ell_2}\bP\{\tilde{A}_5(\ell_1,\ell_2)>u\}$, respectively.
Recall $A_4\leq \max_{\ell_1,\ell_2}A_4(\ell_1,\ell_2)+\max_{\ell_1,\ell_2}\tilde{A}_4(\ell_1,\ell_2)$ and $A_5\leq \max_{\ell_1,\ell_2}A_5(\ell_1,\ell_2)+\max_{\ell_1,\ell_2}\tilde{A}_5(\ell_1,\ell_2)$.
By the Bonferroni inequality,
\begin{align}\label{eq:A4.tail}
\bP(A_4>u)
\lesssim  nr^2l_n^2\exp\{-Cn^{(3-2\rho)/3}l_n^{-1}u^{2/3}\} + nr^2l_n^2\exp\{-Cn^{(3-\rho)/9}l_n^{-1/3}u^{1/9}\}
\end{align}
for any $u>0$, and
\begin{align}\label{eq:A5.tail}
\bP(A_5>u)
\lesssim&~
nr^2l_n^2\exp(-Cn^{1-\rho}l_n^{-1}u)
+nr^2l_n^2\exp\{-Cn^{(2-\rho)/6}l_n^{-1/3}u^{1/6}\}  \notag\\
& + r^2l_n^2\exp\{-Cn^{(2-\rho)/2}u^{1/2}\}+r^2l_n^2\exp\{-Cn^{(4-\rho)/12}u^{1/12}\}
\end{align}
for any $u\gg l_n^2n^{\rho-2}$.
Since $l_n=o(n)$ and $r\geq n^{\kappa}$ for some sufficiently small constant $\kappa>0$, then
$
A_4
= O_{\rm p}\{n^{(2\rho-3)/2}l_n^{3/2}\log^{3/2}(r)\}
$ and
$
A_5 = O_{\rm p}(n^{\rho-1}l_n\log r)
$
provided that $\log r\ll n^{1/5}l_n^{-1/5}$.  $\hfill\Box$

\subsubsection{Convergence rates of $A_7$, $A_8$ and $A_{13}$}
Notice that $A_7\leq \max_{\ell_1,\ell_2}A_7(\ell_1,\ell_2)+\max_{\ell_1,\ell_2}\tilde{A}_7(\ell_1,\ell_2)$, where
\begin{align*} 
A_7(\ell_1,\ell_2)
:=&~ \bigg|\sum_{q=0}^{\tilde{n}-1} \mathcal{K}\bigg(\frac{q}{b_n}\bigg) \bigg\{\frac{1}{\tilde{n}}\sum_{t=q+1}^{\tilde{n}}r^{(1)}_{\ell_1,t}r_{\ell_2,t-q}^{(1)}\bigg\}\bigg|  \,, \\
\tilde{A}_7(\ell_1,\ell_2):=&~ \bigg|\sum_{q=-\tilde{n}+1}^{-1} \mathcal{K}\bigg(\frac{q}{b_n}\bigg) \bigg\{\frac{1}{\tilde{n}}\sum_{t=-q+1}^{\tilde{n}}r^{(1)}_{\ell_1,t+q}r_{\ell_2,t}^{(1)}\bigg\}\bigg|  \,.
\end{align*}
For given $\ell_1$ and $\ell_2$, there exist unique triples $(i_1,j_1,k_1), (i_2,j_2,k_2)\in \cI\times\{-l_n,\ldots,l_n\}$ such that $r_{\ell_1,t}^{(1)}=(2\pi)^{-1}x_{i_1,t+l_n+k_1}\bar{x}_{j_1}$ and $r_{\ell_2,t}^{(1)}=(2\pi)^{-1}x_{i_2,t+l_n+k_2}\bar{x}_{j_2}$, respectively, which implies
\begin{align*}
{A}_7(\ell_1,\ell_2)
= |\bar{x}_{j_1}| |\bar{x}_{j_2}| \bigg|\sum_{q=0}^{\tilde n-1}\cK\bigg(\frac{q}{b_n}\bigg)\bigg\{\frac{1}{4\pi^2\tilde n}\sum_{t=q+1}^{\tilde{n}}x_{i_1,t+l_n+k_1}x_{i_2,t+l_n+k_2-q}\bigg\}\bigg|  \,.
\end{align*}
By Lemma 2 in the supplementary material of \cite{CTW_2013_supp}  
and Condition \ref{as:moment}, we have
$\max_{t,q,k_1,k_2,i_1,i_2}
\bP\{|x_{i_1,t+l_n+k_1}x_{i_2,t+l_n+k_2-q}-\bE(x_{i_1,t+l_n+k_1}x_{i_2,t+l_n+k_2-q})|>u\}\lesssim \exp(-Cu)$
for any $u>0$.	
Given $M_{7n}=o(n)\rightarrow\infty$ satisfying $M_{7n}\geq 2l_n$,
using the technique for deriving \eqref{eq:sumx3.tail},
\begin{align*}
&\max_{i_1,i_2,k_1,k_2}\bP\bigg\{\bigg|\sum_{q=0}^{\tilde n-1}\cK\bigg(\frac{q}{b_n}\bigg)\bigg(\frac{1}{4\pi^2\tilde{n}}\sum_{t=q+1}^{\tilde n} x_{i_1,t+l_n+k_1}x_{i_2,t+l_n+k_2-q}\bigg)\bigg| >2\tilde{C}_7''n^{\rho}\bigg\}  \\
&~~~~~~~~~~ 	\lesssim M_{7n}\exp(-CM_{7n}^{-1/3}n^{1/3}) + n^2\exp( -CM_{7n}^{\vartheta-1}n^{-\rho\vartheta+\rho}  )
\end{align*}
for some sufficiently large constant $\tilde{C}_7''>0$. 	Selecting $M_{7n}=n^{(3\rho\vartheta-3\rho+1)/(3\vartheta-2)}$,
\begin{align*}
&\max_{i_1,i_2,k_1,k_2}\bP\bigg\{\bigg|\sum_{q=0}^{\tilde n-1}\cK\bigg(\frac{q}{b_n}\bigg)\bigg(\frac{1}{4\pi^2\tilde{n}}\sum_{t=q+1}^{\tilde n} x_{i_1,t+l_n+k_1}x_{i_2,t+l_n+k_2-q}\bigg)\bigg| >2\tilde{C}_7''n^{\rho}\bigg\}  \\
&~~~~~~~~~~ 	\lesssim	n^2\exp\{-Cn^{(\vartheta+\rho-\rho\vartheta-1)/(3\vartheta-2)}\}
\lesssim \exp\{-Cn^{(\vartheta+\rho-\rho\vartheta-1)/(3\vartheta-2)}\} \,.
\end{align*}
By the Bonferroni inequality and  \eqref{eq:sumx.tail},  we have
\begin{align}\label{eq:A71tail}
&\max_{\ell_1,\ell_2}\bP\{A_{7}(\ell_1,\ell_2)>u\}   \notag\\
&~~~~~~	\lesssim    \max_{i_1,i_2,k_1,k_2} \bP\bigg\{\bigg|\sum_{q=0}^{\tilde n-1}\cK\bigg(\frac{q}{b_n}\bigg)\bigg(\frac{1}{4\pi^2\tilde{n}}\sum_{t=q+1}^{\tilde n} x_{i_1,t+l_n+k_1}x_{i_2,t+l_n+k_2-q}\bigg)\bigg|>2\tilde{C}_7''n^{\rho} \bigg\}  \notag\\
&~~~~~~~~~ + \max_{j}\bP\bigg\{|\bar{x}_{j}|\geq \frac{u^{1/2}}{(2\tilde{C}_7'')^{1/2}n^{\rho/2}}\bigg\}   \\
&~~~~~~ \lesssim \exp\{-Cn^{(\vartheta+\rho-\rho\vartheta-1)/(3\vartheta-2)}\}
+\exp(-Cn^{1-\rho}u)+\exp\{-Cn^{(2-\rho)/6}u^{1/6}\}  \notag
\end{align}
for any $u>0$.
Analogously, we can  show the upper bound given in \eqref{eq:A71tail}  also holds for $\max_{\ell_1,\ell_2}\bP\{\tilde{A}_7(\ell_1,\ell_2)>u\}$.
Since $A_7\leq \max_{\ell_1,\ell_2}A_7(\ell_1,\ell_2)+\max_{\ell_1,\ell_2}\tilde{A}_7(\ell_1,\ell_2)$,
by Bonferroni inequality, 	for any $u>0$,
\begin{align}\label{eq:A7tail}
\bP(A_7>u)
\lesssim&\, r^2l_n^2\exp\{-Cn^{(\vartheta+\rho-\rho\vartheta-1)/(3\vartheta-2)}\}
+r^2l_n^2\exp(-Cn^{1-\rho}u) \notag\\
&+r^2l_n^2\exp\{-Cn^{(2-\rho)/6}u^{1/6}\} \,.
\end{align}
Since $l_n=o(n)$ and $r\geq n^{\kappa}$ for some sufficiently small constant $\kappa>0$, then
$
A_7
=  O_{\rm p}(n^{\rho-1}\log r)
$
provided that $\log r\ll  \min\{n^{(\vartheta+\rho-\rho\vartheta-1)/(3\vartheta-2)}, n^{1/5}\}$.
By the definition of $A_8$ and $A_{13}$, using the same arguments, we know the upper bound given in \eqref{eq:A7tail} also holds for $\bP(A_8>u)$ and $\bP(A_{13}>u)$. Furthermore,
$
A_8=   O_{\rm p}(n^{\rho-1}\log r)=A_{13} $.  $\hfill\Box$

\subsubsection{Convergence rates of $A_9$, $A_{14}$, $A_{19}$ and $A_{25}$}
Notice that $A_9\leq \max_{\ell_1,\ell_2}A_9(\ell_1,\ell_2)+\max_{\ell_1,\ell_2}\tilde{A}_9(\ell_1,\ell_2)$, $A_{19}\leq \max_{\ell_1,\ell_2}A_{19}(\ell_1,\ell_2)+\max_{\ell_1,\ell_2}\tilde{A}_{19}(\ell_1,\ell_2)$ and $A_{25}\leq \max_{\ell_1,\ell_2}A_{25}(\ell_1,\ell_2)+\max_{\ell_1,\ell_2}\tilde{A}_{25}(\ell_1,\ell_2)$, where
\begin{align*} 
A_9(\ell_1,\ell_2):=&~ \bigg|\sum_{q=0}^{\tilde n-1}\cK\bigg(\frac{q}{b_n}\bigg)\bigg\{{\frac{1}{\tilde n}}\sum_{t=q+1}^{\tilde n} r_{\ell_1,t}^{(1)}r_{\ell_2}^{(3)}\bigg\}\bigg| \,, \\
\tilde{A}_9(\ell_1,\ell_2):=&~ \bigg|\sum_{q=-\tilde{n}+1}^{-1} \mathcal{K}\bigg(\frac{q}{b_n}\bigg) \bigg\{\frac{1}{\tilde{n}}\sum_{t=-q+1}^{\tilde{n}}r^{(1)}_{\ell_1,t+q}r_{\ell_2}^{(3)}\bigg\}\bigg|  \,, \\
A_{19}(\ell_1,\ell_2):=&~ \bigg|\sum_{q=0}^{\tilde n-1}\cK\bigg(\frac{q}{b_n}\bigg)\bigg\{{\frac{1}{\tilde n}}\sum_{t=q+1}^{\tilde n} r_{\ell_1}^{(3)}r_{\ell_2}^{(3)}\bigg\}\bigg| \,, \\
\tilde{A}_{19}(\ell_1,\ell_2):=&~ \bigg|\sum_{q=-\tilde{n}+1}^{-1}\cK\bigg(\frac{q}{b_n}\bigg)\bigg\{{\frac{1}{\tilde n}}\sum_{t=-q+1}^{\tilde n} r_{\ell_1}^{(3)}r_{\ell_2}^{(3)}\bigg\}\bigg| \,, \\
A_{25}(\ell_1,\ell_2)
:=&~ \bigg|\sum_{q=0}^{\tilde{n}-1} \mathcal{K}\bigg(\frac{q}{b_n}\bigg) \bigg(\frac{1}{\tilde{n}}\sum_{t=q+1}^{\tilde{n}}\tilde{\gamma}_{\ell_1}\tilde{\gamma}_{\ell_2}\bigg)\bigg|  \,, \\
\tilde{A}_{25}(\ell_1,\ell_2):=&~ \bigg|\sum_{q=-\tilde{n}+1}^{-1} \mathcal{K}\bigg(\frac{q}{b_n}\bigg) \bigg(\frac{1}{\tilde{n}}\sum_{t=-q+1}^{\tilde{n}} \tilde{\gamma}_{\ell_1}\tilde{\gamma}_{\ell_2}\bigg)\bigg|  \,.
\end{align*}
For given $\ell_1$ and $\ell_2$, there exist unique triples $(i_1,j_1,k_1), (i_2,j_2,k_2)\in \cI\times\{-l_n,\ldots,l_n\}$ such that $r_{\ell_1,t}^{(1)}=(2\pi)^{-1}x_{i_1,t+l_n+k_1}\bar{x}_{j_1}$, $r_{\ell_1}^{(3)}=(2\pi)^{-1}\bar{x}_{i_1}\bar{x}_{j_1}$, $\tilde{\gamma}_{\ell_1}=(2\pi)^{-1}\{\hat{\gamma}_{i_1,j_1}(k_1)-\gamma_{i_1,j_1}(k_1)\}$, $r_{\ell_2}^{(3)}=(2\pi)^{-1}\bar{x}_{i_2}\bar{x}_{j_2}$ and $\tilde{\gamma}_{\ell_2}=(2\pi)^{-1}\{\hat{\gamma}_{i_2,j_2}(k_2)-\gamma_{i_2,j_2}(k_2)\}$, respectively, which implies
\begin{align*}
A_9(\ell_1,\ell_2)
=&~ |\bar{x}_{j_1}| |\bar{x}_{i_2}| |\bar{x}_{j_2}|\bigg|\sum_{q=0}^{\tilde n-1}\cK\bigg(\frac{q}{b_n}\bigg)\bigg(\frac{1}{4\pi^2\tilde n}\sum_{t=q+1}^{\tilde n} x_{i_1,t+l_n+k_1}\bigg)\bigg|  \,, \\
A_{19}(\ell_1,\ell_2)
=&~ |\bar{x}_{i_1}| |\bar{x}_{j_1}| |\bar{x}_{i_2}| |\bar{x}_{j_2}| \bigg|\sum_{q=0}^{\tilde n-1}\cK\bigg(\frac{q}{b_n}\bigg)\bigg(\frac{\tilde{n}-q}{4\pi^2\tilde n} \bigg)\bigg| \,, \\
A_{25}(\ell_1,\ell_2)
=&~ |\hat{\gamma}_{i_1,j_1}(k_1)-\gamma_{i_1,j_1}(k_1)| |\hat{\gamma}_{i_2,j_2}(k_2)-\gamma_{i_2,j_2}(k_2)|  \bigg|\sum_{q=0}^{\tilde n-1}\cK\bigg(\frac{q}{b_n}\bigg)\bigg(\frac{\tilde{n}-q}{4\pi^2\tilde n} \bigg)\bigg| \,.
\end{align*}
Due to $b_n\asymp n^\rho$ and $\sum_{q=0}^{\tilde{n}-1}|\cK(q/b_n)(\tilde{n}-q)/\tilde{n}|\lesssim b_n$, by  \eqref{eq:sumx.tail}, 	for any $u>0$,
\begin{align}
\max_{\ell_1,\ell_2}\bP\{A_9(\ell_1,\ell_2)>u\}
\lesssim&~ n \max_{i_1,k_1}\max_{0\leq q\leq \tilde{n}-1}
\bP\bigg(\bigg|\sum_{t=q+1}^{\tilde n}x_{i_1,t+l_n+k_1}\bigg|>\frac{Cnu^{1/4}}{b_n^{1/4}}\bigg)\notag\\
&+ \max_{j}\bP\bigg(|\bar{x}_{j}|>\frac{Cu^{1/4}}{b_n^{1/4}}\bigg) \notag\\
\lesssim&~   n\exp\{-Cn^{(2-\rho)/2}u^{1/2}\} + n\exp\{-Cn^{(4-\rho)/12}u^{1/12}\} \,, \label{eq:A91tail} \\
\max_{\ell_1,\ell_2}\bP\{A_{19}(\ell_1,\ell_2)>u\}
\lesssim&~ \max_{j}\bP\bigg(|\bar{x}_{j}|>\frac{Cu^{1/4}}{b_n^{1/4}}\bigg)  \notag\\
\lesssim&~   \exp\{-Cn^{(2-\rho)/2}u^{1/2}\} + \exp\{-Cn^{(4-\rho)/12}u^{1/12}\}  \,. \label{eq:A19tail}
\end{align}
By Lemma \ref{la:tildegammatail}, for any $u\gg l_n^2n^{\rho-2}$,
\begin{align}\label{eq:A251tail}
\max_{\ell_1,\ell_2}\bP\{A_{25}(\ell_1,\ell_2)>u\}
\lesssim&~ \max_{i,j,k}\bP\bigg\{|\hat{\gamma}_{i,j}(k)-\gamma_{i,j}(k)|>\frac{Cu^{1/2}}{b_n^{1/2}}\bigg\}  \notag\\
\lesssim&~  \exp(-Cn^{1-\rho}l_n^{-1}u)
+\exp\{-Cn^{(2-\rho)/6}l_n^{-1/3}u^{1/6}\} \notag\\
& +\exp\{-Cn^{(2-\rho)/2}u^{1/2}\} +\exp\{-Cn^{(4-\rho)/12}u^{1/12}\} \,.
\end{align}
Analogously, we can  show the upper bounds given in \eqref{eq:A91tail}, \eqref{eq:A19tail} and \eqref{eq:A251tail}  also hold for $\max_{\ell_1,\ell_2}\bP\{\tilde{A}_9(\ell_1,\ell_2)>u\}$, $\max_{\ell_1,\ell_2}\bP\{\tilde{A}_{19}(\ell_1,\ell_2)>u\}$ and $\max_{\ell_1,\ell_2}\bP\{\tilde{A}_{25}(\ell_1,\ell_2)>u\}$, respectively.
Recall $A_9\leq \max_{\ell_1,\ell_2}A_9(\ell_1,\ell_2)+\max_{\ell_1,\ell_2}\tilde{A}_9(\ell_1,\ell_2)$, $A_{19}\leq \max_{\ell_1,\ell_2}A_{19}(\ell_1,\ell_2)+\max_{\ell_1,\ell_2}\tilde{A}_{19}(\ell_1,\ell_2)$ and $A_{25}\leq \max_{\ell_1,\ell_2}A_{25}(\ell_1,\ell_2)+\max_{\ell_1,\ell_2}\tilde{A}_{25}(\ell_1,\ell_2)$.
By the Bonferroni inequality, 	for any $u>0$,
\begin{align}
\bP(A_9>u) +	\bP(A_{19}>u)
\lesssim&~  nr^2l_n^2\exp\{-Cn^{(2-\rho)/2}u^{1/2}\}\notag\\
&+ nr^2l_n^2\exp\{-Cn^{(4-\rho)/12}u^{1/12}\} \,, \label{eq:A9tail}  
\end{align}
and for any  $u\gg l_n^2n^{\rho-2}$,
\begin{align}\label{eq:A25tail}
\bP(A_{25}>u)
\lesssim&~  r^2l_n^2\exp(-Cn^{1-\rho}l_n^{-1}u)
+r^2l_n^2\exp\{-Cn^{(2-\rho)/6}l_n^{-1/3}u^{1/6}\} \notag\\
& +r^2l_n^2\exp\{-Cn^{(2-\rho)/2}u^{1/2}\} +r^2l_n^2\exp\{-Cn^{(4-\rho)/12}u^{1/12}\} \,.
\end{align}
Since $l_n=o(n)$ and $r\geq n^{\kappa}$ for some sufficiently small constant $\kappa>0$, then
$
A_9
= O_{\rm p}\{n^{\rho-2}\log^2(r)\} =A_{19} $
provided that $\log r\ll n^{1/5}$, and
$
A_{25}	
=  O_{\rm p}(n^{\rho-1}l_n\log r )$
provided that $\log r\ll n^{1/5}l_n^{-1/5}$.
By the definition of $A_{14}$, using the same arguments, we know the upper bound given in \eqref{eq:A9tail} also holds for $\bP(A_{14}>u)$ and
$
A_{14}= O_{\rm p}\{n^{\rho-2}\log^2(r)\}  $.
$\hfill\Box$

\subsubsection{Convergence rates of $A_{10}$, $A_{15}$ and $A_{20}$}	\label{subsec:A10}
Notice that $A_{10}\leq \max_{\ell_1,\ell_2}A_{10}(\ell_1,\ell_2)+\max_{\ell_1,\ell_2}\tilde{A}_{10}(\ell_1,\ell_2)$ and $A_{20}\leq \max_{\ell_1,\ell_2}A_{20}(\ell_1,\ell_2)+\max_{\ell_1,\ell_2}\tilde{A}_{20}(\ell_1,\ell_2)$, where
\begin{align*}
A_{10}(\ell_1,\ell_2)
:=&~ \bigg|\sum_{q=0}^{\tilde{n}-1} \mathcal{K}\bigg(\frac{q}{b_n}\bigg) \bigg\{\frac{1}{\tilde{n}}\sum_{t=q+1}^{\tilde{n}}r_{\ell_1,t}^{(1)}\tilde{\gamma}_{\ell_2}\bigg\}\bigg|  \,, \\
\tilde{A}_{10}(\ell_1,\ell_2):=&~ \bigg|\sum_{q=-\tilde{n}+1}^{-1} \mathcal{K}\bigg(\frac{q}{b_n}\bigg) \bigg\{\frac{1}{\tilde{n}}\sum_{t=-q+1}^{\tilde{n}}r_{\ell_1,t+q}^{(1)}\tilde{\gamma}_{\ell_2}\bigg\}\bigg|  \,, \\
A_{20}(\ell_1,\ell_2)
:=&~ \bigg|\sum_{q=0}^{\tilde{n}-1} \mathcal{K}\bigg(\frac{q}{b_n}\bigg) \bigg\{\frac{1}{\tilde{n}}\sum_{t=q+1}^{\tilde{n}}r_{\ell_1}^{(3)}\tilde{\gamma}_{\ell_2}\bigg\}\bigg|  \,, \\
\tilde{A}_{20}(\ell_1,\ell_2):=&~ \bigg|\sum_{q=-\tilde{n}+1}^{-1} \mathcal{K}\bigg(\frac{q}{b_n}\bigg) \bigg\{\frac{1}{\tilde{n}}\sum_{t=-q+1}^{\tilde{n}}r_{\ell_1}^{(3)}\tilde{\gamma}_{\ell_2}\bigg\}\bigg|  \,.
\end{align*}
For given $\ell_1$ and $\ell_2$, there exist unique triples $(i_1,j_1,k_1), (i_2,j_2,k_2)\in \cI\times\{-l_n,\ldots,l_n\}$ such that $r_{\ell_1,t}^{(1)}=(2\pi)^{-1}x_{i_1,t+l_n+k_1}\bar{x}_{j_1}$, $r_{\ell_1}^{(3)}=(2\pi)^{-1}\bar{x}_{i_1}\bar{x}_{j_1}$ and $\tilde{\gamma}_{\ell_2}=(2\pi)^{-1}\{\hat{\gamma}_{i_2,j_2}(k_2)-\gamma_{i_2,j_2}(k_2)\}$, respectively, which implies
\begin{align*}
A_{10}(\ell_1,\ell_2)
=&~ |\bar{x}_{j_1}| |\hat{\gamma}_{i_2,j_2}(k_2)-\gamma_{i_2,j_2}(k_2)| \bigg|\sum_{q=0}^{\tilde n-1}\cK\bigg(\frac{q}{b_n}\bigg)\bigg(\frac{1}{4\pi^2\tilde{n}}\sum_{t=q+1}^{\tilde n} x_{i_1,t+l_n+k_1}\bigg)\bigg| \,, \\
A_{20}(\ell_1,\ell_2)
=&~ |\bar{x}_{i_1}| |\bar{x}_{j_1}| |\hat{\gamma}_{i_2,j_2}(k_2)-\gamma_{i_2,j_2}(k_2)| \bigg|\sum_{q=0}^{\tilde n-1}\cK\bigg(\frac{q}{b_n}\bigg)\bigg(\frac{\tilde{n}-q}{4\pi^2\tilde{n}} \bigg)\bigg|  \,.
\end{align*}
Due to $b_n\asymp n^\rho$ and $\sum_{q=0}^{\tilde{n}-1}|\cK(q/b_n)(\tilde{n}-q)/\tilde{n}|\lesssim b_n$, by Lemma \ref{la:tildegammatail} and \eqref{eq:sumx.tail}, 	for any $u\gg l_n^3n^{\rho-3}$,
\begin{align}\label{eq:A101tail}
\max_{\ell_1,\ell_2}\bP\{A_{10}(\ell_1,\ell_2)>u\}
\leq&~ n \max_{i_1,k_1}\max_{0\leq q\leq \tilde{n}-1}
\bP\bigg(\bigg|\sum_{t=q+1}^{\tilde n} x_{i_1,t+l_n+k_1}\bigg|>\frac{Cnu^{1/3}}{b_n^{1/3}}\bigg)  \notag\\
&+ \max_{j_1}\bP\bigg(|\bar{x}_{j_1}|>\frac{Cu^{1/3}}{b_n^{1/3}}\bigg)\notag\\
&+ \max_{i_2,j_2,k_2}\bP\bigg\{|\hat{\gamma}_{i_2,j_2}(k_2)-\gamma_{i_2,j_2}(k_2)| |>\frac{Cu^{1/3}}{b_n^{1/3}}\bigg\} \notag\\
\lesssim&~ n\exp\{-Cn^{(3-2\rho)/3}l_n^{-1}u^{2/3}\} + n\exp\{-Cn^{(3-\rho)/9}l_n^{-1/3}u^{1/9}\} \notag\\
&+\exp\{-Cn^{(3-\rho)/3}u^{1/3}\}+\exp\{-Cn^{(6-\rho)/18}u^{1/18}\}  \,,  \\
\max_{\ell_1,\ell_2}\bP\{A_{20}(\ell_1,\ell_2)>u\}
\lesssim&~ \max_{j}\bP\bigg(|\bar{x}_{j}|>\frac{Cu^{1/3}}{b_n^{1/3}}\bigg)\notag\\
&+ \max_{i_2,j_2,k_2}\bP\bigg\{|\hat{\gamma}_{i_2,j_2}(k_2)-\gamma_{i_2,j_2}(k_2)| |>\frac{Cu^{1/3}}{b_n^{1/3}}\bigg\} \notag\\
\lesssim&~ \exp\{-Cn^{(3-2\rho)/3}l_n^{-1}u^{2/3}\} + \exp\{-Cn^{(3-\rho)/9}l_n^{-1/3}u^{1/9}\} \notag\\
&+\exp\{-Cn^{(3-\rho)/3}u^{1/3}\}+\exp\{-Cn^{(6-\rho)/18}u^{1/18}\}  \,.  \label{eq:A201tail}
\end{align}
Analogously, we can  show the upper bounds given in \eqref{eq:A101tail} and \eqref{eq:A201tail}  also hold for $\max_{\ell_1,\ell_2}\bP\{\tilde{A}_{10}(\ell_1,\ell_2)>u\}$ and $\max_{\ell_1,\ell_2}\bP\{\tilde{A}_{20}(\ell_1,\ell_2)>u\}$, respectively.
Since $A_{10}\leq \max_{\ell_1,\ell_2}A_{10}(\ell_1,\ell_2)+\max_{\ell_1,\ell_2}\tilde{A}_{10}(\ell_1,\ell_2)$ and $A_{20}\leq \max_{\ell_1,\ell_2}A_{20}(\ell_1,\ell_2)+\max_{\ell_1,\ell_2}\tilde{A}_{20}(\ell_1,\ell_2)$,
by the Bonferroni inequality, 	for any $u\gg l_n^3n^{\rho-3}$,
\begin{align}\label{eq:A10tail}
\bP(A_{10}>u) + \bP(A_{20}>u)
\lesssim&~  nr^2l_n^2\exp\{-Cn^{(3-2\rho)/3}l_n^{-1}u^{2/3}\}+ nr^2l_n^2\exp\{-Cn^{(3-\rho)/9}l_n^{-1/3}u^{1/9}\} \notag\\
&+r^2l_n^2\exp\{-Cn^{(3-\rho)/3}u^{1/3}\}
+r^2l_n^2\exp\{-Cn^{(6-\rho)/18}u^{1/18}\}  \,.
\end{align}
Since $l_n=o(n)$ and $r\geq n^{\kappa}$ for some sufficiently small constant $\kappa>0$, then we have
$
A_{10}
= O_{\rm p}\{n^{(2\rho-3)/2}l_n^{3/2}\log^{3/2}(r)\} =A_{20}
$
provided that $\log r\ll n^{1/5}l_n^{-1/5}$.
By the definition of $A_{15}$, using the same arguments, we know the upper bound given in \eqref{eq:A10tail} also holds for $\bP(A_{15}>u)$ and
$
A_{15}=  O_{\rm p}\{n^{(2\rho-3)/2}l_n^{3/2}\log^{3/2}(r)\}
$.  $\hfill\Box$

\subsubsection{An upper bound for $\bP(|\widehat{\bXi}-\bXi^*|_\infty>u)$}\label{subsec:upperbound}
By \eqref{eq:A1tail}, \eqref{eq:A2tail}, \eqref{eq:A4.tail}, \eqref{eq:A5.tail}, \eqref{eq:A7tail}, \eqref{eq:A9tail} and \eqref{eq:A25tail}, it holds that
\begin{align*} 
\bP(|\widehat{\bXi}-\bXi^*|_\infty>u)
\lesssim &~ nr^2l_n^2 \exp\{-Cn^{(2\vartheta+4\rho-6\rho\vartheta-2)/(2\vartheta-1)}u^2\}
+nr^2l_n^2\exp(-Cn^{1-\rho}l_n^{-1}u)  \notag\\
&
+ nr^2l_n^2\exp\{-Cn^{(3-2\rho)/3}l_n^{-1}u^{2/3}\}
+ n^2r^2l_n^2 \exp\{-Cn^{(\vartheta+2\rho-3\rho\vartheta-1)/(4\vartheta-2)}u^{1/2}\}  \notag\\
&+r^2l_n^2\exp\{-Cn^{(1-\rho)/3}u^{1/3}\}
+nr^2l_n^2 \exp\{-Cn^{(2\vartheta+3\rho-4\rho\vartheta-2)/(8\vartheta-4)}u^{1/4}\} \notag\\
& +nr^2l_n^2\exp\{-Cn^{(2-\rho)/6}l_n^{-1/3}u^{1/6}\}
+ nr^2l_n^2\exp\{-Cn^{(3-\rho)/9}l_n^{-1/3}u^{1/9}\}  \notag\\
& + nr^2l_n^2\exp\{-Cn^{(4-\rho)/12}u^{1/12}\}
+r^2l_n^2\exp\{-Cn^{(6-\rho)/18}u^{1/18}\}  \notag\\
& 	+r^2l_n^2\exp\{-Cn^{(2\vartheta+2\rho-2\rho\vartheta-2)/(7\vartheta-4)}\}  
\end{align*}
for any  $u\gg l_n^2n^{\rho-2}$.   $\hfill\Box$

\subsection{Proof of Lemma \ref{la:gammatail2}}\label{subsec:pf:gammatail2}
Given $i,j\in [r(2l_n+1)]$ and $0\leq q\leq M$, we write $z_{t}=c_{i,t}c_{j,t-q}-\bE(c_{i,t}c_{j,t-q})$.
Due to $\max_{t\in[\tilde{n}]}\max_{j\in[r(2l_n+1)]}\bP(|c_{j,t}|>u)\leq C\exp(-Cu)$ for any $u>0$,
by Lemma 2 in the supplementary material of \cite{CTW_2013_supp},  
$\bP(|z_{t}|>u)\leq C\exp(-Cu^{1/2})$ for any $u>0$. Notice that $\{z_t\}_{t=q+1}^{\tilde{n}}$ is an $\alpha$-mixing sequence with $\alpha$-mixing coefficient  $\alpha_z(k)\lesssim \exp(-C|k-2l_n-M|_+)$.   Due to $2l_n\leq M$, by Lemma L1 in the supplementary material of \cite{CCW_2021_supp} 
with $(\tilde{L}_{\tilde{n}}, \tilde{B}_{\tilde{n}}, \tilde{j}_{\tilde{n}}, r_1, r_2, r)=(1, 1, 2l_n+M, 1/2, 1, 1/4)$,
\begin{align*}
\mathbb{P}\bigg(\bigg|\sum_{t=q+1}^{\tilde{n}}z_t\bigg|\geq u\bigg)
\lesssim \exp\bigg(-\frac{Cu^2}{nM}\bigg)+\exp\bigg(-\frac{Cu^{1/4}}{M^{1/4}}\bigg)
\end{align*}
for any $u> 0$, where the upper bound in the above inequality does not depend on $i,j$ and $q$. We have completed the proof of Lemma \ref{la:gammatail2}.  $\hfill\Box$

\subsection{Proof of Lemma \ref{la:gammatail}}\label{subsec:pf:gammatail}
Given $i,j\in[p]$ and $-l_n\leq k\leq l_n$, we write $z_{t}=\mathring{x}_{i,t+k}\mathring{x}_{j,t}-\gamma_{i,j}(k)$.
By \eqref{eq:t1},   $\bP(|z_{t}|>u)\leq C\exp(-Cu)$ for any $u>0$. Notice that $\{z_t\}_{t=l_n+s_1}^{l_n+s_2}$ is an $\alpha$-mixing sequence with $\alpha$-mixing coefficient  $\alpha_z(k)\lesssim \exp(-C|k-2l_n|_+)$. Write $s=s_2-s_1$.  By Lemma L1 in the supplementary material of \cite{CCW_2021_supp} 
with
$(\tilde{L}_{\tilde{n}}, \tilde{B}_{\tilde{n}}, \tilde{j}_{\tilde{n}}, r_1, r_2, r)=(1,1,2l_n,1,1,1/3)$,
\begin{align*}
\mathbb{P}\bigg(\bigg|\sum_{t=l_n+s_1}^{l_n+s_2}z_t\bigg|\geq u\bigg)
\lesssim \exp\bigg(-\frac{Cu^2}{sl_n}\bigg)+\exp\bigg(-\frac{Cu^{1/3}}{l_n^{1/3}}\bigg)
\end{align*}
for any $u> 0$, where  the upper bound in above inequality does not depend on $i,j$ and $k$. We
have completed the proof of Lemma \ref{la:gammatail}. $\hfill\Box$

\subsection{Proof of Lemma \ref{la:tildegammatail}}\label{subsec:pftildegamma}
Without loss of generality, we assume $\bmu=\bzero$. Then $\gamma_{i,j}(k)={\rm Cov}(x_{i, t+k}, x_{j,t}) = \bE(x_{i, t+k}x_{j, t})$.  Recall $\hat{\gamma}_{i,j}(k)=n^{-1}\sum_{t=1}^{n-k}(x_{i,t+k}-\bar{x}_{i})(x_{j,t}-\bar{x}_{j})$ if $k\geq0$ and  $	\hat{\gamma}_{i,j}(k)=n^{-1}\sum_{t=-k+1}^n(x_{i,t+k}-\bar{x}_{i})(x_{j,t}-\bar{x}_{j})$ if $k<0$,
with $\bar{\bx}=(\bar{x}_{1},\ldots,\bar{x}_p)^{\T}=n^{-1}\sum_{t=1}^n\bx_t$.
By the triangle inequality and the Bonferroni inequality, 	for any $u>0$,
\begin{align*}
&\max_{(i,j)\in\cI}\max_{0\leq k\leq l_n} \bP\big\{|\hat{\gamma}_{i,j}(k)-\gamma_{i,j}(k)|>u\big\}  \\
&~~~~~~~ \leq \underbrace{\max_{(i,j)\in\cI}\max_{0\leq k\leq l_n}\bP\bigg[\bigg|\frac{1}{n}\sum_{t=1}^{n-k}\{x_{i,t+k}x_{j,t} -\gamma_{i,j}(k)\}\bigg|>\frac{u}{5}\bigg]}_{\rm(a)}
+ \underbrace{\max_{(i,j)\in\cI}\max_{0\leq k\leq l_n}\bP\bigg\{\bigg|\frac{k}{n}\gamma_{i,j}(k)\bigg| >\frac{u}{5}\bigg\}}_{\rm (b)}\\
&~~~~~~~~~\, +\underbrace{\max_{(i,j)\in\cI}\max_{0\leq k\leq l_n}\bP\bigg\{\bigg|\bar{x}_{i}\bigg(\frac{1}{n}\sum_{t=1}^{n-k}x_{j,t} \bigg) \bigg| >\frac{u}{5}\bigg\}}_{\rm(c)}
+\underbrace{\max_{(i,j)\in\cI}\max_{0\leq k\leq l_n}\bP\bigg\{\bigg|\bar{x}_{j}\bigg(\frac{1}{n}\sum_{t=1}^{n-k}x_{i,t+k}\bigg) \bigg| >\frac{u}{5}\bigg\}}_{\rm(d)} \\
&~~~~~~~~~\, + \underbrace{\max_{(i,j)\in\cI}\max_{0\leq k\leq l_n}\bP\bigg(\bigg|\frac{n-k}{n}\bar{x}_{i}\bar{x}_{j}\bigg| >\frac{u}{5}\bigg)}_{\rm(e)}  \,.
\end{align*}
By Lemma \ref{la:gammatail},
\begin{align}\label{eq:term_a}
{\rm(a)} \lesssim&~  \exp\bigg(-\frac{Cnu^2}{l_n}\bigg)
+\exp\bigg(-\frac{Cn^{1/3}u^{1/3}}{l_n^{1/3}}\bigg)
\end{align}
for any $u>0$.
Since $\max_{(i,j)\in\cI}\max_{-l_n\leq k\leq l_n}|\gamma_{i,j}(k)|\leq C$, then
\begin{align}\label{eq:term_b}
{\rm (b)}
=\max_{(i,j)\in\cI}\max_{0\leq k\leq l_n} \bP\bigg\{|\gamma_{i,j}(k)|>\frac{Cnu}{|k|}\bigg\}
\leq \max_{(i,j)\in\cI}\max_{0\leq k\leq l_n}\bP\bigg\{|\gamma_{i,j}(k)|>\frac{Cnu}{l_n}\bigg\}=0
\end{align}
for any $u\gg l_nn^{-1}$.
By the Bonferroni inequality and \eqref{eq:sumx.tail},
\begin{align}\label{eq:term_c}
{\rm (c)} \leq&~ \max_{i}\bP\bigg(\bigg|\frac{1}{n}\sum_{t=1}^n{x}_{i,t}\bigg| >Cu^{1/2}\bigg)
+ \max_j\max_{0\leq k\leq l_n} \bP\bigg(\bigg|\frac{1}{n}\sum_{t=1}^{n-k}x_{j,t}  \bigg| >Cu^{1/2}\bigg) \notag\\
\lesssim&~  \exp(-Cnu)+ \exp(-Cn^{1/3}u^{1/6})
\end{align}
for any $u>0$. Analogously, we  can show  the same upper bound given in \eqref{eq:term_c} also holds  for ${\rm (d)}$ and ${\rm(e)}$. Combining with \eqref{eq:term_a}, \eqref{eq:term_b} and \eqref{eq:term_c}, it holds that
\begin{align}\label{eq:gamma.tail}
\max_{(i,j)\in\cI}\max_{0\leq k\leq l_n} \bP\{|\hat{\gamma}_{i,j}(k)-\gamma_{i,j}(k)|>u\}
\lesssim&~
\exp(-Cnl_n^{-1}u^2)
+\exp(-Cn^{1/3}l_n^{-1/3}u^{1/3}) \notag\\
& + \exp(-Cnu)+ \exp(-Cn^{1/3}u^{1/6})
\end{align}
for any $u\gg l_nn^{-1}$.
Using the same arguments, we can  have the same upper bound given in \eqref{eq:gamma.tail} also holds for $\max_{(i,j)\in\cI}\max_{-l_n\leq k<0}\bP\{|\hat{\gamma}_{i,j}(k)-\gamma_{i,j}(k)|>u\}$ for any $u\gg l_nn^{-1}$. We have completed the proof of Lemma \ref{la:tildegammatail}.    $\hfill\Box$

\section{Statistical inference for coherence matrix}

The coherence matrix at frequency $\omega$ is defined as 
    $$
        {\rm coh}(\omega) = ({\rm coh}_{i,j}(\omega))_{p\times p} = \{\bD(\omega)\}^{-1/2} \bF(\omega) \{\bD(\omega)\}^{-1/2}\,, $$
    where ${\bf F}(\omega)=(f_{i,j}(\omega))_{p\times p}$ and $\bD(\omega)=\diag\{f_{1,1}(\omega),\ldots,f_{p,p}(\omega)\}$. If we would like to consider the following two statistical inference problems for coherences, our current procedures stated in Section \ref{applications} can be applied directly since ${\rm coh}_{i,j}(\omega)=0$ if and only if $f_{i,j}(\omega)=0$:
    \begin{itemize}
    \item (Global hypothesis testing)
    \[
    H_0: {\rm coh}_{i,j}(\omega)=0~\textrm{for any}~(i,j)\in\mathcal{I}~\textrm{and}~\omega\in\mathcal{J}~~~~\textrm{versus}~~~~H_1:H_0~\textrm{is not true}.
    \]
    \item (Multiple testing with FDR control) Given $\{\mathcal{I}^{(q)},\mathcal{J}^{(q)}\}$ with $\mathcal{I}^{(q)}\subset[p]^2$ and $\mathcal{J}^{(q)}\subset[-\pi,\pi)$, consider $Q$ hypothesis testing problems
    \[
    H_{0,q}:{\rm coh}_{i,j}(\omega)=0~\textrm{for any}~(i,j)\in\mathcal{I}^{(q)}~\textrm{and}~\omega\in\mathcal{J}^{(q)}~~~~\textrm{versus}~~~~H_{1,q}:H_{0,q}~\textrm{is not true}.
    \]
    \end{itemize}
  Also, Theorems \ref{tm:H0} and \ref{tm:FDR} can be used for these two problems. 
  
  If we would like to construct the simultaneous inference for the coherence matrix, the procedure will be a little bit different. 
    Based on our proposed  spectral density estimator $\widehat{\bF}(\omega)$, we can estimate ${\rm coh}(\omega)$ by
    \begin{align*}
        \widehat{\rm coh}(\omega)=(\widehat{{\rm coh}}_{i,j}(\omega))_{p\times p} = \{\widehat{\bD}(\omega)\}^{-1/2} \widehat{\bF}(\omega) \{\widehat{\bD}(\omega)\}^{-1/2} \,,
    \end{align*}
    where $\widehat{\bD}(\omega)=\diag\{\hat{f}_{1,1}(\omega),\ldots,\hat{f}_{p,p}(\omega)\}$. To construct the simultaneous inference of the coherence matrix, we need to consider the Gaussian approximation to the distribution of
    \[
    \mathcal{T}_{{\rm coh}}=\sup_{\omega\in\mathcal{J}}\max_{(i,j)\in\mathcal{I}}|\sqrt{n/l_n}\{\widehat{{\rm coh}}_{i,j}(\omega)-{\rm coh}_{i,j}(\omega)\}|^2\,.
    \]
Lemma \ref{la:uniform} in Section \ref{sec:auxlemma} indicates that
$$\sup_{\omega\in[-\pi,\pi]}\max_{(i,j)\in \cI}|\hat{f}_{i,j}(\omega)-f_{i,j}(\omega)| = O_{\rm p}\{n^{-1/2}l_n^{3/2}(\log l_n)\log^{1/2}(r)\}$$
under some regularity conditions. For each given $\omega\in\mathcal{J}$ and $(i,j)\in\mathcal{I}$, it holds that
\begin{align*}
& \sqrt{n/l_n}\{\widehat{{\rm coh}}_{i,j}(\omega)-{\rm coh}_{i,j}(\omega)\}  \\
&~~~~~~~=\frac{1}{f^{1/2}_{i,i}(\omega)f^{1/2}_{j,j}(\omega)}\sqrt{\frac{n}{l_n}}\{\hat{f}_{i,j}(\omega)-f_{i,j}(\omega)\}-\frac{1}{2}\frac{f_{i,j}(\omega)}{f_{i,i}^{3/2}(\omega)f_{j,j}^{1/2}(\omega)}\sqrt{\frac{n}{l_n}}\{\hat{f}_{i,i}(\omega)-f_{i,i}(\omega)\}\\
&~~~~~~~~~~-\frac{1}{2}\frac{f_{i,j}(\omega)}{f_{i,i}^{1/2}(\omega)f_{j,j}^{3/2}(\omega)}\sqrt{\frac{n}{l_n}}\{\hat{f}_{j,j}(\omega)-f_{j,j}(\omega)\}+O_{\rm p}\{n^{-1/2}l_n^{5/2}(\log l_n)^2\log r\}\,,
\end{align*}
where the reminder term $O_{\rm p}\{n^{-1/2}l_n^{5/2}(\log l_n)^2\log r\}$ holds uniformly over $\omega\in\mathcal{J}$ and $(i,j)\in\mathcal{I}$. By \eqref{eq:expansion} and Lemmas \ref{la:remainder}--\ref{la:leading} in Section \ref{sec:auxlemma}, we know 
$$\sup_{\omega\in[-\pi,\pi]}\max_{(i,j)\in\cI}|\hat{f}_{i,j}(\omega) - f_{i,j}(\omega)-\zeta_{i,j}(\omega)|  \lesssim \exp(-Cl_n) + O_{\rm p}(l_n^2 n^{-1} \log r)  $$ with
\begin{align*}
  \zeta_{i,j}(\omega) = \frac{1}{  n}\sum_{t=l_n+1}^{n-l_n}\bigg[ \frac{1}{2\pi}\sum_{k=-l_n}^{l_n} \cW\bigg(\frac{k}{l_n}\bigg)\{\mathring{x}_{i,t+k}\mathring{x}_{j,t} - \gamma_{i,j}(k)\}e^{-\iota k\omega}\bigg] \,, 
\end{align*}
which implies 
\begin{align*}
&\sqrt{n/l_n}\{\widehat{{\rm coh}}_{i,j}(\omega)-{\rm coh}_{i,j}(\omega)\}\\
&~~~~~~=\frac{\sqrt{n/l_n}\zeta_{i,j}(\omega)}{f^{1/2}_{i,i}(\omega)f^{1/2}_{j,j}(\omega)}-\frac{1}{2}\frac{f_{i,j}(\omega)\sqrt{n/l_n}\zeta_{i,i}(\omega)}{f_{i,i}^{3/2}(\omega)f_{j,j}^{1/2}(\omega)}-\frac{1}{2}\frac{f_{i,j}(\omega)\sqrt{n/l_n}\zeta_{j,j}(\omega)}{f_{i,i}^{1/2}(\omega)f_{j,j}^{3/2}(\omega)}\\
&~~~~~~~~~~+O_{\rm p}\{n^{-1/2}l_n^{5/2}(\log l_n)^2\log r\}\\
&~~~~~~=\frac{1}{\sqrt{n}f_{i,i}^{1/2}(\omega)f_{j,j}^{1/2}(\omega)}\sum_{t=l_n+1}^{n-l_n}\frac{1}{2\pi\sqrt{l_n}}\sum_{k=-l_n}^{l_n}\cW\bigg(\frac{k}{l_n}\bigg)\bigg[\{\mathring{x}_{i,t+k}\mathring{x}_{j,t} - \gamma_{i,j}(k)\}e^{-\iota k\omega}\\
&~~~~~~~~~\underbrace{~~~~~~~~~~~~~~~-\frac{1}{2}\{\mathring{x}_{i,t+k}\mathring{x}_{i,t}-\gamma_{i,i}(k)\}\frac{e^{-\iota k\omega}f_{i,j}(\omega)}{f_{i,i}(\omega)}-\frac{1}{2}\{\mathring{x}_{j,t+k}\mathring{x}_{j,t}-\gamma_{j,j}(k)\}\frac{e^{-\iota k\omega}f_{i,j}(\omega)}{f_{j,j}(\omega)}\bigg]}_{{\rm LT}_{i,j}(\omega)}\\
&~~~~~~~~~~+O_{\rm p}\{n^{-1/2}l_n^{5/2}(\log l_n)^2\log r\}
\end{align*}
provided that $l_n\geq C'\log n$ for some sufficiently large constant $C'>0$. Let $\bchi(\cdot)=\{\chi_1(\cdot),\chi_2(\cdot)\}$ be a given bijective mapping from $[r]$ to $\mathcal{I}$ such that for any $(i,j)\in\mathcal{I}$, there exists a unique $\ell\in[r]$ satisfying $(i,j)=\bchi(\ell)$. Recall $\tilde{n}=n-2l_n$. For each $t\in[\tilde{n}]$ and $\ell\in[r]$, we define a $(6l_n+3)$-dimensional vector
\[
\tilde{\bc}_{\ell,t}=\frac{1}{2\pi}
\begin{pmatrix}
\mathring{x}_{\chi_1(\ell),t}\mathring{x}_{\chi_2(\ell),t+l_n}-\gamma_{\bchi(\ell)}(-l_n)\\
\mathring{x}_{\chi_1(\ell),t}\mathring{x}_{\chi_1(\ell),t+l_n}-\gamma_{\chi_1(\ell),\chi_1(\ell)}(-l_n)\\
\mathring{x}_{\chi_2(\ell),t}\mathring{x}_{\chi_2(\ell),t+l_n}-\gamma_{\chi_2(\ell),\chi_2(\ell)}(-l_n)\\
\vdots\\
\mathring{x}_{\chi_1(\ell),t+2l_n}\mathring{x}_{\chi_2(\ell),t+l_n}-\gamma_{\bchi(\ell)}(l_n)\\
\mathring{x}_{\chi_1(\ell),t+2l_n}\mathring{x}_{\chi_1(\ell),t+l_n}-\gamma_{\chi_1(\ell),\chi_1(\ell)}(l_n) \\ 
\mathring{x}_{\chi_2(\ell),t+2l_n}\mathring{x}_{\chi_2(\ell),t+l_n}-\gamma_{\chi_2(\ell),\chi_2(\ell)}(l_n)
\end{pmatrix}
\]
and
\begin{align*}
&\tilde{\bA}(\omega)=\frac{1}{\sqrt{l_n}}
\begin{pmatrix}
\bI_3\otimes\begin{pmatrix}
\cos(-l_n\omega)\\
-\sin(-l_n\omega)
\end{pmatrix},\ldots,\bI_3\otimes\begin{pmatrix}
\cos(l_n\omega)\\
-\sin(l_n\omega)
\end{pmatrix}\end{pmatrix}\\
&~~~~~~~~~~~~~~~~~~\times
\big[{\rm diag}\{
\cW(-l_n/l_n),\ldots,\cW(l_n/l_n)
\}\otimes \bI_3\big]\,.
\end{align*}
For each given $(i,j)\in\mathcal{I}$, let $(i,j)=\bchi(\ell)$ for some $\ell\in[r]$. Write 
\[
\bB_\ell(\omega)=\begin{pmatrix}
{\rm Re}\{f_{\bchi(\ell)}(\omega)\} & -{\rm Im}\{f_{\bchi(\ell)}(\omega)\}\\
{\rm Im}\{f_{\bchi(\ell)}(\omega)\} & {\rm Re}\{f_{\bchi(\ell)}(\omega)\}
\end{pmatrix}\,.
\]
It then holds that  
\begin{align*}
\begin{pmatrix}
{\rm Re}\{{\rm LT}_{i,j}(\omega)\}\\
{\rm Im}\{{\rm LT}_{i,j}(\omega)\}
\end{pmatrix} 
=&~\begin{pmatrix}
\bI_2,
-\frac{1}{2}f_{\chi_1(\ell),\chi_1(\ell)}^{-1}(\omega)\bB_\ell(\omega),
-\frac{1}{2}f_{\chi_2(\ell),\chi_2(\ell)}^{-1}(\omega)\bB_\ell(\omega)
\end{pmatrix} \\
&~~~~~~~~~~~~~~~~~~~
\times \frac{\tilde{\bA}(\omega)}{ f_{\chi_1(\ell),\chi_1(\ell)}^{1/2}(\omega) f_{\chi_2(\ell),\chi_2(\ell)}^{1/2}(\omega)} 
\bigg(\frac{1}{\sqrt{n}}\sum_{t=1}^{\tilde{n}}\tilde{\bc}_{\ell,t}\bigg)\,.
\end{align*}

Define
\begin{align*}
\bW(\omega)=&~{\rm diag}\bigg\{
\begin{pmatrix}
\bI_2,
-\frac{1}{2}f_{\chi_1(1),\chi_1(1)}^{-1}(\omega)\bB_1(\omega), -\frac{1}{2}f_{\chi_2(1),\chi_2(1)}^{-1}(\omega)\bB_1(\omega) 
\end{pmatrix}, \\
&~~~~~~~~~~~~~~~~\ldots, 
\begin{pmatrix}
\bI_2,
-\frac{1}{2}f^{-1}_{\chi_1(r),\chi_1(r)}(\omega)\bB_r(\omega), -\frac{1}{2}f^{-1}_{\chi_2(r),\chi_2(r)}(\omega)\bB_r(\omega)
\end{pmatrix} 
\bigg\}\\
&~~~~~\times \left[ {\rm diag}\big\{f_{\chi_1(1),\chi_1(1)}^{-1/2}(\omega)f_{\chi_2(1),\chi_2(1)}^{-1/2}(\omega),\ldots,f_{\chi_1(r),\chi_1(r)}^{-1/2}(\omega)f_{\chi_2(r),\chi_2(r)}^{-1/2}(\omega)\big\}\otimes \tilde{\bA}(\omega) \right] \,.
\end{align*}
For given $\mathcal{J}=\{\omega_1,\ldots,\omega_K\}$, repeating the Gaussian approximation technique used to establish Proposition \ref{tm:1}(i), we know the distribution of $\mathcal{T}_{{\rm coh}}$ can be approximated by the distribution of 
$
\max_{j\in[Kr]}(\tilde{s}_{2j-1}^2+\tilde{s}_{2j}^2)$,
where $\tilde{\bs}=(\tilde{s}_1,\ldots,\tilde{s}_{2Kr})^{\T}$ is a $(2Kr)$-dimensional normally distributed random vector with mean zero and covariance
\[
\begin{pmatrix}
\bW(\omega_1) \\
\vdots \\
\bW(\omega_K) 
\end{pmatrix}
{\rm Var}\bigg(\frac{1}{\sqrt{\tilde{n}}}\sum_{t=1}^{\tilde{n}}\tilde{\bc}_t\bigg)
\begin{pmatrix}
\bW(\omega_1) \\
\vdots \\
\bW(\omega_K) 
\end{pmatrix}^{\T}
\]
for $\tilde{\bc}_t=(\tilde{\bc}_{1,t}^{\T},\ldots,\tilde{\bc}_{r,t}^{\T})^{\T}$. For given $\cJ = [\omega_L,\omega_U]$, repeating the Gaussian approximation technique used to establish Proposition \ref{tm:1}(ii), we can also obtain the distribution of $\mathcal{T}_{\rm coh}$ which can be approximated by the distribution of $\sup_{\omega\in\cJ}\max_{j\in[r]}\{\tilde{g}_{2j-1}^2(\omega)+\tilde{g}_{2j}^2(\omega)\}$, where $\tilde{\bg}(\omega)=\{\tilde{g}_1(\omega),\ldots,\tilde{g}_{2r}(\omega)\}^{\T}$ is a $(2r)$-dimensional Gaussian process with mean zero and covariance function 
$$
\bW(\omega_1){\rm Var}\bigg(\frac{1}{\sqrt{\tilde{n}}}\sum_{t=1}^{\tilde{n}}\tilde{\bc}_t\bigg) \bW^{\T}(\omega_2) $$ 
for $\tilde{\bc}_t=(\tilde{\bc}_{1,t}^{\T},\ldots,\tilde{\bc}_{r,t}^{\T})^{\T}$. Then the parametric bootstrap procedure given in Section \ref{sec:parboot} can still be applied.

\setlength{\bibsep}{0.2pt plus 1ex}

\begin{spacing}{1.14}

\end{spacing}

\end{document}